\documentclass[12pt,a4paper]{article}

\usepackage[margin=3cm]{geometry}

\usepackage{hyperref}
\usepackage[toc,page]{appendix}
\usepackage{amsmath,amsthm,amssymb}
\usepackage{color}
\usepackage{url}
\usepackage{tikz}
\usepackage{multirow} 
\usepackage{xcolor}
\usepackage{graphicx}
\usepackage{subdepth}
\usepackage{multicol}
\usepackage[space]{grffile}
\usepackage[nocompress]{cite}
\usepackage[outdir=./]{epstopdf}
\usepackage[all,rotate,poly]{xy}
\usepackage{rotating}
\usepackage{lscape}
\usepackage{amsbsy}
\usepackage{geometry}
\usepackage{xspace}
\usepackage{extarrows}
\usepackage{docmute}
\usepackage[textfont=it, labelfont=it]{caption}
\usepackage{xr}
\usepackage[numbers, sort]{natbib}
\usepackage{enumerate}
\usepackage{amsmath, amsthm, amssymb}
\usepackage[latin1]{inputenc}
\usepackage{stmaryrd}
\usepackage{braket}
\usepackage{xspace}
\usepackage{amsmath,amssymb,amsthm}
\usepackage{youngtab}

\usepackage{etoolbox}
\usepackage{tikz}
\usetikzlibrary{positioning}
\usetikzlibrary{arrows}
\usetikzlibrary{decorations.markings}
\usetikzlibrary{decorations.pathreplacing}
\usetikzlibrary{backgrounds}
\usetikzlibrary{shapes.geometric}
\usetikzlibrary{decorations.pathmorphing,intersections,calc}

\allowdisplaybreaks

\def\R{\mathbb{R}}

\def\I{\mathrm{I}}
\def\II{\mathrm{II}}
\def\III{\mathrm{III}}
\def\IV{\mathrm{IV}}
\def\V{\mathrm{V}}
\def\VI{\mathrm{VI}}

\def\lsscalevertex{1}

\DeclareGraphicsExtensions{.jpg,.png,.gif,.svg}

\newcommand{\rewrite}[4]{#1.\mathrm{\Pi}[#2, #4]}

\newcommand{\lift}[2]{#1.\mathrm{\Lambda}[#2]}
\newcommand{\diagid}[1]{\mathrm{Id}(#1)}

\newcommand{\inc}[2]{{\mathrm{Inc}}_{\mathrm r} (#1, #2)}
\newcommand{\increverse}[2]{{\mathrm{Inc}}_{\mathrm l} (#1,#2)}
\newcommand{\comm}[1]{}
\newcommand{\defref}[1]{\text{Def.}~\ref{#1}}
\newcommand{\myeqref}[1]{\text{Eq. }(\ref{#1})}

\renewcommand{\-}[0]{\nobreakdash-\hspace{0pt}}

\definecolor{drawcolour}{rgb}{0.75,0.75,1}
\definecolor{fillcolour}{rgb}{0.8,0.8,1}

\tikzset{region/.style={draw=drawcolour, fill=fillcolour, thick, fill opacity=0.8}}
\tikzset{morphism/.style={draw, circle, fill=white, inner sep=0pt, minimum width=13pt, font=\scriptsize, draw=black, opacity=1}}

\tikzset{blob/.style={circle, draw, fill=white, inner sep=2pt}}
\tikzset{label/.style={font=\tiny, inner sep=1pt}}

\usepackage{etoolbox}
\usepackage{tikz}
\usetikzlibrary{
    decorations.pathmorphing,intersections
}
\usepackage{calc}

\tikzset{dlabel/.style={inner sep=1pt}}
\tikzset{surface/.style={draw=none, fill=blue!50, fill opacity=0.5}}
\tikzset{surface2/.style={draw=none, fill=yellow!50, fill opacity=0.5}}
\tikzset{biunitary/.style={}}
\tikzset{edge/.style={thick}}

\tikzset{vertex/.style={draw, circle, inner sep=1.3pt, fill=black}}
\newcommand*\arrowoffset{0.45pt}
\makeatletter
\pgfarrowsdeclare{new double arrowhead}{new double arrowhead}
{
  \pgfutil@tempdima=-0.84pt%
  \advance\pgfutil@tempdima by-1.3\pgflinewidth%
  \pgfutil@tempdimb=0.21pt%
  \advance\pgfutil@tempdimb by.625\pgflinewidth%
  \pgfarrowsleftextend{+\pgfutil@tempdima}
  \advance\pgfutil@tempdimb by \arrowoffset 
  \pgfarrowsrightextend{+\pgfutil@tempdimb}
}
{
  \pgfsetstrokecolor{black}
  \pgfsetlinewidth{0.9pt}
  \pgfutil@tempdima=0.28pt%
  \advance\pgfutil@tempdima by.3\pgflinewidth%
  \pgfsetlinewidth{0.8\pgflinewidth}
  \pgfsetdash{}{+0pt}
  \pgfsetroundcap
  \pgfsetroundjoin
  \pgfpathmoveto{\pgfpoint{-3\pgfutil@tempdima+\arrowoffset}{4\pgfutil@tempdima}}
  \pgfpathcurveto{\pgfpoint{-2.75\pgfutil@tempdima+\arrowoffset}{2.5\pgfutil@tempdima}}
                {\pgfpoint{\arrowoffset}{0.25\pgfutil@tempdima}}
                {\pgfpoint{0.75\pgfutil@tempdima+\arrowoffset}{0pt}}
  \pgfpathcurveto
                {\pgfpoint{\arrowoffset}{-0.25\pgfutil@tempdima}}
{\pgfpoint{-2.75\pgfutil@tempdima+\arrowoffset}{-2.5\pgfutil@tempdima}}
{\pgfpoint{-3\pgfutil@tempdima+\arrowoffset}{-4\pgfutil@tempdima}}
  \pgfusepathqstroke
}
\makeatother
\makeatletter
\pgfkeys{%
  /tikz/on layer/.code={
    \pgfonlayer{#1}\begingroup
    \aftergroup\endpgfonlayer
    \aftergroup\endgroup
  },
  /tikz/node on layer/.code={
    \gdef\node@@on@layer{%
      \setbox\tikz@tempbox=\hbox\bgroup\pgfonlayer{#1}\unhbox\tikz@tempbox\endpgfonlayer\egroup}
    \aftergroup\node@on@layer
  },
  /tikz/end node on layer/.code={
    \endpgfonlayer\endgroup\endgroup
  }
}
\def\node@on@layer{\aftergroup\node@@on@layer} 
\makeatother

\def\A{{\color{blue!50}\ensuremath{A}}}
\def\B{{\color{red!50}\ensuremath{B}}}
\def\CC{{\color{green!70}\ensuremath{C}}}

\pgfdeclarelayer{background}
\pgfdeclarelayer{foreground}
\pgfsetlayers{background,main,foreground}

\tikzset{comod/.style={rectangle, minimum width=25pt, minimum height=10pt, draw, inner sep=1pt}}

\def\id{\ensuremath{\mathrm{id}}}

\renewcommand{\-}[0]{\nobreakdash-\hspace{0pt}}

\tikzset{arrow/.style={decoration={
    markings,
    mark=at position #1 with \arrow{angle 60}},
    postaction=decorate}
}
\tikzset{reversed arrow/.style={decoration={
    markings,
    mark=at position #1 with \arrow{angle 60 reversed}},
    postaction=decorate}
}

\tikzset{keyvertexcolour/.initial=black}
\tikzset{vertex colour/.style={keyvertexcolour={#1}}}
\tikzset {triangle/.style={
        draw,
        shape border rotate=-90,
        isosceles triangle,
        isosceles triangle apex angle=80,        
        minimum height=2em}}
\tikzset {triangleup/.style={
        draw,
        shape border rotate=90,
        isosceles triangle,
        isosceles triangle apex angle=80,
         minimum height=2em}}
\newlength\vertexradius
\newlength\innerradius
\def\halfanglesep{24}
\def\tripleanglesep{24}
\def\sideangle{30}
\def\tripleanglesep{40}
\makeatletter

\pgfdeclareshape{Vertex}
{
    \savedanchor\centerpoint
    {
        \pgf@x=0pt
        \pgf@y=0pt
    }
    \anchor{n1}
    {
        \pgfextractx{\pgf@x}{\pgfpointpolar{90+\halfanglesep}{\innerradius}}
        \pgfextracty{\pgf@y}{\pgfpointpolar{90+\halfanglesep}{\innerradius}}
    }
    \anchor{n2}
    {
        \pgfextractx{\pgf@x}{\pgfpointpolar{90-\halfanglesep}{\innerradius}}
        \pgfextracty{\pgf@y}{\pgfpointpolar{90-\halfanglesep}{\innerradius}}
    }
    \anchor{n}
    {
        \pgfextractx{\pgf@x}{\pgfpointpolar{90}{\innerradius}}
        \pgfextracty{\pgf@y}{\pgfpointpolar{90}{\innerradius}}
    }
    \anchor{ne1}
    {
        \pgfextractx{\pgf@x}{\pgfpointpolar{\sideangle+\halfanglesep}{\innerradius}}
        \pgfextracty{\pgf@y}{\pgfpointpolar{\sideangle+\halfanglesep}{\innerradius}}
    }
    \anchor{ne2}
    {
        \pgfextractx{\pgf@x}{\pgfpointpolar{\sideangle-\halfanglesep}{\innerradius}}
        \pgfextracty{\pgf@y}{\pgfpointpolar{\sideangle-\halfanglesep}{\innerradius}}
    }
    \anchor{ne}
    {
        \pgfextractx{\pgf@x}{\pgfpointpolar{\sideangle}{\innerradius}}
        \pgfextracty{\pgf@y}{\pgfpointpolar{\sideangle}{\innerradius}}
    }
    \anchor{se1}
    {
        \pgfextractx{\pgf@x}{\pgfpointpolar{-\sideangle+\halfanglesep}{\innerradius}}
        \pgfextracty{\pgf@y}{\pgfpointpolar{-\sideangle+\halfanglesep}{\innerradius}}
    }
    \anchor{se2}
    {
        \pgfextractx{\pgf@x}{\pgfpointpolar{-\sideangle-\halfanglesep}{\innerradius}}
        \pgfextracty{\pgf@y}{\pgfpointpolar{-\sideangle-\halfanglesep}{\innerradius}}
    }
    \anchor{s1}
    {
        \pgfextractx{\pgf@x}{\pgfpointpolar{-90-\halfanglesep}{\innerradius}}
        \pgfextracty{\pgf@y}{\pgfpointpolar{-90-\halfanglesep}{\innerradius}}
    }
    \anchor{s2}
    {
        \pgfextractx{\pgf@x}{\pgfpointpolar{-90+\halfanglesep}{\innerradius}}
        \pgfextracty{\pgf@y}{\pgfpointpolar{-90+\halfanglesep}{\innerradius}}
    }
    \anchor{sA}
    {
        \pgfextractx{\pgf@x}{\pgfpointpolar{-90-\tripleanglesep}{\innerradius}}
        \pgfextracty{\pgf@y}{\pgfpointpolar{-90-\tripleanglesep}{\innerradius}}
    }
    \anchor{sB}
    {
        \pgfextractx{\pgf@x}{\pgfpointpolar{-90}{\innerradius}}
        \pgfextracty{\pgf@y}{\pgfpointpolar{-90}{\innerradius}}
    }
    \anchor{sC}
    {
        \pgfextractx{\pgf@x}{\pgfpointpolar{-90+\tripleanglesep}{\innerradius}}
        \pgfextracty{\pgf@y}{\pgfpointpolar{-90+\tripleanglesep}{\innerradius}}
    }
    \anchor{s}
    {
        \pgfextractx{\pgf@x}{\pgfpointpolar{-90}{\innerradius}}
        \pgfextracty{\pgf@y}{\pgfpointpolar{-90}{\innerradius}}
    }
    \anchor{sw1}
    {
        \pgfextractx{\pgf@x}{\pgfpointpolar{180+\sideangle+\halfanglesep}{\innerradius}}
        \pgfextracty{\pgf@y}{\pgfpointpolar{180+\sideangle+\halfanglesep}{\innerradius}}
    }
    \anchor{sw2}
    {
        \pgfextractx{\pgf@x}{\pgfpointpolar{180+\sideangle-\halfanglesep}{\innerradius}}
        \pgfextracty{\pgf@y}{\pgfpointpolar{180+\sideangle-\halfanglesep}{\innerradius}}
    }
    \anchor{sw}
    {
        \pgfextractx{\pgf@x}{\pgfpointpolar{180+\sideangle}{\innerradius}}
        \pgfextracty{\pgf@y}{\pgfpointpolar{180+\sideangle}{\innerradius}}
    }
    \anchor{nw1}
    {
        \pgfextractx{\pgf@x}{\pgfpointpolar{180-\sideangle+\halfanglesep}{\innerradius}}
        \pgfextracty{\pgf@y}{\pgfpointpolar{180-\sideangle+\halfanglesep}{\innerradius}}
    }
    \anchor{nw2}
    {
        \pgfextractx{\pgf@x}{\pgfpointpolar{180-\sideangle-\halfanglesep}{\innerradius}}
        \pgfextracty{\pgf@y}{\pgfpointpolar{180-\sideangle-\halfanglesep}{\innerradius}}
    }
    \anchor{nw}
    {
        \pgfextractx{\pgf@x}{\pgfpointpolar{180-\sideangle}{\innerradius}}
        \pgfextracty{\pgf@y}{\pgfpointpolar{180-\sideangle}{\innerradius}}
    }
    \anchor{center}{\centerpoint}
    \anchorborder{\centerpoint}
    \backgroundpath
    {
        \pgfkeysgetvalue{/tikz/keyvertexcolour}{\vcol}
        \begin{pgfonlayer}{foreground}
        \pgfsetfillcolor{\vcol}
        \pgfsetfillopacity{1}
        \pgfpathcircle{\pgfpoint{0cm}{0cm}}{\vertexradius}
        \pgfusepath{fill}
        \pgfsetstrokecolor{black}
        \pgfsetlinewidth{0.8pt}
        \pgfpathcircle{\pgfpoint{0cm}{0cm}}{\vertexradius}
        \pgfusepath{stroke}
        \end{pgfonlayer}
    }
}
\makeatother

\usepgflibrary{shapes.geometric}

\theoremstyle{plain}
\newtheorem{theorem}{Theorem}
\newtheorem{lemma}[theorem]{Lemma}

\newtheorem{proposition}[theorem]{Proposition}

\theoremstyle{definition}
\newtheorem{definition}[theorem]{Definition}

\renewcommand{\-}[0]{\nobreakdash-\hspace{0pt}}

\makeatletter
\def\calign@preamble{%
   &\hfil\strut@
    \setboxz@h{\@lign$\m@th\displaystyle{##}$}%
    \ifmeasuring@\savefieldlength@\fi
    \set@field
    \hfil
    \tabskip\alignsep@
}
\let\cmeasure@\measure@
\patchcmd\cmeasure@{\divide\@tempcntb\tw@}{}{}{}
\patchcmd\cmeasure@{\divide\@tempcntb\tw@}{}{}{}
\patchcmd\cmeasure@{\ifodd\maxfields@
  \global\advance\maxfields@\@ne
  \fi}{}{}{}    
\newenvironment{calign}
{%
  \let\align@preamble\calign@preamble
  \let\measure@\cmeasure@
  \align
}
{%
  \endalign
}  
\makeatother
\newcommand\ignore[1]{}

\newcounter{kbcommcounter}
\newcounter{jvcommcounter}
\newcommand\JVcomm[1]{}

\begin{document}

\title{\bf Data structures for\\\bf quasistrict higher categories}
\author{
\begin{tabular}{@{}cc@{}}
\begin{tabular}{@{}c@{}}Krzysztof Bar\\Department of Computer Science\\University of Oxford\\\texttt{krzysztof.bar@cs.ox.ac.uk}\end{tabular}
&
\begin{tabular}{@{}c@{}}Jamie Vicary\\Department of Computer Science\\University of Oxford\\\texttt{jamie.vicary@cs.ox.ac.uk}\end{tabular}
\end{tabular}
}

\maketitle

\begin{abstract}
We present new data structures for quasistrict higher categories, in which associativity and unit laws hold strictly. Our approach has low axiomatic complexity compared to traditional algebraic definitions of higher categories, and we use it to give a practical definition of quasistrict 4\-category. It is also amenable to computer implementation, and we exploit this to give a new machine-verified algebraic proof that every adjunction of 1\-cells in a quasistrict 4\-category can be promoted to a coherent adjunction satisfying the butterfly equations.
\end{abstract}

\newpage
\tableofcontents

\newpage

\section{Introduction}
\label{chapterintro}

\subsection{Our contribution}

\paragraph{Motivation.} Higher category theory now plays an essential role in many areas of mathematics, physics and computer science. Among the most striking applications are the homotopy type theory programme on univalent foundations for mathematics~\cite{hottbook, Awodey_2008, Voevodsky_2006}, motivated by the work of Hofmann and Streicher on an intensional groupoid model for Martin\-L\"of type theory~\cite{Hofmann_1994}; the outline proof by Lurie~\cite{Lurie_2008} of the cobordism hypothesis of Baez and Dolan~\cite{Baez_1995}, and the associated revolution in topological quantum field theory~\cite{Lawrence_1993, Atiyah_1990, CSPthesis}; and Lurie's higher topos theory programme~\cite{Lurie_2009}, developing ideas going back to Grothendieck~\cite{Grothendieck_1983}, with broad implications for geometry. Other applications in computer science include concurrency~\cite{Bruni_2002, Goubault_2003}, rewriting~\cite{Mimram_2014, Lafont_1997, Guiraud_2012}, and quantum computation~\cite{Vicary_2012, Reutter_2016, Jaffe_2016}.

Nonetheless, higher category theory is ``generally regarded as a technical and forbidding subject''~\cite{Lurie_2009}. This may be in part because of the complexity of the definitions from an algebraic perspective. For example, even in the semistrict case where as much structure as possible is suppressed, semistrict 2\-categories comprise 3 sets equipped with 6 functions satisfying 12 equations~\cite{nLab:Strict2Category}; semistrict 3\-categories comprise 4 sets equipped with 19 functions satisfying 58 equations~\cite{DouglasHenriques}; and semistrict 4\-categories comprise 5 sets equipped with 34 functions\JVcomm{Does this include the Gray functor data?} satisfying 118 equations~\cite{Crans_1998}\footnote{These counts are conservative estimates obtained from close readings of the definitions, and are to some degree subjective. There are some issues with this definition of semistrict 4-category (see Section~\ref{fourtasaxioms}), but these do not affect the magnitudes of these numbers. We neglect set-theoretic issues.}. This complexity masks to some extent the naturalness of these structures, and can make them hard to work with directly when constructing or verifying proofs. Traditional proof assistants such as Coq cannot easily help with these difficulties\footnote{Homotopy type theory is perhaps an exception; we discuss this in Section~\ref{sec:relatedwork}.}, not least because they do not support the pasting diagram or string diagram notations which are prevalent in higher category theory.

In this paper we present a new approach to defining and working with globular higher categories, applying in the  \textit{quasistrict} case, meaning that composition is strictly associative and unital.\footnote{This includes the semistrict case mentioned above.} Our approach allows us to give a concise definition of semistrict 4\-category, which corrects some previous errors in the literature. We also give a definition of quasistrict 4\-category, a weaker structure which is more practical for the purpose of constructing proofs. Our proposal is computationally implemented, and we give details of a substantial formalized proof that we have developed in the quasistrict 4\-category setting, giving evidence of the correctness and practicality of our approach.

\paragraph{Signatures and diagrams.}

A group can be defined explicitly in terms of a set of elements, or implicitly in terms of a \textit{presentation} involving generators and relations. Similarly, we work with higher categories in terms of a presentation, encoded by a \textit{signature} that gives the generating objects, 1\-cells, 2\-cells, and so on, up to $n$\-cells for some $n \in \mathbb N$. For $k > 0$, the generating $k$\-cells are equipped with source and target $(k{-}1)$\-diagrams, encoded in terms of \textit{diagram} structures, a new concept that we introduce.

A diagram of dimension $0<k\leq n$ comprises a source $(k{-}1)$-diagram, a list of generating $k$-cells, and a list of attachment coordinates for these generators. It can be thought of as  a construction procedure: to build the composite $k$\-diagram, begin with the source $(k{-}1)$\-diagram, and attach the generators sequentially at their specified attachment coordinates.  
 As a result, each $n$\-cell in an $n$\-diagram is at a unique height; our diagrams therefore satisfy a \textit{generic position} property reminiscent of the central property in higher Morse theory~\cite{CSPthesis}.

A key advantage of our approach is reduced axiomatic complexity, and we can explain the source of this reduction informally. Let $D$ be a  traditional definition of some flavour of quasistrict higher category, such as one of the semistrict definitions listed above. Some of the function data in $D$ will correspond to vertical composition operations, which specify how two $k$-cells compose to produce another $k$-cell; in our approach vertical composition is simply concatenation of the appropriate lists, and no additional composition data is needed. Furthermore, some of the equation data in $D$ will encode strict associativity or unitality of vertical composition; in our approach these equations can be neglected, since it is a trivial property of lists that concatenation is strictly associative and unital, with the unit given by the empty list. Indeed, these missing equations become \textit{theorems} in our approach, and the major technical work of this paper is giving proofs of results of this kind.

For our new definitions, outlined below, we obtain the following: a semistrict 2\-category comprises 4 sets equipped with 7 functions satisfying 4 equations; a semistrict 3\-category comprises 5 sets equipped with 11 functions satisfying 6 equations; and a semistrict 4\-category comprises 6 sets equipped with X functions satisfying 8 equations\JVcomm{Re-evaluate these numbers.}. Note that we must consider more sets than before, since, for example, a presentation of a 2\-category comprises not only generating objects, 1\-cells and 2\-cells, but also generating equations. However, these sets can be finite even in nontrivial cases, since infinite categories can have finite presentations\footnote{That this happens for some higher categories of manifolds follows from the cobordism hypothesis, and indeed is one of the main reasons for interest in the hypothesis.}; this is the source of another substantial reduction in practice of complexity of the required data.

\paragraph{Graphical calculus.}
Based on ideas of Trimble~\cite{TrimbleDiagrams}, we sketch in Section~\ref{secgraphicalformalism} an informal graphical calculus for $n$\-diagrams, which can be made precise for dimensions $n \leq 3$, in which  an $n$\-diagram is represented as a labelled partitioned subspace of $\R^n$. This is consistent with previous proposals in dimension $n\leq 3$~\cite{joyal-street, barrett-graycategories, streettopology, selingersurvey}. However, a unique feature of our diagrams is their similarity to generic-position Morse diagrams~\cite{CSPthesis}. This is a crucial feature of our approach, which we now briefly explore.

In ordinary algebraic approaches to higher category theory, a $p$\-cell and a $q$\-cell can be composed in $\text{min}(p,q)$ ways. Two 3\-cells $\alpha$, $\beta$ can therefore be composed in 3 ways, which we can illustrate graphically as follows:
\begin{calign}
\nonumber
\begin{aligned}
\begin{tikzpicture}[thick, scale=0.8]
\draw [region] (-0.2,0) rectangle +(3,2);
\draw (1.3,0) to +(0,2);
\node [morphism] at (1.3,1) {$\beta$};
\draw [region] (0.0,-0.2) rectangle +(3,2);
\draw (1.5,-0.2) to (1.5,0.5);
\draw (1.5,1.07) to (1.5,1.8);
\node [morphism, fill opacity=0.7] at (1.5,0.8) {$\alpha$};
\end{tikzpicture}
\end{aligned}
&
\begin{aligned}
\begin{tikzpicture}[thick, scale=0.8]
\draw [region] (0,0) rectangle +(3,2);
\draw (1,0) to +(0,2);
\draw (2,0) to +(0,2);
\node [morphism] at (2,1) {$\beta$};
\node [morphism] at (1,1) {$\alpha$};
\end{tikzpicture}
\end{aligned}
&
\begin{aligned}
\begin{tikzpicture}[thick, scale=0.8]
\draw [region] (0,0) rectangle +(2,3);
\draw (1,0) to (1,2) to (1,3);
\node [morphism] at (1,2) {$\beta$};
\node [morphism] at (1,1) {$\alpha$};
\end{tikzpicture}
\end{aligned}
\end{calign}
In the first diagram the two 3\-cells are overlapping, while in the second they are at the same height; only the third diagram is in generic position, which is a essential requirement of our representation, as discussed above. This is a general phenomenon: of all the $\text{min}(p,q)$ ways to compose a $p$\-cell and a $q$\-cell, only the highest-dimensional composition yields a generic position diagram. We internalize this in our formalism, allowing \textit{at most one} composition of $p$\-diagram and a $q$\-diagram, along a common boundary $(\text{min}(p,q){-}1)$\-cell; this smaller number of composition operations helps to further reduce the complexity of our algebraic system. Generic-position perturbations of the missing composites can still be accessed by repeated whiskering and composition, as follows:
\begin{calign}
\nonumber
\begin{aligned}
\begin{tikzpicture}[thick, scale=0.8]
\draw [region] (0,0) rectangle +(3,3);
\draw (1,0) to (1, 3);
\node [morphism] at (1,0.8) {$\beta$};
\draw [region] (0.2,-0.2) rectangle +(3,3);
\draw (2.2,-0.2) to (2.2, 2.8);
\node [morphism] at (2.2,2) {$\alpha$};
\end{tikzpicture}
\end{aligned}
&
\begin{aligned}
\begin{tikzpicture}[thick, scale=0.8]
\draw [region] (0,0) rectangle +(3,3);
\draw (1,0) to (1, 3);
\draw (2,0) to (2, 3);
\node [morphism] at (1,1) {$\alpha$};
\node [morphism] at (2,2) {$\beta$};
\end{tikzpicture}
\end{aligned}
\end{calign}
In this way we retain full compositional expressivity while reaping the advantages of a generic-position representation.

\paragraph{Homotopy generators.}

\begin{figure}[t]
\tikzset{morphism/.style={fill=white, draw=black, circle, inner sep=0pt, font=\scriptsize, minimum width=12pt}}
\[
\tikzset{every picture/.style={yscale=0.8, xscale=0.8}}
\begin{array}{@{}c@{}}
\begin{aligned}
\begin{tikzpicture}[thick, scale=0.8]
\draw [region] (0,0) rectangle +(3,3);
\draw (1,0) to (1, 3);
\draw (2,0) to (2, 3);
\node [morphism] at (1,1) {$\alpha$};
\node [morphism] at (2,2) {$\beta$};
\end{tikzpicture}
\end{aligned}
\stackrel{\displaystyle\I}{\to}
\begin{aligned}
\begin{tikzpicture}[thick, scale=0.8]
\draw [region] (0,0) rectangle +(3,3);
\draw (1,0) to (1, 3);
\draw (2,0) to (2, 3);
\node [morphism] at (1,2) {$\alpha$};
\node [morphism] at (2,1) {$\beta$};
\end{tikzpicture}
\end{aligned}
\\[70pt]
\begin{aligned}
\begin{tikzpicture}[thick]
\draw [region] (0.7,2.1) rectangle +(2.2,3);
\draw (2.6,2.1) to [out=up, in=down] (1.9,2.8) to [out=up, in=down] (1.2,3.5) ;
\draw (1.2,3.5) to (1.2,5.1);
\draw [region] (0.9, 1.9) rectangle +(2.2,3);
\draw (1.2,1.9) to (1.2,2.8) to [out=up, in=down] (1.9,3.5);
\draw (1.9,1.9) to (1.9,2.1)  to [out=up, in=down] (2.6,2.8) to (2.6,3.5);
\draw (1.9,3.5) to [out=up, in=down] (2.25,4.2) to [out=up, in=down](1.9, 4.9) to (1.9,4.9);
\draw (2.6,3.5) to [out=up, in=down] (2.25,4.2) to [out=up, in=down](2.6, 4.9)to (2.6,4.9);
\node [morphism] at (2.25,4.2) {$\alpha$};
\end{tikzpicture}
\end{aligned}
\stackrel{\displaystyle\II}{\rightarrow}
\begin{aligned}
\begin{tikzpicture}[thick]
\draw [region] (0,1.4) rectangle +(2.2,3);
\draw (1.9,2.8) to [out=up, in=down] (1.2,3.5) to [out=up, in=down] (0.5,4.2) to(0.5,4.4);
\draw (1.9, 1.4)to(1.9,2.8);
\draw [region] (0.2, 1.2) rectangle +(2.2,3);
\draw (0.5,2.8) to [out=up, in=down] (0.5,3.5) to [out=up, in=down] (1.2,4.2) ;
\draw (1.2,2.8)  to [out=up, in=down] (1.9,3.5) to [out=up, in=down] (1.9,4.2);
\draw (0.5,1.2) to(0.5,1.4) to [out=up, in=down] (0.85,2.1) to [out=up, in=down](0.5, 2.8);
\draw (1.2,1.2) to(1.2,1.4) to [out=up, in=down] (0.85,2.1) to [out=up, in=down](1.2, 2.8);
\node [morphism] at (0.85,2.1) {$\alpha$};
\end{tikzpicture}
\end{aligned}
\\[80pt]
\begin{aligned}
\begin{tikzpicture}[thick, scale=0.8]
\draw [region] (-0.3,0) rectangle +(2.4,3);
\draw (1.8,0)  to [out=up, in=down] (1.1,0.8) to [out=up, in=down] (0.4,1.6)to (0.4,3);
\draw [region] (-0.1,-0.2) rectangle +(2.4,3);
\draw (1.1,-0.2)to (1.1,0) to [out=up, in=down] (1.8,0.8) to (1.8,1.6) to [out=up, in=down]  (1.1,2.4)to (1.1,2.8);
\draw [region] (0.1,-0.4) rectangle +(2.4,3);
\draw (0.4,-0.4)  to (0.4,0.8)  to [out=up, in=down] (1.1,1.6) to [out=up, in=down] (1.8,2.4)to (1.8,2.6);
\end{tikzpicture}
\end{aligned}
\begin{array}{@{}c@{}}
\stackrel{\II}{\rightarrow}\\
\begin{aligned}
\begin{tikzpicture}[thick, scale=0.8]
\node (1) [scale=1] at (0.6,0) {$\Uparrow$};
\node [right=-7pt] at (1.east) {$\IV$};
\end{tikzpicture}
\end{aligned}\\
\stackrel{\II ^{-1}}{\rightarrow}
\end{array}
\begin{aligned}
\begin{tikzpicture}[thick, scale=0.8]
\draw [region] (-0.3,0) rectangle +(2.4,3);
\draw (1.8,0)  to (1.8,0.8)  to[out=up, in=down] (1.1,1.6) to [out=up, in=down] (0.4,2.4)to (0.4,3);
\draw [region] (-0.1,-0.2) rectangle +(2.4,3);
\draw (1.1,-0.2)to (1.1,0) to [out=up, in=down] (0.4,0.8) to (0.4,1.6) to [out=up, in=down]  (1.1,2.4)to (1.1,2.8);
\draw [region] (0.1,-0.4) rectangle +(2.4,3);
\draw (0.4,-0.4)  to (0.4,0)  to [out=up, in=down] (1.1,0.8) to [out=up, in=down] (1.8,1.6)to (1.8,2.6);
\end{tikzpicture}
\end{aligned}
\end{array}
\hspace{1.4cm}
\begin{array}{@{}c@{}}
\begin{array}{@{}c@{\hspace{-5pt}}c@{\hspace{-5pt}}c@{}}
\begin{aligned}
\begin{tikzpicture}[thick, scale=0.8]
\draw [region] (0.7,2.1) rectangle +(2.2,3);
\draw (2.6,2.1) to [out=up, in=down] (1.9,2.8) to [out=up, in=down] (1.2,3.5);
\draw (1.2,3.5) to (1.2,5.1);
\draw [region] (0.9, 1.9) rectangle +(2.2,3);
\draw (1.2,1.9) to (1.2,2.8) to [out=up, in=down] (1.9,3.5);
\draw (1.9,1.9) to (1.9,2.1)  to [out=up, in=down] (2.6,2.8) to (2.6,3.5);
\draw (1.9,3.5) to [out=up, in=down] (2.25,4.2) to [out=up, in=down](1.9, 4.9) to (1.9,4.9);
\draw (2.6,3.5) to [out=up, in=down] (2.25,4.2) to [out=up, in=down](2.6, 4.9)to (2.6,4.9);
\node [morphism] at (2.25,4.2) {$A$};
\end{tikzpicture}
\end{aligned}
&\stackrel{\II}{\rightarrow}&
\begin{aligned}
\begin{tikzpicture}[thick, scale=0.8]
\draw [region] (0,1.4) rectangle +(2.2,3);
\draw (1.9,2.8) to [out=up, in=down] (1.2,3.5) to [out=up, in=down] (0.5,4.2) to(0.5,4.4);
\draw (1.9, 1.4)to(1.9,2.8);
\draw [region] (0.2, 1.2) rectangle +(2.2,3);
\draw (0.5,2.8) to [out=up, in=down] (0.5,3.5) to [out=up, in=down] (1.2,4.2) ;
\draw (1.2,2.8)  to [out=up, in=down] (1.9,3.5) to [out=up, in=down] (1.9,4.2);
\draw (0.5,1.2) to(0.5,1.4) to [out=up, in=down] (0.85,2.1) to [out=up, in=down](0.5, 2.8);
\draw (1.2,1.2) to(1.2,1.4) to [out=up, in=down] (0.85,2.1) to [out=up, in=down](1.2, 2.8);
\node [morphism] at (0.85,2.1) {$A$};
\end{tikzpicture}
\end{aligned}
\\[-9pt]
\downarrow\scriptstyle\mu
&
\begin{aligned}
\tikz{\node [rotate=45] at (0,0) {$\Rightarrow$} node [below right=0pt] {$\III$};}
\end{aligned}
&
\downarrow\scriptstyle\mu
\\[-5pt]
\begin{aligned}
\begin{tikzpicture}[thick, scale=0.8]
\draw [region] (0.7,2.1) rectangle +(2.2,3);
\draw (2.6,2.1) to [out=up, in=down] (1.9,2.8) to [out=up, in=down] (1.2,3.5);
\draw (1.2,3.5) to (1.2,5.1);
\draw [region] (0.9, 1.9) rectangle +(2.2,3);
\draw (1.2,1.9) to (1.2,2.8) to [out=up, in=down] (1.9,3.5);
\draw (1.9,1.9) to (1.9,2.1)  to [out=up, in=down] (2.6,2.8) to (2.6,3.5);
\draw (1.9,3.5) to [out=up, in=down] (2.25,4.2) to [out=up, in=down](1.9, 4.9) to (1.9,4.9);
\draw (2.6,3.5) to [out=up, in=down] (2.25,4.2) to [out=up, in=down](2.6, 4.9)to (2.6,4.9);
\node [morphism] at (2.25,4.2) {$B$};
\end{tikzpicture}
\end{aligned}
&\stackrel{\II}{\rightarrow}&
\begin{aligned}
\begin{tikzpicture}[thick, scale=0.8]
\draw [region] (0,1.4) rectangle +(2.2,3);
\draw (1.9,2.8) to [out=up, in=down] (1.2,3.5) to [out=up, in=down] (0.5,4.2) to(0.5,4.4);
\draw (1.9, 1.4)to(1.9,2.8);
\draw [region] (0.2, 1.2) rectangle +(2.2,3);
\draw (0.5,2.8) to [out=up, in=down] (0.5,3.5) to [out=up, in=down] (1.2,4.2) ;
\draw (1.2,2.8)  to [out=up, in=down] (1.9,3.5) to [out=up, in=down] (1.9,4.2);
\draw (0.5,1.2) to(0.5,1.4) to [out=up, in=down] (0.85,2.1) to [out=up, in=down](0.5, 2.8);
\draw (1.2,1.2) to(1.2,1.4) to [out=up, in=down] (0.85,2.1) to [out=up, in=down](1.2, 2.8);
\node [morphism] at (0.85,2.1) {$B$};
\end{tikzpicture}
\end{aligned}
\end{array}
\\[80pt]
\tikzset{typeVdiagram/.style={scale=0.9, yscale=0.9}}
\begin{array}{@{}c@{\hspace{3pt}}c@{\hspace{2pt}}c@{\hspace{3pt}}c@{\hspace{2pt}}c@{}}
\begin{aligned}
\begin{tikzpicture}[thick, scale=0.8, typeVdiagram]
\draw [region] (-0.2,-2) rectangle +(2.4,4);
\draw (0.6,-2) to (0.6,-1) to (0.6,0) to [out=up, in=down] (1.6,1.8) to (1.6,2);
\node [morphism, anchor=north] at (0.6,-1) {$\alpha$};
\begin{scope}[]
\draw [region] (0,-2.2) rectangle +(2.4,4);
\end{scope}
\draw (1.6,-2.2) to (1.6,-1.2) to (1.6,0) to [out=up, in=down] (0.6,1.8);
\node [morphism, anchor=north] at (1.6,0) {$\beta$};
\end{tikzpicture}
\end{aligned}
&
\stackrel{\II}{\rightarrow}
&
\begin{aligned}
\begin{tikzpicture}[thick, scale=0.8, typeVdiagram]
\draw [region] (-0.2,-1) rectangle +(2.4,4);
\draw (0.6,-1) to (0.6,-1) to (0.6,0) to [out=up, in=down] (1.6,1.8) to (1.6,3);
\node [morphism, anchor=north] at (0.6,0) {$\alpha$};
\begin{scope}[]
\draw [region] (0,-1.2) rectangle +(2.4,4);
\end{scope}
\draw (1.6,-1.2) to (1.6,0) to [out=up, in=down] (0.6,1.8) to (0.6,2.8);
\node [morphism, anchor=south] at (0.6,1.8) {$\beta$};
\end{tikzpicture}
\end{aligned}
&
\stackrel{\II}{\rightarrow}
&
\begin{aligned}
\begin{tikzpicture}[thick, scale=0.8, typeVdiagram]
\draw [region] (-0.2,-1) rectangle +(2.4,4);
\draw (0.6,-1) to [out=up, in=down] (1.6,0.8) to (1.6,3);
\node [morphism, anchor=south] at (1.6,0.8) {$\alpha$};
\begin{scope}[]
\draw [region] (0,-1.2) rectangle +(2.4,4);
\end{scope}
\draw (1.6,-1.2) to (1.6,-1) to [out=up, in=down] (0.6,0.8) to (0.6,2.8);
\node [morphism, anchor=south] at (0.6,1.8) {$\beta$};
\end{tikzpicture}
\end{aligned}
\\[-6pt]
\downarrow \scriptstyle{\I}
&&
\begin{aligned}
\begin{tikzpicture}[thick, scale=1, baseline=10pt]
\node [rotate=45, scale = 1.2] at (0.3,0.8) {$\Rightarrow $};
\node [rotate=0, scale = 1] at (0.7,0.8) {$\V$};
\end{tikzpicture}
\end{aligned}
&&
\downarrow \scriptstyle{\I}
\\[-2pt]
\begin{aligned}
\begin{tikzpicture}[thick, scale=0.8, typeVdiagram]
\draw [region] (-0.2,-2) rectangle +(2.4,4);
\draw (0.6,-2) to (0.6,-1) to (0.6,0) to [out=up, in=down] (1.6,1.8) to (1.6,2);
\node [morphism, anchor=north] at (0.6,0) {$\alpha$};
\begin{scope}[]
\draw [region] (0,-2.2) rectangle +(2.4,4);
\end{scope}
\draw (1.6,-2.2) to (1.6,-1.2) to (1.6,0) to [out=up, in=down] (0.6,1.8);
\node [morphism, anchor=north] at (1.6,-1) {$\beta$};
\end{tikzpicture}
\end{aligned}
&
\stackrel{\II}{\rightarrow}
&
\begin{aligned}
\begin{tikzpicture}[thick, scale=0.8, typeVdiagram]
\draw [region] (-0.2,-2) rectangle +(2.4,4);
\draw (0.6,-2) to (0.6,-1) to [out=up, in=down] (1.6,0.8) to (1.6,2);
\node [morphism, anchor=south] at (1.6,0.8) {$\alpha$};
\begin{scope}
\draw [region] (0,-2.2) rectangle +(2.4,4);
\end{scope}
\draw (1.6,-2.2) to (1.6,-1) to [out=up, in=down] (0.6,0.8) to (0.6,1.8);
\node [morphism, anchor=north] at (1.6,-1) {$\beta$};
\end{tikzpicture}
\end{aligned}
&
\stackrel{\II}{\rightarrow}
&
\begin{aligned}
\begin{tikzpicture}[thick, scale=0.8, typeVdiagram]
\draw [region] (-0.2,-1) rectangle +(2.4,4);
\draw (0.6,-1) to [out=up, in=down] (1.6,0.8) to (1.6,3);
\node [morphism, anchor=south] at (1.6,1.8) {$\alpha$};
\begin{scope}[]
\draw [region] (0,-1.2) rectangle +(2.4,4);
\end{scope}
\draw (1.6,-1.2)to (1.6,-1) to [out=up, in=down] (0.6,0.8) to (0.6,2.8);
\node [morphism, anchor=south] at (0.6,0.8) {$\beta$};
\end{tikzpicture}
\end{aligned}
\end{array}
\end{array}
\]
\[
\hspace{1cm}
\]

\vspace{-15pt}
\caption{Simple instances of the five types of homotopy generator\label{fig:homotopy}}
\end{figure}
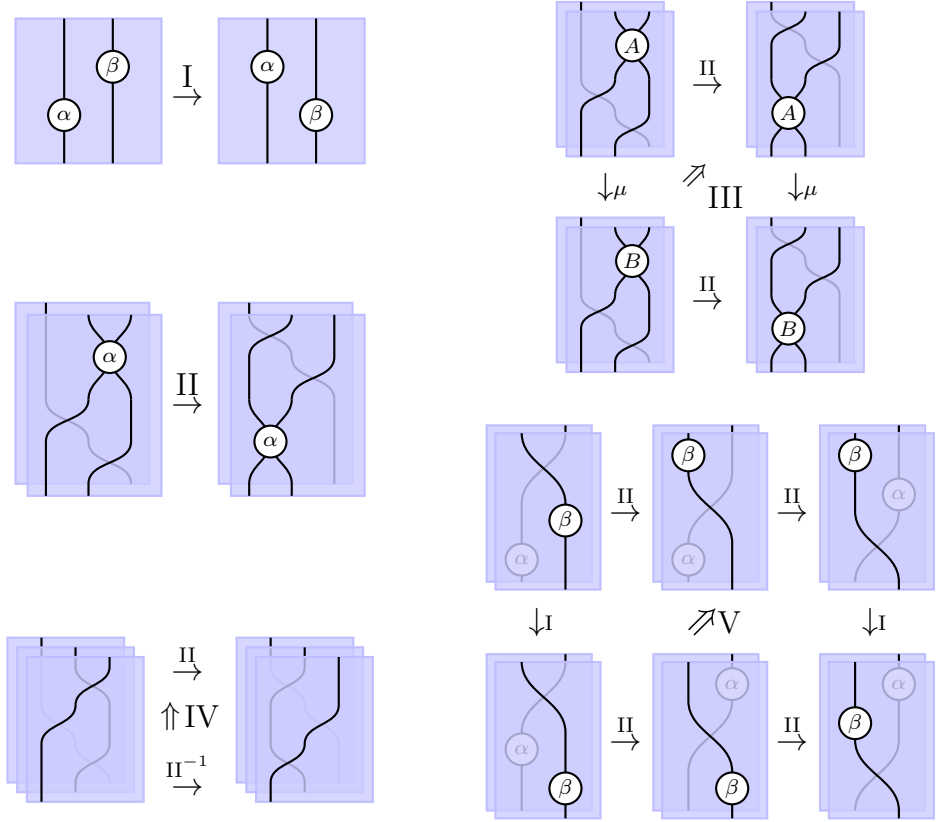

Diagram and signature structures give efficient tools for handling higher-dimensional composition, but lack a notion of \textit{homotopy} which is fundamental to  higher category theory. \JVcomm{Say something about sesquicategories?}\ignore{In this sense they resemble higher versions of classical sesquicategories~\cite{Sesquicategories}, which are 2\-categories without the interchange law. The following definition records this.\footnote{There is an abuse of terminology here, since `sesquicategory' really means `$1 \frac 1 2$ category'.}
\begin{definition}
An \textit{$n$\-sesquicategory} is an $(n{+}1)$\-signature.
\end{definition}

\noindent
Under our scheme this is fully precise for all $n$, including $n = \infty$.

}In Section~\ref{chapterinterchangers} we supply the needed homotopical data in terms of \textit{homotopy generators} of certain types, previewed in simple form in Figure~\ref{fig:homotopy} in terms of our graphical calculus. This must be done separately in each dimension. We use these to give new definitions of semistrict 2\-, 3\- and 4\-category as follows.

\begin{definition}
\label{defsemistricttwocategory}
A \emph{semistrict 2\-category} is a $3$\-signature supporting homotopy generators of type $\I$.
\end{definition}

\begin{definition}
\label{defsemistrictthreecategory}
A \emph{semistrict 3\-category} is a $4$\-signature supporting homotopy generators of types $\I$ and $\II$.
\end{definition}

\begin{definition}
\label{defsemistrictfourcategory}
A \emph{semistrict 4\-category} is a $5$\-signature supporting homotopy generators of types $\I$, $\II$, $\III$, $\IV$ and $\V$.
\end{definition}

\noindent
These are considerably simpler than traditional algebraic definitions. This arises in part from the nature of signature and diagram structures, as discussed above. However, our treatment of homotopy generators yields two further simplifications.

Firstly, we are able to neglect equations  governing redundant encodings of the same homotopy. For example, in Gray's definition of semistrict 3\-category, the following are equivalent descriptions of the same 3\-cell, based on different implicit descriptions of the source diagram:
\begin{calign}
\nonumber
\begin{aligned}
\begin{tikzpicture}[thick, scale=0.7]
\draw [region] (0,0) rectangle +(4,3);
\draw (1,0) to +(0, 3);
\draw (2,0) to +(0, 3);
\draw (3,0) to +(0, 3);
\node [morphism] at (1,1) {$\alpha$};
\node [morphism] at (3,2) {$\beta$};
\draw [blue, dashed] (0.5,0.5) rectangle +(1.8,1);
\draw [blue, dashed] (2.5,1.5) rectangle +(1,1);
\end{tikzpicture}
\end{aligned}
\stackrel {\Phi_1} \to
\begin{aligned}
\begin{tikzpicture}[thick, scale=0.7]
\draw [region] (0,0) rectangle +(4,3);
\draw (1,0) to +(0, 3);
\draw (2,0) to +(0, 3);
\draw (3,0) to +(0, 3);
\node [morphism] at (1,2) {$\alpha$};
\node [morphism] at (3,1) {$\beta$};
\draw [blue, dashed] (0.5,1.5) rectangle +(1.8,1);
\draw [blue, dashed] (2.5,0.5) rectangle +(1,1);
\end{tikzpicture}
\end{aligned}
&
\begin{aligned}
\begin{tikzpicture}[thick, scale=0.7]
\draw [region] (0,0) rectangle +(4,3);
\draw (1,0) to +(0, 3);
\draw (2,0) to +(0, 3);
\draw (3,0) to +(0, 3);
\node [morphism] at (1,1) {$\alpha$};
\node [morphism] at (3,2) {$\beta$};
\draw [blue, dashed] (0.5,0.5) rectangle +(1,1);
\draw [blue, dashed] (1.7,1.5) rectangle +(1.8,1);
\end{tikzpicture}
\end{aligned}
\stackrel {\Phi_2} \to
\begin{aligned}
\begin{tikzpicture}[thick, scale=0.7]
\draw [region] (0,0) rectangle +(4,3);
\draw (1,0) to +(0, 3);
\draw (2,0) to +(0, 3);
\draw (3,0) to +(0, 3);
\node [morphism] at (1,2) {$\alpha$};
\node [morphism] at (3,1) {$\beta$};
\draw [blue, dashed] (0.5,1.5) rectangle +(1,1);
\draw [blue, dashed] (1.7,0.5) rectangle +(1.8,1);
\end{tikzpicture}
\end{aligned}
\end{calign}
In a Gray category there is an axiom that says $\Phi_1 = \Phi_2$. In our approach, this homotopy arises in a unique way, with no redundancy in its description, and therefore with no need for additional equations to control this redundancy.

Secondly, our approach is able to treat simultaneously homotopy generators of the same kind appearing in different dimensions, while traditional algebraic definitions must treat them separately. For example, in a semistrict 3\-category, there are interchangers of 2\-cells that yield 3\-cells such as  $\Phi_1$, $\Phi_2$ shown above, as well as interchangers of 3\-cells that yield equations, specified separately. We treat these uniformly in Definition~\ref{defsemistrictthreecategory}, as they all arise as instances of the type I homotopy generator.

\paragraph{Quasistrict $n$\-categories.} A definition of $n$\-category is semistrict when it is `as strict as possible', while still being able to model arbitrary homotopy $n$\-types. When studying $n$\-categories per se, it can often be attractive to have a semistrict definition available, to get a sense of the minimal algebraic complexity of the theory. However, if our desire is to prove theorems internal to our $n$\-categories, then there is a more important concern: shorter proofs. A simpler definition of $n$\-category can lead to simpler proofs to some degree, but if taken too far it can have the opposite effect: the language can become so meagre that while everything remains possible in principle, some computations which are conceptually simple can become long-winded.

We illustrate these ideas with the following heuristic diagram of the weak-strict spectrum:
\begin{equation}
\nonumber
\begin{tikzpicture}
\tikzset{point/.style={fill=black, inner sep=2pt, circle}};
\draw [thick] (0,0) node (1) [point] {} to
    node (2) [point, pos=0.25] {}
    node (3) [point, pos=0.75] {}
    (8,0)
    node (4) [point] {};
\node [above] at (1.north) {Weak};
\node [above] at (2.north) {\smash{?}};
\node [above] at (3.north) {Semistrict};
\node [above] at (4.north) {Strict};
\draw [decorate,decoration={brace,amplitude=8pt}, thick] (5.75,-0.25) -- +(-3 .5,0);
\node [below] at (4,-0.5) {Quasistrict};
\end{tikzpicture}
\end{equation}
Here, `weak' marks the weakest possible definition of $n$\-categories
; `strict' marks strict $n$\-categories. We now define the new term `quasistrict'.
\begin{definition}
A definition of $n$\-categories is \emph{quasistrict} if it is strictly associative and unital, and able to model all homotopy $n$\-types.
\end{definition}

\noindent
Quasistrict $n$\-categories therefore occupy a region in the centre of the weak-strict spectrum. Semistrict $n$\-categories are the strictest quasistrict $n$-categories\footnote{Here we rule out definitions such as Simpson's snucategories~\cite{}, which are semistrict yet have weak units. Of course, there may be other valid approaches to higher semistrict categories which generalize this notion, and rule out Gray categories.}. We propose that the definition of $n$\-categories most amenable to computation will be the \textit{weakest} quasistrict $n$\-categories; this is marked `?' on the diagram above, since we do not know how to give such a definition.

Here we give a definition of quasistrict 4\-category which is somewhere in the interior of the quasistrict region of the weak-strict spectrum.
\begin{definition}
\label{defquasistrict4cat}
A \emph{quasistrict 4-category} is a $5$-signature supporting homotopy generators of types $\I$, $\II'$, $\III'$, $\IV'$, $\V'$ and $\VI$.
\end{definition}

\noindent
We slightly abuse terminology here; this is merely \emph a definition of quasistrict 4\-category, not \emph{the} definition, since `quasistrict' is a broad notion. These new classes are stronger than those used in the definition of semistrict 4\-category given above, allowing some deductions which would be long-winded in the semistrict case to be given in 1 step in the quasistrict case. \ignore{For example, consider the following sequence:
\[
\text{[PULL THROUGH THE WRONG WAY]}
\]
All we are achieving here is naturality of the inverse braiding, but it takes several steps, and indeed the sequence is non-obvious. It would be laborious to go to such lengths every time one wanted to implement naturality of the inverse braiding in a proof, and would impair the practicality of a proof assistant. The solution is to add this to our definition of $n$\-category, and indeed this move is available atomically in the type $\II'$ homotopy generator family used in Definition~\ref{defquasistrict4cat}. This is a simple instance of the general idea described above: by weakening our definition of $n$\-category, we make it more amenable to calculation. Adding these requires the addition of a further class of homotopy generators, which we call type VI, which controls the interaction between the different type $\II'$ homotopy generators.}

\paragraph{Computer formalization.}

Since quasistrict 4\-categories are designed to aid proof construction, it makes sense to formalize them using in a proof assistant. Based on the ideas of this paper, the authors have developed  (in collaboration with Kissinger) \emph{Globular}~\cite{fscd}\footnote{This citation is to a tool paper which focuses on the user interface and the practical implementation, giving only an outline sketch of the theoretical foundations. The present work can be considered as the corresponding theory paper.}, a new proof assistant which allows the user to construct formal proofs in finitely-presented quasistrict 4\-categories\footnote{The restriction to finite presentability is because the presentation must be entered into the computer; our theoretical approach is compatible with infinitely-generated higher categories.}, in the sense of Definition~\ref{defquasistrict4cat}. The interface is fully graphical, allowing construction and interaction with proof objects by clicking-and-dragging, and the tool is available online through a web page, minimizing barriers to use and allowing formal proofs to be hyperlinked directly from research papers. The proof assistant has been well-received by the community, being accessed 8,755 times by 1,978 distinct users in the first 11 months since launch in December 2015.

This proof assistant has been used to develop a new result, which is the final main result presented in this article: in a quasistrict 4\-category, an adjunction of 1\-morphisms can be promoted to give a coherent adjunction satisfying the butterfly equations. This proof is available online at \href{http://globular.science/1605.002}{globular.science/1605.002} for direct inspection, and presented in detail in Section~\ref{chapterbutterfly}.

By general results of Riehl and Verity~\cite{Riehl_2016}, it is expected that such a theorem ought to hold in any reasonable algebraic definition of 4\-category. Our proof is the first such that has been given; indeed, we believe it to be the first nontrivial proof in the literature of any sort internal to an algebraic 4\-category.\footnote{We distinguish here between algebraic definitions of 4-category (such as those due to Trimble and Crans, and that described here), and homotopical ones (such as quasicategories) which build in the homotopy hypothesis.} That we were able to construct this proof is the strongest evidence we can supply for the correctness and utility of our framework, and for the definition of quasistrict 4\-category that we have built within it. This is just one of many substantial results which have been formalized using our new techniques; a longer list is given at \href{https://ncatlab.org/nlab/show/Globular}{ncatlab.org/nlab/show/Globular}.

\paragraph{Criticisms.} We criticise our work in the following ways. Firstly, we do not prove that our definitions of semistrict 2\-, 3\- and 4\-categories are equivalent to the standard weak notions~\cite{coherencefortricategories, Tetracategories}. In the case of 2\- and 3\-categories we are very confident that this is the case, and in the case of 4\-categories we are reasonably confident; as evidence, we show in Section~\ref{sectionswitchthreecategory} that a Gray category gives rise to a semistrict 3\-category in our sense, and in Section~\ref{fourtasaxioms} that a 4\-tas in the sense of Crans~\cite{Crans_1998} gives rise to a semistrict 4\-category (modulo some issues that we identify with Crans' definition.)  It is worth noting that for traditional semistrict 3\-categories, it took 20 years for such an equivalence theorem to proved, between Gray's work~\cite{Gray_1974, Gray_1976} around 1974 and the coherence theorem of Gordon, Power and Street~\cite{coherencefortricategories} in 1995; furthermore, the definition of semistrict 4\-category due to Crans has never been shown equivalent to the tetracategories of Trimble~\cite{Tetracategories, Hoffnung_2011}, or to any other definition of 4\-category.

Secondly, while we hold that our methods are simpler than existing approaches, we do not propose a definition of semistrict 5\-category. The basic foundation of signature and diagram structures developed here should continue to be usable in any finite dimension, and indeed the many properties proved in Section~\ref{sec:core} are dimension-agnostic. We expect that semistrict 5\-categories would require a total of 9 homotopy generator families\footnote{We conjecture that the homotopy generators required for semistrict $n$-categories correspond to Batanin trees~\cite{Batanin_1998} with $n{+}2$ edges, with no leaves at depth 1 from their nearest branching. By this method we predict 4 more types of homotopy generators would be required for semistrict 5-categories, in addition to those used for semistrict 4-categories.}; identifying these would require careful manual analysis to write down. It would be better to find a more systematic approach to specifying the homotopy generators automatically, preferably in a way which moves closer to the weak end of the weak-strict spectrum as discussed above.

\paragraph{Acknowledgements.} We are extremely grateful to Aleks Kissinger for many detailed discussions, especially in the early stages of this work. We also thank John Baez, Manuel B\"arenz, Bruce Bartlett, Eugenia Cheng, Chris Douglas, Eric Finster, Nick Gurski, Andr\'e Henriques, Samuel Mimram, Christopher Schommer-Pries, Dominic Verdon and Steve Vickers for useful conversations, and the audiences of HDRA 2015, HDRA 2016 and CT 2016 for feedback on early versions of our results.

\subsection{Related work}
\label{sec:relatedwork}

\paragraph{Homotopy type theory.} This broad research programme seeks to develop a new foundation for mathematics, in a setting based on a homotopical interpretation of Martin-L\"of's dependent type theory~\cite{hottbook}. It is amenable to computer formalization, and there has been substantial activity in producing formal proofs~\cite{HoTTCoq, Licata_2013, UniMath}. While the motivations are related, there is little overlap between proofs suitable for formalization in homotopy type theory, and those suitable for formalization in our system. Homotopy type theory is far more expressive, with a rich syntax of term constructors that includes, but goes far beyond, the basic higher composition operations considered here. At the same time, the direct access to homotopy generators that our system provides allows direct construction of sophisticated homotopies---such as the complex proof illustrated in Section~\ref{chapterbutterfly}, or many of the examples listed at \href{https://ncatlab.org/nlab/show/Globular}{ncatlab.org/nlab/show/Globular}---which would be hard to construct directly in a formal system for homotopy type theory. Furthermore, our system does not require invertibility of higher cells, while the natural categorical interpretation of homotopy type theory is in an $(\infty,1)$\-setting where $k$\-cells for $k > 1$ are necessarily invertible. It remains unclear what the true relation is between the approaches, and whether there is a meaningful single approach to formal higher category theory that combines the strengths of both.

A good arena for comparing the approaches more closely may be reasoning about the higher inductive types arising as presentations of CW complexes, which are basic objects of study in homotopy type theory, and which give the data of a signature exactly of the form that we study here. For example, the \textit{Globular} tool could be used to construct elements of higher identity types of such spaces, which could then be exported back to a homotopy type theory setting.

\paragraph{Rewriting theory.} There is a large body of work on polygraphs and their applications to higher-dimensional rewriting~\cite{mimram-talk, Mimram_2014, Guiraud_2012, Lafont_1997}. The setting is strict $n$\-categories, and polygraphs directly correspond to our signature structures. Our work owes a lot to the perspective developed by this community, which emphasizes the idea that a composite $(n{+}1)$\-cell is precisely the same as a rewrite sequence of $n$\-cells. Our research programme can be considered an attempt to apply this perspective in the more general semistrict setting.

Central questions of study in the polygraph community are confluence and normalization, which we entirely neglect here. The basic techniques are also different: while the polygraph community relies on a pushout formalism for constructing rewrites~\cite{Mimram_2009}, our approach is more combinatorial. Furthermore, our diagram structures would not apply in the strict setting, since they contain layout data which breaks the strict notion of equivalence of diagrams appropriate to polygraphs.

\paragraph{Enriched approaches.} It is natural to propose defining semistrict $(n{+}1)$\-categories recursively, as categories enriched in the category of semistrict $n$\-categories. This heuristic idea is successful in leading to a definition of semistrict 3\-category, but fails to give a definition of semistrict 4\-category because the relevant tensor product of Gray categories fails to be monoidal closed, an argument sketched by Crans~\cite{Crans_1999} and refined by Bourke and Gurski~\cite{Bourke_2015}. While this issue remains unresolved, work by Batanin, Cisinski, Garner and Weber~\cite{Weber_2013, Batanin_2013, Weber_2013_II, Garner_2010, Garner_2008} suggests an operadic perspective allowing recursive enriched definitions of higher categories with strict units. However, this approach has not yet led to a concrete definition of semistrict 4\-category, and it is unclear whether our definition of semistrict 4\-category could arise in principle from their approach.

\paragraph{Opetopic higher categories.} Baez and Dolan's opetopic theory of higher categories~\cite{Baez_1998}, given a combinatorial interpretation in terms of trees by Kock and collaborators~\cite{Kock_2010}, has been developed by Finster~\cite{Finster_Talk} as the foundation for a higher-dimensional type theory, and implemented as the proof assistant \emph{Opetopic!}~\cite{opetopic} allowing description of formal proofs in opetopic $(\infty,\infty)$\-categories. The overall structure is at once elegant and powerful, with the the tree interpretation giving a graphical calculus, albeit one which is quite different from the generic-position diagrams used here. Opetopic categories have a more restricted form pasting diagram than globular categories, however, and the approaches are not in general expected to be equivalent. Furthermore, we are not aware of a substantial formalized proof in the opetopic type theory setting, comparable to the proof we present in Section~\ref{chapterbutterfly}.

\section{Data structures}

Here we introduce the basic data structures that underlie our approach: signatures, diagrams, embeddings, slices, rewrites and lifts. These structures are mutually recursive, so we will necessarily have to refer to some of these before they are formally defined. We therefore sketch informal definitions of each structure here.
\begin{itemize}
\item an \textit{$n$-signature} $\Sigma$ is the data of a presentation in dimension $n$, specifying each of the available generating cells $g$ in dimension $k \leq n$, each having (for $k>0$) source and target $(k{-}1)$\-diagrams $g.s$ and $g.t$ respectively;
\item an \textit{$n$-diagram} $D$ over a given $n$-signature $\Sigma$ is a composite of generating cells, with (for $n>0$) source and target $(n{-}1)$\-diagrams $D.s$ and  $D.t$;
\item an \textit{embedding} is the data witnessing one diagram as a subdiagram of another;
\item a \textit{slice} is an intermediate diagram of codimension 1 within a given diagram;
\item a \textit{rewrite} is the diagram $D'$ obtained by removing a subdiagram $S_1$ of a given diagram $D$, and replacing it with another diagram $S_2$;
\item a \emph{lift} is the canonical embedding of $S_2$ within $D'$ induced by a rewrite.
\end{itemize}
The key structure is that of signature: in the introduction we define semistrict and quasistrict $n$\-categories for $n \leq 4$ as $(n{+}1)$\-signatures with particular properties. Signatures are similar to polygraphs as used in the rewriting community\JVcomm{Citation needed}; however, our notion of diagram is unrelated, and can be considered the key innovation.

As well as the basic definitions, we give conditions for instances of these structures be well-defined. We state propositions asserting that slices, rewrites and lifts of well-defined structures are again well-defined. We also define composition and whiskering of diagrams, and state a number of propositions asserting that composition is strictly associative and unital, and that whiskering is distributive; these are among the familiar properties of semistrict $n$-categories. Altogether, these propositions demonstrate that the structures we propose can serve as a foundation for modelling cell composition in higher categories. Proofs are left until Section~\ref{sec:core}, due to their significant length.

\subsection{Signatures and diagrams}

We first give the definitions of signature and diagram.
\begin{definition}[Signature]
\label{defsignature}
For $n \geq 0$, an \emph{$n$-signature} $\Sigma$ is a family of sets $(\Sigma_0, ..., \Sigma_n)$, such that when $n>0$:
\begin{itemize}
\item $(\Sigma_0, \ldots, \Sigma_{n{-}1})$ is an $(n{-}1)$-signature; 
\item each $g \in \Sigma_n$ is equipped with $(n{-}1)$-diagrams $g.s$, $g.t$, such that when $n>1$,  $g.s.s = g.s.t$ and $g.t.s = g.t.t$.
\end{itemize}
\end{definition}

\begin{definition}[Diagram]
\label{defdiagram}
For $n \geq 0$, an \emph{$n$\-diagram} $D$ over an $n$\-signature $\Sigma$ is a list of length $|D|$ such that each element $D[i]$ with $0 < i < |D|$ consists of the following:
\begin{itemize}
\item a generating $n$-cell  $D[i].g \in \Sigma_n$;
\item if $n>0$, an embedding $D[i].e: D[i].g.s \hookrightarrow D[i].d$;
\item if $n>0$, a source $(n{-}1)$-diagram $D.s$ over $\Sigma' := (\Sigma_0, \ldots, \Sigma_{n-1})$.
\end{itemize}
Furthermore, if $n=0$ then we must have $|D|=1$. 
\end{definition}

We now consider what it means for a diagram to be well-defined. This makes use of the notion of \textit{slice}, defined formally later in this section. Informally, the slices of an $n$\-diagram $D$ for $n>0$ are $(n{-}1)$\-diagrams $D[i].d$ appearing at intermediate heights $0 \leq i \leq |D|$ within the diagram. The zeroth slice $D[0].d$ equals the source $D.s$.

\begin{definition}
\label{defwellformed}
An $n$\-diagram $D$ over an $n$\-signature $\Sigma$ is \emph{well-defined} when $n=0$, or when $n>0$ and the source diagram $D.s$ is well-defined, and for every $0 < i \leq |D|$ the slice $D[i].d$ exists and is well-defined.
\end{definition}

\subsection{Embeddings}

For $n$\-diagrams $S,T$, on an intuitive level, an embedding $e:S\hookrightarrow D$ specifies a way in which $S$ appears as a subdiagram of $D$. We illustrate this with the following example of a pair of 2\-diagrams, with two embeddings $e_1,e_2: S \hookrightarrow D$:
\begin{align*}
S\quad=\quad
&
\begin{aligned}
\begin{tikzpicture}[thick, scale = 0.7]
\draw (-0.5,0) to (-0.5,0.5) to [out = up, in = left] (0,1)
to [out = right, in = up](0.5,0.5) to(0.5,0);
\draw (-0.5,3) to (-0.5,2.5) to [out = down, in = left] (0,2)
to [out = right, in = down](0.5,2.5) to(0.5,3);
\draw (0,1) to (0, 2);
\node [blob] at (0,1) {};
\node [blob] at (0,2) {};
\end{tikzpicture}
\end{aligned}
\qquad
\begin{array}{c}
e_1\\
\hookrightarrow\\\\
e_2\\
\hookrightarrow
\end{array}
\qquad
\begin{aligned}
\begin{tikzpicture}[thick, scale = 0.7]
\draw (-0.5,0) to (-0.5,0.5) to [out = up, in = left] (0,1)
to [out = right, in = up](0.5,0.5) to(0.5,0);
\draw (-0.5,2.5) to [out = down, in = left] (0,2)
to [out = right, in = down](0.5,2.5) to(0.5,3);
\draw (0,1) to (0, 2);
\node [blob] at (0,1) {};
\node [blob] at (0,2) {};
\draw (0.5,3) to (0.5,3.5) to [out = up, in = left] (1,4)
to [out = right, in = up](1.5,3.5) to(1.5,3);
\draw (0.5,6) to (0.5,5.5) to [out = down, in = left] (1,5)
to [out = right, in = down](1.5,5.5) to(1.5,6);
\draw (1,4) to (1, 5);
\node [blob] at (1,4) {};
\node [blob] at (1,5) {};
\draw (1.5,0) to (1.5, 3);
\draw (-0.5,2.5) to [out = up, in = right](-1, 3) to [out = left, in = up](-1.5, 2.5) to (-1.5, 0);
\draw (-1,3) to (-1, 6);
\node [blob] at (-1,3) {};
\node [rectangle, minimum width = 30pt, minimum height = 40pt, draw, fill = none, color = blue, dashed] at (0,1.5) {};
\node [rectangle, minimum width = 30pt, minimum height = 40pt, draw, fill = none, color = blue, dashed] at (1,4.5) {};
\end{tikzpicture}
\end{aligned}
\quad=\quad D
\end{align*}
The formal definition goes as follows.
\begin{definition}[Embedding]
\label{defembedding}
Given two $n$\-diagrams $S$ and $D$, an \emph{embedding} $e:S\hookrightarrow D$ consists of no data when $n=0$, and of the following data when $n>0$:
\begin{itemize}
\item an embedding height $e.h \in \mathbb{N}$;
\item a source embedding $e.e: S.s\hookrightarrow D[e.h].d$.
\end{itemize}
\end{definition}

\noindent
Embeddings are therefore specified by sequences of natural numbers. The first embedding above has values $e_1.h=0$ and $e_1.e.h=1$, meaning that its image begins at height 0 and has 1 wire to the left. The second embedding $e_2$ has values $e_2.h=3$ and $e_2.e.h = 1$, meaning that its image begins at height 3 and has 1 wire to the left. Note that there are no vertices to the left or right of the images of $e_1$ and $e_2$, and so in particular the number of wires to the left and right of each image is determinate; this follows from the well-definedness property given below.

Just as for diagram structures, we give a formal statement of what it means for an embedding to be well-defined.
\begin{definition}
\label{defembwelldef}
An $n$-diagram embedding $e:S\hookrightarrow D$ between well-defined diagrams  is \emph{well-defined} if it satisfies the following properties: 
\begin{itemize}
\item if $n=0$, we have
\begin{equation}
\label{eqembwelldef1}
S[0].g=D[0].g;
\end{equation}
\item if $n>0$, the component embedding $e.e$ is well-defined;
\item if $n>0$, for all $0 \leq i < |S|$ we have the following:
\begin{align}
\label{eqembwelldef2}
S[i].g &= D[i+e.h].g
\\
\label{eqembwelldef3}
(\lift{e.e}{S[i].d}) \circ S[i].e &= D[i+e.h].e
\end{align}
\end{itemize}
\end{definition}

\noindent
The symbol $\Lambda$ represents the lift construction, described below. Intuitively, condition~\eqref{eqembwelldef2} says that the generators of $S$ are the same of those of $D$ in the image of the embedding, and condition~\eqref{eqembwelldef3} says that the embedding maps are the compatible.

\subsection{Rewriting}
\label{sectionrewriting}

The intuition behind the procedure of rewriting is that we want to transform a diagram $D$, by removing a part of the diagram $S$ embedded in $D$ by the embedding $e:S\hookrightarrow D$ and replace it with another diagram $T$, somehow preserving connectivity between corresponding elements. For that reason we need $S.s=T.s$, $S.t=T.t$ and we need to update the relevant embeddings $T[i].e$ corresponding to the generators $T[i].g$ being added to $D$.

With the auxiliary structures in place, we can define the first modification of a diagram: a \emph{rewrite}. Notice how for an $n$\-diagram such that $n>0$, the length $|D|$ of the list of generators and embeddings changes to $|D|-|S|+|T|$ as we remove $|S|$ elements of the source of the rewrite and replace them with $|T|$ elements of the target. For this reason the lists in the rewritten diagram consist of three segments: the initial and the final segment that remain unaltered and the middle slice which gets replaced. Depending on the value of $e.h$ and on $|S|$ any of the three segments in the rewritten diagram may be empty.
\begin{definition}[Rewrite]
\label{defrewrite}
Given an $n$-diagram $D$ with a subdiagram $e: S\hookrightarrow D$, and an $n$-diagram $T$ globular with respect to $S$, the \emph{rewrite} $\rewrite{D}{e}{S}{T}$ of $D$ is the following $n$-diagram:

\begin{itemize}

\item If $n=0$, then: 
\begin{align}
\label{defrewritezerosource}
|\rewrite{D}{e}{S}{T}| &= |D| - |S| + |T|
\\
\label{defrewritezerogenerator}
\rewrite{D}{e}{S}{T}[0].g &= T[0].g
\intertext
{
\item If $n>0$, then:
}
\label{defrewritesource}
\rewrite{D}{e}{S}{T}.s &= D.s 
\\
\label{defrewritesize}
|\rewrite{D}{e}{S}{T}| &= |D| - |S| + |T|
\\
\label{defrewritegenerator}
\rewrite{D}{e}{S}{T}[i].g &= 
\begin{cases}
D[i].g & \text{ if } 0\leq i \leq e.h \\
T[i - e.h].g & \text{ if } e.h < i < e.h + |T| \\
D[i - |T| + |S|].g & \text{ if }  e.h + |T| \leq i < |D| {-} |S| {+} |T| \\
\end{cases}
\\
\label{defrewriteembedding}
\rewrite{D}{e}{S}{T}[i].e
&=\begin{cases}
D[i].e & \text{if } 0\leq i \leq e.h \\
\lift{e.e}{T[i-e.h].d}
\\
\quad{}\circ T[i-e.h].e   & \text{if }e.h < i < e.h + |T| \\
D[i - |T| + |S|].e &  \text{if }e.h + |T| \leq i < |D| {-} |S| {+} |T| \\
\end{cases}
\end{align}
\end{itemize}

\end{definition}
This definition can be illustrated with an example. Consider the following diagrams $S, T, D$. 
\begin{align*}
S\quad=\quad
&
\begin{aligned}
\begin{tikzpicture}[thick, scale = 0.7]
\draw (-0.5,0) to (-0.5,0.5) to [out = up, in = left] (0,1)
to [out = right, in = up](0.5,0.5) to(0.5,0);
\draw (-0.5,3) to (-0.5,2.5) to [out = down, in = left] (0,2)
to [out = right, in = down](0.5,2.5) to(0.5,3);
\draw (0,1) to (0, 2);
\node [blob] at (0,1) {};
\node [blob] at (0,2) {};
\end{tikzpicture}
\end{aligned}
&
T\quad=\quad
&
\begin{aligned}
\begin{tikzpicture}[thick, scale = 0.7]
\draw (0,1) to [out = right, in=down] (0.5,1.5) to [out = up, in=left](1,2);
\draw (1,3) to (1,2);
\draw (0,1) to (0,0);
\draw (1,2) to [out = right, in=up](1.5, 1.5)to (1.5, 0);
\draw (0,1) to [out = left, in=down](-0.5, 1.5)to (-0.5, 3);
\node [blob] at (0,1) {};
\node [blob] at (1,2) {};
\end{tikzpicture}
\end{aligned}
\end{align*}
Note that there is an embedding $e: S\hookrightarrow D$ denoted by the blue dashed rectangle. Additionally $S, T$ are globular with respect to each other as their respective sources and targets match.
\begin{align*}
D\quad=\quad
&
\begin{aligned}
\begin{tikzpicture}[thick, scale = 0.7]
\draw (-0.5,0) to (-0.5,0.5) to [out = up, in = left] (0,1)
to [out = right, in = up](0.5,0.5) to(0.5,0);
\draw (-0.5,2.5) to [out = down, in = left] (0,2)
to [out = right, in = down](0.5,2.5) to(0.5,4);
\draw (0,1) to (0, 2);
\node [blob] at (0,1) {};
\node [blob] at (0,2) {};
\draw (-0.5,2.5) to [out = up, in = right] (-1,3)to [out = left, in = up] (-1.5,2.5) to (-1.5, 0);
\draw (-1,3) to (-1, 4);
\draw (-2.5,0) to (-2.5, 4);
\draw (1.5,0) to (1.5, 4);
\node [blob] at (-1,3) {};
\node [rectangle, minimum width = 25pt, minimum height = 40pt, draw, fill = none, color = blue, dashed] at (0,1.5) {};
\end{tikzpicture}
\end{aligned}
\qquad
\xrightarrow{}
\qquad
\begin{aligned}
\begin{tikzpicture}[thick, scale = 0.7]
\draw (0,1) to [out = right, in=down] (0.5,1.5) to [out = up, in=left](1,2);
\draw (1,4) to (1,2);
\draw (0,1) to (0,0);
\draw (1,2) to [out = right, in=up](1.5, 1.5)to (1.5, 0);
\draw (0,1) to [out = left, in=down](-0.5, 1.5)to (-0.5, 2.5);
\node [blob] at (0,1) {};
\node [blob] at (1,2) {};
\draw (-0.5,2.5) to [out = up, in = right] (-1,3)to [out = left, in = up] (-1.5,2.5) to (-1.5, 0);
\draw (-1,3) to (-1, 4);
\draw (-2.5,0) to (-2.5, 4);
\draw (2.5,0) to (2.5, 4);
\node [blob] at (-1,3) {};
\node [rectangle, minimum width = 45pt, minimum height = 40pt, draw, fill = none, color = blue, dashed] at (0.5,1.5) {};
\end{tikzpicture}
\end{aligned}
\quad=\quad \rewrite{D}{e}{S}{T}
\end{align*}
In the resulting diagram $\rewrite{D}{e}{S}{T}$, the generators of $S$ have been removed, and generators of $T$ inserted in their place. Their positions have been determined by the corresponding component embeddings. In this particular example the  final segment is empty, as $e.h+|S|=|D|$.

The definition of rewrite gives rise to the definition of slice.
\begin{definition}[Slice]
\label{defslice}
Given an $n$\-diagram $D$, for $0 < i \leq |D|$, its \emph{slice} $D[i].d$ is the following $(n{-}1)$-diagram:
\begin{itemize}
\item if $i=0$, $D[0].d := D.s$;
\item if $i>0$, $D[i].d := \rewrite{D[i-1].d}{D[i-1].e}{D[i].g.s} {D[i-1].g.t}$
\end{itemize}
\end{definition}

We can use the concept of slice to define  the target of a diagram, and the notion of globularity of a pair of diagrams. Recall that the source $D.s$ of a diagram $D$ is provided as part of the defining data.
\begin{definition}
\label{definitialterminalslice}
For $n>0$, given an $n$\-diagram $D$, its \textit{target} is the $(n{-}1)$-diagram $D.t := D[|D|].d$.
\end{definition}

\begin{definition}
Two $n$\-diagrams $S$, $T$ are \emph{globular} just when $n=0$, or when $n>0$ and $S.s=T.s$ and $S.t=T.t$.
\end{definition}

\subsection{Lifts}

Given a globular pair of diagrams $S,T$, and an embedding $e: S \hookrightarrow D$, we can form the rewritten diagram $\rewrite D e S T$, given intuitively by $D$ with the image of $e$ replaced by $T$. It follows that there is a canonical embedding $T \hookrightarrow \rewrite D e S T$; we call this the lift of $e$ with respect to $T$, and define it as follows.
\begin{definition}
\label{liftedembedding}
For $n \geq 0$, given globular $n$\-diagrams $S,T$ and an embedding $e:S \hookrightarrow D$, the \textit{lifted embedding} $\lift{e}{T}: T\hookrightarrow \rewrite{D}{e}{S}{T}$ is defined as follows:
\begin{itemize}
\item $\lift{e}{T}.h = e.h$
\item $\lift{e}{T}.e = e.e$
\end{itemize}
\end{definition}

\noindent
We illustrate this with the following example:
\begin{align*}
\begin{array}{ccccc}
S\quad=\quad&
\begin{aligned}
\begin{tikzpicture}[thick, scale = 0.7]
\draw (-0.5,0) to (-0.5,0.5) to [out = up, in = left] (0,1)
to [out = right, in = up](0.5,0.5) to(0.5,0);
\draw (-0.5,3) to (-0.5,2.5) to [out = down, in = left] (0,2)
to [out = right, in = down](0.5,2.5) to(0.5,3);
\draw (0,1) to (0, 2);
\node [blob] at (0,1) {};
\node [blob] at (0,2) {};
\end{tikzpicture}
\end{aligned}
&
\qquad\stackrel{e}{\hookrightarrow}\qquad
&
\begin{aligned}
\begin{tikzpicture}[thick, scale = 0.7]
\draw (-0.5,0) to (-0.5,0.5) to [out = up, in = left] (0,1)
to [out = right, in = up](0.5,0.5) to(0.5,0);
\draw (-0.5,2.5) to [out = down, in = left] (0,2)
to [out = right, in = down](0.5,2.5) to(0.5,4);
\draw (0,1) to (0, 2);
\node [blob] at (0,1) {};
\node [blob] at (0,2) {};
\draw (-0.5,2.5) to [out = up, in = right] (-1,3)to [out = left, in = up] (-1.5,2.5) to (-1.5, 0);
\draw (-1,3) to (-1, 4);
\draw (-2.5,0) to (-2.5, 4);
\draw (1.5,0) to (1.5, 4);
\node [blob] at (-1,3) {};
\node [rectangle, minimum width = 25pt, minimum height = 40pt, draw, fill = none, color = blue, dashed] at (0,1.5) {};
\end{tikzpicture}
\end{aligned}
&\quad=\quad D
\\\\
T\quad=\quad&
\begin{aligned}
\begin{tikzpicture}[thick, scale = 0.7]
\draw (0,1) to [out = right, in=down] (0.5,1.5) to [out = up, in=left](1,2);
\draw (1,3) to (1,2);
\draw (0,1) to (0,0);
\draw (1,2) to [out = right, in=up](1.5, 1.5)to (1.5, 0);
\draw (0,1) to [out = left, in=down](-0.5, 1.5)to (-0.5, 3);
\node [blob] at (0,1) {};
\node [blob] at (1,2) {};
\end{tikzpicture}
\end{aligned}
&
\qquad\stackrel{\lift{e}{T}}{\hookrightarrow}\qquad
&
\begin{aligned}
\begin{tikzpicture}[thick, scale = 0.7]
\draw (0,1) to [out = right, in=down] (0.5,1.5) to [out = up, in=left](1,2);
\draw (1,4) to (1,2);
\draw (0,1) to (0,0);
\draw (1,2) to [out = right, in=up](1.5, 1.5)to (1.5, 0);
\draw (0,1) to [out = left, in=down](-0.5, 1.5)to (-0.5, 2.5);
\node [blob] at (0,1) {};
\node [blob] at (1,2) {};
\draw (-0.5,2.5) to [out = up, in = right] (-1,3)to [out = left, in = up] (-1.5,2.5) to (-1.5, 0);
\draw (-1,3) to (-1, 4);
\draw (-2.5,0) to (-2.5, 4);
\draw (2.5,0) to (2.5, 4);
\node [blob] at (-1,3) {};
\node [rectangle, minimum width = 45pt, minimum height = 40pt, draw, fill = none, color = blue, dashed] at (0.5,1.5) {};
\end{tikzpicture}
\end{aligned}
&\quad=\quad \rewrite{D}{e}{S}{T}
\end{array}
\end{align*}

Lifted embeddings are used to construct the embeddings in Definition~\ref{defrewrite} if a rewrite. We do this by extracting the data of a component embedding $e.e$ and use it to give an embedding of the $i$th slice of $T$ into the $(i+e.h)$th slice of the rewrite of $D$. 

We also use the lifted embedding to define composition of two embeddings $e: S\hookrightarrow D$ and $f: D\hookrightarrow A$, given below. This is necessary since there is a mismatch between the source of $f.e$ and the target of $e.e$, and $\lift{f.e}{D[e.h].d}$ is needed to make the transition between them.

\begin{definition}[Composite embedding]
\label{defcomposition}
For $n \geq 0$, given $n$\-diagram embeddings $e: S\hookrightarrow D$ and $f: D\hookrightarrow A$, their \emph{composite embedding} $f\circ e: S\hookrightarrow A$ is defined as follows for $n>0$:
\begin{align*}
(f\circ e).h &:= e.h+f.h \\
(f\circ e).e &:= \lift{(f.e)}{D[e.h]}\circ (e.e)
\end{align*}
For $n=0$, it is defined to have no data.
\end{definition}

\noindent
This has a clear interpretation: that the notion of subdiagram is transitive. Consider the following example chain of embeddings:
\begin{align*}
S\quad:=\quad
&
\begin{aligned}
\begin{tikzpicture}[thick, scale = 0.7]
\draw [white] (0, 0) to (0, 3);
\draw (-0.5,0) to (-0.5,0.5) to [out = up, in = left] (0,1)
to [out = right, in = up](0.5,0.5) to(0.5,0);
\draw (0,1) to (0, 2);
\node [blob] at (0,1) {};
\end{tikzpicture}
\end{aligned}
\qquad
\stackrel{e}{\hookrightarrow}
\qquad
\begin{aligned}
\begin{tikzpicture}[thick, scale = 0.7]
\draw (-0.5,0) to (-0.5,0.5) to [out = up, in = left] (0,1)
to [out = right, in = up](0.5,0.5) to(0.5,0);
\draw (-0.5,3) to (-0.5,2.5) to [out = down, in = left] (0,2)
to [out = right, in = down](0.5,2.5) to(0.5,3);
\draw (0,1) to (0, 2);
\draw (-1.5,0) to (-1.5, 3);
\node [blob] at (0,1) {};
\node [blob] at (0,2) {};
\node [rectangle, minimum width = 25pt, minimum height = 20pt, draw, fill = none, color = orange, dashed] at (0,1) {};
\end{tikzpicture}
\end{aligned}
\qquad
\stackrel{f}{\hookrightarrow}
\qquad
\begin{aligned}
\begin{tikzpicture}[thick, scale = 0.7]
\draw (-0.5,0) to (-0.5,0.5) to [out = up, in = left] (0,1)
to [out = right, in = up](0.5,0.5) to(0.5,0);
\draw (-0.5,2.5) to [out = down, in = left] (0,2)
to [out = right, in = down](0.5,2.5) to(0.5,4);
\draw (0,1) to (0, 2);
\node [blob] at (0,1) {};
\node [blob] at (0,2) {};
\draw (-0.5,2.5) to [out = up, in = right] (-1,3)to [out = left, in = up] (-1.5,2.5) to (-1.5, 0);
\draw (-1,3) to (-1, 4);
\draw (-2.5,0) to (-2.5, 4);
\draw (1.5,0) to (1.5, 4);
\node [blob] at (-1,3) {};
\node [rectangle, minimum width = 55pt, minimum height = 45pt, draw, fill = none, color = blue, dashed] at (-0.5,1.5) {};
\node [rectangle, minimum width = 25pt, minimum height = 20pt, draw, fill = none, color = purple, dashed] at (0,1) {};
\end{tikzpicture}
\end{aligned}
\quad:=\quad A
\end{align*}
The diagram $S$ is directly embedded in $D$ by $e$, which is indicated by the orange rectangle. $D$ in turn is embedded in $A$ by $f$, indicated by the blue rectangle. These two can be combined together to obtain the composed embedding $f\circ e$ of $S$ inside $A$, indicated in the picture by a purple rectangle.

\subsection{Equivalence}

We now introduce simple formal notions of equivalence for diagrams and embeddings.\footnote{A possible alternative approach which we have not investigated in detail would be to define equivalence via embeddings; for example, defining diagrams $D$ and $S$ to be equivalent just when there exist well-defined embeddings $D \hookrightarrow S$ and $S \hookrightarrow D$.} These notions are mutually recursive.
\begin{definition}
\label{embeddingequiv}
Two $n$\-diagram embeddings $e:A\hookrightarrow B$ and $f: C\hookrightarrow D$  are \emph{equivalent}, written $e=f$, when $A=C$ and $B=D$, and furthermore if $n>0$ when $e.h=f.h$ and $e.e=f.e$.
\end{definition}
\begin{definition}
\label{diagequiv}
Two $n$-diagrams $D$ and $S$ are \emph{equivalent}, written $D = S$, when:
\begin{itemize}
\item $|S|=|D|$; 
\item for $0\leq i<|D|$ we have $S[i].g = D[i].g$;
\item if $n>0$, for $0\leq i<|D|$ we have $S[i].e = D[i].e$;
\item if $n>0$, then $D.s=S.s$
\end{itemize}
\end{definition}

\noindent
Note that for $n=0$ there is just a single cell $D[0].g$, and no source or embeddings to compare.


\subsection{Properties}
\label{sec:properties}

Diagrams and signatures have the following properties. These include many properties which traditionally appear as axioms in algebraic approaches to semistrict $n$\-categories; by obtaining them here as propositions we are able to reduce the axiomatic weight of our approach. We give only the statements here; the proofs are given in Section~\ref{sec:core}. 

\paragraph{Well-definedness.} When our basic procedures operate on well-defined structures, the result is again well-defined.

\begin{proposition}[Well-defined rewrites] 
For $n \geq 0$, given well-defined $n$\-diagrams $D, S, T$ with $S.s=T.s$ and $S.t=T.t$, and a well-defined embedding $e: S\hookrightarrow D$, the rewrite $\rewrite{D}{e}{S}{T}$ is well-defined.
\end{proposition}

\begin{proposition}[Well-defined lifts] 
For $n\geq 0$, given well-defined $n$\-diagrams $S, T, A$ with $S.s=T.s$ and $T.s=T.t$, and given a well-defined embedding $e:S\hookrightarrow A$, then the lifted embedding $\lift{e}{T}: T\hookrightarrow \rewrite{A}{e}{S}{T}$ is well-defined.
\end{proposition}

\begin{proposition}[Well-defined composition]  
For $n,m \geq 0$, given an $n$\-diagram $D$ and an $m$-diagram $S$, both well-defined, such that either $n \geq m$ and \mbox{$S.t=s^{n-m+1}(D)$}, or $m > n$ and $t^{m-n+1}(S)=D.s$, then $S\circ D$ is well-defined.
\end{proposition}

\begin{proposition}[Well-defined composite embeddings] 
For $n\geq 0$, given well-defined $n$\-diagram embeddings $e: S\hookrightarrow D$ and $f: D\hookrightarrow M$, then $f\circ e: S\hookrightarrow M$ is well-defined.
\end{proposition}

\paragraph{Further properties.} We show that the constructions we make satisfy straightforward additional properties.

\begin{proposition}[Identity rewrites]
\label{vacuousrewrite}
Given an $n$-diagram $D$ and $e: S\hookrightarrow D$, then $\rewrite{D}{e}{S}{S} = D$.
\end{proposition}

\begin{proposition}[Globularity on slices] 
For $n \geq 2$, given a well-defined $n$\-diagram $D$, we have $D.s.s = D[i].d.s$ and $D.s.t = D[i].d.t$ for any $0\leq i < |D|$.
\end{proposition}

\begin{proposition}[Explicit rewrites] 
For $n\geq 1$, given well-defined $n$\-diagrams $D, S, T$ such that $S.s=T.s$ and $S.t=T.t$, and a well-defined $n$\-diagram embedding $e: S\hookrightarrow A$, the following holds:
\begin{align*}
&\rewrite{A}{e}{S}{T}[j].d =\begin{cases}
A[j].d & \text{if } 0\leq j \leq e.h \\
\rewrite{A[e.h].d}{e.e}{S.s}{T[j-e.h].d}   & \text{if }e.h \leq j \leq e.h + |T| \\
A[j + |S| - |T|].d &  \text{if }e.h + |T| \leq j < |A| - |S| + |T| \\
\end{cases}
\end{align*}
\end{proposition}

\begin{proposition}[Composite lifts] 
For $n \geq 0$, given well-defined $n$\-diagrams $S, T, A, B, C$ with $S.s=T.s$, $S.t=T.t$, $A.s=C.s$, $A.t=C.t$, and given well-defined embeddings $e:S\hookrightarrow A$, $f:C\hookrightarrow B$, the following holds:
\begin{align*}
\lift{(\lift{f}{A}\circ e)}{T} = (\lift{f}{\rewrite{A}{e}{S}{T}})\circ\lift{e}{T}
\end{align*}
\end{proposition}

\begin{proposition}[Composite rewrites] 
For $n\geq 0$, given well-defined $n$\-diagrams  $S, T, A, B, C$ with $S.s=T.s$, $S.t=T.t$, $A.s=C.s$, $A.t=C.t$, and given well-defined embeddings $e:S\hookrightarrow A$, $f:C\hookrightarrow B$, the following holds for $0\leq j\leq e.h$:
\begin{align*}
\rewrite{B}{f}{C}{\rewrite{A}{e}{S}{T}}= \rewrite{\rewrite{B}{f}{C}{A}}{\lift{f}{A}\circ e}{S}{T}
\end{align*}
\end{proposition}

\begin{proposition}[Associative composite embeddings] 
For $n\geq 0$, given three well-defined $n$\-diagram embeddings $e: S\hookrightarrow D$, $f: D\hookrightarrow M$, $g: M\hookrightarrow N$, the following equality holds:
\[
g\circ(f\circ e) = (g\circ f) \circ e
\]
\end{proposition}

\begin{proposition}[Well-defined inclusions] 
For $k\geq 0$, given a well-defined $n$\-diagram $D$ and a well-defined $m$\-diagram $S$, then:
\begin{itemize}
\item if $n\geq m$ and $S.t=s^{n-m+1}(D)$ and the composite $S\circ D$ exists\JVcomm{Surely this existence condition is not needed?}, the inclusion embedding $\inc{S}{D}: D\hookrightarrow S\circ D$   is well-defined; 
\item if $n< m$ and $t^{m-n+1}(S)= D.s$ and the composite $S\circ D$ exists, the inclusion embedding $\increverse{S}{D}: D\hookrightarrow D\circ S$ is well-defined. 
\end{itemize}
\end{proposition}

\begin{proposition}[Well-behaved whiskering]  
For $n,m \geq 0$ with $n \neq m$, given a well-defined $n$\-diagram $D$ and well-defined $m$\-diagram $S$ such that the composite $S\circ D$ exists, then:
\begin{itemize}
\item if $n>m$ then $(S\circ D)[i].d = S\circ (D[i].d)$ for any $0\leq i < |D|$;
\item if $n<m$ then $(S\circ D)[i].d = (S[i].d) \circ D$ for any $0\leq i < |S|$.
\end{itemize}
\end{proposition}

\begin{proposition}[Interaction of lifts and inclusions] 
For $n,m \geq 0$ with $n \neq m$, given a well-defined $n$\-diagram $D$ and a well-defined diagram $m$\-diagram $S$, then:
\begin{itemize}
\item if $n>m$ then $\inc{S}{D[i].d}=\lift{(\inc{S}{D}.e)}{D[i].d}$ for any $0\leq i < |D|$;
\item if $n<m$ then $\increverse{S[i].d}{D}=\lift{(\increverse{S}{D}.e)}{S[i].d}$ for any $0 \leq i < |S|$.
\end{itemize}
\end{proposition}

\begin{proposition}[Associative diagram composition]
\JVcomm{Missing the reflected version?}%
Given two well-defined $n$\-diagrams $D, S$, and a well-defined $m$\-diagram $M$ such that $m>n$, the following holds:\[
S\circ(D\circ M) = (S\circ D) \circ M
\]
\end{proposition}

\begin{proposition}[Composition of inclusions] 
For $k, n\geq 0$, given two well-defined $n$\-diagrams $D, S$, and a well-defined $l$\-diagram $M$ such that $l>n$, the following holds:
\[
\inc{S}{D\circ M}\circ \inc{D}{M} = \inc{S\circ D}{M}
\]
\end{proposition}

\begin{proposition}[Distributive diagram composition] 
Given an $n$\-diagrams $D$, an $m$\-diagram $S$ and an $l$\-diagram $M$, all well-defined and such that $l,n>m>0$, the following holds:\JVcomm{Missing a well-typed condition}\JVcomm{Do we really need the $a,b$?}
\[
S\circ_{b}(D\circ_{a} M)  = (S\circ_{b} D)\circ_{a}(S\circ_{b} M)
\]
Here, we have $a=\mathrm{min}(n,l)-1$, $b=\mathrm{min}(m, \mathrm{max}(n,l)) - 1$.
\end{proposition}

\begin{proposition}[Triple inclusion property] 
Given an $n$\-diagrams $D$, an $m$\-diagram $S$ and an $l$\-diagram $M$, all well-defined and such that $l,n>m>0$, then provided that these composites are well-defined, the following holds:
\[
\inc{S}{D\circ M}\circ \inc{D}{M} = \inc{S\circ D}{S\circ M}\circ \inc{S}{M}
\]
\end{proposition}

\subsection{Graphical formalism}
\label{secgraphicalformalism}

Here we sketch an informal graphical calculus for diagram structures, following the ideas of Trimble~\cite{TrimbleDiagrams}. In particular, we introduce the new idea of $k$\-projected diagram, in which only the top $k$ dimensions of a diagram are depicted. However, we do not present solutions to the technical difficulties that Trimble encounters, and there is substantial work still required to formalize these ideas and prove correctness and completeness of the approach.

\paragraph{General procedure.}
We provide a method of translating a diagram structure into a graphical representation. While we state this in arbitrary dimension, it is only precise up to dimension 3.

In a $k$-projected graphical representation of an $n$\-diagram $D$, each $p$\-cell in $D$ for $k \leq p \leq n$ is represented by an $(n{-}p)$\-dimensional subspace of $\R^k$. The $p$\-cells for $p < k$ are not depicted. For $k < n$, the $k$\-projected representation of an $n$\-diagram does not provide complete information about that diagram. Nonetheless, it is the correct notion for describing the action of homotopy generators, as we explore in Section~\ref{chapterinterchangers}.

\begin{definition}
\label{diagsiggraphicalrep}
For an $n$\-diagram $D$ over a signature $\Sigma$, for $k \leq n$, its \textit{$k$\-projected graphical representation} $G_{D}^k \subset \R^k$ is a labelled partitioned subspace, defined as follows:
\begin{itemize}
\item For $n=0$, to be $G_D^0:=\R^{0}$
\item For $n>0$:
\begin{itemize}
\item at height $i$, to agree with  \mbox{$G_{D[i].d}^{k-1} \subset \R ^{n{-}1} \subset \R^n$;\hspace{-2cm}}
\item between heights $i$ and $i+1$, as a glued double cone, with centre point labelled by the element $D[i].g \in \Sigma_n$.
\end{itemize}
\end{itemize}
\end{definition}

\noindent
The gluing scheme here is straightforward in low dimensions, but not precise in general. We illustrate it here by example.

\paragraph{Example.}
\tikzset{blob/.style={circle, fill, inner sep=2pt, draw=none}}
Consider the following 2\-signature $\Sigma$, with some components coloured to help identification:
\begin{align*}
\text{\bf 0-cells}&&& \A,\B,\CC
\\\\[-15pt]
\text{\bf 1-cells}&& F&: \A \to \B
\\
&& G&: \B \to \CC
\\
&& H&: \CC \to \CC
\\\\[-15pt]
\text{\bf 2-cells} && \mu &: H \Rightarrow G \circ H
\\
&& \nu &: F \circ G \Rightarrow F
\\
&& \phi &: H \Rightarrow \text{id}_C
\end{align*}
Then consider the following $2$\-diagram $D$ defined over $\Sigma$:
\begin{align*}
&D.s.s[0].g = \A \\
&D.s[0].g = F &&D.s[0].e = [] \\
&D.s[1].g = F &&D.s[1].e = [] \\
&D[0].g = \mu &&D[0].e = [1] \\
&D[1].g = \nu && D[1].e = [0] \\
&D[2].g = \phi && D[2].e = [1] \\
\end{align*}
Note that according to Definition~\ref{defembedding}, embeddings between $0$\-diagrams consist of no data. For this reason the lists $D.s[0].e$ and $D.s[1].e$ are empty.

Following the recursive procedure, we begin at the base cases, depicting $A$, $B$ and $C$ as points as follows:
\[
\begin{aligned}
\begin{tikzpicture}[scale=1.2]
\path (0.5,0) node [blob, fill=blue] {} node [above] {$A$};
\path (1.5,0) node [blob, fill=red] {} node [above] {$B$};
\path (2.5,0) node [blob, fill=green] {} node [above] {$C$};
\path (0,1) node [blob, fill=blue] {} node [above] {$A$};
\path (1,1) node [blob, fill=red] {} node [above] {$B$};
\path (2,1) node [blob, fill=green] {} node [above] {$C$};
\path (3,1) node [blob, fill=green] {} node [above] {$C$};
\path (0.5,2) node [blob, fill=blue] {} node [above] {$A$};
\path (1.5,2) node [blob, fill=green] {} node [above] {$C$};
\path (2.5,2) node [blob, fill=green] {} node [above] {$C$};
\path (1,3) node [blob, fill=blue] {} node [above] {$A$};
\path (2,3) node [blob, fill=green] {} node [above] {$C$};
\end{tikzpicture}
\end{aligned}
\]
One step up the call stack, we represent the 1\-diagram slices $D[0].d$, $D[1].d$, $D[2].d$ and $D[3].d$ as 1\-dimensional glued cones over these points, labelling the centre point of the cones with the name of the associated generator: 
\[
\begin{aligned}
\begin{tikzpicture}[scale=1.2]
\node (A) at (1, 0) {};
\node (B) at (2, 0) {};
\node (C) at (2,0.5) {};
\node (D) at (0.5,1) {};
\node (E) at (1.5,1) {};
\node (F) at (2.5,1) {};
\node (G) at (1,1.5) {};
\node (H) at (1, 2) {};
\path (1,0) {} to (0.5,1) {} to (1,1.5) {} to (1,2) {};
\path (2,0) {} to (2,0.5) {} to (1.5,1) {} to (1,1.5) {};
\path (C.center) to (2.5,1) node (F) {} to (2,2) node (I) {};
\node (J) at (2,2.5) {};
\node (K) at (1.5,3) {};
\node (1) at (0,0) {};
\node (3) at (1,3) {};
\node (2) at (2.5,0) {};
\node (4) at (2,3) {};
\draw [draw=none, fill=blue!30, opacity=0.5] (0.5,0) node (1) {} to (A.center) to (D.center) to (G.center) to (H.center) to (K.center) to (3.center) to (0.5,2) to (0,1) to (1.center);
\draw [draw=none, fill=red!30, opacity=0.5] (A.center) to (D.center) to (G.center) to (C.center) to (B.center);
\draw [draw=none, fill=green!50, opacity=0.5] (K.center) to (4.center) to (2.5,2) to (3,1) to (2.center) to (B.center) to (C.center) to (G.center) to (H.center);
\draw [thick] (A.center) to (D.center) to (G.center) to (H.center) to (K.center);
\draw [thick] (B.center) to (C.center) to (E.center) to (G.center);
\path [thick] (C.center) to (F.center) to (I.center);
\draw [thick] (J.center) to (I.center) to (F.center) to (C.center);
\node [circle, draw, fill=white, inner sep=2pt, scale=0.7] at (C.center) {$\mu$};
\node [circle, draw, fill=white, inner sep=2pt, scale=0.7] at (G.center) {$\nu$};
\node [circle, draw, fill=white, inner sep=2pt, scale=0.7] at (J.center) {$\sigma$};
\def\coveropacity{1}
\path [fill=white, opacity=\coveropacity] (0,0) rectangle (3,3);
\draw [thick, blue] (1.center) to (A.center);
\draw [thick, red] (A.center) to (B.center) ;
\draw [thick, green] (B.center) to (2.center);
\draw [thick, blue] (0,1) to (D.center);
\draw [thick, red] (D.center) to (E.center);
\draw [thick, green] (E.center) to (3,1);
\draw [thick, blue] (0.5,2) to (H.center) ;
\draw [thick, green] (H.center) to (I.center) to (2.5,2);
\draw [thick, blue] (1,3) to (K.center) ;
\draw [thick, green] (K.center) to (2,3);
\path (1,0) node [blob] {} node [above] {$F$};
\path (2,0) node [blob] {} node [above] {$G$};
\path (1.5,3) node [blob] {} node [above] {$F$};
\path (0.5,1) node [blob] {} node [above] {$F$};
\path (1.5,1) node [blob] {} node [above] {$G$};
\path (2.5,1) node [blob] {} node [above] {$H$};
\path (2,2) node [blob] {} node [above] {$H$};
\path (1,2) node [blob] {} node [above] {$F$};
\end{tikzpicture}
\end{aligned}
\]
Then at the top level, we glue 2\-dimensional cones between the 1\-dimensional slices to produce the 2\-dimensional image:
\begin{align*}
\begin{aligned}
\begin{tikzpicture}[scale=1.2]
\path [use as bounding box] (0,0) rectangle (3,3);
\node (A) at (1, 0) {};
\node (B) at (2, 0) {};
\node (C) at (2,0.5) {};
\node (D) at (0.5,1) {};
\node (E) at (1.5,1) {};
\node (F) at (2.5,1) {};
\node (G) at (1,1.5) {};
\node (H) at (1, 2) {};
\path (1,0) {} to (0.5,1) {} to (1,1.5) {} to (1,2) {};
\path (2,0) {} to (2,0.5) {} to (1.5,1) {} to (1,1.5) {};
\path (C.center) to (2.5,1) node (F) {} to (2,2) node (I) {};
\node (J) at (2,2.5) {};
\node (K) at (1.5,3) {};
\node (1) at (0,0) {};
\node (3) at (1,3) {};
\node (2) at (2.5,0) {};
\node (4) at (2,3) {};
\draw [draw=none, fill=blue!30, opacity=0.5] (0.5,0) node (1) {} to (A.center) to (D.center) to (G.center) to (H.center) to (K.center) to (3.center) to (0.5,2) to (0,1) to (1.center);
\draw [draw=none, fill=red!30, opacity=0.5] (A.center) to (D.center) to (G.center) to (C.center) to (B.center);
\draw [draw=none, fill=green!50, opacity=0.5] (K.center) to (4.center) to (2.5,2) to (3,1) to (2.center) to (B.center) to (C.center) to (G.center) to (H.center);
\draw [thick] (A.center) to (D.center) to (G.center) to (H.center) to (K.center);
\draw [thick] (B.center) to (C.center) to (E.center) to (G.center);
\path [thick] (C.center) to (F.center) to (I.center);
\draw [thick] (J.center) to (I.center) to (F.center) to (C.center);
\node [circle, draw, fill=white, inner sep=2pt, scale=0.7, thick] at (C.center) {$\mu$};
\node [circle, draw, fill=white, inner sep=2pt, scale=0.7, thick] at (G.center) {$\nu$};
\node [circle, draw, fill=white, inner sep=2pt, scale=0.7, thick] at (J.center) {$\sigma$};
\def\coveropacity{0.75}
\path [fill=white, opacity=\coveropacity] (0,0) rectangle (3,3);
\draw [thick, blue!50] (1.center) to (A.center);
\draw [thick, red!50] (A.center) to (B.center);
\draw [thick, green!50] (B.center) to (2.center);
\draw [thick, blue!50] (0,1) to (D.center);
\draw [thick, red!50] (D.center) to (E.center);
\draw [thick, green!50] (E.center) to (F.center) to (3,1);
\draw [thick, blue!50] (0.5,2) to (H.center);
\draw [thick, green!50] (H.center) to (I.center) to (2.5,2);
\draw [thick, blue!50] (1,3) to (K.center);
\draw [thick, green!50] (K.center) to (2,3);
\path (1,0) node [blob] {};
\path (2,0) node [blob] {};
\path (1.5,3) node [blob] {};
\path (0.5,1) node [blob] {};
\path (1.5,1) node [blob] {};
\path (2.5,1) node [blob] {};
\path (2,2) node [blob] {};
\path (1,2) node [blob] {};
\end{tikzpicture}
\end{aligned}
\qquad\leadsto\qquad
\begin{aligned}
\begin{tikzpicture}[scale=1.2]
\path [use as bounding box] (0,0) rectangle (3,3);
\node (A) at (1, 0) {};
\node (B) at (2, 0) {};
\node (C) at (2,0.5) {};
\node (D) at (0.5,1) {};
\node (E) at (1.5,1) {};
\node (F) at (2.5,1) {};
\node (G) at (1,1.5) {};
\node (H) at (1, 2) {};
\path (1,0) {} to (0.5,1) {} to (1,1.5) {} to (1,2) {};
\path (2,0) {} to (2,0.5) {} to (1.5,1) {} to (1,1.5) {};
\path (C.center) to (2.5,1) node (F) {} to (2,2) node (I) {};
\node (J) at (2,2.5) {};
\node (K) at (1.5,3) {};
\node (1) at (0,0) {};
\node (3) at (1,3) {};
\node (2) at (2.5,0) {};
\node (4) at (2,3) {};
\draw [draw=none, fill=blue!30, opacity=0.5] (0.5,0) node (1) {} to (A.center) to (D.center) to (G.center) to (H.center) to (K.center) to (3.center) to (0.5,2) to (0,1) to (1.center);
\draw [draw=none, fill=red!30, opacity=0.5] (A.center) to (D.center) to (G.center) to (C.center) to (B.center);
\draw [draw=none, fill=green!50, opacity=0.5] (K.center) to (4.center) to (2.5,2) to (3,1) to (2.center) to (B.center) to (C.center) to (G.center) to (H.center);
\draw [thick] (A.center) to (D.center) to (G.center) to (H.center) to (K.center);
\draw [thick] (B.center) to (C.center) to (E.center) to (G.center);
\path [thick] (C.center) to (F.center) to (I.center);
\draw [thick] (J.center) to (I.center) to (F.center) to (C.center);
\node [circle, draw, fill=white, inner sep=2pt, scale=0.7, thick] at (C.center) {$\mu$};
\node [circle, draw, fill=white, inner sep=2pt, scale=0.7, thick] at (G.center) {$\nu$};
\node [circle, draw, fill=white, inner sep=2pt, scale=0.7, thick] at (J.center) {$\sigma$};
\end{tikzpicture}
\end{aligned}
\end{align*}
This idea can be continued straightforwardly into dimension 3, yielding representations of 3\-projected $n$\-diagrams for $n \geq 3$ as subspaces of $\R^3$.

\paragraph{Movies.}
Another way to visualize an $n$\-diagram $D$ is by visualizing each $(n{-}1)$\-dimensional slice $D[i].d$ for $0 \leq i \leq |D|$ sequentially. We call this a \textit{movie}. By iterating this idea, one can visualize diagrams of arbitrary dimension without using projections.  We demonstrate this in Section~\ref{chapterbutterfly} with a movie presentation of a 5\-diagram, giving each slice as a 2\-projection.

\subsection{Whiskered composites}

\def\threecompxscale{0.9}
\def\threecompyscale{0.9}
\def\threecomplabelscale{0.6}

Let us consider the example of two 3\-diagrams $A$ and $B$, each consisting of a single generator $\alpha$ and $\beta$ respectively, represented graphically as follows:
\begin{align*}
\begin{aligned}
\begin{tikzpicture}[thick, yscale = \threecompyscale, xscale = \threecompxscale]
\draw [region] (0,0) rectangle +(2,2);
\draw (1,0) to (1, 2);
\node [morphism] at (1,1) {$\alpha$};
\end{tikzpicture}
\end{aligned}&&&
\begin{aligned}
\begin{tikzpicture}[thick, yscale = \threecompyscale, xscale = \threecompxscale]
\draw [region] (0,0) rectangle +(2,2);
\draw (1,0) to (1, 2);
\node [morphism] at (1,1) {$\beta$};
\end{tikzpicture}
\end{aligned}
\end{align*}
These are both 3D pictures, the 0-cells represented by spaces at the front and the back of the picture are not labelled.
Here we examine graphical representations of all the ways that these can be composed, assuming that there is no type restriction.
\paragraph{Vertical composition.}
By Definition~\ref{diagcompose}, there is only one way in which they could be directly composed. This is realised by diagram composition directly, so we consider the composite $\beta \circ \alpha$. We take the graphical representations $G_A$ and $G_B$ and paste them together by placing $G_D$ {on top} of $G_S$:
\begin{calign}
\nonumber
\begin{aligned}
\begin{tikzpicture}[thick, yscale = \threecompyscale, xscale = \threecompxscale]
\draw [region] (0,0) rectangle +(2,3);
\draw (1,0) to (1, 3);
\node [morphism] at (1,1) {$\alpha$};
\node [morphism] at (1,2) {$\beta$};
\end{tikzpicture}
\end{aligned}
\\\nonumber
B \circ A
\end{calign}

\paragraph{Codimension 1 composition.} By whiskering $A$ and $B$ before we compose them, we can effectively perform a composition along a common boundary 1\-cell. There are 2 variants.
\begin{calign}
\nonumber
\begin{aligned}
\begin{tikzpicture}[thick, yscale = \threecompyscale, xscale = \threecompxscale]
\draw [region] (0,0) rectangle +(3,3);
\draw (1,0) to (1, 3);
\draw (2,0) to (2, 3);
\node [morphism] at (1,1) {$\alpha$};
\node [morphism] at (2,2) {$\beta$};
\end{tikzpicture}
\end{aligned}
&
\begin{aligned}
\begin{tikzpicture}[thick, yscale = \threecompyscale, xscale = \threecompxscale]
\draw [region] (0,0) rectangle +(3,3);
\draw (1,0) to (1, 3);
\draw (2,0) to (2, 3);
\node [morphism] at (1,2) {$\alpha$};
\node [morphism] at (2,1) {$\beta$};
\end{tikzpicture}
\end{aligned}
\\
\nonumber
(A\circ B.s) \circ (A.t \circ B)
&
(A.s \circ B) \circ (A\circ B.t)
\end{calign}

\paragraph{Codimension 2 composition.} By whiskering twice, we can effectively perform a composition along a common boundary 0\-cell. There are 4 variants.
\begin{calign}
\nonumber
\begin{aligned}
\begin{tikzpicture}[thick, yscale = \threecompyscale, xscale = \threecompxscale]
\draw [region] (0,0) rectangle +(3,3);
\draw (1,0) to (1, 3);
\node [morphism] at (1,0.8) {$\alpha$};
\draw [region] (0.2,-0.2) rectangle +(3,3);
\draw (2.2,-0.2) to (2.2, 2.8);
\node [morphism] at (2.2,2) {$\beta$};
\end{tikzpicture}
\end{aligned}
&
\begin{aligned}
\begin{tikzpicture}[thick, yscale = \threecompyscale, xscale = \threecompxscale]
\draw [region] (0,0) rectangle +(3,3);
\draw (1,0) to (1, 3);
\node [morphism] at (1,2) {$\alpha$};
\draw [region] (0.2,-0.2) rectangle +(3,3);
\draw (2.2,-0.2) to (2.2, 2.8);
\node [morphism] at (2.2,0.8) {$\beta$};
\end{tikzpicture}
\end{aligned}
\\\nonumber
\big((A\circ B.s.s) \circ (A.s.t \circ B.s)\big)
&
\big((A.s\circ A.s.s) \circ (A.t.t \circ B) \big)
\\\nonumber
{\,\,\,}\circ \big((A.t\circ B.s.s)\circ (A.t.t\circ B)\big)
&
{\,\,\,}\circ \big( (A\circ B.t.s)\circ (A.t.t\circ B.t) \big)
\\[10pt]\nonumber
\begin{aligned}
\begin{tikzpicture}[thick, yscale = \threecompyscale, xscale = \threecompxscale]
\draw [region] (0,0) rectangle +(3,3);
\draw (2.2,0) to (2.2, 3);
\node [morphism] at (2.2,2) {$\alpha$};
\draw [region] (0.2,-0.2) rectangle +(3,3);
\draw (1,-0.2) to (1, 2.8);
\node [morphism] at (1,0.8) {$\beta$};
\end{tikzpicture}
\end{aligned}
&
\begin{aligned}
\begin{tikzpicture}[thick, yscale = \threecompyscale, xscale = \threecompxscale]
\draw [region] (0,0) rectangle +(3,3);
\draw (2.2,0) to (2.2, 3);
\node [morphism] at (2.2,0.8) {$\alpha$};
\draw [region] (0.2,-0.2) rectangle +(3,3);
\draw (1,-0.2) to (1, 2.8);
\node [morphism] at (1,2) {$\beta$};
\end{tikzpicture}
\end{aligned}
\\\nonumber
\big( (A.s.s \circ B) \circ (A.s \circ B.t.t) \big)
&
\big( (A.s.s\circ B.s) \circ (A \circ B.t.s) \big)
\\\nonumber
{\,\,\,} \circ \big( (A.s.s\circ B.t)\circ (A\circ B.t.t) \big)
&
{\,\,\,} \circ \big( (A.s.t \circ B) \circ (A.t \circ B.t.t) \big)
\end{calign}

\section{Homotopy generators}
\label{chapterinterchangers}

In this section we introduce the homotopy generators which must be supported by signatures as part of the definitions of semistrict $n$\-category, for $2 \leq n \leq 4$, which we give in the introduction. We also describe their extended versions, required for the definition of quasistrict 4\-category. 

\paragraph{Invertibility.} We use a familiar coinductive definition of invertible cell in a higher category\JVcomm{Should really cite Shulman's cafe entry.}. For the purposes of this definition, for a generator $f$, we write $[f]$ for the diagram of height 1 consisting only of that generator.
\begin{definition}
For $n>0$, given an $n$\-signature $\Sigma$, a $k$\-generator $f \in \Sigma_k$ for $0 < k \leq n$ is \textit{invertible} when it is equipped with:
\begin{itemize}
\item an \textit{inverse} $f^{-1} \in \Sigma_k$ with $f^{-1}.t = f.s$ and $f ^{-1}.s = f.t$;
\item when $k<n$, invertible $(k{+}1)$\-cells $f', f'' \in \Sigma_{k+1}$ as follows:
\begin{align*}
f'.s &= [f] \circ [f^{-1}] & f''.s &= f.t.\text{Id}
\\
f'.t &= {f.s.\text{Id}} & f''.t &= [f ^{-1}] \circ [f]
\end{align*}
\end{itemize}
\end{definition}

\subsection{Type I homotopy generators}

This move allows height rearrangement of adjacent cells. Note that these generators are always of dimension 3 or higher.
\begin{definition}[Type I homotopy generator]
\label{defswap}
For $n \geq 0$, an $n$\-signature $\Sigma$ supports \emph{type~I homotopy generators} if for any $2 \leq k < n$, for any $k$\-diagram with 2\-projection given by the the left-hand diagram below, there is a chosen invertible $(k{+}1)$\-cell as follows:
\begin{equation}
\begin{aligned}
\begin{tikzpicture}[thick, yscale=0.9]
\draw (0.3,1.3) to [out=up, in=down] (0,2) to [out=up, in=down](0.3, 2.7);
\draw (-0.3,1.3) to [out=up, in=down] (0,2) to [out=up, in=down](-0.3, 2.7);
\draw (1.5, 1.3) to (1.5, 2.7);
\draw (2.1, 1.3) to (2.1, 2.7);
\draw (0.3, 1.3) to (0.3, -0.1);
\draw (-0.3, 1.3) to (-0.3, -0.1);
\draw (1.5,-0.1) to [out=up, in=down] (1.8,0.6) to [out=up, in=down](1.5, 1.3);
\draw (2.1,-0.1) to [out=up, in=down] (1.8,0.6) to [out=up, in=down](2.1, 1.3);
\draw (0.6, 2.7) to (0.6, -0.1);
\draw (1.2, 2.7) to (1.2, -0.1);
\node [morphism] at (0,2) {$f$};
\node [morphism] at (1.8,0.6) {$g$};
\node [scale=0.6] at (1.8,0) {$\ldots$};
\node [scale=0.6] at (1.8,1.2) {$\ldots$};
\node [scale=0.6] at (0,2.6) {$\ldots$};
\node [scale=0.6] at (0,1.4) {$\ldots$};
\node [scale=0.6] at (0.9,0) {$\ldots$};
\node [scale=0.6] at (0.9,2.6) {$\ldots$};
\end{tikzpicture}
\end{aligned}
\quad\stackrel{\I_k}{\rightarrow}\quad
\begin{aligned}
\begin{tikzpicture}[thick, yscale=0.9]
\draw (2.1,1.3) to [out=up, in=down] (1.8,2) to [out=up, in=down](2.1, 2.7);
\draw (1.5,1.3) to [out=up, in=down] (1.8,2) to [out=up, in=down](1.5, 2.7);
\draw (0.3, 1.3) to (0.3, 2.7);
\draw (-0.3, 1.3) to (-0.3, 2.7);
\draw (1.5, 1.3) to (1.5, -0.1);
\draw (2.1, 1.3) to (2.1, -0.1);
\draw (0.3,-0.1) to [out=up, in=down] (0,0.6) to [out=up, in=down](0.3, 1.3);
\draw (-0.3,-0.1) to [out=up, in=down] (0,0.6) to [out=up, in=down](-0.3, 1.3);
\draw (0.6, 2.7) to (0.6, -0.1);
\draw (1.2, 2.7) to (1.2, -0.1);
\node [morphism] at (1.8,2) {$g$};
\node [morphism] at (0,0.6) {$f$};
\node [scale=0.6] at (0,0) {$\ldots$};
\node [scale=0.6] at (0,1.2) {$\ldots$};
\node [scale=0.6] at (1.8,2.6) {$\ldots$};
\node [scale=0.6] at (1.8,1.4) {$\ldots$};
\node [scale=0.6] at (0.9,0) {$\ldots$};
\node [scale=0.6] at (0.9,2.6) {$\ldots$};
\end{tikzpicture}
\end{aligned}
\end{equation}
\end{definition}

\noindent
By abuse of notation, we simply refer to these $(k{+}1)$\-cells as $\I_k$. They are indexed formally by their entire source diagram, drawn here on the left; while the left-hand diagram may be a $k$\-diagram for $k > 2$, only the features visible in its 2\-projection are relevant for constructing the associated type I homotopy generator. Since they are invertible, the reverse rewrite is also allowed. We emphasize some key features of the left-hand diagram:
\begin{itemize}
\item there are no additional wires to the left of $f$ or right of $g$;
\item there are arbitrary interposing wires;
\item $f$ and $g$ have arbitrary input and output wires;
\item $f$ and $g$ are generators, not composite diagrams.
\end{itemize}

\def\dotscale{0.8}
\def\opacityscale{0.4}

At times, we also depict type I homotopy generators with 3\-projected diagrams as follows, using a braiding convention:
\begin{align*}
\begin{aligned}
\begin{tikzpicture}[thick, scale = 1.6]
\draw  [region, fill opacity = \opacityscale, draw=black](2.2, 1.95) to(2.1, 1.95) to [out=left, in=up] (1.7, 1.8) to (1.7, 1.2) to [out=down, in=up](0.9, 0) to [out=up, in=left] (1.3, 0.15) to (2.2, 0.15);
\node [rotate=120, scale=\dotscale] at (2.3,0) {$\ldots$};
\draw [region, draw=black](2.4, 1.65) to (2.1, 1.65) to [out=left, in=down] (1.7, 1.8)to (1.7, 1.2) to [out=down, in=up](0.9, 0) to [out=down, in=left] (1.3, -0.15) to (2.4, -0.15);
\draw  [region, fill opacity = \opacityscale, draw=black](-0.1, 1.95) to(1.3, 1.95) to [out=right, in=up] (1.7, 1.8) to (1.7, 1.2) to [out=down, in=up](0.9, 0) to [out=up, in=right] (0.5, 0.15) to (-0.1, 0.15);
\draw [region, fill opacity = \opacityscale, draw=black](0.1, 1.65) to (1.3, 1.65) to [out=right, in=down] (1.7, 1.8)to (1.7, 1.2) to [out=down, in=up](0.9, 0) to [out=down, in=right] (0.5, -0.15) to (0.1, -0.15);
\node [scale=0.6] at (0.9,0.4) {$g$};
\node [rotate=120, scale=\dotscale] at (2.5,-0.3) {$\ldots$};
\draw [region, draw=black, fill opacity = \opacityscale](2.6, 1.35) to (1.3, 1.35) to [out=left, in=up] (0.9,1.2) to [out=down, in=up](1.7,0) to (1.7,-0.6) to [out=up, in=left] (2.1, -0.45) to (2.6, -0.45);
\node [rotate=120, scale=\dotscale] at (2.7,-0.6) {$\ldots$};
\draw [region,  fill opacity = \opacityscale,draw=black](2.8, 1.05) to (1.3, 1.05) to [out=left, in=down] (0.9,1.2) to [out=down, in=up](1.7,0) to (1.7,-0.6)  to [out=down, in=left] (2.1, -0.75) to (2.8, -0.75);
\draw [region,  fill opacity = \opacityscale,draw=black](0.3, 1.35) to (0.5, 1.35) to [out=right, in=up] (0.9,1.2) to [out=down, in=up](1.7,0) to (1.7,-0.6) to [out=up, in=right] (1.3, -0.45) to (0.3, -0.45);
\draw [region, fill opacity = \opacityscale, draw=black](0.5, 1.05)  to [out=right, in=down] (0.9,1.2) to [out=down, in=up](1.7,0) to (1.7,-0.6)  to [out=down, in=right] (1.3, -0.75) to (0.5, -0.75);
\node [rotate=120, scale=\dotscale] at (2.3,1.8) {$\ldots$};
\node [rotate=120, scale=\dotscale] at (2.5,1.5) {$\ldots$};
\node [rotate=120, scale=\dotscale] at (2.7,1.2) {$\ldots$};
\node [rotate=120, scale=\dotscale] at (0,1.8) {$\ldots$};
\node [rotate=120, scale=\dotscale] at (0.2,1.5) {$\ldots$};
\node [rotate=120, scale=\dotscale] at (0.4,1.2) {$\ldots$};
\node [rotate=120, scale=\dotscale] at (0,0) {$\ldots$};
\node [rotate=120, scale=\dotscale] at (0.2,-0.3) {$\ldots$};
\node [rotate=120, scale=\dotscale] at (0.4,-0.6) {$\ldots$};
\node [scale=0.6] at (1.8,0) {$f$};
\end{tikzpicture}
\end{aligned}
\end{align*}
This braiding convention is an artistic style, and we do not attempt to formalize it. Note also that for clarity we omit the intervening sheets between $f$ and $g$; we will often adopt this convention to reduce diagram complexity.

\def\expansionscale{0.6}

\paragraph{Expansion scheme.}
We require an expansion scheme for type I homotopy generators, defining their action on composite diagrams. We provide this scheme recursively.

\begin{definition}[Type I homotopy composite]
\label{defgammainterchanger}
In an $n$\-signature that supports type~I homotopy generators, a \emph{type I homotopy composite} is a diagram formed from a sequence of type I homotopy generators, with values defined recursively as follows:
\def\intexpansionscale{0.7}
\def\intexpansionscaleY{0.6}
\begin{align*}
&\begin{aligned}
\begin{tikzpicture}[thick, yscale=\intexpansionscale, scale=\intexpansionscale]
\draw (0.3,1.3) to [out=up, in=down] (0,2) to [out=up, in=down](0.3, 2.7);
\draw (-0.3,1.3) to [out=up, in=down] (0,2) to [out=up, in=down](-0.3, 2.7);
\node [morphism] at (0,2) {$f$};
\draw (1.5, 1.3) to (1.5, 3.4);
\draw (2.1, 1.3) to (2.1, 3.4);
\draw (0.3, 1.3) to (0.3, -2.9);
\draw (-0.3, 1.3) to (-0.3, -2.9);
\draw (1.5,-0.1) to [out=up, in=down] (1.8,0.6) to [out=up, in=down](1.5, 1.3);
\draw (2.1,-0.1) to [out=up, in=down] (1.8,0.6) to [out=up, in=down](2.1, 1.3);
\node [morphism] at (1.8,0.6) {$g_1$};
\draw (1.5,-2.9) to (1.5,-2.2) to [out=up, in=down] (1.8,-1.5) to [out=up, in=down](1.5, -0.8);
\draw (2.1,-2.9) to (2.1,-2.2) to [out=up, in=down] (1.8,-1.5) to [out=up, in=down](2.1, -0.8);
\node [morphism] at (1.8,-1.5) {$g_s$};
\draw (0.6, 3.4) to (0.6, -2.9);
\draw (1.2, 3.4) to (1.2, -2.9);
\draw [dotted](1.5, -0.1) to (1.5, -0.8);
\draw [dotted](2.1, -0.1) to (2.1, -0.8);
\draw (-0.3, 2.7) to (-0.3, 3.4);
\draw (0.3, 2.7) to (0.3, 3.4);
\node [scale=0.6] at (1.8,-2.8) {$\ldots$};
\node [scale=0.6] at (1.8,-0.9) {$\ldots$};
\node [scale=0.6] at (1.8,0) {$\ldots$};
\node [scale=0.6] at (1.8,1.2) {$\ldots$};
\node [scale=0.6] at (1.8,3.3) {$\ldots$};
\node [scale=0.6] at (0,3.3) {$\ldots$};
\node [scale=0.6] at (0,-2.8) {$\ldots$};
\node [scale=0.6] at (0,1.4) {$\ldots$};
\node [scale=0.6] at (0.9,-2.8) {$\ldots$};
\node [scale=0.6] at (0.9,3.3) {$\ldots$};
\end{tikzpicture}
\end{aligned}
\stackrel{\widetilde{\I}_k}{\rightarrow}
\begin{aligned}
\begin{tikzpicture}[thick, yscale=\intexpansionscale, scale=\intexpansionscale]
\draw (0.3,-0.1) to [out=up, in=down] (0,0.6) to [out=up, in=down](0.3, 1.3);
\draw (-0.3,-0.1) to [out=up, in=down] (0,0.6) to [out=up, in=down](-0.3, 1.3);
\node [morphism] at (0,0.6) {$f$};
\draw (1.5, 1.3) to (1.5, -0.8);
\draw (2.1, 1.3) to (2.1, -0.8);
\draw (0.3, 1.3) to (0.3, 5.5);
\draw (-0.3, 1.3) to (-0.3, 5.5);
\draw (1.5,1.3) to [out=up, in=down] (1.8,2) to [out=up, in=down](1.5, 2.7);
\draw (2.1,1.3) to [out=up, in=down] (1.8,2) to [out=up, in=down](2.1, 2.7);
\node [morphism] at (1.8,2) {$g_1$};
\draw (1.5,3.4) to [out=up, in=down] (1.8,4.1) to [out=up, in=down](1.5, 4.8) to (1.5, 5.5);
\draw (2.1,3.4) to [out=up, in=down] (1.8,4.1) to [out=up, in=down](2.1, 4.8) to (2.1, 5.5);
\node [morphism] at (1.8,4.1) {$g_s$};
\draw (0.6, 5.5) to (0.6, -0.8);
\draw (1.2, 5.5) to (1.2, -0.8);
\draw [dotted](1.5, 3.4) to (1.5, 2.7);
\draw [dotted](2.1, 3.4) to (2.1, 2.7);
\draw (-0.3, -0.1) to (-0.3, -0.8);
\draw (0.3, -0.1) to (0.3, -0.8);
\node [scale=0.6] at (1.8,-0.7) {$\ldots$};
\node [scale=0.6] at (1.8,1.4) {$\ldots$};
\node [scale=0.6] at (1.8,2.6) {$\ldots$};
\node [scale=0.6] at (1.8,3.5) {$\ldots$};
\node [scale=0.6] at (1.8,5.4) {$\ldots$};
\node [scale=0.6] at (0,-0.7) {$\ldots$};
\node [scale=0.6] at (0,5.4) {$\ldots$};
\node [scale=0.6] at (0,1.2) {$\ldots$};
\node [scale=0.6] at (0.9,-0.7) {$\ldots$};
\node [scale=0.6] at (0.9,5.4) {$\ldots$};
\end{tikzpicture}
\end{aligned}
\quad:=\quad
\begin{aligned}
\begin{tikzpicture}[thick, yscale=\intexpansionscale, scale=\intexpansionscale]
\draw (0.3,1.3) to [out=up, in=down] (0,2) to [out=up, in=down](0.3, 2.7);
\draw (-0.3,1.3) to [out=up, in=down] (0,2) to [out=up, in=down](-0.3, 2.7);
\node [morphism] at (0,2) {$f$};
\draw (1.5, 1.3) to (1.5, 3.4);
\draw (2.1, 1.3) to (2.1, 3.4);
\draw (0.3, 1.3) to (0.3, -2.9);
\draw (-0.3, 1.3) to (-0.3, -2.9);
\draw (1.5,-0.1) to [out=up, in=down] (1.8,0.6) to [out=up, in=down](1.5, 1.3);
\draw (2.1,-0.1) to [out=up, in=down] (1.8,0.6) to [out=up, in=down](2.1, 1.3);
\node [morphism] at (1.8,0.6) {$g_1$};
\draw (1.5,-2.9) to (1.5,-2.2) to [out=up, in=down] (1.8,-1.5) to [out=up, in=down](1.5, -0.8);
\draw (2.1,-2.9) to (2.1,-2.2) to [out=up, in=down] (1.8,-1.5) to [out=up, in=down](2.1, -0.8);
\node [morphism] at (1.8,-1.5) {$g_s$};
\draw (0.6, 3.4) to (0.6, -2.9);
\draw (1.2, 3.4) to (1.2, -2.9);
\draw [dotted](1.5, -0.1) to (1.5, -0.8);
\draw [dotted](2.1, -0.1) to (2.1, -0.8);
\draw (-0.3, 2.7) to (-0.3, 3.4);
\draw (0.3, 2.7) to (0.3, 3.4);
\node [scale=0.6] at (1.8,-2.8) {$\ldots$};
\node [scale=0.6] at (1.8,-0.9) {$\ldots$};
\node [scale=0.6] at (1.8,0) {$\ldots$};
\node [scale=0.6] at (1.8,1.2) {$\ldots$};
\node [scale=0.6] at (1.8,3.3) {$\ldots$};
\node [scale=0.6] at (0,3.3) {$\ldots$};
\node [scale=0.6] at (0,-2.8) {$\ldots$};
\node [scale=0.6] at (0,1.4) {$\ldots$};
\node [scale=0.6] at (0.9,-2.8) {$\ldots$};
\node [scale=0.6] at (0.9,3.3) {$\ldots$};
\end{tikzpicture}
\end{aligned}
\stackrel{\I_k}{\rightarrow}
\begin{aligned}
\begin{tikzpicture}[thick, yscale=\intexpansionscale, scale=\intexpansionscale]
\draw (0.3,0.6) to [out=up, in=down] (0,1.3) to [out=up, in=down](0.3, 2);
\draw (-0.3,0.6) to [out=up, in=down] (0,1.3) to [out=up, in=down](-0.3, 2);
\node [morphism] at (0,1.3) {$f$};
\draw (1.5, 0.6) to (1.5, 2);
\draw (2.1, 0.6) to (2.1, 2);
\draw (0.3, 0.6) to (0.3, -2.9);
\draw (-0.3, 0.6) to (-0.3, -2.9);
\draw (1.5,2) to [out=up, in=down] (1.8,2.7) to [out=up, in=down](1.5, 3.4);
\draw (2.1,2) to [out=up, in=down] (1.8,2.7) to [out=up, in=down](2.1, 3.4);
\node [morphism] at (1.8,2.7) {$g_1$};
\draw (1.5,-0.8) to [out=up, in=down] (1.8,-0.1) to [out=up, in=down](1.5, 0.6);
\draw (2.1,-0.8) to [out=up, in=down] (1.8,-0.1) to [out=up, in=down](2.1, 0.6);
\node [morphism] at (1.8,-0.1) {$g_2$};
\draw (1.5,-2.9) to [out=up, in=down] (1.8,-2.2) to [out=up, in=down](1.5, -1.5);
\draw (2.1,-2.9) to [out=up, in=down] (1.8,-2.2) to [out=up, in=down](2.1, -1.5);
\node [morphism] at (1.8,-2.2) {$g_s$};
\draw (0.6, 3.4) to (0.6, -2.9);
\draw (1.2, 3.4) to (1.2, -2.9);
\draw [dotted](1.5, -0.8) to (1.5, -1.5);
\draw [dotted](2.1, -0.8) to (2.1, -1.5);
\draw (-0.3, 2) to (-0.3, 3.4);
\draw (0.3, 2) to (0.3, 3.4);
\node [scale=0.6] at (1.8,-2.8) {$\ldots$};
\node [scale=0.6] at (1.8,-0.7) {$\ldots$};
\node [scale=0.6] at (1.8,-1.6) {$\ldots$};
\node [scale=0.6] at (1.8,2) {$\ldots$};
\node [scale=0.6] at (1.8,0.5) {$\ldots$};
\node [scale=0.6] at (1.8,3.3) {$\ldots$};
\node [scale=0.6] at (0,3.3) {$\ldots$};
\node [scale=0.6] at (0,-2.8) {$\ldots$};
\node [scale=0.6] at (0,0.7) {$\ldots$};
\node [scale=0.6] at (0,1.9) {$\ldots$};
\node [scale=0.6] at (0.9,-2.8) {$\ldots$};
\node [scale=0.6] at (0.9,3.3) {$\ldots$};
\end{tikzpicture}
\end{aligned}
\stackrel{\widetilde{\I}_k}{\rightarrow}
\begin{aligned}
\begin{tikzpicture}[thick, yscale=\intexpansionscale, scale=\intexpansionscale]
\draw (0.3,-0.1) to [out=up, in=down] (0,0.6) to [out=up, in=down](0.3, 1.3);
\draw (-0.3,-0.1) to [out=up, in=down] (0,0.6) to [out=up, in=down](-0.3, 1.3);
\node [morphism] at (0,0.6) {$f$};
\draw (1.5, 1.3) to (1.5, -0.8);
\draw (2.1, 1.3) to (2.1, -0.8);
\draw (0.3, 1.3) to (0.3, 5.5);
\draw (-0.3, 1.3) to (-0.3, 5.5);
\draw (1.5,1.3) to [out=up, in=down] (1.8,2) to [out=up, in=down](1.5, 2.7);
\draw (2.1,1.3) to [out=up, in=down] (1.8,2) to [out=up, in=down](2.1, 2.7);
\node [morphism] at (1.8,2) {$g_1$};
\draw (1.5,3.4) to [out=up, in=down] (1.8,4.1) to [out=up, in=down](1.5, 4.8) to (1.5, 5.5);
\draw (2.1,3.4) to [out=up, in=down] (1.8,4.1) to [out=up, in=down](2.1, 4.8) to (2.1, 5.5);
\node [morphism] at (1.8,4.1) {$g_s$};
\draw (0.6, 5.5) to (0.6, -0.8);
\draw (1.2, 5.5) to (1.2, -0.8);
\draw [dotted](1.5, 3.4) to (1.5, 2.7);
\draw [dotted](2.1, 3.4) to (2.1, 2.7);
\draw (-0.3, -0.1) to (-0.3, -0.8);
\draw (0.3, -0.1) to (0.3, -0.8);
\node [scale=0.6] at (1.8,-0.7) {$\ldots$};
\node [scale=0.6] at (1.8,1.4) {$\ldots$};
\node [scale=0.6] at (1.8,2.6) {$\ldots$};
\node [scale=0.6] at (1.8,3.5) {$\ldots$};
\node [scale=0.6] at (1.8,5.4) {$\ldots$};
\node [scale=0.6] at (0,-0.7) {$\ldots$};
\node [scale=0.6] at (0,5.4) {$\ldots$};
\node [scale=0.6] at (0,1.2) {$\ldots$};
\node [scale=0.6] at (0.9,-0.7) {$\ldots$};
\node [scale=0.6] at (0.9,5.4) {$\ldots$};
\end{tikzpicture}
\end{aligned}
\\
&\begin{aligned}
\begin{tikzpicture}[thick, yscale=\intexpansionscale, scale=\intexpansionscale]
\draw (0.3,3.4) to [out=up, in=down] (0,4.1) to [out=up, in=down](0.3, 4.8);
\draw (-0.3,3.4) to [out=up, in=down] (0,4.1) to [out=up, in=down](-0.3, 4.8);
\node [morphism] at (0,4.1) {$f_r$};
\draw (0.3,1.3) to [out=up, in=down] (0,2) to [out=up, in=down](0.3, 2.7);
\draw (-0.3,1.3) to [out=up, in=down] (0,2) to [out=up, in=down](-0.3, 2.7);
\node [morphism] at (0,2) {$f_1$};
\draw (-0.3, 4.8) to (-0.3, 5.5);
\draw (0.3, 4.8) to (0.3, 5.5);
\draw (1.5, 1.3) to (1.5, 5.5);
\draw (2.1, 1.3) to (2.1, 5.5);
\draw (0.3, 1.3) to (0.3, -2.9);
\draw (-0.3, 1.3) to (-0.3, -2.9);
\draw (1.5, -2.2) to (1.5, -2.9);
\draw (2.1, -2.2) to (2.1, -2.9);
\draw (1.5,-0.1) to [out=up, in=down] (1.8,0.6) to [out=up, in=down](1.5, 1.3);
\draw (2.1,-0.1) to [out=up, in=down] (1.8,0.6) to [out=up, in=down](2.1, 1.3);
\node [morphism] at (1.8,0.6) {$g_1$};
\draw (1.5,-2.2) to [out=up, in=down] (1.8,-1.5) to [out=up, in=down](1.5, -0.8);
\draw (2.1,-2.2) to [out=up, in=down] (1.8,-1.5) to [out=up, in=down](2.1, -0.8);
\node [morphism] at (1.8,-1.5) {$g_s$};
\draw (0.6, 5.5) to (0.6, -2.9);
\draw (1.2, 5.5) to (1.2, -2.9);
\draw [dotted](1.5, -0.1) to (1.5, -0.8);
\draw [dotted](2.1, -0.1) to (2.1, -0.8);
\draw [dotted](-0.3, 2.7) to (-0.3, 3.4);
\draw [dotted](0.3, 2.7) to (0.3, 3.4);
\node [scale=0.6] at (1.8,-2.8) {$\ldots$};
\node [scale=0.6] at (1.8,-0.9) {$\ldots$};
\node [scale=0.6] at (1.8,0) {$\ldots$};
\node [scale=0.6] at (1.8,1.2) {$\ldots$};
\node [scale=0.6] at (1.8,5.4) {$\ldots$};
\node [scale=0.6] at (0,2.6) {$\ldots$};
\node [scale=0.6] at (0,3.5) {$\ldots$};
\node [scale=0.6] at (0,5.4) {$\ldots$};
\node [scale=0.6] at (0,-2.8) {$\ldots$};
\node [scale=0.6] at (0,1.4) {$\ldots$};
\node [scale=0.6] at (0.9,-2.8) {$\ldots$};
\node [scale=0.6] at (0.9,5.4) {$\ldots$};
\end{tikzpicture}
\end{aligned}
\stackrel{\widetilde{\I}_k}{\rightarrow}
\begin{aligned}
\begin{tikzpicture}[thick, yscale=\intexpansionscale, scale=\intexpansionscale]
\draw (0.3,-0.1) to [out=up, in=down] (0,0.6) to [out=up, in=down](0.3, 1.3);
\draw (-0.3,-0.1) to [out=up, in=down] (0,0.6) to [out=up, in=down](-0.3, 1.3);
\node [morphism] at (0,0.6) {$f_r$};
\draw (0.3,-2.2) to [out=up, in=down] (0,-1.5) to [out=up, in=down](0.3, -0.8);
\draw (-0.3,-2.2) to [out=up, in=down] (0,-1.5) to [out=up, in=down](-0.3, -0.8);
\node [morphism] at (0,-1.5) {$f_1$};
\draw (1.5, 1.3) to (1.5, -2.9);
\draw (2.1, 1.3) to (2.1, -2.9);
\draw (0.3, 1.3) to (0.3, 5.5);
\draw (-0.3, 1.3) to (-0.3, 5.5);
\draw (2.1, 4.8) to (2.1, 5.5);
\draw (1.5, 4.8) to (1.5, 5.5);
\draw (1.5,1.3) to [out=up, in=down] (1.8,2) to [out=up, in=down](1.5, 2.7);
\draw (2.1,1.3) to [out=up, in=down] (1.8,2) to [out=up, in=down](2.1, 2.7);
\node [morphism] at (1.8,2) {$g_1$};
\draw (1.5,3.4) to [out=up, in=down] (1.8,4.1) to [out=up, in=down](1.5, 4.8);
\draw (2.1,3.4) to [out=up, in=down] (1.8,4.1) to [out=up, in=down](2.1, 4.8);
\node [morphism] at (1.8,4.1) {$g_s$};
\draw (0.6, 5.5) to (0.6, -2.9);
\draw (1.2, 5.5) to (1.2, -2.9);
\draw (-0.3, -2.2) to (-0.3, -2.9);
\draw (0.3, -2.2) to (0.3, -2.9);
\draw [dotted](1.5, 3.4) to (1.5, 2.7);
\draw [dotted](2.1, 3.4) to (2.1, 2.7);
\draw [dotted](-0.3, -0.1) to (-0.3, -0.8);
\draw [dotted](0.3, -0.1) to (0.3, -0.8);
\node [scale=0.6] at (1.8,-2.8) {$\ldots$};
\node [scale=0.6] at (1.8,1.4) {$\ldots$};
\node [scale=0.6] at (1.8,2.6) {$\ldots$};
\node [scale=0.6] at (1.8,3.5) {$\ldots$};
\node [scale=0.6] at (1.8,5.4) {$\ldots$};
\node [scale=0.6] at (0,0) {$\ldots$};
\node [scale=0.6] at (0,-0.8) {$\ldots$};
\node [scale=0.6] at (0,5.4) {$\ldots$};
\node [scale=0.6] at (0,-2.8) {$\ldots$};
\node [scale=0.6] at (0,1.2) {$\ldots$};
\node [scale=0.6] at (0.9,-2.8) {$\ldots$};
\node [scale=0.6] at (0.9,5.4) {$\ldots$};
\end{tikzpicture}
\end{aligned}
\quad:=\quad
\begin{aligned}
\begin{tikzpicture}[thick, yscale=\intexpansionscale, scale=\intexpansionscale]
\draw (0.3,3.4) to [out=up, in=down] (0,4.1) to [out=up, in=down](0.3, 4.8);
\draw (-0.3,3.4) to [out=up, in=down] (0,4.1) to [out=up, in=down](-0.3, 4.8);
\node [morphism] at (0,4.1) {$f_r$};
\draw (0.3,1.3) to [out=up, in=down] (0,2) to [out=up, in=down](0.3, 2.7);
\draw (-0.3,1.3) to [out=up, in=down] (0,2) to [out=up, in=down](-0.3, 2.7);
\node [morphism] at (0,2) {$f_1$};
\draw (-0.3, 4.8) to (-0.3, 5.5);
\draw (0.3, 4.8) to (0.3, 5.5);
\draw (1.5, 1.3) to (1.5, 5.5);
\draw (2.1, 1.3) to (2.1, 5.5);
\draw (0.3, 1.3) to (0.3, -2.9);
\draw (-0.3, 1.3) to (-0.3, -2.9);
\draw (1.5, -2.2) to (1.5, -2.9);
\draw (2.1, -2.2) to (2.1, -2.9);
\draw (1.5,-0.1) to [out=up, in=down] (1.8,0.6) to [out=up, in=down](1.5, 1.3);
\draw (2.1,-0.1) to [out=up, in=down] (1.8,0.6) to [out=up, in=down](2.1, 1.3);
\node [morphism] at (1.8,0.6) {$g_1$};
\draw (1.5,-2.2) to [out=up, in=down] (1.8,-1.5) to [out=up, in=down](1.5, -0.8);
\draw (2.1,-2.2) to [out=up, in=down] (1.8,-1.5) to [out=up, in=down](2.1, -0.8);
\node [morphism] at (1.8,-1.5) {$g_s$};
\draw (0.6, 5.5) to (0.6, -2.9);
\draw (1.2, 5.5) to (1.2, -2.9);
\draw [dotted](1.5, -0.1) to (1.5, -0.8);
\draw [dotted](2.1, -0.1) to (2.1, -0.8);
\draw [dotted](-0.3, 2.7) to (-0.3, 3.4);
\draw [dotted](0.3, 2.7) to (0.3, 3.4);
\node [scale=0.6] at (1.8,-2.8) {$\ldots$};
\node [scale=0.6] at (1.8,-0.9) {$\ldots$};
\node [scale=0.6] at (1.8,0) {$\ldots$};
\node [scale=0.6] at (1.8,1.2) {$\ldots$};
\node [scale=0.6] at (1.8,5.4) {$\ldots$};
\node [scale=0.6] at (0,2.6) {$\ldots$};
\node [scale=0.6] at (0,3.5) {$\ldots$};
\node [scale=0.6] at (0,5.4) {$\ldots$};
\node [scale=0.6] at (0,-2.8) {$\ldots$};
\node [scale=0.6] at (0,1.4) {$\ldots$};
\node [scale=0.6] at (0.9,-2.8) {$\ldots$};
\node [scale=0.6] at (0.9,5.4) {$\ldots$};
\end{tikzpicture}
\end{aligned}
\stackrel{\widetilde{\I}_k}{\rightarrow}
\begin{aligned}
\begin{tikzpicture}[thick, yscale=\intexpansionscale, scale=\intexpansionscale]
\draw (0.3,3.4) to [out=up, in=down] (0,4.1) to [out=up, in=down](0.3, 4.8);
\draw (-0.3,3.4) to [out=up, in=down] (0,4.1) to [out=up, in=down](-0.3, 4.8);
\node [morphism] at (0,4.1) {$f_r$};
\draw (0.3,1.3) to [out=up, in=down] (0,2) to [out=up, in=down](0.3, 2.7);
\draw (-0.3,1.3) to [out=up, in=down] (0,2) to [out=up, in=down](-0.3, 2.7);
\node [morphism] at (0,2) {$f_2$};
\draw (0.3,-3.6) to [out=up, in=down] (0,-2.9) to [out=up, in=down](0.3, -2.2);
\draw (-0.3,-3.6) to [out=up, in=down] (0,-2.9) to [out=up, in=down](-0.3, -2.2);
\node [morphism] at (0,-2.9) {$f_1$};
\draw (1.5, 1.3) to (1.5, 4.8);
\draw (2.1, 1.3) to (2.1, 4.8);
\draw (0.3, 1.3) to (0.3, -2.2);
\draw (-0.3, 1.3) to (-0.3, -2.2);
\draw (1.5,-0.1) to [out=up, in=down] (1.8,0.6) to [out=up, in=down](1.5, 1.3);
\draw (2.1,-0.1) to [out=up, in=down] (1.8,0.6) to [out=up, in=down](2.1, 1.3);
\node [morphism] at (1.8,0.6) {$g_1$};
\draw (1.5,-2.2) to [out=up, in=down] (1.8,-1.5) to [out=up, in=down](1.5, -0.8);
\draw  (2.1,-2.2) to [out=up, in=down] (1.8,-1.5) to [out=up, in=down](2.1, -0.8);
\node [morphism] at (1.8,-1.5) {$g_s$};
\draw (2.1, -3.6) to (2.1,-2.2);
\draw (1.5, -3.6) to (1.5,-2.2);
\draw (0.6, 4.8) to (0.6, -3.6);
\draw (1.2, 4.8) to (1.2, -3.6);
\draw [dotted](1.5, -0.1) to (1.5, -0.8);
\draw [dotted](2.1, -0.1) to (2.1, -0.8);
\draw [dotted](-0.3, 2.7) to (-0.3, 3.4);
\draw [dotted](0.3, 2.7) to (0.3, 3.4);
\node [scale=0.6] at (1.8,-3.5) {$\ldots$};
\node [scale=0.6] at (1.8,-2.1) {$\ldots$};
\node [scale=0.6] at (1.8,-0.9) {$\ldots$};
\node [scale=0.6] at (1.8,0) {$\ldots$};
\node [scale=0.6] at (1.8,1.2) {$\ldots$};
\node [scale=0.6] at (1.8,4.7) {$\ldots$};
\node [scale=0.6] at (0,2.6) {$\ldots$};
\node [scale=0.6] at (0,3.5) {$\ldots$};
\node [scale=0.6] at (0,4.7) {$\ldots$};
\node [scale=0.6] at (0,-2.1) {$\ldots$};
\node [scale=0.6] at (0,-3.5) {$\ldots$};
\node [scale=0.6] at (0,1.4) {$\ldots$};
\node [scale=0.6] at (0.9,-3.5) {$\ldots$};
\node [scale=0.6] at (0.9,4.7) {$\ldots$};
\end{tikzpicture}
\end{aligned}
\stackrel{\widetilde{\I}_k}{\rightarrow}
\begin{aligned}
\begin{tikzpicture}[thick, yscale=\intexpansionscale, scale=\intexpansionscale]
\draw (0.3,-0.1) to [out=up, in=down] (0,0.6) to [out=up, in=down](0.3, 1.3);
\draw (-0.3,-0.1) to [out=up, in=down] (0,0.6) to [out=up, in=down](-0.3, 1.3);
\node [morphism] at (0,0.6) {$f_r$};
\draw (0.3,-2.2) to [out=up, in=down] (0,-1.5) to [out=up, in=down](0.3, -0.8);
\draw (-0.3,-2.2) to [out=up, in=down] (0,-1.5) to [out=up, in=down](-0.3, -0.8);
\node [morphism] at (0,-1.5) {$f_1$};
\draw (1.5, 1.3) to (1.5, -2.9);
\draw (2.1, 1.3) to (2.1, -2.9);
\draw (0.3, 1.3) to (0.3, 5.5);
\draw (-0.3, 1.3) to (-0.3, 5.5);
\draw (2.1, 4.8) to (2.1, 5.5);
\draw (1.5, 4.8) to (1.5, 5.5);
\draw (1.5,1.3) to [out=up, in=down] (1.8,2) to [out=up, in=down](1.5, 2.7);
\draw (2.1,1.3) to [out=up, in=down] (1.8,2) to [out=up, in=down](2.1, 2.7);
\node [morphism] at (1.8,2) {$g_1$};
\draw (1.5,3.4) to [out=up, in=down] (1.8,4.1) to [out=up, in=down](1.5, 4.8);
\draw (2.1,3.4) to [out=up, in=down] (1.8,4.1) to [out=up, in=down](2.1, 4.8);
\node [morphism] at (1.8,4.1) {$g_s$};
\draw (0.6, 5.5) to (0.6, -2.9);
\draw (1.2, 5.5) to (1.2, -2.9);
\draw (-0.3, -2.2) to (-0.3, -2.9);
\draw (0.3, -2.2) to (0.3, -2.9);
\draw [dotted](1.5, 3.4) to (1.5, 2.7);
\draw [dotted](2.1, 3.4) to (2.1, 2.7);
\draw [dotted](-0.3, -0.1) to (-0.3, -0.8);
\draw [dotted](0.3, -0.1) to (0.3, -0.8);
\node [scale=0.6] at (1.8,-2.8) {$\ldots$};
\node [scale=0.6] at (1.8,1.4) {$\ldots$};
\node [scale=0.6] at (1.8,2.6) {$\ldots$};
\node [scale=0.6] at (1.8,3.5) {$\ldots$};
\node [scale=0.6] at (1.8,5.4) {$\ldots$};
\node [scale=0.6] at (0,0) {$\ldots$};
\node [scale=0.6] at (0,-0.8) {$\ldots$};
\node [scale=0.6] at (0,5.4) {$\ldots$};
\node [scale=0.6] at (0,-2.8) {$\ldots$};
\node [scale=0.6] at (0,1.2) {$\ldots$};
\node [scale=0.6] at (0.9,-2.8) {$\ldots$};
\node [scale=0.6] at (0.9,5.4) {$\ldots$};
\end{tikzpicture}
\end{aligned}
\end{align*}
\end{definition}

\noindent
Note that any signature supporting type I homotopy generators always supports type I homotopy composites; this is not extra structure.

\subsection{Type II homotopy generators} 

Homotopy generators of type II are defined in the following way. Note that these generators are always of dimension 4 or higher.
\def\pullthroughscale{0.8}
\def\pullthroughscalevertex{0.8}
\begin{definition}[Type II\ homotopy generator]
\label{defpullthrough}
For $n \geq 0$, given an $n$\-signature $\Sigma$ supporting type I homotopy generators, $\Sigma$ supports \emph{type~II homotopy generators} if for any $3 \leq k < n$, for $k$\-diagrams with 3\-projections given by the source diagrams below for $\alpha,\beta \in \Sigma_k$, there are chosen invertible $(k{+}1)$\-cells $\II_k$ as follows:
\begin{calign}
\label{eq:basicpullthroughs}
\begin{aligned}
\begin{tikzpicture}[thick, scale=\pullthroughscale]
\draw [region] (0.7,2.1) rectangle +(2.2,3);
\draw (2.6,2.1) to [out=up, in=down] (1.9,2.8) to [out=up, in=down] (1.2,3.5) ;
\draw (1.2,3.5) to (1.2,5.1);
\node [scale=0.4, rotate = 135, thick] at (0.8,5) {$\ldots$};
\node [scale=0.4, rotate = 135, thick] at (0.8,2) {$\ldots$};
\node [scale=0.4, rotate = 135, thick] at (3,5) {$\ldots$};
\node [scale=0.4, rotate = 135, thick] at (3,2) {$\ldots$};
\draw [region] (0.9, 1.9) rectangle +(2.2,3);
\draw (1.2,1.9) to (1.2,2.8) to [out=up, in=down] (1.9,3.5);
\draw (1.9,1.9) to (1.9,2.1)  to [out=up, in=down] (2.6,2.8) to (2.6,3.5);
\draw (1.9,3.5) to [out=up, in=down] (2.25,4.2) to [out=up, in=down](1.9, 4.9) to (1.9,4.9);
\draw (2.6,3.5) to [out=up, in=down] (2.25,4.2) to [out=up, in=down](2.6, 4.9)to (2.6,4.9);
\node [morphism, scale=\pullthroughscalevertex] at (2.25,4.2) {$\alpha$};
\node [scale=0.6] at (1.55,2) {$\ldots$};
\node [scale=0.6] at (2.25,3.6) {$\ldots$};
\node [scale=0.6] at (2.25,4.8) {$\ldots$};
\end{tikzpicture}
\end{aligned}
\stackrel{\II_k}{\rightarrow}
\begin{aligned}
\begin{tikzpicture}[thick, scale=\pullthroughscale]
\draw [region] (0,1.4) rectangle +(2.2,3);
\draw (1.9,2.8) to [out=up, in=down] (1.2,3.5) to [out=up, in=down] (0.5,4.2) to(0.5,4.4);
\draw (1.9, 1.4)to(1.9,2.8);
\node [scale=0.4, rotate = 135, thick] at (0.1,4.3) {$\ldots$};
\node [scale=0.4, rotate = 135, thick] at (0.1,1.3) {$\ldots$};
\node [scale=0.4, rotate = 135, thick] at (2.3,4.3) {$\ldots$};
\node [scale=0.4, rotate = 135, thick] at (2.3,1.3) {$\ldots$};
\draw [region] (0.2, 1.2) rectangle +(2.2,3);
\draw (0.5,2.8) to [out=up, in=down] (0.5,3.5) to [out=up, in=down] (1.2,4.2) ;
\draw (1.2,2.8)  to [out=up, in=down] (1.9,3.5) to [out=up, in=down] (1.9,4.2);
\draw (0.5,1.2) to(0.5,1.4) to [out=up, in=down] (0.85,2.1) to [out=up, in=down](0.5, 2.8);
\draw (1.2,1.2) to(1.2,1.4) to [out=up, in=down] (0.85,2.1) to [out=up, in=down](1.2, 2.8);
\node [morphism, scale=\pullthroughscalevertex] at (0.85,2.1) {$\alpha$};
\node [scale=0.6] at (0.85,1.3) {$\ldots$};
\node [scale=0.6] at (0.85,2.7) {$\ldots$};
\node [scale=0.6] at (1.55,4.1) {$\ldots$};
\end{tikzpicture}
\end{aligned}
&
\begin{aligned}
\begin{tikzpicture}[thick, scale=\pullthroughscale]
\draw [region] (0,2.1) rectangle +(2.2,3);
\draw (1.2,2.1)  to [out=up, in=down] (0.5,2.8) to (0.5,3.5);
\draw (1.9,2.1) to [out=up, in=down] (1.9,2.8) to [out=up, in=down] (1.2,3.5) ;
\draw (0.5,3.5) to [out=up, in=down] (0.85,4.2) to [out=up, in=down](0.5, 4.9)to (0.5,5.1);
\draw (1.2,3.5) to [out=up, in=down] (0.85,4.2) to [out=up, in=down](1.2, 4.9)to (1.2,5.1);
\node [morphism, scale=\pullthroughscalevertex] at (0.85,4.2) {$\alpha$};
\node [scale=0.6] at (1.55,2.2) {$\ldots$};
\node [scale=0.6] at (0.85,3.6) {$\ldots$};
\node [scale=0.6] at (0.85,5) {$\ldots$};
\node [scale=0.4, rotate = 135, thick] at (0.1,5) {$\ldots$};
\node [scale=0.4, rotate = 135, thick] at (0.1,2) {$\ldots$};
\node [scale=0.4, rotate = 135, thick] at (2.3,5) {$\ldots$};
\node [scale=0.4, rotate = 135, thick] at (2.3,2) {$\ldots$};
\draw [region] (0.2, 1.9) rectangle +(2.2,3);
\draw (0.5,1.9) to (0.5,2.1) to [out=up, in=down] (1.2,2.8) to [out=up, in=down] (1.9,3.5) to (1.9,4.9);
\end{tikzpicture}
\end{aligned}
\stackrel{\II_k}{\rightarrow}
\begin{aligned}
\begin{tikzpicture}[thick, scale=\pullthroughscale]
\draw [region] (0.7,-1.4) rectangle +(2.2,3);
\draw (1.9,0) to(1.9,0) to [out=up, in=down] (1.2,0.7) to [out=up, in=down] (1.2,1.4)to(1.2,1.6) ;
\draw (2.6,0)  to  [out=up, in=down] (2.6,0.7) to[out=up, in=down](1.9,1.4)to(1.9,1.6);
\draw (1.9,-1.4) to(1.9,-1.4) to [out=up, in=down] (2.25,-0.7) to [out=up, in=down](1.9, 0);
\draw (2.6,-1.4) to(2.6,-1.4) to [out=up, in=down] (2.25,-0.7) to [out=up, in=down](2.6, 0);
\node [morphism, scale=\pullthroughscalevertex] at (2.25,-0.7) {$\alpha$};
\node [scale=0.6] at (1.55,1.5) {$\ldots$};
\node [scale=0.6] at (2.25,-0.1) {$\ldots$};
\node [scale=0.4, rotate = 135, thick] at (0.8,1.5) {$\ldots$};
\node [scale=0.4, rotate = 135, thick] at (0.8,-1.5) {$\ldots$};
\node [scale=0.4, rotate = 135, thick] at (3,1.5) {$\ldots$};
\node [scale=0.4, rotate = 135, thick] at (3,-1.5) {$\ldots$};
\node [scale=0.6] at (2.25,-1.3) {$\ldots$};
\draw [region] (0.9, -1.6) rectangle +(2.2,3);
\draw (1.2,-1.6) to(1.2,0.0);
\draw (1.2,0.0)  to [out=up, in=down] (1.9,0.7) to [out=up, in=down] (2.6,1.4);
\end{tikzpicture}
\end{aligned}
\end{calign}
\end{definition}

\noindent
Note that these diagrams involve homotopy generators of type I and their expansion scheme. As per the convention set out above, these diagrams may have arbitrary interposing sheets, but we omit them for notational clarity. However, it is deliberate that there are no additional sheets in front or behind the displayed diagrams.

\def\expansionpullthroughscale{0.65}
\def\expansionpullthroughscaleprime{0.7}

\paragraph{Expansion scheme.}
We require expansion schemes for type II homotopy generators, defining their action on composite diagrams.
\begin{definition}
\label{pullthroughexpansion}
In an $n$\-signature that supports type~II homotopy generators, a \emph{type II homotopy composite} is a diagram formed from a sequence of type~II homotopy generators, denoted $\widetilde \II$, with values defined recursively as follows.
\begin{align*}
\begin{aligned}
\begin{tikzpicture}[thick, scale=\expansionpullthroughscale]
\draw [region] (-0.7,-2.1) rectangle +(3.6,7.9);
\draw (-0.2,-2.1)  to (-0.2,-1.4)  to [out=up, in=down](0.5,-0.7)  to [out=up, in=down](1.2,0.0)  to [out=up, in=down] (1.9,0.7) to [out=up, in=down] (2.6,1.4) to (2.6,5.8);
\node [scale=0.4, rotate = 135, thick] at (3,-2.2) {$\ldots$};
\draw [region] (-0.5, -2.3) rectangle +(3.6,7.9);
\draw (1.9,-2.3) to(1.9,0) to [out=up, in=down] (1.2,0.7) to [out=up, in=down] (1.2,1.4) ;
\draw (2.6,-2.3) to  (2.6,0.7) to [out=up, in=down] (1.9,1.4) to (1.9,5.6);
\draw (1.2,-2.3) to (1.2,-0.7) to [out=up, in=down]  (0.5,0) to (0.5,1.4);
\draw (0.5,-2.3) to (0.5,-1.4) to [out=up, in=down]  (-0.2,-0.7) to (-0.2,5.6);
\draw [dotted](1.2,2.8)  to (1.2,3.5);
\draw [dotted](0.5,2.8)  to (0.5,3.5);
\draw (0.5,3.5) to [out=up, in=down] (0.85,4.2) to [out=up, in=down](0.5, 4.9)to (0.5,5.6);
\draw (1.2,3.5) to [out=up, in=down] (0.85,4.2) to [out=up, in=down](1.2, 4.9)to (1.2,5.6);
\node [morphism, scale=\lsscalevertex] at (0.85,4.2) {$\alpha_m$};
\draw (0.5,1.4) to [out=up, in=down] (0.85,2.1) to [out=up, in=down](0.5, 2.8);
\draw (1.2,1.4) to [out=up, in=down] (0.85,2.1) to [out=up, in=down](1.2, 2.8);
\node [morphism, scale=\lsscalevertex] at (0.85,2.1) {$\alpha_1$};
\node [scale=0.6] at (1.55,-2.2) {$\ldots$};
\node [scale=0.6] at (0.85,-2.2) {$\ldots$};
\node [scale=0.6] at (2.25,-2.2) {$\ldots$};
\node [scale=0.6] at (0.85,3.6) {$\ldots$};
\node [scale=0.6] at (0.85,2.7) {$\ldots$};
\node [scale=0.6] at (0.85,1.5) {$\ldots$};
\node [scale=0.6] at (0.15,1.5) {$\ldots$};
\node [scale=0.6] at (0.15,2.7) {$\ldots$};
\node [scale=0.6] at (0.15,3.6) {$\ldots$};
\node [scale=0.6] at (0.15,5.5) {$\ldots$};
\node [scale=0.6] at (1.55,1.5) {$\ldots$};
\node [scale=0.6] at (1.55,2.7) {$\ldots$};
\node [scale=0.6] at (1.55,3.6) {$\ldots$};
\node [scale=0.6] at (0.85,5.5) {$\ldots$};
\node [scale=0.6] at (1.55,5.5) {$\ldots$};
\node [scale=0.4, rotate = 135, thick] at (-0.6,5.7) {$\ldots$};
\node [scale=0.4, rotate = 135, thick] at (-0.6,-2.2) {$\ldots$};
\node [scale=0.4, rotate = 135, thick] at (3,5.7) {$\ldots$};
\end{tikzpicture}
\end{aligned}
\stackrel{\widetilde{\II}_k}{\rightarrow}
\begin{aligned}
\begin{tikzpicture}[thick, scale=\expansionpullthroughscale]
\draw [region] (-1.4,0.7) rectangle +(3.6,7.9);
\draw (-0.9,0.7)  to (-0.9,4.9)  to [out=up, in=down](-0.2,5.6) to [out=up, in=down](0.5,6.3) to [out=up, in=down](1.2,7) to [out=up, in=down](1.9,7.7)to [out=up, in=down](1.9,7.9)to(1.9,8.6) ;
\node [scale=0.4, rotate = 135, thick] at (2.3,0.6) {$\ldots$};
\draw [region] (-1.2, 0.5) rectangle +(3.6,7.9);
\draw (1.9,0.5) to(1.9,1.4) to (1.9,4.9);
\draw (-0.2,0.5) to (-0.2,4.9);
\draw [dotted](1.2,2.8)  to (1.2,3.5);
\draw [dotted](0.5,2.8)  to (0.5,3.5);
\draw (0.5,3.5) to [out=up, in=down] (0.85,4.2) to [out=up, in=down](0.5, 4.9);
\draw (1.2,3.5) to [out=up, in=down] (0.85,4.2) to [out=up, in=down](1.2, 4.9);
\node [morphism, scale=\lsscalevertex] at (0.85,4.2) {$\alpha_m$};
\draw (0.5,0.5) to(0.5,1.4) to [out=up, in=down] (0.85,2.1) to [out=up, in=down](0.5, 2.8);
\draw (1.2,0.5) to (1.2,1.4) to [out=up, in=down] (0.85,2.1) to [out=up, in=down](1.2, 2.8);
\draw (-0.2,4.9)  to [out=up, in=down] (-0.9,5.6)  to [out=up, in=down](-0.9,7.7) to (-0.9,8.4);
\draw (0.5,4.9)  to [out=up, in=down] (0.5,5.6)  to [out=up, in=down](-0.2,6.3)  to [out=up, in=down](-0.2,8.4);
\draw (1.2,4.9)  to [out=up, in=down] (1.2,6.3)  to [out=up, in=down](0.5,7)  to [out=up, in=down](0.5,8.4);
\draw (1.9,4.9)  to [out=up, in=down] (1.9,7)  to [out=up, in=down](1.2,7.7) to (1.2,8.4);
\node [morphism, scale=\lsscalevertex] at (0.85,2.1) {$\alpha_1$};
\node [scale=0.6] at (0.15,8.3) {$\ldots$};
\node [scale=0.6] at (-0.55,8.3) {$\ldots$};
\node [scale=0.6] at (0.85,8.3) {$\ldots$};
\node [scale=0.6] at (0.85,3.6) {$\ldots$};
\node [scale=0.6] at (0.85,2.7) {$\ldots$};
\node [scale=0.6] at (0.85,0.6) {$\ldots$};
\node [scale=0.6] at (0.15,0.6) {$\ldots$};
\node [scale=0.6] at (0.15,2.7) {$\ldots$};
\node [scale=0.6] at (0.15,3.6) {$\ldots$};
\node [scale=0.6] at (0.15,4.8) {$\ldots$};
\node [scale=0.6] at (1.55,0.6) {$\ldots$};
\node [scale=0.6] at (1.55,2.7) {$\ldots$};
\node [scale=0.6] at (1.55,3.6) {$\ldots$};
\node [scale=0.6] at (0.85,4.8) {$\ldots$};
\node [scale=0.6] at (1.55,4.8) {$\ldots$};
\node [scale=0.4, rotate = 135, thick] at (-1.3,8.5) {$\ldots$};
\node [scale=0.4, rotate = 135, thick] at (-1.3,0.6) {$\ldots$};
\node [scale=0.4, rotate = 135, thick] at (2.3,8.5) {$\ldots$};
\end{tikzpicture}
\end{aligned}
:=&
\begin{aligned}
\begin{tikzpicture}[thick, scale=\expansionpullthroughscale]
\draw [region] (-0.7,-2.1) rectangle +(3.6,7.9);
\draw (-0.2,-2.1)  to (-0.2,-1.4)  to [out=up, in=down](0.5,-0.7)  to [out=up, in=down](1.2,0.0)  to [out=up, in=down] (1.9,0.7) to [out=up, in=down] (2.6,1.4) to (2.6,5.8);
\node [scale=0.4, rotate = 135, thick] at (3,-2.2) {$\ldots$};
\draw [region] (-0.5, -2.3) rectangle +(3.6,7.9);
\draw (1.9,-2.3) to(1.9,0) to [out=up, in=down] (1.2,0.7) to [out=up, in=down] (1.2,1.4) ;
\draw (2.6,-2.3) to  (2.6,0.7) to [out=up, in=down] (1.9,1.4) to (1.9,5.6);
\draw (1.2,-2.3) to (1.2,-0.7) to [out=up, in=down]  (0.5,0) to (0.5,1.4);
\draw (0.5,-2.3) to (0.5,-1.4) to [out=up, in=down]  (-0.2,-0.7) to (-0.2,5.6);
\draw [dotted](1.2,2.8)  to (1.2,3.5);
\draw [dotted](0.5,2.8)  to (0.5,3.5);
\draw (0.5,3.5) to [out=up, in=down] (0.85,4.2) to [out=up, in=down](0.5, 4.9)to (0.5,5.6);
\draw (1.2,3.5) to [out=up, in=down] (0.85,4.2) to [out=up, in=down](1.2, 4.9)to (1.2,5.6);
\node [morphism, scale=\lsscalevertex] at (0.85,4.2) {$\alpha_m$};
\draw (0.5,1.4) to [out=up, in=down] (0.85,2.1) to [out=up, in=down](0.5, 2.8);
\draw (1.2,1.4) to [out=up, in=down] (0.85,2.1) to [out=up, in=down](1.2, 2.8);
\node [morphism, scale=\lsscalevertex] at (0.85,2.1) {$\alpha_1$};
\node [scale=0.6] at (1.55,-2.2) {$\ldots$};
\node [scale=0.6] at (0.85,-2.2) {$\ldots$};
\node [scale=0.6] at (2.25,-2.2) {$\ldots$};
\node [scale=0.6] at (0.85,3.6) {$\ldots$};
\node [scale=0.6] at (0.85,2.7) {$\ldots$};
\node [scale=0.6] at (0.85,1.5) {$\ldots$};
\node [scale=0.6] at (0.15,1.5) {$\ldots$};
\node [scale=0.6] at (0.15,2.7) {$\ldots$};
\node [scale=0.6] at (0.15,3.6) {$\ldots$};
\node [scale=0.6] at (0.15,5.5) {$\ldots$};
\node [scale=0.6] at (1.55,1.5) {$\ldots$};
\node [scale=0.6] at (1.55,2.7) {$\ldots$};
\node [scale=0.6] at (1.55,3.6) {$\ldots$};
\node [scale=0.6] at (0.85,5.5) {$\ldots$};
\node [scale=0.6] at (1.55,5.5) {$\ldots$};
\node [scale=0.4, rotate = 135, thick] at (-0.6,5.7) {$\ldots$};
\node [scale=0.4, rotate = 135, thick] at (-0.6,-2.2) {$\ldots$};
\node [scale=0.4, rotate = 135, thick] at (3,5.7) {$\ldots$};
\end{tikzpicture}
\end{aligned}
\stackrel{\II_k}{\rightarrow}
\begin{aligned}
\begin{tikzpicture}[thick, scale=\expansionpullthroughscale]
\draw [region] (-0.7,-2.1) rectangle +(3.6,7.9);
\draw (-0.2,-2.1)  to(-0.2,-1.4)  to [out=up, in=down](0.5,-0.7)  to [out=up, in=down](1.2,0.0)  to [out=up, in=down] (1.9,0.7)to (1.9,2.1) to [out=up, in=down] (2.6,2.8) to (2.6,5.8);
\node [scale=0.4, rotate = 135, thick] at (3,-2.2) {$\ldots$};
\draw [region] (-0.5, -2.3) rectangle +(3.6,7.9);
\draw (1.9,-2.3) to(1.9,0) to [out=up, in=down] (1.2,0.7) to [out=up, in=down] (1.2,0.7) ;
\draw (2.6,-2.3) to  (2.6,2.1) to [out=up, in=down] (1.9,2.8) to (1.9,5.6);
\draw (1.2,-2.3) to (1.2,-0.7) to [out=up, in=down]  (0.5,0) to (0.5,0.7);
\draw (0.5,-2.3) to (0.5,-1.4) to [out=up, in=down]  (-0.2,-0.7) to (-0.2,5.6);
\draw [dotted](1.2,2.8)  to (1.2,3.5);
\draw [dotted](0.5,2.8)  to (0.5,3.5);
\draw (0.5,3.5) to [out=up, in=down] (0.85,4.2) to [out=up, in=down](0.5, 4.9)to (0.5,5.6);
\draw (1.2,3.5) to [out=up, in=down] (0.85,4.2) to [out=up, in=down](1.2, 4.9)to (1.2,5.6);
\node [morphism, scale=\lsscalevertex] at (0.85,4.2) {$\alpha_m$};
\draw (0.5,0.7) to [out=up, in=down] (0.85,1.4) to [out=up, in=down](0.5, 2.1)to (0.5,2.8);
\draw (1.2,0.7) to [out=up, in=down] (0.85,1.4) to [out=up, in=down](1.2, 2.1)to (1.2,2.8);
\node [morphism, scale=\lsscalevertex] at (0.85,1.4) {$\alpha_1$};
\node [scale=0.6] at (1.55,-2.2) {$\ldots$};
\node [scale=0.6] at (0.85,-2.2) {$\ldots$};
\node [scale=0.6] at (2.25,-2.2) {$\ldots$};
\node [scale=0.6] at (0.85,3.6) {$\ldots$};
\node [scale=0.6] at (0.85,2.7) {$\ldots$};
\node [scale=0.6] at (0.85,0.8) {$\ldots$};
\node [scale=0.6] at (0.15,2.7) {$\ldots$};
\node [scale=0.6] at (0.15,3.6) {$\ldots$};
\node [scale=0.6] at (0.15,5.5) {$\ldots$};
\node [scale=0.6] at (1.55,2.7) {$\ldots$};
\node [scale=0.6] at (1.55,3.6) {$\ldots$};
\node [scale=0.6] at (0.85,5.5) {$\ldots$};
\node [scale=0.6] at (1.55,5.5) {$\ldots$};
\node [scale=0.4, rotate = 135, thick] at (-0.6,5.7) {$\ldots$};
\node [scale=0.4, rotate = 135, thick] at (-0.6,-2.2) {$\ldots$};
\node [scale=0.4, rotate = 135, thick] at (3,5.7) {$\ldots$};
\end{tikzpicture}
\end{aligned}
\stackrel{\widetilde{\II}_k}{\rightarrow}
\begin{aligned}
\begin{tikzpicture}[thick, scale=\expansionpullthroughscale]
\draw [region] (-0.7,-2.8) rectangle +(3.6,7.9);
\draw (-0.2,-2.8)   to [out=up, in=down](0.5,-2.1) to (0.5,-0.7) to [out=up, in=down](1.2,0.0)  to [out=up, in=down] (1.9,0.7) to [out=up, in=down] (2.6,1.4) to (2.6,5.1);
\node [scale=0.4, rotate = 135, thick] at (3,-2.9) {$\ldots$};
\draw [region] (-0.5, -3) rectangle +(3.6,7.9);
\draw (1.9,-1.4) to(1.9,0) to [out=up, in=down] (1.2,0.7) to [out=up, in=down] (1.2,7.0) ;
\draw (2.6,-3) to  (2.6,0.7) to [out=up, in=down] (1.9,1.4) to (1.9,4.9);
\draw (1.2,-1.4) to (1.2,-0.7) to [out=up, in=down]  (0.5,0) to (0.5,0.7);
\draw (0.5,-3) to (0.5,-2.8) to [out=up, in=down]  (-0.2,-2.1) to (-0.2,4.9);
\draw [dotted](1.2,2.8)  to (1.2,3.5);
\draw [dotted](0.5,2.8)  to (0.5,3.5);
\draw (0.5,3.5) to [out=up, in=down] (0.85,4.2) to [out=up, in=down](0.5, 4.9);
\draw (1.2,3.5) to [out=up, in=down] (0.85,4.2) to [out=up, in=down](1.2, 4.9);
\node [morphism, scale=\lsscalevertex] at (0.85,4.2) {$\alpha_m$};
\draw (0.5,1.4) to [out=up, in=down] (0.85,2.1) to [out=up, in=down](0.5, 2.8);
\draw (1.2,1.4) to [out=up, in=down] (0.85,2.1) to [out=up, in=down](1.2, 2.8);
\node [morphism, scale=\lsscalevertex] at (0.85,2.1) {$\alpha_1$};
\draw (1.2,-0.7) to [out=down, in=up] (1.55,-1.4) to [out=down, in=up](1.2, -2.1);
\draw (1.9,-0.7) to [out=down, in=up] (1.55,-1.4) to [out=down, in=up](1.9, -2.1);
\node [morphism, scale=\lsscalevertex] at (1.55,-1.4) {$\alpha_1$};
\draw (1.2,-3) to  (1.2,-2.1);
\draw (1.9,-3) to  (1.9,-2.1);
\node [scale=0.6] at (1.55,-2.9) {$\ldots$};
\node [scale=0.6] at (0.85,-2.9) {$\ldots$};
\node [scale=0.6] at (2.25,-2.9) {$\ldots$};
\node [scale=0.6] at (1.55,-0.8) {$\ldots$};
\node [scale=0.6] at (0.85,3.6) {$\ldots$};
\node [scale=0.6] at (0.85,2.7) {$\ldots$};
\node [scale=0.6] at (0.85,1.5) {$\ldots$};
\node [scale=0.6] at (0.15,1.5) {$\ldots$};
\node [scale=0.6] at (0.15,2.7) {$\ldots$};
\node [scale=0.6] at (0.15,3.6) {$\ldots$};
\node [scale=0.6] at (0.15,4.8) {$\ldots$};
\node [scale=0.6] at (1.55,1.5) {$\ldots$};
\node [scale=0.6] at (1.55,2.7) {$\ldots$};
\node [scale=0.6] at (1.55,3.6) {$\ldots$};
\node [scale=0.6] at (0.85,4.8) {$\ldots$};
\node [scale=0.6] at (1.55,4.8) {$\ldots$};
\node [scale=0.4, rotate = 135, thick] at (-0.6,5) {$\ldots$};
\node [scale=0.4, rotate = 135, thick] at (-0.6,-2.9) {$\ldots$};
\node [scale=0.4, rotate = 135, thick] at (3,5) {$\ldots$};
\end{tikzpicture}
\end{aligned}
\\
&\stackrel{\widetilde{\II}_k}{\rightarrow}
\begin{aligned}
\begin{tikzpicture}[thick, scale=\expansionpullthroughscale]
\draw [region] (-0.7,-2.8) rectangle +(3.6,7.9);
\draw (-0.2,-2.8)  to (-0.2,-1.4)  to [out=up, in=down](0.5,-0.7)  to [out=up, in=down](1.2,0.0)  to [out=up, in=down] (1.9,0.7) to [out=up, in=down] (2.6,1.4) to (2.6,5.1);
\node [scale=0.4, rotate = 135, thick] at (3,-2.9) {$\ldots$};
\draw [region] (-0.5, -3) rectangle +(3.6,7.9);
\draw (1.9,-1.4) to(1.9,0) to [out=up, in=down] (1.2,0.7) to [out=up, in=down] (1.2,1.4) ;
\draw (2.6,-3) to  (2.6,0.7) to [out=up, in=down] (1.9,1.4) to (1.9,4.9);
\draw (1.2,-1.4) to (1.2,-0.7) to [out=up, in=down]  (0.5,0) to (0.5,1.4);
\draw (0.5,-3) to (0.5,-1.4) to [out=up, in=down]  (-0.2,-0.7) to (-0.2,4.9);
\draw [dotted](1.2,2.8)  to (1.2,3.5);
\draw [dotted](0.5,2.8)  to (0.5,3.5);
\draw (0.5,3.5) to [out=up, in=down] (0.85,4.2) to [out=up, in=down](0.5, 4.9);
\draw (1.2,3.5) to [out=up, in=down] (0.85,4.2) to [out=up, in=down](1.2, 4.9);
\node [morphism, scale=\lsscalevertex] at (0.85,4.2) {$\alpha_m$};
\draw (0.5,1.4) to [out=up, in=down] (0.85,2.1) to [out=up, in=down](0.5, 2.8);
\draw (1.2,1.4) to [out=up, in=down] (0.85,2.1) to [out=up, in=down](1.2, 2.8);
\node [morphism, scale=\lsscalevertex] at (0.85,2.1) {$\alpha_2$};
\draw (1.2,-1.4) to [out=down, in=up] (1.55,-2.1) to [out=down, in=up](1.2, -2.8);
\draw (1.9,-1.4) to [out=down, in=up] (1.55,-2.1) to [out=down, in=up](1.9, -2.8);
\node [morphism, scale=\lsscalevertex] at (1.55,-2.1) {$\alpha_1$};
\draw (1.2,-3) to  (1.2,-2.8);
\draw (1.9,-3) to  (1.9,-2.8);
\node [scale=0.6] at (1.55,-2.9) {$\ldots$};
\node [scale=0.6] at (0.85,-2.9) {$\ldots$};
\node [scale=0.6] at (2.25,-2.9) {$\ldots$};
\node [scale=0.6] at (1.55,-1.5) {$\ldots$};
\node [scale=0.6] at (0.85,-1.5) {$\ldots$};
\node [scale=0.6] at (2.25,-1.5) {$\ldots$};
\node [scale=0.6] at (0.85,3.6) {$\ldots$};
\node [scale=0.6] at (0.85,2.7) {$\ldots$};
\node [scale=0.6] at (0.85,1.5) {$\ldots$};
\node [scale=0.6] at (0.15,1.5) {$\ldots$};
\node [scale=0.6] at (0.15,2.7) {$\ldots$};
\node [scale=0.6] at (0.15,3.6) {$\ldots$};
\node [scale=0.6] at (0.15,4.8) {$\ldots$};
\node [scale=0.6] at (1.55,1.5) {$\ldots$};
\node [scale=0.6] at (1.55,2.7) {$\ldots$};
\node [scale=0.6] at (1.55,3.6) {$\ldots$};
\node [scale=0.6] at (0.85,4.8) {$\ldots$};
\node [scale=0.6] at (1.55,4.8) {$\ldots$};
\node [scale=0.4, rotate = 135, thick] at (-0.6,5) {$\ldots$};
\node [scale=0.4, rotate = 135, thick] at (-0.6,-2.9) {$\ldots$};
\node [scale=0.4, rotate = 135, thick] at (3,5) {$\ldots$};
\end{tikzpicture}
\end{aligned}
\stackrel{\widetilde{\II}_k}{\rightarrow}
\begin{aligned}
\begin{tikzpicture}[thick, scale=\expansionpullthroughscale]
\draw [region] (-1.4,0.7) rectangle +(3.6,7.9);
\draw (-0.9,0.7)  to (-0.9,4.9)  to [out=up, in=down](-0.2,5.6) to [out=up, in=down](0.5,6.3) to [out=up, in=down](1.2,7) to [out=up, in=down](1.9,7.7)to [out=up, in=down](1.9,7.9)to(1.9,8.6) ;
\node [scale=0.4, rotate = 135, thick] at (2.3,0.6) {$\ldots$};
\draw [region] (-1.2, 0.5) rectangle +(3.6,7.9);
\draw (1.9,0.5) to(1.9,1.4) to (1.9,4.9);
\draw (-0.2,0.5) to (-0.2,4.9);
\draw [dotted](1.2,2.8)  to (1.2,3.5);
\draw [dotted](0.5,2.8)  to (0.5,3.5);
\draw (0.5,3.5) to [out=up, in=down] (0.85,4.2) to [out=up, in=down](0.5, 4.9);
\draw (1.2,3.5) to [out=up, in=down] (0.85,4.2) to [out=up, in=down](1.2, 4.9);
\node [morphism, scale=\lsscalevertex] at (0.85,4.2) {$\alpha_m$};
\draw (0.5,0.5) to(0.5,1.4) to [out=up, in=down] (0.85,2.1) to [out=up, in=down](0.5, 2.8);
\draw (1.2,0.5) to (1.2,1.4) to [out=up, in=down] (0.85,2.1) to [out=up, in=down](1.2, 2.8);
\draw (-0.2,4.9)  to [out=up, in=down] (-0.9,5.6)  to [out=up, in=down](-0.9,7.7) to (-0.9,8.4);
\draw (0.5,4.9)  to [out=up, in=down] (0.5,5.6)  to [out=up, in=down](-0.2,6.3)  to [out=up, in=down](-0.2,8.4);
\draw (1.2,4.9)  to [out=up, in=down] (1.2,6.3)  to [out=up, in=down](0.5,7)  to [out=up, in=down](0.5,8.4);
\draw (1.9,4.9)  to [out=up, in=down] (1.9,7)  to [out=up, in=down](1.2,7.7) to (1.2,8.4);
\node [morphism, scale=\lsscalevertex] at (0.85,2.1) {$\alpha_1$};
\node [scale=0.6] at (0.15,8.3) {$\ldots$};
\node [scale=0.6] at (-0.55,8.3) {$\ldots$};
\node [scale=0.6] at (0.85,8.3) {$\ldots$};
\node [scale=0.6] at (0.85,3.6) {$\ldots$};
\node [scale=0.6] at (0.85,2.7) {$\ldots$};
\node [scale=0.6] at (0.85,0.6) {$\ldots$};
\node [scale=0.6] at (0.15,0.6) {$\ldots$};
\node [scale=0.6] at (0.15,2.7) {$\ldots$};
\node [scale=0.6] at (0.15,3.6) {$\ldots$};
\node [scale=0.6] at (0.15,4.8) {$\ldots$};
\node [scale=0.6] at (1.55,0.6) {$\ldots$};
\node [scale=0.6] at (1.55,2.7) {$\ldots$};
\node [scale=0.6] at (1.55,3.6) {$\ldots$};
\node [scale=0.6] at (0.85,4.8) {$\ldots$};
\node [scale=0.6] at (1.55,4.8) {$\ldots$};
\node [scale=0.4, rotate = 135, thick] at (-1.3,8.5) {$\ldots$};
\node [scale=0.4, rotate = 135, thick] at (-1.3,0.6) {$\ldots$};
\node [scale=0.4, rotate = 135, thick] at (2.3,8.5) {$\ldots$};
\end{tikzpicture}
\end{aligned}
\\
\begin{aligned}
\begin{tikzpicture}[thick, scale=\expansionpullthroughscaleprime]
\draw [region] (0,-0.7) rectangle +(2.9,6.5);
\draw (0.5,-0.7) to  (0.5,1.4) ;
\draw (1.2,-0.7) to (1.2,0.0)  to [out=up, in=down] (1.9,0.7) to [out=up, in=down] (2.6,1.4);
\draw [dotted](2.6,1.4)  to (2.6,2.1);
\draw [dotted](0.5,1.4)  to (0.5,2.1);
\draw (0.5,2.1) to [out=up, in=down](1.2,2.8) to [out=up, in=down] (1.9,3.5) to (1.9,5.8);
\draw (2.6,2.1) to (2.6,5.8);
\node [scale=0.6] at (0.85,-0.6) {$\ldots$};
\node [scale=0.6] at (2.25,5.7) {$\ldots$};
\node [scale=0.4, rotate = 135, thick] at (3,-0.8) {$\ldots$};
\draw [region] (0.2, -0.9) rectangle +(2.9,6.5);
\draw (1.9,-0.9) to(1.9,0) to [out=up, in=down] (1.2,0.7) to [out=up, in=down] (1.2,1.4) ;
\draw (2.6,-0.9) to  (2.6,0.7) to [out=up, in=down] (1.9,1.4);
\draw [dotted](1.2,1.4)  to (1.2,2.1);
\draw [dotted](1.9,1.4)  to (1.9,2.1);
\draw (1.2,2.1)  to [out=up, in=down] (0.5,2.8) to (0.5,3.5);
\draw (1.9,2.1) to [out=up, in=down] (1.9,2.8) to [out=up, in=down] (1.2,3.5) ;
\draw (0.5,3.5) to [out=up, in=down] (0.85,4.2) to [out=up, in=down](0.5, 4.9)to (0.5,5.6);
\draw (1.2,3.5) to [out=up, in=down] (0.85,4.2) to [out=up, in=down](1.2, 4.9)to (1.2,5.6);
\node [morphism, scale=\pullthroughscalevertex] at (0.85,4.2) {$\alpha$};
\node [scale=0.6] at (2.25,-0.8) {$\ldots$};
\node [scale=0.6] at (0.85,3.6) {$\ldots$};
\node [scale=0.6] at (0.85,5.5) {$\ldots$};
\node [scale=0.4, rotate = 135, thick] at (0.1,5.7) {$\ldots$};
\node [scale=0.4, rotate = 135, thick] at (0.1,-0.8) {$\ldots$};
\node [scale=0.4, rotate = 135, thick] at (3,5.7) {$\ldots$};
\end{tikzpicture}
\end{aligned}
\stackrel{\widetilde{\II}_k}{\rightarrow}
\begin{aligned}
\begin{tikzpicture}[thick, scale=\expansionpullthroughscaleprime]
\draw [region] (0,-2.1) rectangle +(2.9,6.5);
\draw (0.5,-2.1) to(0.5,0);
\draw (1.2,-2.1) to(1.2,0.0);
\draw (0.5,0) to (0.5,0.7) to [out=up, in=down] (0.5,1.4) ;
\draw (1.2,0.0)  to [out=up, in=down] (1.9,0.7) to [out=up, in=down] (2.6,1.4);
\draw [dotted](2.6,1.4)  to (2.6,2.1);
\draw [dotted](0.5,1.4)  to (0.5,2.1);
\draw (0.5,2.1) to [out=up, in=down] (1.2,2.8) to [out=up, in=down] (1.9,3.5) to (1.9,4.4);
\draw (2.6,2.1)  to [out=up, in=down] (2.6,2.8) to (2.6,4.4);
\node [scale=0.6] at (0.85,-2) {$\ldots$};
\node [scale=0.6] at (2.25,4.3) {$\ldots$};
\node [scale=0.4, rotate = 135, thick] at (3,-2.2) {$\ldots$};
\draw [region] (0.2, -2.3) rectangle +(2.9,6.5);
\draw (1.9,0) to(1.9,0) to [out=up, in=down] (1.2,0.7) to [out=up, in=down] (1.2,1.4) ;
\draw (2.6,0)  to  [out=up, in=down] (2.6,0.7) to[out=up, in=down](1.9,1.4);
\draw [dotted](1.2,1.4)  to (1.2,2.1);
\draw [dotted](1.9,1.4)  to (1.9,2.1);
\draw (1.2,2.1)  to [out=up, in=down] (0.5,2.8) to (0.5,4.2);
\draw (1.9,2.1) to [out=up, in=down] (1.9,2.8) to [out=up, in=down] (1.2,3.5) to (1.2,4.2);
\draw (1.9,-2.3) to(1.9,-1.4) to [out=up, in=down] (2.25,-0.7) to [out=up, in=down](1.9, 0);
\draw (2.6,-2.3) to(2.6,-1.4) to [out=up, in=down] (2.25,-0.7) to [out=up, in=down](2.6, 0);
\node [morphism, scale=\pullthroughscalevertex] at (2.25,-0.7) {$\alpha$};
\node [scale=0.6] at (0.85,4.1) {$\ldots$};
\node [scale=0.6] at (2.25,-0.1) {$\ldots$};
\node [scale=0.6] at (2.25,-2.2) {$\ldots$};
\node [scale=0.4, rotate = 135, thick] at (0.1,4.3) {$\ldots$};
\node [scale=0.4, rotate = 135, thick] at (0.1,-2.2) {$\ldots$};
\node [scale=0.4, rotate = 135, thick] at (3,4.3) {$\ldots$};
\end{tikzpicture}
\end{aligned}
:=&
\begin{aligned}
\begin{tikzpicture}[thick, scale=\expansionpullthroughscaleprime]
\draw [region] (0,-0.7) rectangle +(2.9,6.5);
\draw (0.5,-0.7) to  (0.5,1.4) ;
\draw (1.2,-0.7) to (1.2,0.0)  to [out=up, in=down] (1.9,0.7) to [out=up, in=down] (2.6,1.4);
\draw [dotted](2.6,1.4)  to (2.6,2.1);
\draw [dotted](0.5,1.4)  to (0.5,2.1);
\draw (0.5,2.1) to [out=up, in=down](1.2,2.8) to [out=up, in=down] (1.9,3.5) to (1.9,5.8);
\draw (2.6,2.1) to (2.6,5.8);
\node [scale=0.6] at (0.85,-0.6) {$\ldots$};
\node [scale=0.6] at (2.25,5.7) {$\ldots$};
\node [scale=0.4, rotate = 135, thick] at (3,-0.8) {$\ldots$};
\draw [region] (0.2, -0.9) rectangle +(2.9,6.5);
\draw (1.9,-0.9) to(1.9,0) to [out=up, in=down] (1.2,0.7) to [out=up, in=down]
(1.2,1.4) ;
\draw (2.6,-0.9) to  (2.6,0.7) to [out=up, in=down] (1.9,1.4);
\draw [dotted](1.2,1.4)  to (1.2,2.1);
\draw [dotted](1.9,1.4)  to (1.9,2.1);
\draw (1.2,2.1)  to [out=up, in=down] (0.5,2.8) to (0.5,3.5);
\draw (1.9,2.1) to [out=up, in=down] (1.9,2.8) to [out=up, in=down] (1.2,3.5) ;
\draw (0.5,3.5) to [out=up, in=down] (0.85,4.2) to [out=up, in=down](0.5, 4.9)to (0.5,5.6);
\draw (1.2,3.5) to [out=up, in=down] (0.85,4.2) to [out=up, in=down](1.2, 4.9)to (1.2,5.6);
\node [morphism, scale=\pullthroughscalevertex] at (0.85,4.2) {$\alpha$};
\node [scale=0.6] at (2.25,-0.8) {$\ldots$};
\node [scale=0.6] at (0.85,3.6) {$\ldots$};
\node [scale=0.6] at (0.85,5.5) {$\ldots$};
\node [scale=0.4, rotate = 135, thick] at (0.1,5.7) {$\ldots$};
\node [scale=0.4, rotate = 135, thick] at (0.1,-0.8) {$\ldots$};
\node [scale=0.4, rotate = 135, thick] at (3,5.7) {$\ldots$};
\end{tikzpicture}
\end{aligned}
\stackrel{\II_k}{\rightarrow}
\begin{aligned}
\begin{tikzpicture}[thick, scale=\expansionpullthroughscaleprime]
\draw [region] (-0.7,0) rectangle +(3.6,6.5);
\draw (0.5,0) to  (0.5,1.4) ;
\draw (1.2,0.0)  to [out=up, in=down] (1.9,0.7) to [out=up, in=down] (2.6,1.4);
\draw [dotted](2.6,1.4)  to (2.6,2.1);
\draw [dotted](0.5,1.4)  to (0.5,2.1);
\draw (0.5,2.1) to [out=up, in=down](1.2,2.8) to [out=up, in=down] (1.9,3.5) to (1.9,6.5);
\draw (-0.2,-0.2) to (-0.2,4.9) to [out=up, in=down](0.5,5.6) to [out=up, in=down] (1.2,6.3) to (1.2,6.5);
\draw (2.6,2.1) to (2.6,6.5);
\node [scale=0.6] at (0.85,0.1) {$\ldots$};
\node [scale=0.6] at (2.25,6.4) {$\ldots$};
\node [scale=0.4, rotate = 135, thick] at (3,-0.1) {$\ldots$};
\draw [region] (-0.5, -0.2) rectangle +(3.6,6.5);
\draw (1.9,-0.2) to(1.9,0) to [out=up, in=down] (1.2,0.7) to [out=up, in=down] (1.2,1.4) ;
\draw (2.6,-0.2) to  (2.6,0.7) to [out=up, in=down] (1.9,1.4);
\draw [dotted](1.2,1.4)  to (1.2,2.1);
\draw [dotted](1.9,1.4)  to (1.9,2.1);
\draw (1.2,2.1)  to [out=up, in=down] (0.5,2.8) to (0.5,3.5);
\draw (1.9,2.1) to [out=up, in=down] (1.9,2.8) to [out=up, in=down] (1.2,3.5) ;
\draw (0.5,3.5) to [out=up, in=down] (0.85,4.2) to [out=up, in=down](0.5, 4.9);
\draw (1.2,3.5) to [out=up, in=down] (0.85,4.2) to [out=up, in=down](1.2, 4.9);
\draw (0.5,4.9)  to [out=up, in=down](-0.2,5.6)  to  (-0.2,6.3);
\draw (1.2,4.9)  to (1.2,5.6)  to [out=up, in=down] (0.5,6.3);
\node [morphism, scale=\pullthroughscalevertex] at (0.85,4.2) {$\alpha$};
\node [scale=0.6] at (2.25,-0.1) {$\ldots$};
\node [scale=0.6] at (0.85,3.6) {$\ldots$};
\node [scale=0.6] at (0.85,4.8) {$\ldots$};
\node [scale=0.4, rotate = 135, thick] at (-0.6,6.4) {$\ldots$};
\node [scale=0.4, rotate = 135, thick] at (-0.6,-0.1) {$\ldots$};
\node [scale=0.4, rotate = 135, thick] at (3,6.4) {$\ldots$};
\node [scale=0.6] at (0.15,6.2) {$\ldots$};
\end{tikzpicture}
\end{aligned}
\stackrel{\widetilde{\II}_k}{\rightarrow}
\begin{aligned}
\begin{tikzpicture}[thick, scale=\expansionpullthroughscaleprime]
\draw [region] (0,-2.1) rectangle +(2.9,6.5);
\draw (0.5,-2.1) to(0.5,0);
\draw (1.2,-2.1) to(1.2,0.0);
\draw (0.5,0) to (0.5,0.7) to [out=up, in=down] (0.5,1.4) ;
\draw (1.2,0.0)  to [out=up, in=down] (1.9,0.7) to [out=up, in=down] (2.6,1.4);
\draw [dotted](2.6,1.4)  to (2.6,2.1);
\draw [dotted](0.5,1.4)  to (0.5,2.1);
\draw (0.5,2.1) to [out=up, in=down] (1.2,2.8) to [out=up, in=down] (1.9,3.5) to (1.9,4.4);
\draw (2.6,2.1)  to [out=up, in=down] (2.6,2.8) to (2.6,4.4);
\node [scale=0.6] at (0.85,-2) {$\ldots$};
\node [scale=0.6] at (2.25,4.3) {$\ldots$};
\node [scale=0.4, rotate = 135, thick] at (3,-2.2) {$\ldots$};
\draw [region] (0.2, -2.3) rectangle +(2.9,6.5);
\draw (1.9,0) to(1.9,0) to [out=up, in=down] (1.2,0.7) to [out=up, in=down] (1.2,1.4) ;
\draw (2.6,0)  to  [out=up, in=down] (2.6,0.7) to[out=up, in=down](1.9,1.4);
\draw [dotted](1.2,1.4)  to (1.2,2.1);
\draw [dotted](1.9,1.4)  to (1.9,2.1);
\draw (1.2,2.1)  to [out=up, in=down] (0.5,2.8) to (0.5,4.2);
\draw (1.9,2.1) to [out=up, in=down] (1.9,2.8) to [out=up, in=down] (1.2,3.5) to (1.2,4.2);
\draw (1.9,-2.3) to(1.9,-1.4) to [out=up, in=down] (2.25,-0.7) to [out=up, in=down](1.9, 0);
\draw (2.6,-2.3) to(2.6,-1.4) to [out=up, in=down] (2.25,-0.7) to [out=up, in=down](2.6, 0);
\node [morphism, scale=\pullthroughscalevertex] at (2.25,-0.7) {$\alpha$};
\node [scale=0.6] at (0.85,4.1) {$\ldots$};
\node [scale=0.6] at (2.25,-0.1) {$\ldots$};
\node [scale=0.6] at (2.25,-2.2) {$\ldots$};
\node [scale=0.4, rotate = 135, thick] at (0.1,4.3) {$\ldots$};
\node [scale=0.4, rotate = 135, thick] at (0.1,-2.2) {$\ldots$};
\node [scale=0.4, rotate = 135, thick] at (3,4.3) {$\ldots$};
\end{tikzpicture}
\end{aligned}
\end{align*}
\end{definition}

\noindent
A similar composition scheme is needed for the other style of type II homotopy generator, where the vertex is on the rear sheet. As for type I homotopy generators, any signature supporting type II homotopy generators always supports type II homotopy composites.

\paragraph{Extended variant.} The type $\II'$ homotopy generators include 2 further variants of the moves shown in \eqref{eq:basicpullthroughs}, corresponding to pulling a vertex through an inverse type I homotopy generator, as well as the corresponding composition schemes.

\subsection{Type III homotopy generators}

We define type III homotopy generators as follows. Note that these generators are always of dimension 5 or higher.

\begin{definition}[Type III homotopy generator]
\label{defpullthrough}
For $n \geq 0$, given an $n$\-signature $\Sigma$ supporting type II homotopy generators, $\Sigma$ supports \emph{type~III homotopy generators} if for any $4 \leq k < n$, for $k$\-diagrams with 4\-projections given by the source diagrams below, there are chosen invertible $(k{+}1)$\-cells $\III_k$ as follows:
\def\pullthroughscalevertex{1}
\begin{calign}
\label{eq:basictypeIII}
\begin{array}{@{}c@{\,\,\,}c@{\,\,}c@{}}
\begin{aligned}
\begin{tikzpicture}[thick, scale=\pullthroughscale]
\draw [region] (0.7,2.1) rectangle +(2.2,3);
\draw (2.6,2.1) to [out=up, in=down] (1.9,2.8) to [out=up, in=down] (1.2,3.5);
\draw (1.2,3.5) to (1.2,5.1);
\node [scale=0.4, rotate = 135, thick] at (0.8,5) {$\ldots$};
\node [scale=0.4, rotate = 135, thick] at (0.8,2) {$\ldots$};
\node [scale=0.4, rotate = 135, thick] at (3,5) {$\ldots$};
\node [scale=0.4, rotate = 135, thick] at (3,2) {$\ldots$};
\draw [region] (0.9, 1.9) rectangle +(2.2,3);
\draw (1.2,1.9) to (1.2,2.8) to [out=up, in=down] (1.9,3.5);
\draw (1.9,1.9) to (1.9,2.1)  to [out=up, in=down] (2.6,2.8) to (2.6,3.5);
\draw (1.9,3.5) to [out=up, in=down] (2.25,4.2) to [out=up, in=down](1.9, 4.9) to (1.9,4.9);
\draw (2.6,3.5) to [out=up, in=down] (2.25,4.2) to [out=up, in=down](2.6, 4.9)to (2.6,4.9);
\node [morphism, scale=\pullthroughscalevertex] at (2.25,4.2) {$\mu.s$};
\node [scale=0.6] at (1.55,2) {$\ldots$};
\node [scale=0.6] at (2.25,3.6) {$\ldots$};
\node [scale=0.6] at (2.25,4.8) {$\ldots$};
\end{tikzpicture}
\end{aligned}
&\stackrel{\widetilde \II_k}{\rightarrow}&
\begin{aligned}
\begin{tikzpicture}[thick, scale=\pullthroughscale]
\draw [region] (0,1.4) rectangle +(2.2,3);
\draw (1.9,2.8) to [out=up, in=down] (1.2,3.5) to [out=up, in=down] (0.5,4.2) to(0.5,4.4);
\draw (1.9, 1.4)to(1.9,2.8);
\node [scale=0.4, rotate = 135, thick] at (0.1,4.3) {$\ldots$};
\node [scale=0.4, rotate = 135, thick] at (0.1,1.3) {$\ldots$};
\node [scale=0.4, rotate = 135, thick] at (2.3,4.3) {$\ldots$};
\node [scale=0.4, rotate = 135, thick] at (2.3,1.3) {$\ldots$};
\draw [region] (0.2, 1.2) rectangle +(2.2,3);
\draw (0.5,2.8) to [out=up, in=down] (0.5,3.5) to [out=up, in=down] (1.2,4.2) ;
\draw (1.2,2.8)  to [out=up, in=down] (1.9,3.5) to [out=up, in=down] (1.9,4.2);
\draw (0.5,1.2) to(0.5,1.4) to [out=up, in=down] (0.85,2.1) to [out=up, in=down](0.5, 2.8);
\draw (1.2,1.2) to(1.2,1.4) to [out=up, in=down] (0.85,2.1) to [out=up, in=down](1.2, 2.8);
\node [morphism, scale=\pullthroughscalevertex] at (0.85,2.1) {$\mu.s$};
\node [scale=0.6] at (0.85,1.3) {$\ldots$};
\node [scale=0.6] at (0.85,2.7) {$\ldots$};
\node [scale=0.6] at (1.55,4.1) {$\ldots$};
\end{tikzpicture}
\end{aligned}
\\
\downarrow\mu
&
\begin{aligned}
\tikz{\node [rotate=45] at (0,0) {$\Rightarrow$} node [below right, font=\small] {$\III_k$};}
\end{aligned}
&
\downarrow\mu
\\[4pt]
\begin{aligned}
\begin{tikzpicture}[thick, scale=\pullthroughscale]
\draw [region] (0.7,2.1) rectangle +(2.2,3);
\draw (2.6,2.1) to [out=up, in=down] (1.9,2.8) to [out=up, in=down] (1.2,3.5);
\draw (1.2,3.5) to (1.2,5.1);
\node [scale=0.4, rotate = 135, thick] at (0.8,5) {$\ldots$};
\node [scale=0.4, rotate = 135, thick] at (0.8,2) {$\ldots$};
\node [scale=0.4, rotate = 135, thick] at (3,5) {$\ldots$};
\node [scale=0.4, rotate = 135, thick] at (3,2) {$\ldots$};
\draw [region] (0.9, 1.9) rectangle +(2.2,3);
\draw (1.2,1.9) to (1.2,2.8) to [out=up, in=down] (1.9,3.5);
\draw (1.9,1.9) to (1.9,2.1)  to [out=up, in=down] (2.6,2.8) to (2.6,3.5);
\draw (1.9,3.5) to [out=up, in=down] (2.25,4.2) to [out=up, in=down](1.9, 4.9) to (1.9,4.9);
\draw (2.6,3.5) to [out=up, in=down] (2.25,4.2) to [out=up, in=down](2.6, 4.9)to (2.6,4.9);
\node [morphism, scale=\pullthroughscalevertex] at (2.25,4.2) {$\mu.t$};
\node [scale=0.6] at (1.55,2) {$\ldots$};
\node [scale=0.6] at (2.25,3.6) {$\ldots$};
\node [scale=0.6] at (2.25,4.8) {$\ldots$};
\end{tikzpicture}
\end{aligned}
&\stackrel{\widetilde \II_k}{\rightarrow}&
\begin{aligned}
\begin{tikzpicture}[thick, scale=\pullthroughscale]
\draw [region] (0,1.4) rectangle +(2.2,3);
\draw (1.9,2.8) to [out=up, in=down] (1.2,3.5) to [out=up, in=down] (0.5,4.2) to(0.5,4.4);
\draw (1.9, 1.4)to(1.9,2.8);
\node [scale=0.4, rotate = 135, thick] at (0.1,4.3) {$\ldots$};
\node [scale=0.4, rotate = 135, thick] at (0.1,1.3) {$\ldots$};
\node [scale=0.4, rotate = 135, thick] at (2.3,4.3) {$\ldots$};
\node [scale=0.4, rotate = 135, thick] at (2.3,1.3) {$\ldots$};
\draw [region] (0.2, 1.2) rectangle +(2.2,3);
\draw (0.5,2.8) to [out=up, in=down] (0.5,3.5) to [out=up, in=down] (1.2,4.2) ;
\draw (1.2,2.8)  to [out=up, in=down] (1.9,3.5) to [out=up, in=down] (1.9,4.2);
\draw (0.5,1.2) to(0.5,1.4) to [out=up, in=down] (0.85,2.1) to [out=up, in=down](0.5, 2.8);
\draw (1.2,1.2) to(1.2,1.4) to [out=up, in=down] (0.85,2.1) to [out=up, in=down](1.2, 2.8);
\node [morphism, scale=\pullthroughscalevertex] at (0.85,2.1) {$\mu.t$};
\node [scale=0.6] at (0.85,1.3) {$\ldots$};
\node [scale=0.6] at (0.85,2.7) {$\ldots$};
\node [scale=0.6] at (1.55,4.1) {$\ldots$};
\end{tikzpicture}
\end{aligned}
\end{array}
&
\begin{array}{@{}c@{\,\,\,}c@{\,\,}c@{}}
\begin{aligned}
\begin{tikzpicture}[thick, scale=\pullthroughscale]
\draw [region] (0,2.1) rectangle +(2.2,3);
\draw (1.2,2.1)  to [out=up, in=down] (0.5,2.8) to (0.5,3.5);
\draw (1.9,2.1) to [out=up, in=down] (1.9,2.8) to [out=up, in=down] (1.2,3.5) ;
\draw (0.5,3.5) to [out=up, in=down] (0.85,4.2) to [out=up, in=down](0.5, 4.9)to (0.5,5.1);
\draw (1.2,3.5) to [out=up, in=down] (0.85,4.2) to [out=up, in=down](1.2, 4.9)to (1.2,5.1);
\node [morphism, scale=\pullthroughscalevertex] at (0.85,4.2) {$\mu.s$};
\node [scale=0.6] at (1.55,2.2) {$\ldots$};
\node [scale=0.6] at (0.85,3.6) {$\ldots$};
\node [scale=0.6] at (0.85,5) {$\ldots$};
\node [scale=0.4, rotate = 135, thick] at (0.1,5) {$\ldots$};
\node [scale=0.4, rotate = 135, thick] at (0.1,2) {$\ldots$};
\node [scale=0.4, rotate = 135, thick] at (2.3,5) {$\ldots$};
\node [scale=0.4, rotate = 135, thick] at (2.3,2) {$\ldots$};
\draw [region] (0.2, 1.9) rectangle +(2.2,3);
\draw (0.5,1.9) to (0.5,2.1) to [out=up, in=down] (1.2,2.8) to [out=up, in=down] (1.9,3.5) to (1.9,4.9);
\end{tikzpicture}
\end{aligned}
&\stackrel{\widetilde \II_k}{\rightarrow}&
\begin{aligned}
\begin{tikzpicture}[thick, scale=\pullthroughscale]
\draw [region] (0.7,-1.4) rectangle +(2.2,3);
\draw (1.9,0) to(1.9,0) to [out=up, in=down] (1.2,0.7) to [out=up, in=down] (1.2,1.4)to(1.2,1.6) ;
\draw (2.6,0)  to  [out=up, in=down] (2.6,0.7) to[out=up, in=down](1.9,1.4)to(1.9,1.6);
\draw (1.9,-1.4) to(1.9,-1.4) to [out=up, in=down] (2.25,-0.7) to [out=up, in=down](1.9, 0);
\draw (2.6,-1.4) to(2.6,-1.4) to [out=up, in=down] (2.25,-0.7) to [out=up, in=down](2.6, 0);
\node [morphism, scale=\pullthroughscalevertex] at (2.25,-0.7) {$\mu.s$};
\node [scale=0.6] at (1.55,1.5) {$\ldots$};
\node [scale=0.6] at (2.25,-0.1) {$\ldots$};
\node [scale=0.4, rotate = 135, thick] at (0.8,1.5) {$\ldots$};
\node [scale=0.4, rotate = 135, thick] at (0.8,-1.5) {$\ldots$};
\node [scale=0.4, rotate = 135, thick] at (3,1.5) {$\ldots$};
\node [scale=0.4, rotate = 135, thick] at (3,-1.5) {$\ldots$};
\node [scale=0.6] at (2.25,-1.3) {$\ldots$};
\draw [region] (0.9, -1.6) rectangle +(2.2,3);
\draw (1.2,-1.6) to(1.2,0.0);
\draw (1.2,0.0)  to [out=up, in=down] (1.9,0.7) to [out=up, in=down] (2.6,1.4);
\end{tikzpicture}
\end{aligned}
\\
\downarrow\mu
&
\begin{aligned}
\tikz{\node [rotate=45] at (0,0) {$\Rightarrow$} node [below right, font=\small] {$\III_k$};}
\end{aligned}
&
\downarrow\mu
\\[4pt]
\begin{aligned}
\begin{tikzpicture}[thick, scale=\pullthroughscale]
\draw [region] (0,2.1) rectangle +(2.2,3);
\draw (1.2,2.1)  to [out=up, in=down] (0.5,2.8) to (0.5,3.5);
\draw (1.9,2.1) to [out=up, in=down] (1.9,2.8) to [out=up, in=down] (1.2,3.5) ;
\draw (0.5,3.5) to [out=up, in=down] (0.85,4.2) to [out=up, in=down](0.5, 4.9)to (0.5,5.1);
\draw (1.2,3.5) to [out=up, in=down] (0.85,4.2) to [out=up, in=down](1.2, 4.9)to (1.2,5.1);
\node [morphism, scale=\pullthroughscalevertex] at (0.85,4.2) {$\mu.t$};
\node [scale=0.6] at (1.55,2.2) {$\ldots$};
\node [scale=0.6] at (0.85,3.6) {$\ldots$};
\node [scale=0.6] at (0.85,5) {$\ldots$};
\node [scale=0.4, rotate = 135, thick] at (0.1,5) {$\ldots$};
\node [scale=0.4, rotate = 135, thick] at (0.1,2) {$\ldots$};
\node [scale=0.4, rotate = 135, thick] at (2.3,5) {$\ldots$};
\node [scale=0.4, rotate = 135, thick] at (2.3,2) {$\ldots$};
\draw [region] (0.2, 1.9) rectangle +(2.2,3);
\draw (0.5,1.9) to (0.5,2.1) to [out=up, in=down] (1.2,2.8) to [out=up, in=down] (1.9,3.5) to (1.9,4.9);
\end{tikzpicture}
\end{aligned}
&\stackrel{\widetilde \II_k}{\rightarrow}&
\begin{aligned}
\begin{tikzpicture}[thick, scale=\pullthroughscale]
\draw [region] (0.7,-1.4) rectangle +(2.2,3);
\draw (1.9,0) to(1.9,0) to [out=up, in=down] (1.2,0.7) to [out=up, in=down] (1.2,1.4)to(1.2,1.6) ;
\draw (2.6,0)  to  [out=up, in=down] (2.6,0.7) to[out=up, in=down](1.9,1.4)to(1.9,1.6);
\draw (1.9,-1.4) to(1.9,-1.4) to [out=up, in=down] (2.25,-0.7) to [out=up, in=down](1.9, 0);
\draw (2.6,-1.4) to(2.6,-1.4) to [out=up, in=down] (2.25,-0.7) to [out=up, in=down](2.6, 0);
\node [morphism, scale=\pullthroughscalevertex] at (2.25,-0.7) {$\mu.t$};
\node [scale=0.6] at (1.55,1.5) {$\ldots$};
\node [scale=0.6] at (2.25,-0.1) {$\ldots$};
\node [scale=0.4, rotate = 135, thick] at (0.8,1.5) {$\ldots$};
\node [scale=0.4, rotate = 135, thick] at (0.8,-1.5) {$\ldots$};
\node [scale=0.4, rotate = 135, thick] at (3,1.5) {$\ldots$};
\node [scale=0.4, rotate = 135, thick] at (3,-1.5) {$\ldots$};
\node [scale=0.6] at (2.25,-1.3) {$\ldots$};
\draw [region] (0.9, -1.6) rectangle +(2.2,3);
\draw (1.2,-1.6) to(1.2,0.0);
\draw (1.2,0.0)  to [out=up, in=down] (1.9,0.7) to [out=up, in=down] (2.6,1.4);
\end{tikzpicture}
\end{aligned}
\end{array}
\end{calign}
\end{definition}

\noindent
Note that these diagrams make use of the expansion scheme for type II\ homotopy generators, since $\mu.s$ and $\mu.t$ are diagrams, not generators. 

\paragraph{Expansion scheme.} We have no need for an expansion scheme for type III homotopy generators. Such a scheme would be needed for the definition of semistrict 5\-category.

\paragraph{Extended variant.} The type $\III'$ homotopy generators include 2 further variants of the moves shown in \eqref{eq:basictypeIII}, corresponding to the extended type $\II'$ homotopy generators.

\subsection{Type IV homotopy generators}

We define type IV homotopy generators as follows. Note that these generators are always of dimension 5 or higher.

\begin{definition}[Type IV homotopy generator]
\label{defpullthrough}
For $n \geq 0$, given an $n$\-signature $\Sigma$ supporting type II homotopy generators, $\Sigma$ supports \emph{type~IV homotopy generators} if for any $4 \leq k < n$, for a $k$\-diagram with 4\-projections given by the counterclockwise diagram below, there are chosen invertible $(k{+}1)$\-cells $\IV_k$ as follows:
\begin{align}
\label{eq:basictypeIV}
\begin{aligned}
\begin{tikzpicture}[thick, scale=0.8]
\draw [region] (-0.3,0) rectangle +(2.4,3);
\draw (1.8,0)  to [out=up, in=down] (1.1,0.8) to [out=up, in=down] (0.4,1.6)to (0.4,3);
\draw [region] (-0.1,-0.2) rectangle +(2.4,3);
\draw (1.1,-0.2)to (1.1,0) to [out=up, in=down] (1.8,0.8) to (1.8,1.6) to [out=up, in=down]  (1.1,2.4)to (1.1,2.8);
\draw [region] (0.1,-0.4) rectangle +(2.4,3);
\draw (0.4,-0.4)  to (0.4,0.8)  to [out=up, in=down] (1.1,1.6) to [out=up, in=down] (1.8,2.4)to (1.8,2.6);
\end{tikzpicture}
\end{aligned}
\quad\begin{array}{@{}c@{}}\stackrel{\II_{k-1}}{\rightarrow}\\
\begin{aligned}
\begin{tikzpicture}[thick, scale=0.8]
\node (1) [scale=1] at (0.6,0) {$\Uparrow$};
\node [right=-7pt] at (1.east) {$\IV_k$};
\end{tikzpicture}
\end{aligned}\\
\stackrel{\II_{k-1} ^{-1}}{\rightarrow}\end{array}\quad
\begin{aligned}
\begin{tikzpicture}[thick, scale=0.8]
\draw [region] (-0.3,0) rectangle +(2.4,3);
\draw (1.8,0)  to (1.8,0.8)  to[out=up, in=down] (1.1,1.6) to [out=up, in=down] (0.4,2.4)to (0.4,3);
\draw [region] (-0.1,-0.2) rectangle +(2.4,3);
\draw (1.1,-0.2)to (1.1,0) to [out=up, in=down] (0.4,0.8) to (0.4,1.6) to [out=up, in=down]  (1.1,2.4)to (1.1,2.8);
\draw [region] (0.1,-0.4) rectangle +(2.4,3);
\draw (0.4,-0.4)  to (0.4,0)  to [out=up, in=down] (1.1,0.8) to [out=up, in=down] (1.8,1.6)to (1.8,2.6);
\end{tikzpicture}
\end{aligned}
\end{align}
In the clockwise path, the upper crossing is `pulled down' by the $\II_{k-1}$ move. In the counterclockwise path, the lower crossing is `pulled up' by the $\II_{k-1}^{-1}$ move. Both of these can be seen as ways to interpret the third Reidemeister move. The type IV homotopy generator says that these two diagrams are related by an invertible cell.
\end{definition}

\noindent
The need for this homotopy was identified by Breen~\cite{BreentoHopkins} and the cell has been referred to in the literature as the `Breenator'.

\paragraph{Expansion scheme.} We have no need for an expansion scheme for type IV homotopy generators.

\paragraph{Extended variant.} The type $\IV'$ homotopy generators include 2 further variants of the moves shown in \eqref{eq:basictypeIV}, corresponding to the inverse type I homotopy generators.

\subsection{Type V homotopy generators}

We now describe type V homotopy generators. These generators are always of dimension 5 or higher. They require a rearrangement scheme for type I homotopy generators, labelled $\widehat{ \I}$; this is described below.

\begin{definition}[Type V homotopy generator]
\label{defpullthrough}
For $n \geq 0$, given an $n$\-signature $\Sigma$ supporting type II homotopy generators, $\Sigma$ supports \emph{type~V homotopy generators} if for any $4 \leq k < n$, a $k$\-diagram with 4\-projections given by the counterclockwise diagram below for $\alpha,\beta \in \Sigma_{k-1}$, there are chosen invertible $(k{+}1)$\-cells $\V_k$ as follows:
\def\lnscalevertex{0.5}
\def\lnxscale{0.5}
\def\lnyscale{0.5}
\def\lnscaledot{0.4}
\def\lnscaledotperspective{0.2}
\begin{equation*}
\tikzset{every picture/.style={yscale=0.85, scale=0.6}}
\label{eq:binaturalitypicture}
\begin{array}{ccccccccc}
\begin{aligned}
\begin{tikzpicture}[thick, scale=0.8]
\draw [region] (0,0) rectangle +(2.9,6.5);
\draw (1.9,0) to(1.9,0) to [out=up, in=down] (1.2,0.7) to [out=up, in=down] (0.5,1.4) ;
\draw (2.6,0)  to (2.6,1.4);
\draw [dotted](2.6,1.4)  to (2.6,2.1);
\draw [dotted](0.5,1.4)  to (0.5,2.1);
\draw (0.5,2.1)  to (0.5,3.5);
\draw (2.6,2.1) to [out=up, in=down] (1.9,2.8) to [out=up, in=down] (1.2,3.5) ;
\draw (0.5,3.5) to (0.5,4.9) to [out=up, in=down] (0.85,5.6) to [out=up, in=down](0.5, 6.3)to (0.5,6.5);
\draw (1.2,3.5) to (1.2,4.9) to [out=up, in=down] (0.85,5.6) to [out=up, in=down](1.2, 6.3)to (1.2,6.5);
\node [morphism, scale=\lnscalevertex] at (0.85,5.6) {$\alpha$};
\node [scale=\lnscaledot] at (2.25,0.1) {$\ldots$};
\node [scale=\lnscaledot] at (0.85,5) {$\ldots$};
\node [scale=\lnscaledot] at (0.85,6.4) {$\ldots$};
\node [scale=\lnscaledotperspective, rotate = 135, thick] at (0.1,6.4) {$\ldots$};
\node [scale=\lnscaledotperspective, rotate = 135, thick] at (0.1,-0.1) {$\ldots$};
\node [scale=\lnscaledotperspective, rotate = 135, thick] at (3,6.4) {$\ldots$};
\draw [region] (0.2, -0.2) rectangle +(2.9,6.5);
\draw (0.5,-0.2) to (0.5,0) to (0.5,0.7) to [out=up, in=down] (1.2,1.4) ;
\draw (1.2,-0.2) to(1.2,0.0)  to [out=up, in=down] (1.9,0.7) to (1.9,1.4);
\draw [dotted](1.2,1.4)  to (1.2,2.1);
\draw [dotted](1.9,1.4)  to (1.9,2.1);
\draw (1.2,2.1) to (1.2,2.8) to [out=up, in=down] (1.9,3.5);
\draw (1.9,2.1)  to [out=up, in=down] (2.6,2.8) to (2.6,3.5);
\draw (1.9,3.5) to [out=up, in=down] (2.25,4.2) to [out=up, in=down](1.9, 4.9) to (1.9,6.3);
\draw (2.6,3.5) to [out=up, in=down] (2.25,4.2) to [out=up, in=down](2.6, 4.9)to (2.6,6.3);
\node [morphism, scale=\lnscalevertex] at (2.25,4.2) {$\beta$};
\node [scale=\lnscaledot] at (0.85,-0.1) {$\ldots$};
\node [scale=\lnscaledot] at (2.25,3.6) {$\ldots$};
\node [scale=\lnscaledot] at (2.25,4.8) {$\ldots$};
\node [scale=\lnscaledot] at (2.25,6.2) {$\ldots$};
\node [scale=\lnscaledotperspective, rotate = 135, thick] at (3,-0.1) {$\ldots$};
\end{tikzpicture}
\end{aligned}
&
\stackrel{\II_{k-1}}{\rightarrow}
&
\begin{aligned}
\begin{tikzpicture}[thick, scale=0.8]
\draw [region] (0,-1.4) rectangle +(2.9,6.5);
\draw (2.6,-1.4)  to (2.6,0);
\draw (1.9,-1.4)  to (1.9,0);
\draw (1.9,-0.2) to(1.9,0) to [out=up, in=down] (1.2,0.7) to [out=up, in=down] (0.5,1.4) ;
\draw (2.6,-0.2)  to (2.6,1.4);
\draw [dotted](2.6,1.4)  to (2.6,2.1);
\draw [dotted](0.5,1.4)  to (0.5,2.1);
\draw (0.5,2.1)  to (0.5,3.5);
\draw (2.6,2.1) to [out=up, in=down] (1.9,2.8) to [out=up, in=down] (1.2,3.5) ;
\draw (0.5,3.5) to [out=up, in=down] (0.85,4.2) to [out=up, in=down](0.5, 4.9) to (0.5,5.1);
\draw (1.2,3.5) to [out=up, in=down] (0.85,4.2) to [out=up, in=down](1.2, 4.9) to (1.2,5.1);
\node [morphism, scale=\lnscalevertex] at (0.85,4.2) {$\alpha$};
\node [scale=\lnscaledot] at (2.25,-1.3) {$\ldots$};
\node [scale=\lnscaledot] at (0.85,3.6) {$\ldots$};
\node [scale=\lnscaledot] at (0.85,5) {$\ldots$};
\node [scale=\lnscaledotperspective, rotate = 135, thick] at (0.1,5) {$\ldots$};
\node [scale=\lnscaledotperspective, rotate = 135, thick] at (0.1,-1.3) {$\ldots$};
\node [scale=\lnscaledotperspective, rotate = 135, thick] at (3,5) {$\ldots$};
\draw [region] (0.2, -1.6) rectangle +(2.9,6.5);
\draw (0.5,0) to (0.5,0.7) to [out=up, in=down] (1.2,1.4) ;
\draw (1.2,0.0)  to [out=up, in=down] (1.9,0.7) to (1.9,1.4);
\draw [dotted](1.2,1.4)  to (1.2,2.1);
\draw [dotted](1.9,1.4)  to (1.9,2.1);
\draw (1.2,2.1) to (1.2,2.8) to [out=up, in=down] (1.9,3.5) to (1.9,4.9);
\draw (1.9,2.1)  to [out=up, in=down] (2.6,2.8) to (2.6,3.5) to (2.6,4.9);
\draw (0.5,-1.6) to (0.5,-1.4) to [out=up, in=down] (0.85,-0.7) to [out=up, in=down](0.5, 0);
\draw (1.2,-1.6) to (1.2,-1.4) to [out=up, in=down] (0.85,-0.7) to [out=up, in=down](1.2, 0);
\node [morphism, scale=\lnscalevertex] at (0.85,-0.7) {$\beta$};
\node [scale=\lnscaledot] at (0.85,-0.1) {$\ldots$};
\node [scale=\lnscaledot] at (0.85,-1.5) {$\ldots$};
\node [scale=\lnscaledot] at (2.25,4.8) {$\ldots$};
\node [scale=\lnscaledotperspective, rotate = 135, thick] at (3,-1.5) {$\ldots$};
\end{tikzpicture}
\end{aligned}
&
\stackrel{\widehat{{\I}}_{k-1}}{\rightarrow}
&
\begin{aligned}
\begin{tikzpicture}[thick, scale=0.8]
\draw [region] (0,-1.4) rectangle +(2.9,6.5);
\draw (2.6,-1.4)  to (2.6,0);
\draw (1.9,-1.4)  to (1.9,0);
\draw (1.9,0) to(1.9,0) to [out=up, in=down] (1.2,0.7) to [out=up, in=down] (1.2,1.4) ;
\draw (2.6,0)  to  [out=up, in=down] (2.6,0.7) to[out=up, in=down](1.9,1.4);
\draw [dotted](1.2,1.4)  to (1.2,2.1);
\draw [dotted](1.9,1.4)  to (1.9,2.1);
\draw (1.2,2.1)  to [out=up, in=down] (0.5,2.8) to (0.5,3.5);
\draw (1.9,2.1) to [out=up, in=down] (1.9,2.8) to [out=up, in=down] (1.2,3.5) ;
\draw (0.5,3.5) to [out=up, in=down] (0.85,4.2) to [out=up, in=down](0.5, 4.9) to (0.5,5.1);
\draw (1.2,3.5) to [out=up, in=down] (0.85,4.2) to [out=up, in=down](1.2, 4.9) to (1.2,5.1);
\node [morphism, scale=\lnscalevertex] at (0.85,4.2) {$\alpha$};
\node [scale=\lnscaledot] at (2.25,-1.3) {$\ldots$};
\node [scale=\lnscaledot] at (0.85,3.6) {$\ldots$};
\node [scale=\lnscaledot] at (0.85,5) {$\ldots$};
\node [scale=\lnscaledotperspective, rotate = 135, thick] at (0.1,5) {$\ldots$};
\node [scale=\lnscaledotperspective, rotate = 135, thick] at (0.1,-1.5) {$\ldots$};
\node [scale=\lnscaledotperspective, rotate = 135, thick] at (3,5) {$\ldots$};
\draw [region] (0.2, -1.6) rectangle +(2.9,6.5);
\draw (2.6,3.5)  to (2.6,4.9);
\draw (1.9,3.5)  to (1.9,4.9);
\draw (0.5,0) to (0.5,0.7) to [out=up, in=down] (0.5,1.4) ;
\draw (1.2,0.0)  to [out=up, in=down] (1.9,0.7) to [out=up, in=down] (2.6,1.4);
\draw [dotted](2.6,1.4)  to (2.6,2.1);
\draw [dotted](0.5,1.4)  to (0.5,2.1);
\draw (0.5,2.1) to [out=up, in=down] (1.2,2.8) to [out=up, in=down] (1.9,3.5) ;
\draw (2.6,2.1)  to [out=up, in=down] (2.6,2.8) to (2.6,3.5);
\draw (0.5,-1.6) to(0.5,-1.4) to [out=up, in=down] (0.85,-0.7) to [out=up, in=down](0.5, 0);
\draw (1.2,-1.6) to(1.2,-1.4) to [out=up, in=down] (0.85,-0.7) to [out=up, in=down](1.2, 0);
\node [morphism, scale=\lnscalevertex] at (0.85,-0.7) {$\beta$};
\node [scale=\lnscaledot] at (0.85,-0.1) {$\ldots$};
\node [scale=\lnscaledot] at (0.85,-1.5) {$\ldots$};
\node [scale=\lnscaledot] at (2.25,4.8) {$\ldots$};
\node [scale=\lnscaledotperspective, rotate = 135, thick] at (3,-1.5) {$\ldots$};
\end{tikzpicture}
\end{aligned}
&
\stackrel{\II_{k-1}}{\rightarrow}
&
\begin{aligned}
\begin{tikzpicture}[thick, scale=0.8]
\draw [region] (0,0) rectangle +(2.9,6.5);
\draw (1.9,2.8) to [out=up, in=down] (1.2,3.5) to [out=up, in=down] (1.2,4.2) ;
\draw (2.6,2.8)  to  [out=up, in=down] (2.6,3.5) to[out=up, in=down](1.9,4.2);
\draw [dotted](1.2,4.2)  to (1.2,4.9);
\draw [dotted](1.9,4.2)  to (1.9,4.9);
\draw (1.2,4.9)  to [out=up, in=down] (0.5,5.6) to (0.5,6.3)to (0.5,6.5);
\draw (1.9,4.9) to [out=up, in=down] (1.9,5.6) to [out=up, in=down] (1.2,6.3) to (1.2,6.5);
\draw (1.9,0) to(1.9,1.4) to [out=up, in=down] (2.25,2.1) to [out=up, in=down](1.9,2.8);
\draw (2.6,0) to(2.6,1.4) to [out=up, in=down] (2.25,2.1) to [out=up, in=down](2.6,2.8);
\node [morphism, scale=\lnscalevertex] at (2.25,2.1) {$\alpha$};
\node [scale=\lnscaledot] at (2.25,0.1) {$\ldots$};
\node [scale=\lnscaledot] at (2.25,2.7) {$\ldots$};
\node [scale=\lnscaledot] at (2.25,1.5) {$\ldots$};
\node [scale=\lnscaledot] at (0.85,6.4) {$\ldots$};
\node [scale=\lnscaledotperspective, rotate = 135, thick] at (0.1,6.4) {$\ldots$};
\node [scale=\lnscaledotperspective, rotate = 135, thick] at (0.1,-0.1) {$\ldots$};
\node [scale=\lnscaledotperspective, rotate = 135, thick] at (3,6.4) {$\ldots$};
\draw [region] (0.2, -0.2) rectangle +(2.9,6.5);
\draw (0.5,2.8) to (0.5,3.5) to [out=up, in=down] (0.5,4.2) ;
\draw (1.2,2.8)  to [out=up, in=down] (1.9,3.5) to [out=up, in=down] (2.6,4.2);
\draw [dotted](2.6,4.2)  to (2.6,4.9);
\draw [dotted](0.5,4.2)  to (0.5,4.9);
\draw (0.5,4.9) to [out=up, in=down] (1.2,5.6) to [out=up, in=down] (1.9,6.3);
\draw (2.6,4.9)  to [out=up, in=down] (2.6,5.6) to (2.6,6.3);
\draw (0.5,-0.2)  to(0.5,0)  to [out=up, in=down] (0.85,0.7) to [out=up, in=down](0.5, 1.4) to (0.5,2.8);
\draw (1.2,-0.2)  to(1.2,0) to [out=up, in=down] (0.85,0.7) to [out=up, in=down](1.2, 1.4)to (1.2,2.8);;
\node [morphism, scale=\lnscalevertex] at (0.85,0.7) {$\beta$};
\node [scale=\lnscaledot] at (0.85,-0.1) {$\ldots$};
\node [scale=\lnscaledot] at (0.85,1.3) {$\ldots$};
\node [scale=\lnscaledot] at (2.25,6.2) {$\ldots$};
\node [scale=\lnscaledotperspective, rotate = 135, thick] at (3,-0.1) {$\ldots$};
\end{tikzpicture}
\end{aligned}
&
\stackrel{\widehat{{\I}}_{k-1}}{\rightarrow}
&
\begin{aligned}
\begin{tikzpicture}[thick, scale=0.8]
\draw [region] (0,0) rectangle +(2.9,6.5);
\draw (1.9,2.8) to [out=up, in=down] (1.2,3.5) to [out=up, in=down] (0.5,4.2) ;
\draw (2.6,2.8)  to  [out=up, in=down] (2.6,3.5) to[out=up, in=down](2.6,4.2);
\draw [dotted](2.6,4.2)  to (2.6,4.9);
\draw [dotted](0.5,4.2)  to (0.5,4.9);
\draw (0.5,4.9)  to [out=up, in=down] (0.5,5.6) to (0.5,6.3) to (0.5,6.5);
\draw (2.6,4.9) to [out=up, in=down] (1.9,5.6) to [out=up, in=down] (1.2,6.3) to (1.2,6.5);
\draw (1.9,0) to(1.9,1.4) to [out=up, in=down] (2.25,2.1) to [out=up, in=down](1.9, 2.8)to(1.9,2.8);
\draw (2.6,0) to(2.6,1.4) to [out=up, in=down] (2.25,2.1) to [out=up, in=down](2.6,2.8);
\node [morphism, scale=\lnscalevertex] at (2.25,2.1) {$\alpha$};
\node [scale=\lnscaledot] at (2.25,0.1) {$\ldots$};
\node [scale=\lnscaledot] at (2.25,1.5) {$\ldots$};
\node [scale=\lnscaledot] at (0.85,6.4) {$\ldots$};
\node [scale=\lnscaledot] at (2.25,2.7) {$\ldots$};
\node [scale=\lnscaledotperspective, rotate = 135, thick] at (0.1,6.4) {$\ldots$};
\node [scale=\lnscaledotperspective, rotate = 135, thick] at (0.1,-0.1) {$\ldots$};
\node [scale=\lnscaledotperspective, rotate = 135, thick] at (3,6.4) {$\ldots$};
\draw [region] (0.2, -0.2) rectangle +(2.9,6.5);
\draw (0.5,2.8) to [out=up, in=down] (0.5,3.5) to [out=up, in=down] (1.2,4.2) ;
\draw (1.2,2.8)  to [out=up, in=down] (1.9,3.5) to [out=up, in=down] (1.9,4.2);
\draw [dotted](1.2,4.2)  to (1.2,4.9);
\draw [dotted](1.9,4.2)  to (1.9,4.9);
\draw (1.2,4.9) to [out=up, in=down] (1.2,5.6) to [out=up, in=down] (1.9,6.3) to (1.9,6.3);
\draw (1.9,4.9)  to [out=up, in=down] (2.6,5.6) to (2.6,6.3)to (2.6,6.3);
\draw (0.5,-0.2) to(0.5,0) to [out=up, in=down] (0.85,0.7) to [out=up, in=down](0.5, 1.4)to (0.5, 2.8);
\draw (1.2,-0.2) to(1.2,0) to [out=up, in=down] (0.85,0.7) to [out=up, in=down](1.2, 1.4) to (1.2, 2.8);
\node [morphism, scale=\lnscalevertex] at (0.85,0.7) {$\beta$};
\node [scale=\lnscaledot] at (0.85,-0.1) {$\ldots$};
\node [scale=\lnscaledot] at (0.85,1.3) {$\ldots$};
\node [scale=\lnscaledot] at (2.25,6.2) {$\ldots$};
\node [scale=\lnscaledotperspective, rotate = 135, thick] at (3,-0.1) {$\ldots$};
\end{tikzpicture}
\end{aligned}
\\
\downarrow \scriptstyle \I_{k-1}&&&&
\begin{aligned}
\begin{tikzpicture}[thick, scale=1]
\node [rotate=0, scale = 0.7] at (-0.7,1.0) {$\V_k$};
\node [rotate=60, scale = 1.2] at (0,1) {$\Rightarrow $};
\end{tikzpicture}
\end{aligned}
&&&&\downarrow \scriptstyle{\I_{k-1}}
\\
\begin{aligned}
\begin{tikzpicture}[thick, scale=0.8]
\draw [region] (0,0) rectangle +(2.9,6.5);
\draw (1.9,0) to(1.9,0) to [out=up, in=down] (1.2,0.7) to [out=up, in=down] (0.5,1.4) ;
\draw (2.6,0)  to (2.6,1.4);
\draw [dotted](2.6,1.4)  to (2.6,2.1);
\draw [dotted](0.5,1.4)  to (0.5,2.1);
\draw (0.5,2.1)  to (0.5,3.5);
\draw (2.6,2.1) to [out=up, in=down] (1.9,2.8) to [out=up, in=down] (1.2,3.5) ;
\draw (0.5,3.5) to [out=up, in=down] (0.85,4.2) to [out=up, in=down](0.5, 4.9)to (0.5,6.5);
\draw (1.2,3.5) to  [out=up, in=down] (0.85,4.2) to [out=up, in=down](1.2, 4.9)to (1.2,6.5);
\node [morphism, scale=\lnscalevertex] at (0.85,4.2) {$\alpha$};
\node [scale=\lnscaledot] at (2.25,0.1) {$\ldots$};
\node [scale=\lnscaledot] at (0.85,4.8) {$\ldots$};
\node [scale=\lnscaledot] at (0.85,6.4) {$\ldots$};
\node [scale=\lnscaledot] at (0.85,3.6) {$\ldots$};
\node [scale=\lnscaledotperspective, rotate = 135, thick] at (0.1,6.4) {$\ldots$};
\node [scale=\lnscaledotperspective, rotate = 135, thick] at (0.1,-0.1) {$\ldots$};
\node [scale=\lnscaledotperspective, rotate = 135, thick] at (3,6.4) {$\ldots$};
\draw [region] (0.2, -0.2) rectangle +(2.9,6.5);
\draw (0.5,-0.2) to (0.5,0) to (0.5,0.7) to [out=up, in=down] (1.2,1.4) ;
\draw (1.2,-0.2) to(1.2,0.0)  to [out=up, in=down] (1.9,0.7) to (1.9,1.4);
\draw [dotted](1.2,1.4)  to (1.2,2.1);
\draw [dotted](1.9,1.4)  to (1.9,2.1);
\draw (1.2,2.1) to (1.2,2.8) to [out=up, in=down] (1.9,3.5);
\draw (1.9,2.1)  to [out=up, in=down] (2.6,2.8) to (2.6,3.5);
\draw (1.9,3.5) to (1.9,4.9) to [out=up, in=down] (2.25,5.6) to [out=up, in=down](1.9,6.3);
\draw (2.6,3.5) to (2.6,4.9) to [out=up, in=down] (2.25,5.6) to [out=up, in=down](2.6,6.3);
\node [morphism, scale=\lnscalevertex] at (2.25,5.6) {$\beta$};
\node [scale=\lnscaledot] at (0.85,-0.1) {$\ldots$};
\node [scale=\lnscaledot] at (2.25,5) {$\ldots$};
\node [scale=\lnscaledot] at (2.25,6.2) {$\ldots$};
\node [scale=\lnscaledotperspective, rotate = 135, thick] at (3,-0.1) {$\ldots$};
\end{tikzpicture}
\end{aligned}
&
\stackrel{\widehat{{\I}}_{k-1}}{\rightarrow}
&
\begin{aligned}
\begin{tikzpicture}[thick, scale=0.8]
\draw [region] (0,0) rectangle +(2.9,6.5);
\draw (1.9,0) to(1.9,0) to [out=up, in=down] (1.2,0.7) to [out=up, in=down] (1.2,1.4) ;
\draw (2.6,0)  to  [out=up, in=down] (2.6,0.7) to[out=up, in=down](1.9,1.4);
\draw [dotted](1.2,1.4)  to (1.2,2.1);
\draw [dotted](1.9,1.4)  to (1.9,2.1);
\draw (1.2,2.1)  to [out=up, in=down] (0.5,2.8) to (0.5,3.5);
\draw (1.9,2.1) to [out=up, in=down] (1.9,2.8) to [out=up, in=down] (1.2,3.5) ;
\draw (0.5,3.5) to [out=up, in=down] (0.85,4.2) to [out=up, in=down](0.5, 4.9)to (0.5,6.3)to (0.5,6.5);
\draw (1.2,3.5) to [out=up, in=down] (0.85,4.2) to [out=up, in=down](1.2, 4.9)to (1.2,6.3) to (1.2,6.5);
\node [morphism, scale=\lnscalevertex] at (0.85,4.2) {$\alpha$};
\node [scale=\lnscaledot] at (2.25,0.1) {$\ldots$};
\node [scale=\lnscaledot] at (0.85,3.6) {$\ldots$};
\node [scale=\lnscaledot] at (0.85,4.8) {$\ldots$};
\node [scale=\lnscaledot] at (0.85,6.4) {$\ldots$};
\node [scale=\lnscaledotperspective, rotate = 135, thick] at (0.1,6.4) {$\ldots$};
\node [scale=\lnscaledotperspective, rotate = 135, thick] at (0.1,-0.1) {$\ldots$};
\node [scale=\lnscaledotperspective, rotate = 135, thick] at (3,6.4) {$\ldots$};
\draw [region] (0.2, -0.2) rectangle +(2.9,6.5);
\draw (0.5,-0.2) to (0.5,0) to (0.5,0.7) to [out=up, in=down] (0.5,1.4) ;
\draw (1.2,-0.2)  to(1.2,0.0)  to [out=up, in=down] (1.9,0.7) to [out=up, in=down] (2.6,1.4);
\draw [dotted](2.6,1.4)  to (2.6,2.1);
\draw [dotted](0.5,1.4)  to (0.5,2.1);
\draw (0.5,2.1) to [out=up, in=down] (1.2,2.8) to [out=up, in=down] (1.9,3.5) to (1.9,4.9);
\draw (2.6,2.1)  to [out=up, in=down] (2.6,2.8) to (2.6,3.5) to (2.6,4.9);
\draw (1.9,4.9) to [out=up, in=down] (2.25,5.6) to [out=up, in=down](1.9, 6.3) to (1.9,6.3);
\draw (2.6,4.9) to [out=up, in=down] (2.25,5.6) to [out=up, in=down](2.6, 6.3)to (2.6,6.3);
\node [morphism, scale=\lnscalevertex] at (2.25,5.6) {$\beta$};
\node [scale=\lnscaledot] at (0.85,-0.1) {$\ldots$};
\node [scale=\lnscaledot] at (2.25,5) {$\ldots$};
\node [scale=\lnscaledot] at (2.25,6.2) {$\ldots$};
\node [scale=\lnscaledotperspective, rotate = 135, thick] at (3,-0.1) {$\ldots$};
\end{tikzpicture}
\end{aligned}
&
\stackrel{\II_{k-1}}{\rightarrow}
&
\begin{aligned}
\begin{tikzpicture}[thick, scale=0.8]
\draw [region] (0,-1.4) rectangle +(2.9,6.5);
\draw (0.5,3.5) to (0.5,5.1);
\draw (1.2,3.5) to (1.2,5.1);
\draw (1.9,0) to(1.9,0) to [out=up, in=down] (1.2,0.7) to [out=up, in=down] (1.2,1.4) ;
\draw (2.6,0)  to  [out=up, in=down] (2.6,0.7) to[out=up, in=down](1.9,1.4);
\draw [dotted](1.2,1.4)  to (1.2,2.1);
\draw [dotted](1.9,1.4)  to (1.9,2.1);
\draw (1.2,2.1)  to [out=up, in=down] (0.5,2.8) to (0.5,3.5);
\draw (1.9,2.1) to [out=up, in=down] (1.9,2.8) to [out=up, in=down] (1.2,3.5) ;
\draw (1.9,-1.4) to(1.9,-1.4) to [out=up, in=down] (2.25,-0.7) to [out=up, in=down](1.9, 0);
\draw (2.6,-1.4) to(2.6,-1.4) to [out=up, in=down] (2.25,-0.7) to [out=up, in=down](2.6, 0);
\node [morphism, scale=\lnscalevertex] at (2.25,-0.7) {$\alpha$};
\node [scale=\lnscaledot] at (0.85,5) {$\ldots$};
\node [scale=\lnscaledot] at (2.25,-0.1) {$\ldots$};
\node [scale=\lnscaledotperspective, rotate = 135, thick] at (0.1,5) {$\ldots$};
\node [scale=\lnscaledotperspective, rotate = 135, thick] at (0.1,-1.5) {$\ldots$};
\node [scale=\lnscaledotperspective, rotate = 135, thick] at (3,5) {$\ldots$};
\node [scale=\lnscaledot] at (2.25,-1.3) {$\ldots$};
\draw [region] (0.2, -1.6) rectangle +(2.9,6.5);
\draw (0.5,-1.6) to(0.5,0);
\draw (1.2,-1.6) to(1.2,0.0);
\draw (0.5,0) to (0.5,0.7) to [out=up, in=down] (0.5,1.4) ;
\draw (1.2,0.0)  to [out=up, in=down] (1.9,0.7) to [out=up, in=down] (2.6,1.4);
\draw [dotted](2.6,1.4)  to (2.6,2.1);
\draw [dotted](0.5,1.4)  to (0.5,2.1);
\draw (0.5,2.1) to [out=up, in=down] (1.2,2.8) to [out=up, in=down] (1.9,3.5) ;
\draw (2.6,2.1)  to [out=up, in=down] (2.6,2.8) to (2.6,3.5);
\draw (1.9,3.5) to [out=up, in=down] (2.25,4.2) to [out=up, in=down](1.9, 4.9) to (1.9,4.9);
\draw (2.6,3.5) to [out=up, in=down] (2.25,4.2) to [out=up, in=down](2.6, 4.9)to (2.6,4.9);
\node [morphism, scale=\lnscalevertex] at (2.25,4.2) {$\beta$};
\node [scale=\lnscaledot] at (0.85,-1.5) {$\ldots$};
\node [scale=\lnscaledot] at (2.25,3.6) {$\ldots$};
\node [scale=\lnscaledot] at (2.25,4.8) {$\ldots$};
\node [scale=\lnscaledotperspective, rotate = 135, thick] at (3,-1.5) {$\ldots$};
\end{tikzpicture}
\end{aligned}
&
\stackrel{\widehat{{\I}}_{k-1}}{\rightarrow}
&
\begin{aligned}
\begin{tikzpicture}[thick, scale=0.8]
\draw [region] (0,-1.4) rectangle +(2.9,6.5);
\draw (1.9,0) to [out=up, in=down] (1.2,0.7) to [out=up, in=down] (0.5,1.4) ;
\draw (2.6,0)  to (2.6,1.4);
\draw [dotted](2.6,1.4)  to (2.6,2.1);
\draw [dotted](0.5,1.4)  to (0.5,2.1);
\draw (0.5,2.1)  to (0.5,3.5);
\draw (2.6,2.1) to [out=up, in=down] (1.9,2.8) to [out=up, in=down] (1.2,3.5) ;
\draw (0.5,3.5) to (0.5,5.1);
\draw (1.2,3.5) to (1.2,5.1);
\draw (1.9,-1.4) to(1.9,-1.4) to [out=up, in=down] (2.25,-0.7) to [out=up, in=down](1.9, 0);
\draw (2.6,-1.4) to(2.6,-1.4) to [out=up, in=down] (2.25,-0.7) to [out=up, in=down](2.6, 0);
\node [morphism, scale=\lnscalevertex] at (2.25,-0.7) {$\alpha$};
\node [scale=\lnscaledot] at (0.85,5) {$\ldots$};
\node [scale=\lnscaledot] at (2.25,-0.1) {$\ldots$};
\node [scale=\lnscaledotperspective, rotate = 135, thick] at (0.1,5) {$\ldots$};
\node [scale=\lnscaledotperspective, rotate = 135, thick] at (0.1,-1.5) {$\ldots$};
\node [scale=\lnscaledotperspective, rotate = 135, thick] at (3,5) {$\ldots$};
\node [scale=\lnscaledot] at (2.25,-1.3) {$\ldots$};
\draw [region] (0.2, -1.6) rectangle +(2.9,6.5);
\draw (0.5,-1.6) to(0.5,0);
\draw (1.2,-1.6) to(1.2,0.0);
\draw (0.5,0) to (0.5,0.7) to [out=up, in=down] (1.2,1.4) ;
\draw (1.2,0.0)  to [out=up, in=down] (1.9,0.7) to (1.9,1.4);
\draw [dotted](1.2,1.4)  to (1.2,2.1);
\draw [dotted](1.9,1.4)  to (1.9,2.1);
\draw (1.2,2.1) to (1.2,2.8) to [out=up, in=down] (1.9,3.5);
\draw (1.9,2.1)  to [out=up, in=down] (2.6,2.8) to (2.6,3.5);
\draw (1.9,3.5) to [out=up, in=down] (2.25,4.2) to [out=up, in=down](1.9, 4.9) to (1.9,4.9);
\draw (2.6,3.5) to [out=up, in=down] (2.25,4.2) to [out=up, in=down](2.6, 4.9)to (2.6,4.9);
\node [morphism, scale=\lnscalevertex] at (2.25,4.2) {$\beta$};
\node [scale=\lnscaledot] at (0.85,-1.5) {$\ldots$};
\node [scale=\lnscaledot] at (2.25,3.6) {$\ldots$};
\node [scale=\lnscaledot] at (2.25,4.8) {$\ldots$};
\node [scale=\lnscaledotperspective, rotate = 135, thick] at (3,-1.5) {$\ldots$};
\end{tikzpicture}
\end{aligned}
&
\stackrel{\II_{k-1}}{\rightarrow}
&
\begin{aligned}
\begin{tikzpicture}[thick, scale=0.8]
\draw [region] (0,0) rectangle +(2.9,6.5);
\draw (1.9,2.8) to [out=up, in=down] (1.2,3.5) to [out=up, in=down] (0.5,4.2) ;
\draw (2.6,2.8)  to  [out=up, in=down] (2.6,3.5) to[out=up, in=down](2.6,4.2);
\draw [dotted](2.6,4.2)  to (2.6,4.9);
\draw [dotted](0.5,4.2)  to (0.5,4.9);
\draw (0.5,4.9)  to [out=up, in=down] (0.5,5.6) to (0.5,6.3) to (0.5,6.5);
\draw (2.6,4.9) to [out=up, in=down] (1.9,5.6) to [out=up, in=down] (1.2,6.3) to (1.2,6.5);
\draw (1.9,0) to(1.9,0) to [out=up, in=down] (2.25,0.7) to [out=up, in=down](1.9, 1.4)to(1.9,2.8);
\draw (2.6,0) to(2.6,0) to [out=up, in=down] (2.25,0.7) to [out=up, in=down](2.6, 1.4)to(2.6,2.8);
\node [morphism, scale=\lnscalevertex] at (2.25,0.7) {$\alpha$};
\node [scale=\lnscaledot] at (2.25,0.1) {$\ldots$};
\node [scale=\lnscaledot] at (2.25,1.3) {$\ldots$};
\node [scale=\lnscaledot] at (0.85,6.4) {$\ldots$};
\node [scale=\lnscaledotperspective, rotate = 135, thick] at (0.1,6.4) {$\ldots$};
\node [scale=\lnscaledotperspective, rotate = 135, thick] at (0.1,-0.1) {$\ldots$};
\node [scale=\lnscaledotperspective, rotate = 135, thick] at (3,6.4) {$\ldots$};
\draw [region] (0.2, -0.2) rectangle +(2.9,6.5);
\draw (0.5,2.8) to [out=up, in=down] (0.5,3.5) to [out=up, in=down] (1.2,4.2) ;
\draw (1.2,2.8)  to [out=up, in=down] (1.9,3.5) to [out=up, in=down] (1.9,4.2);
\draw [dotted](1.2,4.2)  to (1.2,4.9);
\draw [dotted](1.9,4.2)  to (1.9,4.9);
\draw (1.2,4.9) to [out=up, in=down] (1.2,5.6) to [out=up, in=down] (1.9,6.3) to (1.9,6.3);
\draw (1.9,4.9)  to [out=up, in=down] (2.6,5.6) to (2.6,6.3)to (2.6,6.3);
\draw (0.5,-0.2) to(0.5,1.4) to [out=up, in=down] (0.85,2.1) to [out=up, in=down](0.5, 2.8);
\draw (1.2,-0.2) to(1.2,1.4) to [out=up, in=down] (0.85,2.1) to [out=up, in=down](1.2, 2.8);
\node [morphism, scale=\lnscalevertex] at (0.85,2.1) {$\beta$};
\node [scale=\lnscaledot] at (0.85,-0.1) {$\ldots$};
\node [scale=\lnscaledot] at (0.85,1.5) {$\ldots$};
\node [scale=\lnscaledot] at (0.85,2.7) {$\ldots$};
\node [scale=\lnscaledot] at (2.25,6.2) {$\ldots$};
\node [scale=\lnscaledotperspective, rotate = 135, thick] at (3,-0.1) {$\ldots$};
\end{tikzpicture}
\end{aligned}
\end{array}
\end{equation*}
\end{definition}

\paragraph{Type I homotopy generator rearrangement scheme.} The type V homotopy generators require a rearrangement scheme for the type I homotopy generators. To see why this is required, consider the picture below:
\begin{align*}
\begin{aligned}
\begin{tikzpicture}[thick, scale=0.6]
\draw [region] (0,-1.4) rectangle +(2.9,6.5);
\draw (0.5,-1.4) to(0.5,0);
\draw (1.2,-1.4) to(1.2,0.0);
\draw (0.5,0) to (0.5,0.7) to [out=up, in=down] (0.5,1.4) ;
\draw (1.2,0.0)  to [out=up, in=down] (1.9,0.7) to [out=up, in=down] (2.6,1.4);
\draw [dotted](2.6,1.4)  to (2.6,2.1);
\draw [dotted](0.5,1.4)  to (0.5,2.1);
\draw (0.5,2.1) to [out=up, in=down] (1.2,2.8) to [out=up, in=down] (1.9,3.5) ;
\draw (2.6,2.1)  to [out=up, in=down] (2.6,2.8) to (2.6,3.5);
\draw (1.9,3.5) to [out=up, in=down] (2.25,4.2) to [out=up, in=down](1.9, 4.9) to (1.9,5.1);
\draw (2.6,3.5) to [out=up, in=down] (2.25,4.2) to [out=up, in=down](2.6, 4.9)to (2.6,5.1);
\node [morphism, scale=1] at (2.25,4.2) {$\alpha$};
\node [scale=0.6] at (0.85,-1.3) {$\ldots$};
\node [scale=0.6] at (2.25,3.6) {$\ldots$};
\node [scale=0.6] at (2.25,5) {$\ldots$};
\node [scale=0.6] at (2.25,-1.5) {$\ldots$};
\node [scale=0.4, rotate = 135, thick] at (3,-1.5) {$\ldots$};
\draw [region] (0.2, -1.6) rectangle +(2.9,6.5);
\draw (0.5,3.5) to (0.5,4.9);
\draw (1.2,3.5) to (1.2,4.9);
\draw (1.9,0) to(1.9,0) to [out=up, in=down] (1.2,0.7) to [out=up, in=down] (1.2,1.4) ;
\draw (2.6,0)  to  [out=up, in=down] (2.6,0.7) to[out=up, in=down](1.9,1.4);
\draw [dotted](1.2,1.4)  to (1.2,2.1);
\draw [dotted](1.9,1.4)  to (1.9,2.1);
\draw (1.2,2.1)  to [out=up, in=down] (0.5,2.8) to (0.5,3.5);
\draw (1.9,2.1) to [out=up, in=down] (1.9,2.8) to [out=up, in=down] (1.2,3.5) ;
\draw (1.9,-1.6) to(1.9,-1.4) to [out=up, in=down] (2.25,-0.7) to [out=up, in=down](1.9, 0);
\draw (2.6,-1.6) to(2.6,-1.4) to [out=up, in=down] (2.25,-0.7) to [out=up, in=down](2.6, 0);
\node [morphism, scale=1] at (2.25,-0.7) {$\beta$};
\node [scale=0.6] at (0.85,4.8) {$\ldots$};
\node [scale=0.6] at (2.25,-0.1) {$\ldots$};
\node [scale=0.4, rotate = 135, thick] at (0.1,5) {$\ldots$};
\node [scale=0.4, rotate = 135, thick] at (0.1,-1.5) {$\ldots$};
\node [scale=0.4, rotate = 135, thick] at (3,5) {$\ldots$};
\end{tikzpicture}
\end{aligned}
\end{align*}
The generator $\beta$ can be `pulled  up' by a $\II^{-1}$ move, but the generator $\alpha$ cannot be `pulled down' by a $\II$ move, since the type I generators, drawn as crossings, are not correctly arranged. To remedy this, we need a rearrangement scheme for type I homotopy generators. We denote this scheme  ${\widehat \I}$, and construct it recursively as follows.
\def\thetaexpansionscale{0.5}
\begin{definition}
\label{defthetainterchanger}
In an $(n+1)$\-signature $\sigma$ that supports interchangers of type $\I_k$, a \emph{rearrangement $\widehat{{\I_k}}$} consists of a sequence of applications of the composite interchanger of type $\widetilde{\I_k}$, with values defined recursively as follows:
\begin{align*}
&\begin{aligned}
\begin{tikzpicture}[thick, scale=\thetaexpansionscale]
\draw [region] (-0.7,0) rectangle +(4.3,7.2);
\draw (1.2,0.0)  to [out=up, in=down] (1.9,0.7) to [out=up, in=down] (2.6,1.4) to [out=up, in=down] (3.3,2.1) to  (3.3,7.2);
\draw (0.5,0) to (0.5,2.1) to [out=up, in=down](1.2,2.8) to [out=up, in=down] (1.9,3.5)  to [out=up, in=down] (2.6,4.2) to  (2.6,7.2);
\draw (-0.2,0) to (-0.2,4.9) to [out=up, in=down](0.5,5.6) to [out=up, in=down] (1.2,6.3)  to [out=up, in=down] (1.9,7) to  (1.9,7.2) ;
\node [scale=0.6] at (0.15,0.1) {$\ldots$};
\node [scale=0.6] at (2.25,7.1) {$\ldots$};
\node [scale=0.4, rotate = 135, thick] at (3.7,-0.1) {$\ldots$};
\draw [region] (-0.5, -0.2) rectangle +(4.3,7.2);
\draw (1.9,-0.2) to(1.9,0) to [out=up, in=down] (1.2,0.7) to [out=up, in=down] (1.2,1.4)  to  (1.2,2.1) to [out=up, in=down] (0.5,2.8)to(0.5,4.2) ;
\draw (2.6,-0.2) to  (2.6,0.7) to [out=up, in=down] (1.9,1.4) to  (1.9,2.8) to [out=up, in=down] (1.2,3.5)to(1.2,4.2) ;
\draw (3.3,-0.2) to  (3.3,1.4) to [out=up, in=down] (2.6,2.1) to  (2.6,3.5) to [out=up, in=down] (1.9,4.2);
\draw [dotted](0.5,4.2)  to (0.5,4.9);
\draw [dotted](1.2,4.2)  to (1.2,4.9);
\draw [dotted](1.9,4.2)  to (1.9,4.9);
\draw (1.2,4.9) to  (1.2,5.6) to [out=up, in=down] (0.5,6.3) to  (0.5,7);
\draw (0.5,4.9)  to [out=up, in=down] (-0.2,5.6) to  (-0.2,7);
\draw (1.9,4.9) to  (1.9,6.3) to [out=up, in=down] (1.2,7);
\node [scale=0.6] at (2.95,-0.1) {$\ldots$};
\node [scale=0.6] at (0.15,6.9) {$\ldots$};
\node [scale=0.4, rotate = 135, thick] at (-0.6,7.1) {$\ldots$};
\node [scale=0.4, rotate = 135, thick] at (-0.6,-0.1) {$\ldots$};
\node [scale=0.4, rotate = 135, thick] at (3.7,7.1) {$\ldots$};
\end{tikzpicture}
\end{aligned}
\stackrel{\widehat{{\I}}_k}{\rightarrow}
\begin{aligned}
\begin{tikzpicture}[thick, scale=\thetaexpansionscale]
\draw [region] (-0.7,0) rectangle +(4.3,7.2);
\draw (2.6,4.9) to  [out=up, in=down] (3.3,5.6) to  (3.3,7.2);
\draw (1.2,0.0)  to [out=up, in=down] (1.9,0.7) to (1.9,2.1) to[out=up, in=down] (2.6,2.8) to (2.6,4.2) ;
\draw  (1.9,4.9)  to(1.9,5.6)  to[out=up, in=down] (2.6,6.3) to  (2.6,7.2);
\draw (0.5,0) to (0.5,0.7) to [out=up, in=down](1.2,1.4)to (1.2,2.8) to [out=up, in=down] (1.9,3.5)  to (1.9,4.2) ;
\draw (1.2,4.9)  to (1.2,6.3)  to [out=up, in=down] (1.9,7) to  (1.9,7.2);
\draw (-0.2,0) to (-0.2,1.4) to [out=up, in=down](0.5,2.1) to(0.5,3.5) to [out=up, in=down] (1.2,4.2); 
\draw [dotted](2.6,4.2)  to (2.6,4.9);
\draw [dotted](1.9,4.2)  to (1.9,4.9);
\draw [dotted](1.2,4.2)  to (1.2,4.9);
\node [scale=0.6] at (0.15,0.1) {$\ldots$};
\node [scale=0.6] at (2.25,7.1) {$\ldots$};
\node [scale=0.4, rotate = 135, thick] at (3.7,-0.1) {$\ldots$};
\draw [region] (-0.5, -0.2) rectangle +(4.3,7.2);
\draw (1.9,-0.2) to(1.9,0) to [out=up, in=down] (1.2,0.7)   to [out=up, in=down] (0.5,1.4) to [out=up, in=down] (-0.2,2.1)to(-0.2,7) ;
\draw (2.6,-0.2) to  (2.6,2.1) to [out=up, in=down] (1.9,2.8) to [out=up, in=down] (1.2,3.5) to [out=up, in=down] (0.5,4.2)to(0.5, 7) ;
\draw (3.3,-0.2) to  (3.3,4.9) to [out=up, in=down] (2.6,5.6)  to [out=up, in=down] (1.9,6.3) to [out=up, in=down]  (1.2,7);
\node [scale=0.6] at (2.25,-0.1) {$\ldots$};
\node [scale=0.6] at (0.15,6.9) {$\ldots$};
\node [scale=0.4, rotate = 135, thick] at (-0.6,7.1) {$\ldots$};
\node [scale=0.4, rotate = 135, thick] at (-0.6,-0.1) {$\ldots$};
\node [scale=0.4, rotate = 135, thick] at (3.7,7.1) {$\ldots$};
\end{tikzpicture}
\end{aligned}
\\
&\quad:=
\begin{aligned}
\begin{tikzpicture}[thick, scale=\thetaexpansionscale]
\draw [region] (-0.7,0) rectangle +(4.3,7.2);
\draw (1.2,0.0)  to [out=up, in=down] (1.9,0.7) to [out=up, in=down] (2.6,1.4) to [out=up, in=down] (3.3,2.1) to  (3.3,7.2);
\draw (0.5,0) to (0.5,2.1) to [out=up, in=down](1.2,2.8) to [out=up, in=down] (1.9,3.5)  to [out=up, in=down] (2.6,4.2) to  (2.6,7.2);
\draw (-0.2,0) to (-0.2,4.9) to [out=up, in=down](0.5,5.6) to [out=up, in=down] (1.2,6.3)  to [out=up, in=down] (1.9,7) to  (1.9,7.2) ;
\node [scale=0.6] at (0.15,0.1) {$\ldots$};
\node [scale=0.6] at (2.25,7.1) {$\ldots$};
\node [scale=0.4, rotate = 135, thick] at (3.7,-0.1) {$\ldots$};
\draw [region] (-0.5, -0.2) rectangle +(4.3,7.2);
\draw (1.9,-0.2) to(1.9,0) to [out=up, in=down] (1.2,0.7) to [out=up, in=down] (1.2,1.4)  to  (1.2,2.1) to [out=up, in=down] (0.5,2.8)to(0.5,4.2) ;
\draw (2.6,-0.2) to  (2.6,0.7) to [out=up, in=down] (1.9,1.4) to  (1.9,2.8) to [out=up, in=down] (1.2,3.5)to(1.2,4.2) ;
\draw (3.3,-0.2) to  (3.3,1.4) to [out=up, in=down] (2.6,2.1) to  (2.6,3.5) to [out=up, in=down] (1.9,4.2);
\draw [dotted](0.5,4.2)  to (0.5,4.9);
\draw [dotted](1.2,4.2)  to (1.2,4.9);
\draw [dotted](1.9,4.2)  to (1.9,4.9);
\draw (1.2,4.9) to  (1.2,5.6) to [out=up, in=down] (0.5,6.3) to  (0.5,7);
\draw (0.5,4.9)  to [out=up, in=down] (-0.2,5.6) to  (-0.2,7);
\draw (1.9,4.9) to  (1.9,6.3) to [out=up, in=down] (1.2,7);
\node [scale=0.6] at (2.95,-0.1) {$\ldots$};
\node [scale=0.6] at (0.85,6.9) {$\ldots$};
\node [scale=0.4, rotate = 135, thick] at (-0.6,7.1) {$\ldots$};
\node [scale=0.4, rotate = 135, thick] at (-0.6,-0.1) {$\ldots$};
\node [scale=0.4, rotate = 135, thick] at (3.7,7.1) {$\ldots$};
\end{tikzpicture}
\end{aligned}
\stackrel{\widetilde{\I}_k}{\rightarrow}
\begin{aligned}
\begin{tikzpicture}[thick, scale=\thetaexpansionscale]
\draw [region] (-0.7,0) rectangle +(4.3,7.2);
\draw (1.2,0.0)  to [out=up, in=down] (1.9,0.7) to (1.9,1.4) to[out=up, in=down] (2.6,2.1) to [out=up, in=down] (3.3,2.8) to  (3.3,7.2);
\draw (0.5,0) to (0.5,0.7) to [out=up, in=down](1.2,1.4)to (1.2,2.8) to [out=up, in=down] (1.9,3.5)  to [out=up, in=down] (2.6,4.2) to  (2.6,7.2);
\draw (-0.2,0) to (-0.2,4.9) to [out=up, in=down](0.5,5.6) to [out=up, in=down] (1.2,6.3)  to [out=up, in=down] (1.9,7) to  (1.9,7.2) ;
\node [scale=0.6] at (0.15,0.1) {$\ldots$};
\node [scale=0.6] at (2.25,7.1) {$\ldots$};
\node [scale=0.4, rotate = 135, thick] at (3.7,-0.1) {$\ldots$};
\draw [region] (-0.5, -0.2) rectangle +(4.3,7.2);
\draw (1.9,-0.2) to(1.9,0) to [out=up, in=down] (1.2,0.7)   to [out=up, in=down] (0.5,1.4)to(0.5,4.2) ;
\draw (2.6,-0.2) to  (2.6,1.4) to [out=up, in=down] (1.9,2.1) to  (1.9,2.8) to [out=up, in=down] (1.2,3.5)to(1.2,4.2) ;
\draw (3.3,-0.2) to  (3.3,2.1) to [out=up, in=down] (2.6,2.8) to  (2.6,3.5) to [out=up, in=down] (1.9,4.2);
\draw [dotted](0.5,4.2)  to (0.5,4.9);
\draw [dotted](1.2,4.2)  to (1.2,4.9);
\draw [dotted](1.9,4.2)  to (1.9,4.9);
\draw (1.2,4.9) to  (1.2,5.6) to [out=up, in=down] (0.5,6.3) to  (0.5,7);
\draw (0.5,4.9)  to [out=up, in=down] (-0.2,5.6) to  (-0.2,7);
\draw (1.9,4.9) to  (1.9,6.3) to [out=up, in=down] (1.2,7);
\node [scale=0.6] at (2.25,-0.1) {$\ldots$};
\node [scale=0.6] at (0.15,6.9) {$\ldots$};
\node [scale=0.4, rotate = 135, thick] at (-0.6,7.1) {$\ldots$};
\node [scale=0.4, rotate = 135, thick] at (-0.6,-0.1) {$\ldots$};
\node [scale=0.4, rotate = 135, thick] at (3.7,7.1) {$\ldots$};
\end{tikzpicture}
\end{aligned}
\stackrel{\widetilde{\I}_k}{\rightarrow}
\begin{aligned}
\begin{tikzpicture}[thick, scale=\thetaexpansionscale]
\draw [region] (-0.7,0) rectangle +(4.3,7.2);
\draw (1.2,0.0)  to [out=up, in=down] (1.9,0.7) to (1.9,2.1) to[out=up, in=down] (2.6,2.8) to [out=up, in=down] (3.3,3.5) to  (3.3,7.2);
\draw (0.5,0) to (0.5,0.7) to [out=up, in=down](1.2,1.4)to (1.2,3.5) to [out=up, in=down] (1.9,4.2)  to [out=up, in=down] (2.6,4.9) to  (2.6,7.2);
\draw (-0.2,0) to (-0.2,1.4) to [out=up, in=down](0.5,2.1) to(0.5,5.6) to [out=up, in=down] (1.2,6.3)  to [out=up, in=down] (1.9,7) to  (1.9,7.2) ;
\node [scale=0.6] at (0.15,0.1) {$\ldots$};
\node [scale=0.6] at (2.25,7.1) {$\ldots$};
\node [scale=0.4, rotate = 135, thick] at (3.7,-0.1) {$\ldots$};
\draw [region] (-0.5, -0.2) rectangle +(4.3,7.2);
\draw (1.9,-0.2) to(1.9,0) to [out=up, in=down] (1.2,0.7)   to [out=up, in=down] (0.5,1.4) to [out=up, in=down] (-0.2,2.1)to(-0.2,7) ;
\draw (2.6,-0.2) to  (2.6,2.1) to [out=up, in=down] (1.9,2.8) to  (1.9,3.5) to [out=up, in=down] (1.2,4.2)to(1.2,4.9) ;
\draw (3.3,-0.2) to  (3.3,2.8) to [out=up, in=down] (2.6,3.5) to  (2.6,4.2) to [out=up, in=down] (1.9,4.9);
\draw [dotted](1.2,4.9)  to (1.2,5.6);
\draw [dotted](1.9,4.9)  to (1.9,5.6);
\draw (1.2,5.6) to  (1.2,5.6) to [out=up, in=down] (0.5,6.3) to  (0.5,7);
\draw (1.9,5.6) to  (1.9,6.3) to [out=up, in=down] (1.2,7);
\node [scale=0.6] at (2.25,-0.1) {$\ldots$};
\node [scale=0.6] at (0.15,6.9) {$\ldots$};
\node [scale=0.4, rotate = 135, thick] at (-0.6,7.1) {$\ldots$};
\node [scale=0.4, rotate = 135, thick] at (-0.6,-0.1) {$\ldots$};
\node [scale=0.4, rotate = 135, thick] at (3.7,7.1) {$\ldots$};
\end{tikzpicture}
\end{aligned}
\stackrel{\widehat{{\I}}_k}{\rightarrow}
\begin{aligned}
\begin{tikzpicture}[thick, scale=\thetaexpansionscale]
\draw [region] (-0.7,0) rectangle +(4.3,7.2);
\draw (2.6,4.9) to  [out=up, in=down] (3.3,5.6) to  (3.3,7.2);
\draw (1.2,0.0)  to [out=up, in=down] (1.9,0.7) to (1.9,2.1) to[out=up, in=down] (2.6,2.8) to (2.6,4.2) ;
\draw  (1.9,4.9)  to(1.9,5.6)  to[out=up, in=down] (2.6,6.3) to  (2.6,7.2);
\draw (0.5,0) to (0.5,0.7) to [out=up, in=down](1.2,1.4)to (1.2,2.8) to [out=up, in=down] (1.9,3.5)  to (1.9,4.2) ;
\draw (1.2,4.9)  to (1.2,6.3)  to [out=up, in=down] (1.9,7) to  (1.9,7.2) ;
\draw (-0.2,0) to (-0.2,1.4) to [out=up, in=down](0.5,2.1) to(0.5,3.5) to [out=up, in=down] (1.2,4.2); 
\draw [dotted](2.6,4.2)  to (2.6,4.9);
\draw [dotted](1.9,4.2)  to (1.9,4.9);
\draw [dotted](1.2,4.2)  to (1.2,4.9);
\node [scale=0.6] at (0.15,0.1) {$\ldots$};
\node [scale=0.6] at (2.25,7.1) {$\ldots$};
\node [scale=0.4, rotate = 135, thick] at (3.7,-0.1) {$\ldots$};
\draw [region] (-0.5, -0.2) rectangle +(4.3,7.2);
\draw (1.9,-0.2) to(1.9,0) to [out=up, in=down] (1.2,0.7)   to [out=up, in=down] (0.5,1.4) to [out=up, in=down] (-0.2,2.1)to(-0.2,7) ;
\draw (2.6,-0.2) to  (2.6,2.1) to [out=up, in=down] (1.9,2.8) to [out=up, in=down] (1.2,3.5) to [out=up, in=down] (0.5,4.2)to(0.5, 7) ;
\draw (3.3,-0.2) to  (3.3,4.9) to [out=up, in=down] (2.6,5.6)  to [out=up, in=down] (1.9,6.3) to [out=up, in=down]  (1.2,7);
\node [scale=0.6] at (2.25,-0.1) {$\ldots$};
\node [scale=0.6] at (0.15,6.9) {$\ldots$};
\node [scale=0.4, rotate = 135, thick] at (-0.6,7.1) {$\ldots$};
\node [scale=0.4, rotate = 135, thick] at (-0.6,-0.1) {$\ldots$};
\node [scale=0.4, rotate = 135, thick] at (3.7,7.1) {$\ldots$};
\end{tikzpicture}
\end{aligned}
\end{align*}
\end{definition}

\paragraph{Expansion scheme.} We do not  need an expansion scheme for type V homotopy generators.

\paragraph{Extended variant.} The type $\V'$ homotopy generators include 2 further variants of the moves shown in \eqref{eq:basicpullthroughs}, corresponding to pulling a vertex through an inverse type I homotopy generator, as well as the corresponding composition schemes.

\subsection{Type VI homotopy generators}

These are required for the definition of quasistrict 4\-category. They involve reflected versions of the following move:
\begin{equation*}
\begin{array}{ccc}
\begin{aligned}
\begin{tikzpicture}[thick, scale=0.8, yscale=0.8]
\draw [region] (-0.2,-1) rectangle +(2.4,3);
\draw (1.6,-1) to [out=up, in=down] (0.6,1) to (0.6, 2);
\begin{scope}[yshift=-0.2cm, xshift = 0.2cm]
\draw [region] (-0.2,-1) rectangle +(2.4,3);
\draw (0.4,-1) to (0.4,-0.8) to [out=up, in=down] (1.4,1) to (1.4, 2);
\end{scope}
\node [morphism, anchor=north] at (1.6,1.4) {$\alpha$};
\end{tikzpicture}
\end{aligned}
&\stackrel{\II}{\rightarrow}&
\begin{aligned}
\begin{tikzpicture}[thick, scale=0.8, yscale=0.8]
\draw [region] (-0.2,-1) rectangle +(2.4,3);
\draw (1.6,-1) to (1.6,0) to [out=up, in=down] (0.6,1.8) to (0.6, 2);
\begin{scope}[yshift=-0.2cm, xshift = 0.2cm]
\draw [region] (-0.2,-1) rectangle +(2.4,3);
\draw (0.4,-1) to (0.4,0.2) to [out=up, in=down] (1.4,2);
\end{scope}
\node [morphism, anchor=north] at (0.6,0) {$\alpha$};
\end{tikzpicture}
\end{aligned}
\\
\downarrow{\I'}&
\begin{aligned}
\begin{tikzpicture}[thick, scale=0.5]
\node [rotate=0, scale = 0.7] at (-1,0) {$$};
\node [rotate=90, scale = 1.5] at (0,0) {$\Leftarrow $};
\node [rotate=90, scale = 1.5] at (0.8,0) {$\Rightarrow $};
\node [rotate=0, scale = 0.7] at (1.8,0) {$$};
\end{tikzpicture}
\end{aligned}
&\uparrow{ \I' }
\\
\begin{aligned}
\begin{tikzpicture}[thick, scale=0.8, yscale=0.8]
\draw [region] (-0.2,-3) rectangle +(2.4,6.8);
\draw (1.6,-3) to(1.6,-2.8) to [out=up, in=down](0.6,-1)to (0.6, 0) to [out=up, in=down] (1.6,1.8)to [out=up, in=down](0.6,3.6) to (0.6,3.8);
\begin{scope}[yshift=-0.2cm, xshift = 0.2cm]
\draw [region] (-0.2,-3) rectangle +(2.4,6.8);
\draw (0.4,-3) to(0.4,-2.6) to [out=up, in=down](1.4,-0.8) to (1.4,0.2) to [out=up, in=down] (0.4,2) to (0.4, 2)to [out=up, in=down](1.4,3.8);
\end{scope}
\node [morphism, anchor=north] at (1.6,-0.3) {$\alpha$};
\end{tikzpicture}
\end{aligned}
&\stackrel{ \II^{-1} }{\rightarrow}&
\begin{aligned}
\begin{tikzpicture}[thick, scale=0.8, yscale=0.8]
\draw [region] (-0.2,-3) rectangle +(2.4,6.8);
\draw (1.6,-3) to(1.6,-2.8) to [out=up, in=down](0.6,-1)to [out=up, in=down](1.6, 0.8) to [out=up, in=down] (1.6,1.8)to [out=up, in=down](0.6,3.6) to (0.6,3.8);
\begin{scope}[yshift=-0.2cm, xshift = 0.2cm]
\draw [region] (-0.2,-3) rectangle +(2.4,6.8);
\draw (0.4,-3) to(0.4,-2.6) to [out=up, in=down](1.4,-0.8)[out=up, in=down] to (0.4,1) to [out=up, in=down] (0.4,2) to (0.4, 2)to [out=up, in=down](1.4,3.8);
\end{scope}
\node [morphism, anchor=north] at (0.6,1.5) {$\alpha$};
\end{tikzpicture}
\end{aligned}
\end{array}
\end{equation*}

\subsection{Comparison with 4-teisi}
\label{fourtasaxioms}

A previous definition of semistrict 4\-category has been given by Crans~\cite{Crans_1998}, using the terminology of 4\-teisi. Here we show in detail that a semistrict 4\-category in his sense gives rise to a semistrict 4\-category in our sense, modulo some issues which we believe are inaccuracies in his definition:
\begin{itemize}
\item there is no equivalent of our type IV homotopy generators (the `Breenator');
\item there is no recognition for the need for interchanger rearrangements, written $\widehat \I$ in our formalism, which play a key part in the description of our type $\V$ homotopy generators.
\end{itemize}
Aside from these two major points, the definitions are in correspondence.

\paragraph{Detailed comparison.} We list the individual structures and axioms that Crans proposes, and comment on how they relate to the signature structure and homotopy generators that form our definition.

A $4$\-dimensional tas consists of collections $C_0$ of objects, $C_1$ of arrows, $C_2$ of $2$\-arrows, $C_3$ of $3$\-arrows and $C_4$ of $4$\-arrows, together with:
\begin{enumerate}[(S1)]

\item \textit{Functions $s_n, t_n: C_i\to C_n$ for all $0\leq n<i\leq 4$ also denoted $d_n$ and ${d_n}^+$ and called $n$\-source and $n$\-target.}

The collection of sets of arrows (cells) and the source and target functions correspond exactly to those in Definition~\ref{defsignature}.

\item \textit{Functions $\#_n: C_{n+1}\times_{s_n, t_n}C_{n+1}\to C_{n+1}$ for all $0\leq n< 4$ called vertical composition.}

This is realised by (vertical) composition of two $n$\-diagrams $S,D$, as described in Definition~\ref{diagcompose}.

\item \textit{Functions $\#_n: C_{i}\times_{s_n, t_n}C_{n+1}\to C_{i}$ and $\#_m: C_{n+1}\times_{s_n, t_n}C_{i}\to C_{i}$  for all $0\leq n\leq 2$, $n+1<i<4$, called whiskering.}

This is realised by composition of an $n$\-diagram $D$ and an $m$\-diagram $S$ such that $m\neq n$, as described in Definition ~\ref{diagcompose}.

\item \textit{Functions $\#_n: C_{q}\times_{s_n, t_n}C_{p}\to C_{p+q}$ for all $0\leq n\leq 1$, $p,q>n+1$, $p+q-n-1\leq 4$  called horizontal composition.}

In our approach this is not a built-in operation, horizontal $n$\-cell composition is instead achieved by first whiskering and then vertical composition.

\item \textit{Functions $\id: C_{i}\to C_{i+1}$ for all $0\leq n\leq 3$  called identity.}

This is realised by the identity operation $\diagid{D}$ on an $n$\-diagram $D$, as described in Definition~\ref{defidiagid}.
\end{enumerate}

\noindent
Crans also postulates eleven axioms for these structures:
\begin{enumerate}[(R1)]

\item \textit{$\mathbb{C}$ is a $4$\-truncated globular set.}

We understand that $\mathbb{C}=\{C_0, C_1, C_2, C_3, C_4\}$. The corresponding requirement is realised by the equalities imposed on the source and target maps in Definition~\ref{defsignature}.

\item \textit{For every $C, C'\in C_0$ the collection of elements of $\mathbb{C}$ with $0$\-source $C$ and $0$\-target $C'$ forms a 3D tas $\mathbb{C}(C, C')$, with $n$\-composition in $\mathbb{C}(C, C')$ given by $\#_{n+1}$ and identities given by $\id$.}

Crans requires that for every $C,C'\in C_0$, the collection of elements of $\mathbb{C}$ whose lowest lever source is $C$ and lowest level target $C'$ is a 3D tas correspond to the requirement that $\sigma$ supports higher\-level instances of all singularities supported by a signature presenting a quasistrict $3$\-category. The difference here is that we list the interchangers explicitly, whereas for Crans they are built into the definition of the lower level structure.
 
\item \textit{For every $g:C'\to C''$ in $C_1$ and every $C$ and $C'''\in C_0$, $-\#_0 g$ is a functor $\mathbb{C}(C'',C''')\to \mathbb{C}(C',C''')$ and $g\#_0 - $ is a functor $\mathbb{C}(C'',C''')\to \mathbb{C}(C',C'')$.}

This is a statement about composition of $1$\-cells with $0$\-cells, Crans requires this operation to give rise to functors. In our approach to compose a 0\-cell $A$ with 1\-cell $g$, the identity 1\-cell on $A$ needs to be produced first. Then, these two are composed by the ordinary composition operation according to Definition~\ref{diagcompose}. We also explicitly prove that this gives rise to functors in Theorem~\ref{assocdisttheorem}.

\item \textit{For every $C\in C_0$ we have: $s_0(\id_C)=C=t_0(\id_c)$.}

This requirement means that sources and targets of a cell that has undergone a transformation by the identity map are equal to the cell itself. In our approach this is a direct consequence of Definition~\ref{identityembedding}.

\item \textit{For every $C'\in C_0$ and every $C$ and $C''\in C_0$, $-\#_0\id_{C'}$ is equal to the identity functor $\mathbb{C}(C',C'')\to \mathbb{C}(C',C'')$ and $\id_{C'}\#_0 -$ is equal to the identity functor $\mathbb{C}(C,C')\to \mathbb{C}(C,C')$}.

This is a statement about the relation between the identity operation and $0$\-composition of $1$\-cells and that it gives rise to certain functors. In our approach, this is built into the definition of a signature and is captured by Lemma~\ref{compositionidempotent}, where we prove that this construction gives rise to functors.
\item
\begin{enumerate}[(a)]
\item \textit{For every $\gamma:f\to f', \delta:g\to g' \in C_2$
\begin{align*}
s_2(\delta\#_0\gamma) &= (g'\#_0\gamma)\#_1(\delta\#_0 f)\\
t_2(\delta\#_0\gamma) &= (\delta\#_0f')\#_1(g\#_0 \gamma)
\end{align*}
And $\delta\#_0\gamma$ is an iso\-3\-arrow.}

This interchanger is a $3$\-cell that corresponds to type $\I_{2}$, as described in Definition~\ref{defswap}.

\item \textit{For every $\phi: \gamma\to\gamma'\in C_3$ such that $\gamma, \gamma': f\to f'$ and $\delta:g\to g'\in C_2$
\begin{align*}
s_3(\delta\#_0\phi) &= ((\delta\#_0f')\#_1(g\#_0 \phi))\#_2(\delta\#_0\gamma)\\
t_3(\delta\#_0\phi) &= (\delta\#_0\gamma')\#_2((g'\#_0 \phi)\#_1(\delta\#_0 f))
\end{align*}
And $\delta\#_0\gamma$ is an iso\-3\-arrow.}

This interchanger is a $4$\-cell that corresponds to type ${\II}_{3}$, as described in Definition~\ref{defpullthrough}.

\item \textit{For every $\gamma:f\to f'\in C_2$ and $\psi: \delta\to\delta'\in C_3$ such that $\delta, \delta': g\to g'$
\begin{align*}
s_3(\psi\#_0\gamma) &= (\delta\#_0\gamma')\#_2((g'\#_0 \gamma)\#_1(\psi\#_0 f))\\
t_3(\psi\#_0\gamma) &= ((\psi\#_0 f')\#_1(g\#_0 \gamma))\#_2(\delta\#_0\gamma)
\end{align*}
And $\delta\#_0\gamma$ is an iso\-3\-arrow.}

This interchanger is a $4$\-cell that corresponds to type ${\II}_{3}$, as described in Definition~\ref{defpullthrough}. It is a different variant than in the clause above.
\end{enumerate}
\item 
\begin{enumerate} [(a)]
\item \textit{For every $\Gamma: \phi\to \phi'\in C_4$ such that $\phi,\phi':\gamma\to\gamma'$ and $\gamma,\gamma':f\to f'$ and $\delta: g\to g'\in C_2$:}
\begin{align*}
&(((g'\#_0\Gamma)\#_1(\delta\#_0f))\#_2(\delta\#_0\gamma'))\#_3(\delta\#_0\phi)\\
&\quad=(\delta\#_0\phi')\#_3((\delta\#_0\gamma)\#_2((\delta\#_0f')\#_1(g\#_0\Gamma)))
\end{align*}
This interchanger is a $5$\-cell that corresponds to type ${\III}_{4}$.
\item \textit{For every  $\gamma: f\to f'\in C_2$ and $\Delta: \psi\to \psi'\in C_4$ such that $\psi,\psi':\delta\to\delta'$ and $\delta,\delta':g\to g'$}
\begin{align*}
&(\psi'\#_0\gamma)\#_3((\delta'\#_0\gamma)\#_2((g'\#_0\gamma)\#_1(\Delta\#_0f))) \\
&\quad=(((\Delta\#_0f')\#_1(g\#_0\gamma))\#_2(\delta\#_0\gamma))\#_3(\psi\#_0\gamma)\\
\end{align*} 
This interchanger is a $5$\-cell that corresponds to type ${\III}_{4}$.
\item  \textit{For every $\phi: \gamma\to\gamma'\in C_3$ such that $\gamma, \gamma': f\to f'$ and $\psi: \delta\to\delta'\in C_3$ such that $\delta, \delta': g\to g'$}
\begin{align*}
&(((\psi\#_0f')\#_1(g\#_0\gamma))\#_2(\delta\#_0\phi))\#_3 \\
&(((\psi\#_0f')\#_1(g\#_0\phi))\#_2(\delta\#_0\gamma))\#_3 \\
&(((\delta'\#_0f')\#_1(g\#_0\phi))\#_2(\psi\#_0\gamma))\#_3 \\
&\qquad\qquad\qquad\qquad=\\
&(\psi\#_0\gamma')\#_2((g'\#_0\phi)\#_1(\delta\#_0 f)))\#_3\\
&((\delta'\#_0\gamma')\#_2{((g'\#_0\phi)\#_1(\psi\#_0 f))}^{-1})\#_3 \\
&((\delta'\#_0\phi)\#_2((g'\#_0\gamma)\#_1(\psi\#_0 f)))
\end{align*}
This interchanger is a $5$\-cell that corresponds to type ${\V}_{4}$.
\end{enumerate}

\item In this clause, Crans discusses expansion morphisms for different singularities. It is not clear to us whether the intention is to make them actual equalities or higher level morphisms. If the former, then some details pertaining reorganisation of crossings for interchangers of type $\V_4$ are missing from the definition. If the latter, then this would introduce unnecessary weakness, in violation of the goal of a semistrict definition.
\begin{enumerate} [(a)]
\item \textit{For every $\gamma: f\to f', \gamma':f'\to f''$ and $\delta:g\to g'$ in $\mathbb{C}$:}
\begin{align*}
&\delta\#_0(\gamma'\#_1\gamma) = \\ &\qquad((\delta\#_0\gamma')\#_1(g\#_0\gamma))\#_2((g'\#_0\gamma')\#_1(\delta\#_0\gamma))
\end{align*} 
This corresponds to a 4\-cell decomposing an application of composite interchanger of type $\I_2$ into applications of atomic interchangers. This is not a native structure in our approach. However we could still efficiently reason about decompositions of this type, due to the presence of a shorthand $\widetilde{\I_2}$ as described in Definition~\ref{defgammainterchanger}. 

\item \textit{For every $\gamma:f\to f'$ and $\delta: g\to g', \delta':g'\to g''$ in $\mathbb{C}$:}
 \begin{align*}
(\delta'\#_1\delta)\#_0\gamma = ((\delta'\#_0f')\#_1(\delta\#_0\gamma))\#_2((\delta'\#_0\gamma)\#_1(\delta\#_0 f))
\end{align*} 
Similar as above, but for the other variant of the interchanger of type $\I_2$.

\item \textit{For every $\phi:\gamma\to\gamma', \phi':\gamma'\to\gamma''$ such that $\gamma, \gamma', \gamma'': f\to f'$ and $\delta: g\to g'$ in $\mathbb{C}$:}
\begin{align*}
&\delta\#_0(\gamma''\#_2\phi) = \\
&\qquad((\delta\#_0\phi')\#_2((g'\#_0\phi)\#_1(\delta\#_0 f)))\\
&\qquad\qquad\qquad\qquad\#_3 \\
&\qquad (((\delta\#_0 f)\#_1(g\#_0 \phi'))\#_2(\delta\#_0\phi))
\end{align*} 
This is an expansion axiom for interchangers of type ${\II}_{3}$. Intuitively, this says that we can pull-through individual subsequent 3\-cells through a crossing or pull them through all at once in one move. This is not a built-in structure in our approach, however this does not result in reduced expressivity. Analogously to $\widetilde{\I_2}$ for interchangers of type $\I$, this cell is realised by $\widetilde{\II_3}$ as a composite of applications of atomic interchangers of type $\II_3$, as defined in~\ref{pullthroughexpansion}.

\item For every $\phi:\gamma\to\gamma'$ such that $\gamma,\gamma':f\to f'$, for every $\gamma'':f'\to f''$ and $\delta:g\to g'$  in $\mathbb{C}$:
\begin{align*}
&\delta\#_0(\gamma'''\#_2\phi) \\ &\quad=(((\delta\#_0\gamma'')\#_1(g\#_0\gamma'))\#_2((g'\#_0\gamma'')\#_1(\delta\#_0\phi)))\\
&\qquad\#_3(((\delta\#_0\gamma'')\#_1(g\#_0\phi))\#_2((g'\#_0\gamma'')\#_1(\delta\#_0\gamma)))
\end{align*} 
This corresponds to an application of a composite interchanger of type $\II_3$, where some additional crossings have to be interchanged out of the way first using interchangers of type $\I_3$. This is not a native cell in our definition, but similarly as above, we could express it as a sequence of atomic interchangers using the shorthand $\widetilde{\II_3}$.

\item For every $\phi':\gamma'\to\gamma''$ such that $\gamma',\gamma'':f'\to f''$, for every $\gamma:f\to f'$ and $\delta:g\to g'$ in $\mathbb{C}$:
\begin{align*}
&\delta\#_0(\phi'\#_2\gamma) \\ &\quad=(((\delta\#_0\gamma'')\#_1(g\#_0\gamma))\#_2((g\#_0\phi')\#_1(\delta\#_0\gamma)))\\
&\qquad\#_3(((\delta\#_0\phi')\#_1(g\#_0\gamma))\#_2((g'\#_0\gamma')\#_1(\delta\#_0\gamma)))
\end{align*}  
Similarly as above, but for another variant of the interchanger of type $\II_3$.

\item For every $\gamma, \gamma':f\to f'$, for every $\delta:g\to g;$ and $\delta':g'\to g''$ in $\mathbb{C}$:
\begin{align*}
&(\delta'\#_1\delta)\#_0\phi \\
&\quad=(((\delta'\#_0 f')\#_1(\delta'\#_0\gamma'))\#_2((\delta'\#_0\phi)\#_1(\delta\#_0 f)))\\
&\qquad\#_3(((\delta'\#_0 f')\#_1(\delta'\#_0\phi))\#_2((\delta'\#_0\gamma')\#_1(\delta\#_0 f)))
\end{align*}  
This is another expansion axiom for interchangers of type ${\II}_{3}$. By this, we could either pull a single 3\-cell through multiple crossings individually one by one, or we could do this all at once. The same as other expansions above, this is not a native cell in our approach. Again, this does not result in reduced expressivity, as this cell can be realised as a composite of atomic interchangers of type $\II_3$, as defined in~\ref{pullthroughexpansion}. 

\item Crans remarks that there should be additional clauses analogous to the last four, most likely this is to deal with the inverses of the interchanger variants discussed above.
\end{enumerate}

\item For every $\gamma:f\to f'$, $\gamma':f'\to f''$ and for every $\delta:g\to g;$, $\delta':g'\to g''$ in $\mathbb{C}$:
\begin{align*}
&(\delta'\#_1\delta)\#_0(\gamma'\#_1\gamma) \\
&\quad=((\delta'\#_0 f'')\#_1(\delta\#_0\gamma')\#_1(g\#_0\gamma)) \\
&\qquad \#_2 ((\delta'\#_0\gamma')\#_1(\delta\#_0\gamma)) \\
&\quad\qquad \#_2((g''\#_0\gamma')\#_1(\delta'\#_0\gamma)\#_1(\delta\#_0 f)) 
\end{align*}
This is a higher\-level coherence which is the result of the expansion axiom for an interchanger of type $\I_2$. There are two different, equivalent methods of decomposition for swapping heights of four adjacent cells and they have to be related to each other by a higher\-level morphism. This is a direct evidence for how explicit expansion morphisms give rise to further singularities. The presence of an explicit expansion cell for $\I_2$ suggests that perhaps the intention in the definition was to include higher\-level expansion cells explicitly as well.

\item \textit{For every $c\in \mathbb{C}{(C, C')}_p, c'\in \mathbb{C}{(C', C'')}_q, c''\in \mathbb{C}{(C'', C''')}_r$ with $p+q+r\leq 3$ we have: $(c''\#_0 c')\#_0 c = c''\#_0(c'\#_0 c)$.}

This axiom is on associativity of cell composition. The corresponding result is summarised by equality~\myeqref{eqassoc}.

\item \textit{For every $c\in \mathbb{C}{(C, C')}_p, c'\in \mathbb{C}{(C', C'')}_q$ such that $p,q>0$ and $p+q\leq3$ if $ q\leq2$ we have: $c'\#_0\id_c = \id_{c'\#_0 c}$ and if $p\leq2$ we have: $c'\#_0\id_c = \id_{c'\#_0 c}$.}

This axiom defines how to compose a cell with an identity of a cell of a lower dimension. The corresponding result is proved in Lemma~\ref{boostingcomposition}.

\end{enumerate}




\subsection{Comparison with Gray categories}
\label{sectionswitchthreecategory}

In this section we discuss the definition of a switch 3\-category given by Douglas and Henriques~\cite{douglashenriques},  which is an alternative presentation of a \textbf{Gray}-category that has the same generic-position properties as our definition. We analyze each of their structures and axioms turn, showing in detail how a switch 3\-category gives rise to a semistrict 3\-category by our definition.

\begin{itemize}
\item\textbf{S1-1} For a $0$\-diagram $D$, the map $i_x$ is realised by:  $i_x(D):= \diagid{D}$

\item\textbf{S1-2} For two $1$\-diagrams $S$, $D$ the map $m_x$ is realised by $m_x(S,D):= S\circ D$

\item\textbf{S1-3} For a $1$\-diagram $D$, the map $i_x$ is realised by:  $i_x(D):= \diagid{D}$

\item\textbf{S1-4} For two $2$\-diagrams $S$, $D$ the map $m_y$ is realised by $m_y(S,D):= S\circ D$

\item\textbf{S1\-5} For a $1$\-diagram $D$ and a $2$\-diagram $S$, the map $w_r$ is realised by $w_r(S,D):= S\circ D$

\item\textbf{S1\-6} For a $2$\-diagram $D$ and a $1$\-diagram $S$, the map $w_l$ is realised by $w_l(S,D):= S\circ D$

\item\textbf{S1\-7} This is realised by interchangers of type $\I_2$ supported by $\sigma$.

\item\textbf{S1\-8} For a $2$\-diagram $D$, the map $i_z$ is realised by:  $i_z(D):= \diagid{D}$

\item\textbf{S1\-9} For two $3$\-diagrams $S$, $D$ the map $m_z$ is realised by $m_z(S,D):= S\circ D$

\item\textbf{S1\-10} For a $3$\-diagram $D$ and a $2$\-diagram $S$, the map $f_b$ is realised by $f_b(S,D):= S\circ D$

\item\textbf{S1\-11} For a $2$\-diagram $D$ and a $3$\-diagram $S$, the map $f_t$ is realised by $f_t(S,D):= S\circ D$

\item\textbf{S1\-12} For a $3$\-diagram $D$ and a $1$\-diagram $S$, the map $h_r$ is realised by $h_r(S,D):= S\circ D$

\item\textbf{S1\-13} For a $1$\-diagram $D$ and a $3$\-diagram $S$, the map $h_l$ is realised by $h_l(S,D):= S\circ D$

\item\textbf{S2-1} Given that map $m_x$ is realised by diagram composition and map $i_x$ by the identity operation, we need to show that  $\diagid{S}\circ D = D$ for any $1$\- diagram $D$ and any $0$\-diagram $S$ such that $\diagid{S}\circ D$ exists. This follows by Lemma~\ref{compositionidempotent}.

\item\textbf{S2-2} The argument is analogous to \textbf{S2-1}.

\item\textbf{S2\-3} Given that map $m_x$ is realised by diagram composition, we need to show that $S\circ (D\circ M)= (S\circ D)\circ M$ for any three $1$\-diagrams $S, D, M$ that are composable. This follows by equality~\myeqref{eqassoc}.

\item\textbf{S2-4} Given that map $m_y$ is realised by diagram composition and map $i_y$ by the identity operation, we need to show that  $\diagid{S}\circ D = D$ for any $2$\- diagram $D$ and any $1$\-diagram $S$ such that $\diagid{S}\circ D$ exists. This follows by Lemma~\ref{compositionidempotent}.

\item\textbf{S2-5} The argument is analogous to \textbf{S2-4}.

\item\textbf{S2\-6} Given that $m_y$ is realised by diagram composition, we need to show that $S\circ (D\circ M)= (S\circ D)\circ M$ for any three $2$\-diagrams $S, D, M$ that are composable. This follows by equality~\myeqref{eqassoc}.

\item\textbf{S2\-7} Given that map $w_r$ is realised by diagram composition, we need to show that for any two 
$1$\-diagrams $S$, $D$ that are composable, we have: $\diagid{S}\circ D = \diagid{S\circ D}$. This follows by Lemma~\ref{boostingcomposition}.

\item\textbf{S2-8} The argument is analogous to \textbf{S2-7}.

\item\textbf{S2\-9} Given that maps $m_y$ and $w_r$ are realised by diagram composition,  we need to show that $(S\circ D)\circ M= (S\circ M)\circ (D\circ M)$ for any two $2$\-diagrams $S, D$ and a $1$\-diagram $M$ that are composable. This follows by equality~\myeqref{eqdistright}.


\item\textbf{S2-10} The argument is analogous to \textbf{S2-9} and the result follows by equality~\myeqref{eqdistleft}.

\item\textbf{S2\-11} Given that maps $m_x$ and $w_l$ are realised by diagram composition, we need to show that $(S\circ D)\circ M = S\circ (D\circ M)$ for any two $1$\-diagrams $D, M$ and a $2$\-diagram $S$ that are composable. This follows by Theorem~\ref{eqassoc}.

\item\textbf{S2-12} The argument is analogous to \textbf{S2-11}.

\item\textbf{S2\-13} Given that maps $w_l$ and $w_r$ are realised by diagram composition, we need to show that $(S\circ D)\circ M = S\circ (D\circ M)$ for any two $1$\-diagrams $S, M$ and a $2$\-diagram $D$ that are composable. This follows by Theorem~\ref{eqassoc}.

\item\textbf{S2\-14} Given that map $m_y$ is realised by diagram composition and map $i_x$ by the identity operation, we need to show that  $\diagid{S}\circ D = D$ for any $2$\- diagram $D$ and any $0$\-diagram $S$ such that $\diagid{S}\circ D$ exists. This follows by Lemma~\ref{compositionidempotent}.

\item\textbf{S2-15} The argument is analogous to \textbf{S2-14}.

\item\textbf{S2-16} Given that the `switch' map corresponds to interchangers of type $\I_2$, this holds by Definition~\ref{defswap}. As application of an interchanger of type $\I$ at height $h=1$ for $|D|=1$ has no effect on $D$.

\item\textbf{S2-17} This is an expansion axiom for the `switch' morphism, it is not present in Definition~\ref{defsemistrictthreecategory}. For reasons outlined earlier in this section, an explicit expansion morphism for composite interchangers does not add any additional expressivity. In $\sigma$, a composite interchanger of type $\I_2$ can be expressed as a sequence of atomic interchangers of the same type using the construction $\widetilde{\I_2}$ described in Definition~\ref{defgammainterchanger}.

\item\textbf{S2-18} This follows by the result on associativity of composition~\ref{eqassoc} and by the fact that $\sigma$ supports interchangers of type $\I_2$. Since we only interpret the interchanger in the context of the digram $D$ at the particular height $h$, the order in which any additional 1\-cells on either side are composed in does not alter the overall interchange, as required.

\item\textbf{S2-19} The argument is analogous to \textbf{S2-18}, but for 1\-cells that are in between the 2\-cells being interchnaged. By Definition~\ref{defswap}, there may be an arbitrary number of them and the order in which they are composed with 2\-cells has no effect on the interchange, as required.

\item\textbf{S2-20} Given that map $m_z$ is realised by diagram composition and map $i_z$ by the identity operation, we need to show that  $S\circ \diagid{D} = S$ for any $3$\- diagram $S$ and any $2$\-diagram $D$ such that $S\circ \diagid{D}$ exists. This follows by Lemma~\ref{compositionidempotent}.

\item\textbf{S2\-21} Given that map $m_z$ is realised by diagram composition, we need to show that $S\circ (D\circ M)= (S\circ D)\circ M$ for any three $3$\-diagrams $S, D, M$ that are composable. This follows by equality~\myeqref{eqassoc}.

\item\textbf{S2\-22} Given that map $f_t$ is realised by diagram composition, we need to show that for any two $2$\-diagrams $S$, $D$ that are composable, we have: $\diagid{S}\circ D = \diagid{S\circ D}$. This follows by Lemma~\ref{boostingcomposition}.

\item\textbf{S2-23} Given that map $f_t$ is realised by diagram composition and map $i_y$ by the identity operation, we need to show that  $S\circ \diagid{D} = S$ for any $3$\- diagram $S$ and any $1$\-diagram $D$ such that $S\circ \diagid{D}$ exists. This follows by Lemma~\ref{compositionidempotent}.

\item\textbf{S2\-24} Given that maps $m_z$ and $f_t$ are realised by diagram composition,  we need to show that $(S\circ D)\circ M= (S\circ M)\circ (D\circ M)$ for any two $3$\-diagrams $S, D$ and a $2$\-diagram $M$ that are composable. This follows by equality~\myeqref{eqdistright}.

\item\textbf{S2\-25} This is an instance of two 3\-cells being subject to the interchange law, similarly as \textbf{S1-7} for 2\-cells. This axiom holds, since $\sigma$ supports interchangers of type $\I_3$.

\item\textbf{S2\-26} Given that maps $m_y$ and $f_t$ are realised by diagram composition, we need to show that $ S\circ (D\circ M) = (S\circ D)\circ M$ for any two $2$\-diagrams $D, M$ and a $3$\-diagram $S$ that are composable. This follows by equality~\myeqref{eqassoc}.

\item\textbf{S2\-27} Given that maps $f_t$ and $f_b$ are realised by diagram composition, we need to show that $(S\circ D)\circ M = S\circ (D\circ M)$ for any two $2$\-diagrams $S, M$ and a $3$\-diagram $D$ that are composable. This follows by equality~\myeqref{eqassoc}.

\item\textbf{S2\-28} Given that map $h_r$ is realised by diagram composition, we need to show that for a $2$\-diagram $S$, and a $1$\-diagram $D$  that are composable, we have: $\diagid{S}\circ D = \diagid{S\circ D}$. This follows by Lemma~\ref{boostingcomposition}.

\item\textbf{S2-29} Given that map $h_r$ is realised by diagram composition and map $i_x$ by the identity operation, we need to show that  $S\circ \diagid{D} = S$ for any $3$\- diagram $S$ and any $0$\-diagram $D$ such that $S\circ \diagid{D}$ exists. This follows by Lemma~\ref{compositionidempotent}.

\item\textbf{S2\-30} This is an instance of naturality of the switch morphism (\textbf{S1-7}) in one of its inputs and corresponds to interchangers of type $\II$. This axiom holds, since $\sigma$ supports interchangers of type $\II_3$.

\item\textbf{S2\-31} Given that maps $m_z$ and $h_r$ are realised by diagram composition,  we need to show that $(S\circ D)\circ M= (S\circ M)\circ (D\circ M)$ for any two $3$\-diagrams $S, D$ and a $1$\-diagram $M$ that are composable. This follows by equality~\myeqref{eqdistright}.

\item\textbf{S2\-32} Given that maps $f_t$ and $h_r$ are realised by diagram composition,  we need to show that $(S\circ D)\circ M= (S\circ M)\circ (D\circ M)$ for a $3$\-diagrams $S$, a $2$\-diagram $D$ and a $1$\-diagram $M$ that are composable. This follows by equality~\myeqref{eqdistright}.

\item\textbf{S2\-33} Given that maps $m_x$ and $h_r$ are realised by diagram composition, we need to show that $ S\circ (D\circ M) = (S\circ D)\circ M$ for any two $1$\-diagrams $D, M$ and a $3$\-diagram $S$ that are composable. This follows by Theorem~\ref{eqassoc}.

\item\textbf{S2\-34} Given that maps $h_l$ and $h_r$ are realised by diagram composition, we need to show that $(S\circ D)\circ M = S\circ (D\circ M)$ for any two $1$\-diagrams $S, M$ and a $3$\-diagram $D$ that are composable. This follows by Theorem~\ref{eqassoc}.

\end{itemize}
All additional axis flips required by Definition~\cite{DouglasHenriques} follow in an analogous way.


\section{Core proofs}
\label{sec:core}

This is the technical core of the paper. Here we prove the propositions listed in \ref{sec:properties}, and some further propositions. We split the material as follows: Section~\ref{sec:rewriting} handles correctness of diagram rewriting and auxiliary notions, Section~\ref{sec:identityembeddings} handles identity embeddings, Section~\ref{sec:diagramcomposition} handles correctness of diagram composition; Section~\ref{sec:identitydiagrams} handles identity diagrams, and Section~\ref{secassocdist} handles associativity and distributivity of diagram composition.

\paragraph{Example proof.} As an example of our methods, we give here an example proof: rewriting a diagram $D$ by removing a subdiagram $e : S \hookrightarrow D$, and replacing an identical subdiagram $S$, yields $D$ again. In this derivation, as well as those elsewhere in the article, for each equality we refer explicitly to the statement that justifies it. If the reference is absent and we do not make an exceptional statement immediately before the derivation, then the equality follows by arithmetic manipulations.
\begin{lemma}
\label{vacuousrewrite}
Given an $n$-diagram $D$ and $e: S\hookrightarrow D$, the following equality holds:
\begin{align*}
\rewrite{D}{e}{S}{S} = D
\end{align*}
\end{lemma}
\begin{proof}
By Definition~\ref{diagequiv}, we need to check four conditions.

\paragraph{Sources are equivalent diagrams.}
\begin{align*}
&\rewrite{D}{e}{S}{S}.s= D.s &[\myeqref{defrewritesource}]
\end{align*}

\paragraph{Sizes of lists of generators are equal.}
\begin{align*}
&|\rewrite{D}{e}{S}{S}|\\
&\quad= |D| - |S| + |S|  &[\myeqref{defrewritesize}] \\
&\quad= |D| 
\end{align*}

\paragraph{Corresponding generators are equal.} We must consider this for $0\leq j\leq |D|$. We consider this separately for three ranges:

\begin{itemize}
\item For $0\leq j\leq e.h$:
\begin{align*}
&\rewrite{D}{e}{S}{S}[j].g = D[j].g&[\myeqref{defrewritegenerator}]
\end{align*}
\item For $e.h\leq j\leq e.h + |S|$:
\begin{align*}
&\rewrite{D}{e}{S}{S}[j].g \\
&\quad= S[j - e.h].g && [\myeqref{defrewritegenerator}]\\
&\quad= D[(j - e.h) + e.h].g && [\myeqref{eqembwelldef3}]\\
&\quad= D[j].g
\end{align*}

\item For $e.h + |S|\leq j\leq |D|$:
\begin{align*}
&\rewrite{D}{e}{S}{S}[j].g \\
&\quad= D[j - |S| + |S|].g & [\myeqref{defrewritegenerator}]\\
&\quad= D[j].g
\end{align*}
\end{itemize}

\paragraph{Corresponding embeddings are equivalent.} We must consider this for $0\leq j\leq |D|$. We consider this separately for three ranges:

\begin{itemize}
\item For $0\leq j\leq e.h$:
\begin{align*}
&\rewrite{D}{e}{S}{S}[j].e = D[j].e&[\myeqref{defrewriteembedding}]
\end{align*}

\item For $e.h\leq j\leq e.h + |S|$:
\begin{align*}
&\rewrite{D}{e}{S}{S}[j].e \\
&\quad= (\lift{e.e}{S[j - e.h].d}) \circ S[j-e.h].e && [\myeqref{defrewriteembedding}]\\
&\quad= D[(j-e.h) + e.h].e&& [\myeqref{eqembwelldef3}] \\
&\quad= D[j].e
\end{align*}

\item For $e.h + |S|\leq j\leq |D|$:
\begin{align*}
&\rewrite{D}{e}{S}{S}[j].e \\
&\quad= D[j - |S| + |S|].e & [\myeqref{defrewriteembedding}]\\
&\quad= D[j].e
\end{align*}
\end{itemize}
Hence $\rewrite{D}{e}{S}{S}$ and $D$ are indeed equivalent.
\end{proof}

\subsection{Rewriting}
\label{sec:rewriting}

The aim here is to prove that the rewrite construction produces a well-defined diagram when its data is well-defined. Along the way, we also show that the lift construction and embedding composition also produce well-defined structures. Due to the heavily recursive structure of our definitions, we introduce several logical statements that are dependent on the diagram dimension $n$: $P(n), R(n), S(n), T(n), Q(n), A(n), B(n), C(n)$. All these statements concern properties of $n$\-diagrams and embeddings between $n$\-diagrams. 

Of particular importance is the statement $R(n)$, which formalises the fact that a rewrite of a well-defined $n$\-diagram is also well-defined. $B(n)$ and $C(n)$ state respectively that a lifted embedding between well-defined $n$\-diagrams is well-defined and that a composite of two well-defined embeddings between $n$\-diagrams is also well-defined. The remaining statements $P(n), S(n), T(n), Q(n), A(n)$ are auxiliary statements that are necessary to carry out an inductive proof that $R(n)$ holds for all $n\geq 0$, \emph{i.e.} that the process of rewriting preserves the diagram property of being well-defined for any diagram. 

The intuition behind the proof is that in order to say that a rewritten $n$\-diagram $\rewrite{D}{e}{S}{T}$ is well-defined ($R(n)$), we need all its slices, which are $(n{-}1)$\-diagrams, to be well-defined. To achieve that, we express them as slices of $D$ or as rewrites of a well-defined $(n{-}1)$\-diagram $D[e.h].d$ ($S(n)$, $T(n)$) and use the inductive hypothesis $R(n{-}1)$. For the inductive hypothesis to be applicable, the component embeddings $\rewrite{D}{e}{S}{T}[j].e: \rewrite{D}{e}{S}{T}[j].g.s\hookrightarrow \rewrite{D}{e}{S}{T}[j].d$, which are embeddings of $(n{-}1)$\-diagrams, must be well-defined. In the segment of $\rewrite{D}{e}{S}{T}$ between $e.h$ and $e.h+|T|$, these embeddings are built using lifted embeddings (well-defined by $B(n{-}1)$) and embedding composition (well-defined by $C(n{-}1)$), as per Definition~\ref{defrewrite}. The rest follows by the inductive hypothesis $R(n{-}1)$. The base case is straightforward.

\paragraph{The statements.}
The main statement concerns well-definedness of a rewritten $n$\-diagram. 
\begin{definition}
\label{defR}
For $n\geq 0$, $R(n)$ states that given well-defined $n$-diagrams $D, S, T$ such that $S$, $T$ are globular, and a well-defined embedding $e: S\hookrightarrow D$, then $\rewrite{D}{e}{S}{T}$ is a well-defined diagram.
\end{definition}

\noindent
The second statement concerns globularity of individual slices in a well-defined diagram. This is a necessary condition to be able to consider them as rewrites of a well-defined $(n{-}1)$\-diagram.
\begin{definition}
\label{defT}
For $n\geq 2$, $T(n)$ states that given a well-defined $n$\-diagram $D$, we have $D.s.s = D[j].d.s$ and $D.s.t = D[j].d.t$ for any $0\leq j < |D|$.
\end{definition}

\noindent
The next statement allows us to express slices of a rewritten diagram in an explicit way instead of depending on the recursive definition.
\begin{definition}
\label{defS}
For $n\geq 1$, $S(n)$ states that given well-defined $n$\-diagrams $D, S, T$ such that $S,T$ are globular, and a well-defined embedding $e: S\hookrightarrow A$, the following hold:
\begin{align*}
&\rewrite{A}{e}{S}{T}[j].d=\begin{cases}
A[j].d & \text{if } 0\leq j \leq e.h \\
\rewrite{A[e.h].d}{e.e}{S.s}{T[j-e.h].d}   & \text{if }e.h \leq j \leq e.h + |T| \\
A[j + |S| - |T|].d &  \text{if }e.h + |T| \leq j < |A| - |S| + |T| \\
\end{cases}
\end{align*}
\end{definition}

The next two auxiliary statements let us respectively, decompose a lifted embedding into a composite of two individual lifted embeddings and express a composite rewrite in two equivalent manners. These are both indirectly used to establish the statement $S(n)$ on expressing slices of a rewritten $n$\-diagram as rewrites of a well-defined $(n{-}1)$\-diagram.
\begin{definition}
\label{defQ}
For $n \geq 0$, $Q(n)$ states that given well-defined $n$\-diagrams   $A$, $B$, $C$, $S$, $T$ such that the pairs $S,T$ and $A,C$ are globular, and given well-defined embeddings $e:S\hookrightarrow A$, $f:C\hookrightarrow B$, the following holds:
\[
\lift{(\lift{f}{A}\circ e)}{T} = \lift{f}{\rewrite{A}{e}{S}{T}}\circ\lift{e}{T}
\]
\end{definition}

\begin{definition}
\label{defP}
For $n\geq 0$, $P(n)$ states that given well-defined $n$\-diagrams  $A$, $B$, $C$, $S$, $T$ such that pairs $S,T$ and $A,C$ are globular, and given well-defined embeddings $e:S\hookrightarrow A$ and $f:C\hookrightarrow B$, the following holds:
\[
\rewrite{B}{f}{C}{\rewrite{A}{e}{S}{T}}= \rewrite{\rewrite{B}{f}{C}{A}}{\lift{f}{A}\circ e}{S}{T}
\]
\end{definition}

The auxiliary concept of a lifted embedding is used to construct the rewritten diagram, and must itself be a well-defined embedding.

\begin{definition}
\label{defB}
For $n\geq 0$, $B(n)$ states that given well-defined $n$\-diagrams   $S, T, A$ such that  $S,T$ are globular, and given a well-defined embedding $e:S\hookrightarrow A$, the lifted embedding $\lift{e}{T}: T\hookrightarrow \rewrite{A}{e}{S}{T}$ is well-defined.
\end{definition}

\noindent
The composite of two well-defined embeddings must also be well-defined.

\begin{definition}
\label{defC}
For $n\geq 0$,  $C(n)$ states that given well\-defined $n$\-diagrams $S,D,M$, and given well-defined $n$\-diagram embeddings $e: S\hookrightarrow D$ and $f: D\hookrightarrow M$, their composite $f\circ e: S\hookrightarrow M$ is well-defined.
\end{definition}

\noindent
The final statement is associativity of embedding composition.
\begin{definition}
\label{defA}
For $n\geq 0$, $A(n)$ states that given well\-defined $n$\-diagram embeddings $e: S\hookrightarrow D$, $f: D\hookrightarrow M$, $g: M\hookrightarrow N$ between well-defined $n$\-diagrams $S,D,M, N$ the following equality holds:
\[
g\circ(f\circ e) = (g\circ f) \circ e
\]
\end{definition}

\paragraph{Structure of the argument.}
To make the following exposition easier to follow, we summarise the main results that are used in the inductive step of the proof of Theorem~\ref{bigtheorem}, which states the conjunction of all the above statements holds for all $n\geq 0$. These implications are as follows:
\begin{itemize}
\item $S(n)\wedge T(n)\wedge R(n{-}1) \implies R(n)$, proved in Lemma~\ref{SimpliesR}.
\item $T(n) \wedge P(n{-}1)\wedge C(n{-}1)\wedge B(n{-}1)\implies S(n)$, proved in Lemma~\ref{TPimpliesS}.
\item $S(n{-}1)\implies T(n)$, proved in Lemma~\ref{SimpliesT}.
\item $S(n) \wedge P(n)\wedge A(n{-}1)\implies Q(n)$, proved in Lemma~\ref{SPAimpliesQ}.
\item $S(n) \wedge Q(n{-}1)\implies P(n)$, proved in Lemma~\ref{SQimpliesP}.
\item $R(n)\implies B(n)$, proved in Lemma~\ref{RimpliesL}.
\item $C(n{-}1)\wedge Q(n{-}1)\wedge A(n{-}1)\implies C(n)$, proved in Lemma~\ref{CQAimpliesC}.
\item $S(n) \wedge Q(n{-}1) \wedge A(n{-}1)\implies A(n)$, proved in Lemma~\ref{SQAimpliesA}.
\end{itemize}
We visualize these implications with the following diagram, which can be used to verify that the overall mutual induction is well-founded:
\tikzset{conjunction/.style={}}
\newcommand\onefold[2]{
\draw [->] (#1) to (#2);
}
\newcommand\twofold[3]{
\node [conjunction] (1) at ($($(#1)!0.5!(#2)$)!0.3333!(#3)$) {.};
\draw [->,] (1.center) to (#3);
\draw [-] (#1) to (1.center);
\draw [-] (#2) to (1.center);
}
\newcommand\twofoldnew[5]{
\ignore{
\node [conjunction] (1) at (#4) {#5};
\draw [->,] (1.center) to (#3);
\draw [-] (#1) to (1.center);
\draw [-] (#2) to (1.center);
}
\draw [->] (#1) to (#3);
\draw [->] (#2) to (#3);
}
\newcommand\threefold[4]{
\node [conjunction] (1) at ($($(#1)!0.5!(#2)$)!0.5!($(#3)!0.5!(#4)$)$) {};
\draw [->] (1.center) to (#4);
\draw [-] (#1) to (1.center);
\draw [-] (#2) to (1.center);
\draw [-] (#3) to (1.center);
}
\newcommand\threefoldnew[6]{
\ignore{
\node [conjunction] (1) at (#5) {#6};
\draw [->] (1.center) to (#4);
\draw [-] (#1) to (1.center);
\draw [-] (#2) to (1.center);
\draw [-] (#3) to (1.center);
}
\draw [->] (#1) to (#4);
\draw [->] (#2) to (#4);
\draw [->] (#3) to (#4);
}
\newcommand\fourfold[5]{
\node [conjunction] (1) at ($($($(#1)!0.5!(#2)$)!0.5!($(#3)!0.5!(#4)$)$)!0.2!(#5)$) {};
\draw [->] (1.center) to (#5);
\draw [-] (#1) to (1.center);
\draw [-] (#2) to (1.center);
\draw [-] (#3) to (1.center);
\draw [-] (#4) to (1.center);
}
\newcommand\fourfoldnew[7]{
\ignore{
\node [conjunction] (1) at (#6) {#7};
\draw [->] (1.center) to (#5);
\draw [-] (#1) to (1.center);
\draw [-] (#2) to (1.center);
\draw [-] (#3) to (1.center);
\draw [-] (#4) to (1.center);}
\draw [->] (#1) to (#5);
\draw [->] (#2) to (#5);
\draw [->] (#3) to (#5);
\draw [->] (#4) to (#5);
}
\[
\begin{tikzpicture}[xscale=2, yscale=2]
\node (Pn) at (2.5,-0.5) {$P(n)$};
\node (Rn) at (0.5,0.5) {$R(n)$};
\node (Sn) at (1.5,0.5) {$S(n)$};
\node (Tn) at (0.5,0) {$T(n)$};
\node (Qn) at (2.5,0.5) {$Q(n)$};
\node (An) at (3.5,-0) {$A(n)$};
\node (Bn) at (0.5,1) {$B(n)$};
\node (Cn) at (4.5,0.5) {$C(n)$};
\node (Pm) at (1.5,1) {$P(n{-}1)$};
\node (Rm) at (-0.5,0.5) {$R(n{-}1)$};
\node (Sm) at (0.5,-0.5) {$S(n{-}1)$};
\node (Qm) at (3.5,-0.5) {$Q(n{-}1)$};
\node (Am) at (3.5,0.5) {$A(n{-}1)$};
\node (Bm) at (1.5,-0.5) {$B(n{-}1)$};
\node (Cm) at (3.5,1) {$C(n{-}1)$};
\onefold{Sm}{Tn}
\onefold{Rn}{Bn}
\twofoldnew{Sn}{Qm}{Pn}{4,1}{.}
\threefoldnew{Sn}{Rm}{Tn}{Rn}{1.5,0}{.}
\threefoldnew{Sn}{Pn}{Am}{Qn}{5,2}{.}
\threefoldnew{Cm}{Qm}{Am}{Cn}{6,-0.2}{.}
\threefoldnew{Sn}{Qm}{Am}{An}{3.5,0}{.}
\fourfoldnew{Tn}{Pm}{Cm}{Bm}{Sn}{3,-1}{.}
\end{tikzpicture}
\]

\paragraph{Base cases.}
We now prove several lemmas establishing the base cases.
\begin{lemma}
\label{recursionbaseC}
C(0).
\end{lemma}
\begin{proof}
We need to show that given two $0$\-diagram embeddings $e: S\hookrightarrow D$ and $f: D\hookrightarrow M$ between well-defined $0$\-diagrams $S,D,M$, their composite $f\circ e: S\hookrightarrow M$ is well-defined.

The domain and codomain of $f\circ e$ are well-defined by assumption. Both $f$ and $e$ consist of no data, so there is nothing to check and the result is vacuously true.
\end{proof}

\begin{lemma}
\label{recursionbaseR}
R(0).
\end{lemma}
\begin{proof}
We need to show that given well-defined $0$-diagrams $D, S, T$ such that $S$ and $T$ are globular with respect to each other, and a well-defined embedding $e: S\hookrightarrow D$, the rewrite $\rewrite{D}{e}{S}{T}$ of $D$ by $e$ is a well-defined diagram.

The $0$\-diagram embedding $e$ consists of no data. $\rewrite{D}{e}{S}{T}$ is a $0$\-diagram, its list of generators consists of a single cell and there is no source and no slices, so the result is vacuously true.
\end{proof}

\begin{lemma}
\label{recursionbaseP}
P(0).
\end{lemma}
\begin{proof}
We must verify the following:
\[
\rewrite{B}{f}{C}{\rewrite{A}{e}{S}{T}}= \rewrite{(\rewrite{B}{f}{C}{A})}{\lift{f}{A}\circ e}{S}{T}
\]
Since $S, T, A, B, C$  are all $0$\-diagrams, their rewrites exist as $R(0)$ holds by Lemma~\ref{recursionbaseR}. Since $n=0$, for the equation to hold, by Definition~\ref{diagequiv}, we must verify that the diagrams on either side of the equality have precisely 1 generator, and that these generators are the same.

We verify that the diagrams have the same number of generators by following argument, relying on repeated applications of~\myeqref{defrewritezerosource}:
\begin{align*}
&|\rewrite{B}{f}{C}{\rewrite{A}{e}{S}{T}}| \\
&\quad= |B| - |C| + |\rewrite{A}{e}{S}{T}| \\
&\quad= |B| - |C| + |A| - |S| + |T| \\
&\quad= |\rewrite{B}{f}{C}{A}| - |S| + |T| \\
&\quad= |\rewrite{(\rewrite{B}{f}{C}{A})}{\lift{f}{A}\circ e}{S}{T}|
\end{align*}
Furthermore, since $S, T, A, B, C$ are well-defined $0$\-diagrams, we must have $|A|=|B|=|C|=|S|=|T|=1$, and hence $|B| - |C| + |A| - |S| + |T|=1$. By the above chain of equalities, we therefore verify the lengths are 1 are required:
\[
|\rewrite{(\rewrite{B}{f}{C}{A})}{\lift{f}{A}\circ e}{S}{T}|
=|\rewrite{B}{f}{C}{\rewrite{A}{e}{S}{T}}| =1
\]
We now verify that these single generators are identical, employing~\myeqref{defrewritezerogenerator} :
\begin{align*}
&\rewrite{\rewrite{B}{f}{C}{A}}{\lift{f}{A}\circ e}{S}{T}[0].g
= T[0].g
= \rewrite{A}{e}{S}{T}[0].g
=\rewrite{B}{f}{C}{\rewrite{A}{e}{S}{T}}[0].g
\end{align*}
We therefore conclude $P(0)$.
\end{proof}

\paragraph{Inductive steps.}
With these base cases established, we now proceed to prove a series of implications between the logical statements defined earlier in this chapter. Each time we only take the minimal subset of expressions for $(n{-}1)$ that implies the given statement for $n$.
\begin{lemma}
\label{SimpliesR}
For $n\geq 1$, the following holds: $S(n)\wedge T(n)\wedge R(n{-}1) \implies R(n)$. Additionally, if $n=1$, then $S(n)\wedge R(n{-}1) \implies R(n)$
\end{lemma}
\begin{proof}
Let us assume that $S(n)$ and $R(n{-}1)$ hold, then we need to show $\rewrite{D}{e}{S}{T}$ is well-defined. By Definition~\ref{defwellformed}, this means that we require:

\begin{itemize}
\item $\rewrite{D}{e}{S}{T}.s$ is well-defined
\item $\rewrite{D}{e}{S}{T}[j].d$ for $0 < j < |\rewrite{D}{e}{S}{T}|$ exists and is well-defined
\end{itemize}
The first statement follows since:
\begin{align*}
\rewrite{D}{e}{S}{T}.s=D.s
\end{align*}
$D.s$ is well-defined as the source of a well-defined diagram $D$. We prove the second statement separately for three ranges within $0\leq j\leq |\rewrite{D}{e}{S}{T}|$

\begin{itemize}
\item
For $0< j \leq e.h$, by $S(n)$ we obtain that:
\[\rewrite{D}{e}{S}{T}[j].d = D[j].d\]
Hence all slices in this section exist and are well-defined as $D$ is well-defined.

\item
For $e.h\leq j < e.h + |T|$, by $S(n)$ we obtain that:
\[\rewrite{D}{e}{S}{T}[j].d = \rewrite{(D[e.h].d)}{e.e}{S.s}{T[j-e.h].d}\]
\begin{itemize}
\item$D[e.h].d$ is well-defined as the slice of a well-defined diagram $D$
\item The embedding $e.e$ is well-defined as $e$ is well-defined 
\item If $n=1$, $S.s=S(T)$ and $T[j-e.h].d$ are $0$-diagrams and globular with respect to each other
\item If $n>1$, we need $S.s.s=T[j-e.h].d.s$ and $S.s.t=T[j-e.h].d.t$. We note that $S.s=T.s=T[0].d$, then the result follows by $T(n)$.
\end{itemize}

By applying $R(n{-}1)$ to the rewrite of the $(n{-}1)$\-diagram $D[e.h].d$, we get that the diagram $\rewrite{(D[e.h].d)}{e.e}{S.s}{T[j-e.h].d}$ exists and is well-defined, hence so is $\rewrite{D}{e}{S}{T}[j].d$, as required.

\item
For $e.h + |T| \leq j < |D| - |S| + |T|$,  by $S(n)$ we obtain that:
\[\rewrite{D}{e}{S}{T}[j].d = D[j - |T| + |S|].d\]
It follows that all slices in this section exist and are well-defined, since all $D[j].d$ are well-defined as slices of the well-defined diagram $D$.
\end{itemize}
Hence, all the slices $(\rewrite{D}{e}{S}{T})[j].d$ exists and are also well-defined, hence this diagram is itself well-defined. By this $R(n)$ holds, hence the implication is true. Note that $T(n)$ is only used for $n>1$, hence $S(1)\wedge R(0) \implies R(1)$ holds.
\end{proof}


\begin{lemma}
\label{TPimpliesS}
For $n>1$ the following holds: $T(n) \wedge P(n{-}1)\wedge C(n{-}1)\wedge B(n{-}1)\implies S(n)$. Additionally, for $n=1$: $P(0)\wedge C(0)\wedge B(0)\implies S(1)$
\end{lemma}

\begin{proof}
Let us assume that all $T(n)$, $P(n{-}1)$, $C(n{-}1)$, $B(n{-}1)$ hold. We show that $S(n)$ holds in each of the individual three ranges separately:

\begin{itemize}

\item For $0\leq j \leq e.h $:

In this range, we show this result by induction on $j$.
\begin{itemize}
\item \emph{Base case:} For $j=0$
\[\rewrite{A}{e}{S}{T}[0].d = \rewrite{A}{e}{S}{T}.s = A.s = A[0].s\]

\item \emph{Inductive step:} For $0< j \leq e.h$ assume by induction that:
\begin{align*}
\rewrite{A}{e}{S}{T}[j].d = A[j].d &\qquad (\text{\emph{IH}})
\end{align*}

Let us consider $\rewrite{A}{e}{S}{T}[j+1].d$: Note that by Definition~\ref{defrewrite} for $0< j \leq e.h$ we have the following:

\begin{itemize}
\item $\rewrite{A}{e}{S}{T}[j].g = A[j].g$
\item $\rewrite{A}{e}{S}{T}[j].e = A[j].e$
\end{itemize}

Hence:
\begin{align*}
&\rewrite{A}{e}{S}{T}[j+1].d \\ 
&\quad=\rewrite{(\rewrite{A}{e}{S}{T}[j].d)}{\rewrite{A}{e}{S}{T}[j].e}{\rewrite{A}{e}{S}{T}[j].g.s}{\\&\qquad\qquad \rewrite{A}{e}{S}{T}[j].g.t} &&[\defref{defslice}]\\
&\quad=\rewrite{(\rewrite{A}{e}{S}{T}[j].d)}{A[j].e}{A[j].g.s}{A[j].g.t} &&[\defref{defrewrite}]\\
&\quad=\rewrite{(A[j].d)}{A[j].e}{A[j].g.s}{A[j].g.t} &&[\text{\emph{IH}}]\\
&\quad=A[j+1].d &&[\defref{defslice}]
\end{align*}
As required.
\end{itemize}
By induction we have that: $\rewrite{A}{e}{S}{T}[j].d = A[j].d$ holds for $0\leq j \leq e.h$.

\item For $e.h \leq j \leq e.h + |T| $:

First, for the rewrite $\rewrite{(A[e.h].d)}{e.e}{S.s}{T[j-e.h].d}$ to be well-defined, we need to show that the source $S.s$ and the target $T[j-e.h].d$ of the rewrite are globular with respect to each other.

\begin{itemize}

\item For $n=1$, both $S.s$ and $T[j-e.h].d$ are $0$\-diagrams for $e.h \leq j \leq |e.h + |T|$ so there are no conditions to check.

\item For $n>1$, we need $S.s.s = T[j-e.h].d.s$ and $S.s.t = T[j-e.h].d.t$ for $e.h \leq j \leq |e.h + |T|$. 

Note that since $S$ and $T$ are globular, we have: 
\begin{align*}
S.s = T.s = T[0].d = T.s
\end{align*}
We then apply proposition $T(n)$ to the diagram $T$, to obtain the necessary result.

\end{itemize}

We can now prove this result by induction on $e.h \leq j \leq e.h + |T|$:

\begin{itemize}
\item \emph{Base case:} For $j=e.h$
\begin{align*}
&\rewrite{(A[e.h].d)}{e.e}{S.s}{T[j-e.h].d} \\
&\quad=  \rewrite{(A[e.h].d)}{e.e}{S.s}{T[e.h-e.h].d} && [j=e.h]\\
&\quad=  \rewrite{(A[e.h].d)}{e.e}{S.s}{T[0].d} \\
&\quad=  \rewrite{(A[e.h].d)}{e.e}{S.s}{T.s} &&[\defref{definitialterminalslice}] \\
&\quad=  \rewrite{(A[e.h].d)}{e.e}{S.s}{S.s} &&[S.t=T.t]\\
&\quad=  A[e.h].d &&[\text{Lemma }\ref{vacuousrewrite}] \\
&\quad= \rewrite{A}{e}{S}{T}[e.h].d
\end{align*}
The final transformation follows by this result for the initial rage $0\leq j \leq e.h$.

\item \emph{Inductive step:} For $e.h< j \leq e.h + |T|$ assume by induction that:
\begin{align*}
\rewrite{A}{e}{S}{T}[j].d = \rewrite{A[e.h].d}{e.e}{S.s}{T[j-e.h].d}&\qquad \text{[I.H.]}
\end{align*}

Let us consider $\rewrite{A}{e}{S}{T}[j+1].d$:
\begin{align*}
&\rewrite{A}{e}{S}{T}[j+1].d \\
&\quad= \rewrite{\rewrite{A}{e}{S}{T}[j].d}{\rewrite{A}{e}{S}{T}[j].e}{\rewrite{A}{e}{S}{T}[j].g.s}{\\&\qquad\qquad \rewrite{A}{e}{S}{T}[j].g.t} &&[\defref{defslice}]\\
&\quad= \rewrite{\rewrite{A}{e}{S}{T}[j].d}{\lift{(e.e)}{T[j-e.h].d} \\
&\qquad\qquad \circ T[j-e.h].e}{T[j-e.h].g.s}{T[j-e.h].g.t} &&[\defref{defrewrite}]\\
&\quad= \rewrite{\rewrite{A[e.h].d}{e.e}{S.s}{T[j-e.h].d}}{\lift{e.e}{T[j-e.h].d}\\
&\qquad\qquad \circ T[j-e.h].e}{T[j-e.h].g.s}{T[j-e.h].g.t} &&[\mathrm{I.H.}]\\
&\quad= \rewrite{A[e.h].d}{e.e}{S.s}{\rewrite{T[j-e.h].d}{T[j-e.h].e}{T[j-e.h].g.s}{\\&\qquad\qquad T[j-e.h].g.t}} &&[P(n{-}1)] \\
&\quad= \rewrite{A[e.h].d}{e.e}{S.s}{T[(j+1)-e.h].d} &&[\defref{defslice}]
\end{align*}
$P(n{-}1)$ may be applied here as all the embeddings involved in the rewrites are well-defined. The lifted embedding between $(n{-}1)$\-diagrams is well-defined as $B(n{-}1)$ holds and the composed embedding of $(n{-}1)$\-diagram embeddings is well-defined, as $C(n{-}1)$ holds.

\end{itemize}

\item For $e.h + |T| \leq j < |A| - |S| + |T|$:

We show the result in this range by induction on $j$.
\begin{itemize}
\item \emph{Base case:} For $j=e.h + |T|$, we have:
\begin{align*}
&\rewrite{A}{e}{S}{T}[e.h + |T|].d  \\
&\quad=\rewrite{A[e.h].d}{e.e}{S.s}{T[(e.h + |T|)-e.h]} \\
&\quad=\rewrite{A[e.h].d}{e.e}{S.s}{T[|T|]} \\
&\quad=\rewrite{A[e.h].d}{e.e}{S.s}{T.t]} &&[\defref{definitialterminalslice}] \\
&\quad=\rewrite{A[e.h].d}{e.e}{S.s}{S.t]} &&[S.t=T.t]\\
&\quad=\rewrite{A[e.h].d}{e.e}{S.s}{S[|S|]} &&[\defref{definitialterminalslice}] \\
&\quad=\rewrite{A}{e}{S}{S}[e.h+|S|].d &&[\defref{defslice}]\\
&\quad=A[e.h+|S|].d &&[\text{Lemma }~\ref{vacuousrewrite}] 
\end{align*}
The penultimate transformation follows by the result for the range $e.h \leq j \leq e.h + |S|$ applied to the identity rewrite of $A$ $e:S\hookrightarrow A$.

\item \emph{Inductive step:} For $e.h + |T| < j \leq |A| + |T| - |S|$, note that, similarly as in the range $0\leq j \leq e.h $, we have that by Definition~\ref{defrewrite}:

\begin{itemize}
\item $\rewrite{A}{e}{S}{T}[j].g = A[j + |S| - |T].g$
\item $\rewrite{A}{e}{S}{T}[j].e = A[j + |S| - |T].e$
\end{itemize}
The rest of the argument is the same as for the range $0\leq j \leq e.h$.
\end{itemize}
By induction we have proved that: $\rewrite{A}{e}{S}{T}[j].d = A[j].d$ holds for $e.h + |T| \leq j \leq |A| + |T| - |S|$.

\end{itemize}

By this, $S(n)$ holds for each of the three ranges, hence the implication is true. Note that the assumption $T(n)$ is only used to show globularity for $n>1$, hence we also have $P(0)\wedge C(0)\wedge B(0)\implies S(1)$.

\end{proof}


\begin{lemma}
\label{SimpliesT}
For $n\leq 2$ the following holds: $S(n{-}1)\implies T(n)$
\end{lemma}
\begin{proof}
Let us assume that $S(n{-}1)$ holds. We need to show that  for any well-defined $n$\-diagram $D$, we have $(D.s).s = (D[j].d).s$ and $(D.s).t = (D[j].d).t$ for any $0\leq j < |D|$. First, we prove this result by induction on $i$:
\begin{itemize}
\item \emph{Base case}: For $j=0$, by Definition~\ref{defrewrite}, we have: $D.s = D[0].d$, so the result trivially holds. 

\item \emph{Inductive step}: 
For $j>0$, let us assume by induction that: 
\begin{align*}
&(D.s).s = (D[j].d).s & \qquad (\text{\emph{IH}(a)}) \\
&(D.s).t = (D[j].d).t & \qquad (\text{\emph{IH}(b)})
\end{align*}
Firstly, let us consider $(D[i+1].d).s$:
\begin{align*}
&(D[j+1].d).s = \\
&\quad=(\rewrite{(D[j].d)}{D[j].e}{D[j].g.s}{D[j].g.t}).s &&[\defref{defslice}] \\
&\quad=(D[j].d).s &&[\defref{defrewrite}] \\
&\quad= (D.s).s && [\text{\emph{IH}(a)}]
\end{align*}
Secondly, let us consider $(D[i+1].d).t$:
\begin{align*}
&(D[i+1].d).t = \\
&\quad=(\rewrite{(D[j].d)}{D[j].e)}{D[j].g.s}{D[j].g.t}).t &&[\defref{defslice}] \\
&\quad=(\rewrite{(D[j].d)}{D[j].e)}{D[j].g.s}{\\&\qquad\qquad D[j].g.t})[|D[j].d| - |D[j].g.s| + |D[j].g.t|].d  &&[\defref{definitialterminalslice}] \\
&\quad=(D[j].d)[(|D[j].d| - |D[j].g.s| + |D[j].g.t|)\\&\qquad\qquad + |D[j].g.s| - |D[j].g.s|] &&[S(n{-}1)]\\
&\quad=(D[j].d)[|D[j].d|]\\
&\quad=(D[j].d).t &&[\defref{definitialterminalslice}] \\
&\quad= (D.s).t && [\text{\emph{IH(b)}}]
\end{align*}
\end{itemize}
By this, the result holds for all $0\leq i\leq |D|$, hence $T(n)$ holds and the implication is true.
\end{proof}


\begin{lemma}
\label{SPAimpliesQ}
For $n>0$ the following holds: $S(n) \wedge P(n)\wedge A(n{-}1)\implies Q(n)$. Additionally, for $n=0$: $P(0)\implies Q(0)$
\end{lemma}
\begin{proof}
Let us assume that $S(n)$, $P(n)$ and $A(n{-}1)$ all hold. By Lemma~\ref{SimpliesR} the statement $R(n)$ also holds and all the following $n$\-diagram rewrites are well-defined. We need to show that the following two $n$\-diagram embeddings are equivalent:
\begin{align*}
\lift{(\lift{f}{A}\circ e)}{T} = \lift{f}{\rewrite{A}{e}{S}{T}}\circ\lift{e}{T}
\end{align*}
For this, by Definition~\ref{embeddingequiv} domains and codomains must be equivalent diagrams. By Definition~\ref{liftedembedding}, types of the individual embeddings are as follows:
\begin{align*}
\lift{f}{A}\circ e &: S\hookrightarrow  \rewrite{B}{f}{C}{A}\\
\lift{(\lift{f}{A}\circ e)}{T} &: T\hookrightarrow \rewrite{\rewrite{B}{f}{C}{A}}{\lift{f}{A}\circ e}{S}{T} \\
\lift{f}{\rewrite{A}{e}{S}{T}} &: \rewrite{A}{e}{S}{T} \hookrightarrow \rewrite{B}{f}{A}{\rewrite{A}{e}{S}{T}}\\
\lift{e}{T} &: T\hookrightarrow \rewrite{A}{e}{S}{T}
\end{align*}
As the domain of $\lift{f}{\rewrite{A}{e}{S}{T}}$ matches the codomain of $\lift{e}{T}$, we can conclude\JVcomm{How?} that the composite $\lift{f}{\rewrite{A}{e}{S}{T}}\circ\lift{e}{T}$ exists. The domain of this composite matches the domain of $\lift{\lift{f}{A}\circ e}{T}$. 

The codomain of $\lift{f}{\rewrite{A}{e}{S}{T}}\circ\lift{e}{T}$ is $\rewrite{B}{f}{A}{\rewrite{A}{e}{S}{T}}$, the codomain of $\lift{(\lift{f}{A}\circ e)}{T}$ is $\rewrite{\rewrite{B}{f}{C}{A}}{\lift{f}{A}\circ e}{S}{T}$. These two diagrams are equivalent by $P(n)$, hence the codomains of both embeddings also match, as required.

If $n>0$, then we need to check two additional conditions:
\begin{itemize}
\item Component embeddings are equivalent:
\begin{align*}
&(\lift{f}{\rewrite{A}{e}{S}{T}}\circ\lift{e}{T}).e\\
&\quad=\lift{\lift{f}{\rewrite{A}{e}{S}{T}}.e}{\rewrite{A}{e}{S}{T}[\lift{e}{T}.h].d}\circ \lift{e}{T}.e  &&[\myeqref{defrewriteembedding}] \\
&\quad=\lift{\lift{f}{\rewrite{A}{e}{S}{T}}.e}{\rewrite{A}{e}{S}{T}[\lift{e}{T}.h].d}\circ e.e &&[\defref{liftedembedding}] \\
&\quad=\lift{\lift{f}{\rewrite{A}{e}{S}{T}}.e}{\rewrite{A}{e}{S}{T}[e.h].d}\circ e.e &&[\defref{liftedembedding}] \\
&\quad=\lift{\lift{f}{\rewrite{A}{e}{S}{T}}.e}{A[e.h].d}\circ e.e &&[S(n)] \\
&\quad=\lift{f.e}{A[e.h].d}\circ e.e &&[\defref{liftedembedding}] \\
&\quad=\lift{\lift{f}{A}.e}{A[e.h].d}\circ e.e &&[\defref{liftedembedding}] \\
&\quad=(\lift{f}{A}\circ e).e &&[\defref{defcomposition}] \\
&\quad=\lift{(\lift{f}{A}\circ e)}{T}.e &&[\defref{liftedembedding}] 
\end{align*}

\item Heights are equal:
\begin{align*}
&(\lift{f}{\rewrite{A}{e}{S}{T}}\circ\lift{e}{T}).h\\
&\quad=\lift{f}{\rewrite{A}{e}{S}{T}}.h+\lift{e}{T}.h&&[\defref{defcomposition}] \\
&\quad=f.h+ e.h &&[\defref{liftedembedding}] \\
&\quad=\lift{f}{A}.h+ e.h &&[\defref{liftedembedding}] \\
&\quad=\lift{f}{A}\circ e.h &&[\defref{defcomposition}] \\
&\quad=\lift{(\lift{f}{A}\circ e)}{T}.h &&[\defref{liftedembedding}] 
\end{align*}

\end{itemize}
Both conditions are satisfied, hence these two embeddings are equivalent and $Q(n)$ holds. By this, the implication $S(n) \wedge P(n)\implies Q(n)$ is true for $n>0$. Note that the assumptions $S(n)$, $A(n{-}1)$ are only used when $n>0$, hence $P(0)\implies Q(0)$ also holds.

\end{proof}


\begin{lemma}
\label{SQimpliesP}
For $n\geq 1$ the following holds: $S(n) \wedge Q(n{-}1)\implies P(n)$. 
\end{lemma}
\begin{proof}
Let us assume that both $S(n)$ and $Q(n{-}1)$ hold. By Lemma~\ref{SimpliesR} the statement $R(n)$ also holds and all the following $n$\-diagram rewrites are well-defined.

We need to show that for well-defined diagrams $S, T, A, B, C$ such that $S.s=T.s$, $S.t=T.t$, $A.s=C.s$, $A.t=C.t$ and well-defined embeddings $e:S\hookrightarrow A$, $f:C\hookrightarrow B$, the following holds':
\begin{align*}
\rewrite{B}{f}{C}{\rewrite{A}{e}{S}{T}}= \rewrite{(\rewrite{B}{f}{C}{A})}{\lift{f}{A}\circ e}{S}{T}
\end{align*}

By Definition~\ref{diagequiv}, we need to check the following four conditions:

\begin{itemize}

\item Sources are equivalent diagrams: All equalities in this derivation follow by the source clause~\myeqref{defrewritesource} in the definition of a rewrite.
\begin{align*}
&(\rewrite{(\rewrite{B}{f}{C}{A})}{\lift{f}{A}\circ e}{S}{T}).s = \\
&\quad= (\rewrite{B}{f}{C}{A}).s \\
&\quad=B.s\\
&\quad=(\rewrite{B}{f}{C}{\rewrite{A}{e}{S}{T}}).s
\end{align*}

\item Sizes of lists of generators are equal: All equalities in this derivation follow by the size clause~\myeqref{defrewritesize} in the definition of a rewrite.
\begin{align*}
&|\rewrite{B}{f}{C}{\rewrite{A}{e}{S}{T}}| = \\
&\quad= |B| - |C| + |\rewrite{A}{e}{S}{T}| \\
&\quad= |B| - |C| + |A| - |S| + |T| \\
&\quad= |\rewrite{B}{f}{C}{A}| - |S| + |T| \\
&\quad= |\rewrite{(\rewrite{B}{f}{C}{A})}{\lift{f}{A}\circ e}{S}{T}|
\end{align*}

\item Generators are equal, for $0\leq i\leq |\rewrite{B}{f}{C}{\rewrite{A}{e}{S}{T}}|$. All derivations in this section follow by the generator clause~\myeqref{defrewritegenerator} in the definition of a rewrite.

We consider $i$ in five separate ranges:

\begin{itemize}

\item $0\leq i \leq f.h$
\begin{align*}
&(\rewrite{B}{f}{C}{\rewrite{A}{e}{S}{T}})[i].g = \\
&\quad= B[i].g \\
&\quad= (\rewrite{B}{f}{C}{A})[i].g \\
&\quad=(\rewrite{(\rewrite{B}{f}{C}{A})}{\lift{f}{A}\circ e}{S}{T}).g
\end{align*}

\item $f.h\leq i \leq f.h + e.h$
\begin{align*}
&(\rewrite{B}{f}{C}{\rewrite{A}{e}{S}{T}})[i].g = \\
&\quad=(\rewrite{A}{e}{S}{T})[i-f.h].g \\
&\quad= A[i-f.h].g \\
&\quad= (\rewrite{B}{f}{C}{A})[i].g \\
&\quad=(\rewrite{(\rewrite{B}{f}{C}{A})}{\lift{f}{A}\circ e}{S}{T}).g
\end{align*}

\item $f.h + e.h\leq i \leq f.h + e.h + |T|$
\begin{align*}
&(\rewrite{B}{f}{C}{\rewrite{A}{e}{S}{T}})[i].g = \\
&\quad=(\rewrite{A}{e}{S}{T})[i-f.h].g \\
&\quad=T[(i-f.h)-e.h].g \\
&\quad=T[(i-(f.h+e.h)].g \\
&\quad=T[i - (\lift{f}{A}.h+ e.h)].g \\
&\quad=T[i - (\lift{f}{A}\circ e).h].g \\
&\quad=(\rewrite{(\rewrite{B}{f}{C}{A})}{\lift{f}{A}\circ e}{S}{T}).g
\end{align*}

\item $f.h + e.h + |T|\leq i \leq f.h + |A| - |S| + |T|$
\begin{align*}
&(\rewrite{B}{f}{C}{\rewrite{A}{e}{S}{T}})[i].g = \\
&\quad=\rewrite{A}{e}{S}{T})[i- f.h].g \\
&\quad=A[(i- f.h) - |T| + |S|].g \\
&\quad=A[(i - |T| + |S|) - f.h].g \\
&\quad=\rewrite{B}{f}{C}{A}[i - |T| + |S|].g \\
&\quad=(\rewrite{(\rewrite{B}{f}{C}{A})}{\lift{f}{A}\circ e}{S}{T}).g
\end{align*}

\item $f.h + |A| - |S| + |T| \leq i \leq |B| - |C| + |A| - |S| + |T|  $
\begin{align*}
&(\rewrite{B}{f}{C}{\rewrite{A}{e}{S}{T}})[i].g  \\
&\quad= B[i - |\rewrite{A}{e}{S}{T}| + |C|].g \\
&\quad= B[i - |A| + |C| - |T| + |S|].g \\
&\quad=(\rewrite{B}{f}{C}{A})[i - |T| + |S|].g \\
&\quad=(\rewrite{(\rewrite{B}{f}{C}{A})}{\lift{f}{A}\circ e}{S}{T}).g
\end{align*}

\end{itemize}

\item Embeddings are equivalent, similarly as for generators, for $0\leq i\leq |\rewrite{B}{f}{C}{\rewrite{A}{e}{S}{T}}|$. We consider $i$ in five separate ranges:

\begin{itemize}

\item $0\leq i \leq f.h$
\begin{align*}
&(\rewrite{B}{f}{C}{\rewrite{A}{e}{S}{T}})[i].e  \\
&\quad= B[i].e &&[\myeqref{defrewriteembedding}] \\
&\quad= (\rewrite{B}{f}{C}{A})[i].e &&[\myeqref{defrewriteembedding}] \\
&\quad=(\rewrite{(\rewrite{B}{f}{C}{A})}{\lift{f}{A}\circ e}{S}{T}).e &&[\myeqref{defrewriteembedding}]
\end{align*}

\item $f.h\leq i \leq f.h + e.h$
\begin{align*}
&(\rewrite{B}{f}{C}{\rewrite{A}{e}{S}{T}})[i].e  \\
&\quad=\lift{(f.e)}{\rewrite{A}{e}{S}{T}[i-f.h].d}\circ \rewrite{A}{e}{S}{T}[i-f.h].e &&[\myeqref{defrewriteembedding}] \\
&\quad=\lift{(f.e)}{\rewrite{A}{e}{S}{T}[i-f.h].d}\circ A [i-f.h].e &&[\myeqref{defrewriteembedding}] \\
&\quad=\lift{(f.e)}{A[i-f.h].d}\circ A [i-f.h].e && [S(n)] \\
&\quad= (\rewrite{B}{f}{C}{A})[i].e &&[\myeqref{defrewriteembedding}] \\
&\quad=(\rewrite{(\rewrite{B}{f}{C}{A})}{\lift{f}{A}\circ e}{S}{T}).e &&[\myeqref{defrewriteembedding}] 
\end{align*}

\item $f.h + e.h\leq i \leq f.h + e.h + |T|$
\begin{align*}
&(\rewrite{B}{f}{C}{\rewrite{A}{e}{S}{T}})[i].e  \\
&\quad=\lift{(f.e)}{\rewrite{A}{e}{S}{T}[i-f.h].d}\circ \rewrite{A}{e}{S}{T}[i-f.h].e &&[\myeqref{defrewriteembedding}] \\
&\quad=\lift{(f.e)}{\rewrite{A}{e}{S}{T}[i-f.h].d}\\&\qquad\qquad\circ (\lift{(e.e)}{T[(i-f.h)-e.h].d} \circ T[(i-f.h)-e.h].e) &&[\myeqref{defrewriteembedding}] \\
&\quad=(\lift{(f.e)}{\rewrite{A}{e}{S}{T}[i-f.h].d}\\&\qquad\qquad\circ \lift{(e.e)}{T[(i-f.h)-e.h].d}) \circ T[(i-f.h)-e.h].e &&[A(n{-}1)] \\
&\quad=(\lift{(f.e)}{\rewrite{A[e.h]}{e.e}{S.s}{T[i-f.h]}}\\&\qquad\qquad \circ \lift{(e.e)}{T[(i-f.h)-e.h].d}) \circ T[(i-f.h)-e.h].e &&[S(n)]\\
&\quad= (\lift{((\lift{f.e)}{A[e.h].d}\circ e.e)}{T[(i-f.h)-e.h]} \\&\qquad\qquad\circ T[(i-f.h)-e.h].e &&[Q(n{-}1)]\\
&\quad= \lift{((\lift{f.e)}{A[e.h].d}\circ e.e)}{T[i-(f.h+ e.h)]} \\&\qquad\qquad\circ T[i-(f.h+ e.h)].e \\
&\quad= \lift{((\lift{\lift{f}{A}.e)}{A[e.h].d}\circ e.e)}{T[i-(f.h+ e.h)]} \\&\qquad\qquad\circ  T[i-(f.h+ e.h)].e &&[\defref{liftedembedding}] \\
&\quad= \lift{((\lift{f}{A}\circ e).e)}{T[i-(f.h+ e.h)]} \\&\qquad\qquad\circ T[i-(f.h+ e.h)].e &&[\defref{defcomposition}] \\
&\quad= \lift{((\lift{f}{A}\circ e).e)}{T[i-(\lift{f}{A}.h+ e.h)]} \\&\qquad\qquad\circ T[i-(\lift{f}{A}.h+ e.h)].e &&[\defref{liftedembedding}] \\
&\quad= \lift{((\lift{f}{A}\circ e).e)}{T[i-(\lift{f}{A}\circ e).h]} \\&\qquad\qquad\circ T[i-(\lift{f}{A}\circ e).h].e &&[\defref{defcomposition}] \\
&\quad=(\rewrite{(\rewrite{B}{f}{C}{A})}{\lift{f}{A}\circ e}{S}{T}).g &&[\myeqref{defrewriteembedding}] 
\end{align*}

\item $f.h + e.h + |T|\leq i \leq f.h + |A| - |S| + |T|$
\begin{align*}
&(\rewrite{B}{f}{C}{\rewrite{A}{e}{S}{T}})[i].g  \\
&\quad=\lift{(f.e)}{\rewrite{A}{e}{S}{T}[i-f.h].d}\\&\qquad\qquad\circ \rewrite{A}{e}{S}{T}[i-f.h].e &&[\myeqref{defrewriteembedding}] \\
&\quad=\lift{(f.e)}{\rewrite{A}{e}{S}{T}[i-f.h].d}\\&\qquad\qquad\circ A [(i - |T| + |S|)-f.h].e &&[\myeqref{defrewriteembedding}] \\
&\quad=\lift{(f.e)}{A[(i - |T| + |S|)-f.h].d}\\&\qquad\qquad\circ A [(i - |T| + |S|)-f.h].e &&[S(n)] \\
&\quad=\rewrite{B}{f}{C}{A}[i - |T| + |S|].e &&[\myeqref{defrewriteembedding}] \\
&\quad=(\rewrite{(\rewrite{B}{f}{C}{A})}{\lift{f}{A}\circ e}{S}{T}).e &&[\myeqref{defrewriteembedding}]
\end{align*}

\item $f.h + |A| - |S| + |T| \leq i \leq |B| - |C| + |A| - |S| + |T|  $
\begin{align*}
&(\rewrite{B}{f}{C}{\rewrite{A}{e}{S}{T}})[i].e \\
&\quad= B[i - |\rewrite{A}{e}{S}{T}| + |C|].e &&[\myeqref{defrewriteembedding}] \\
&\quad= B[i - |A| + |C| - |T| + |S|].e &&[\defref{defrewrite}] \\
&\quad=(\rewrite{B}{f}{C}{A})[i - |T| + |S|].e &&[\myeqref{defrewriteembedding}] \\
&\quad=(\rewrite{(\rewrite{B}{f}{C}{A})}{\lift{f}{A}\circ e}{S}{T}).e &&[\myeqref{defrewriteembedding}] 
\end{align*}

\end{itemize}

\end{itemize}

\end{proof}

\begin{lemma}
\label{RimpliesL}
For $n\geq 1$ the following holds: $R(n)\implies B(n)$. 

\end{lemma}
\begin{proof}
Let us assume that $R(n)$ holds. We need to show that $\lift{e}{T}: T\hookrightarrow \rewrite{A}{e}{S}{T}$ is well-defined.
For that, the domain and codomain of $\lift{e}{T}$ have to be well-defined. $T$ is well-defined by assumption, $\rewrite{A}{e}{S}{T}$ is well-defined by $R(n)$.

\begin{itemize}
\item If $n=0$, we need $ \rewrite{A}{e}{S}{T}[0].g = T[0].g$, but this follows directly from the definition of a rewrite for $n=0$.

\item If $n>0$, by Definition~\ref{defembwelldef} we need to show three separate conditions:
\begin{itemize}
\item $\lift{e}{T}.e$ is well-defined. 

By Definition~\ref{liftedembedding}, we have $\lift{e}{T}.e = e.e$. As $e$ is well-defined, so is its component embedding $e.e$. Hence, $\lift{e}{T}.e$ is equivalent to a well-defined embedding and itself is well-defined.

\item For every $0\leq i< |T|$ we need to show that the generators of both diagrams satisfy the following: 
\begin{align*}
&\rewrite{A}{e}{S}{T}[i+\lift{e}{T}.h].g  \\
&\quad=\rewrite{A}{e}{S}{T}[i+e.h].g &&[\defref{liftedembedding}] \\
&\quad=T[(i+e.h)-e.h].g &&[\myeqref{defrewriteembedding}] \\
&\quad=T[i].g
\end{align*}

\item For every $0\leq i< |T|$ we need to show that the embeddings of both diagrams satisfy the following:
\begin{align*}
&\rewrite{A}{e}{S}{T}[i+\lift{e}{T}.h].e \\
&\quad=\rewrite{A}{e}{S}{T}[i+e.h].e &&[\defref{liftedembedding}] \\
&\quad=\lift{(e.e)}{T[i].d}\circ T[i].e  &&[\defref{defrewriteembedding}] \\
&\quad=\lift{(\lift{e}{T}.e)}{T[i].d}\circ T[i].e &&[\defref{liftedembedding}] 
\end{align*}

\end{itemize}
\end{itemize}
As all three conditions are satisfied, $\lift{e}{T}$ is a well-defined embedding and $B(n)$ holds, hence the implication is true.
\end{proof}


\begin{lemma}
\label{CQAimpliesC}
For $n\geq 1$ the following holds: $C(n{-}1)\wedge Q(n{-}1)\wedge A(n{-}1)\implies C(n)$. 
\end{lemma}
\begin{proof}

We assume that all $C(n{-}1)$ , $A(n{-}1)$ and $S(n{-}1)$ hold. By combining Lemmas~\ref{RimpliesL} and~\ref{SimpliesR}, we have that since $S(n{-}1)$ holds, both $R(n{-}1)$ and $B(n{-}1)$ also hold.

We need to show that given two $n$\-diagram embeddings $e: S\hookrightarrow D$ and $f: D\hookrightarrow M$ between well-defined $n$\-diagrams $S,D,M$, their composite $f\circ e: S\hookrightarrow M$ is well-defined.

The domain and codomain of $f\circ e$ are well-defined by assumption. To show that $f\circ e$ is well-defined, by Definition~\ref{defembwelldef} we need to show three separate conditions:
\begin{itemize}
\item $(f\circ e).e$ is well-defined. 
\begin{align*}
(f\circ e).e = \lift{(f.e)}{D[e.h].d}\circ e.e
\end{align*}
Both $e.e$ and $f.e$ are well-defined $(n{-}1)$\-embeddings, as they are component embeddings of well-defined embeddings $f, e$. $\lift{(f.e)}{D[e.h].d}$ is well-defined $(n{-}1)$\-embedding by $B(n{-}1)$.

We then apply the inductive hypothesis to conclude that $(f\circ e).e$ is well-defined.

\item For $0\leq i < |S|$ we have: 
\begin{align*}
&S[i].g \\
&\quad=D[e.h + i].g &&[\myeqref{eqembwelldef2} \text{ for } S, D]\\ 
&\quad=M[e.h + f.h + i].g &&[\myeqref{eqembwelldef2} \text{ for } D, M] 
\end{align*}

\item For $0\leq i < |S|$ we have:
\begin{align*}
&M[(f\circ e).h + i].e  \\
&\quad=M[(f.h + e.h) + i].e &&[\defref{defcomposition}] \\
&\quad= \lift{(f.e)}{D[e.h+i].d}\circ D[e.h + i].e &&[\myeqref{eqembwelldef3} \text{ for } D,M]\\ 
&\quad= \lift{(f.e)}{D[e.h+i]}\circ (\lift{(e.e)}{S[i].d}\circ S[i].e) &&[\myeqref{eqembwelldef3} \text{ for } S, D]\\
 &\quad= (\lift{(f.e)}{D[e.h+i]}\circ \lift{(e.e)}{S[i].d}\circ) S[i].e &&[A(n{-}1)]\\
&\quad= \lift{(\lift{(f.e)}{D[e.h]}\circ e.e)}{S[i].d}\circ S[i].e &&[Q(n{-}1)]\\
&\quad= \lift{((f\circ e).e)}{S[i].d}\circ S[i].e &&[\defref{defcomposition}] 
\end{align*}

\end{itemize}


As all conditions of Definition~\ref{diagequiv} are satisfied, the two diagrams are equivalent and the result holds.
\end{proof}

\begin{lemma}
\label{SQAimpliesA}
For $n\geq 1$ the following holds: $S(n) \wedge Q(n{-}1) \wedge A(n{-}1)\implies A(n)$. Additionally, if $n=0$, then A(0) holds with no further assumptions.
\end{lemma}
\begin{proof}
We assume that all $S(n)$, $Q(n{-}1)$ and $A(n{-}1)$ hold.

We need to show that the following two $n$\-diagram embeddings are equivalent:
\begin{align*}
g\circ(f\circ e) = (g\circ f) \circ e
\end{align*}
For that by Definition~\ref{embeddingequiv} domains and codomains must be equivalent diagrams. By Definition~\ref{defcomposition}, types of the individual embeddings are as follows:
\begin{align*}
f\circ e&: S\hookrightarrow M \\
g\circ f&: D\hookrightarrow N \\
g\circ(f\circ e)&: S\hookrightarrow N \\
(g\circ f) \circ e&: S\hookrightarrow N 
\end{align*}
By this we can conclude that the domains and codomains of both embeddings match. If $n>0$, then we must check two additional conditions:
\begin{itemize}
\item Heights are equal:
\begin{align*}
&(g\circ(f\circ e)).h \\
&\quad= g.h+(f\circ e).h & [\defref{defcomposition}]\\
&\quad= g.h+f.h+ e.h & [\defref{defcomposition}]\\
&\quad (g\circ f).h + e.h & [\defref{defcomposition}]\\
&\quad ((g\circ f) \circ e).h & [\defref{defcomposition}]
\end{align*}



\item Component embeddings are equivalent:\JVcomm{Errors here}
\begin{align*}
&(g\circ(f\circ e)).e \\
&\quad= \lift{(g.e)}{M[(f\circ e).h].d}\circ (f\circ e).e && [\defref{defcomposition}]\\
&\quad= \lift{(g.e)}{M[(f\circ e).h].d}\circ (\lift{(f.e)}{D[e.h].d}\circ e.e) && [\defref{defcomposition}]\\
&\quad= (\lift{(g.e)}{M[(f\circ e).h].d}\circ \lift{(f.e)}{D[e.h].d})\circ e.) && [A(n{-}1)]\\
&\quad= (\lift{(g.e)}{M[f.h+ e.h].d}\circ \lift{(f.e)}{D[e.h].d})\circ e.) && [\defref{defcomposition}]\\
&\quad= (\lift{(g.e)}{\rewrite{M}{f}{D}{D}[f.h+ e.h].d}\circ \lift{(f.e)}{D[e.h].d})\circ e.) && [\defref{vacuousrewrite}]\\
&\quad= (\lift{(g.e)}{\rewrite{M[f.h].d}{f.e}{D[0].d}{D[e.h].d}}\circ \lift{(f.e)}{D[e.h].d})\circ e.) && [S(n)]\\
&\quad= \lift{(\lift{g.e}{M[f.h].d}\circ f.e)}{D[e.h].d}\circ e.e&& [Q(n{-}1)]\\
&\quad= \lift{((g\circ f).e)}{D[e.h].d}\circ e.e && [\defref{defcomposition}]\\
&\quad= ((g\circ f) \circ e).e && [\defref{defcomposition}]
\end{align*}

\end{itemize}
As all these conditions hold, both embeddings are equivalent, hence $A(n)$ holds and the implication is true. In addition, we only used the assumptions for $n>0$, hence A(0) also holds with no further requirements.
\end{proof}


With all these implications proved we can now bring them together to prove that the conjunction of selected logical statements holds for all $n\geq 0$.

\begin{theorem}
\label{bigtheorem}
The statement $S(n)\wedge R(n)\wedge P(n)\wedge C(n)\wedge A(n)\wedge Q(n)$ holds for all $n\geq 1$.
\end{theorem}
\begin{proof}
We prove the result by induction on $n$.

\begin{itemize}

\item \emph{Base case}: For $n=1$:

\begin{itemize}
\item To establish $S(1)$, by Lemma~\ref{TPimpliesS}, we need $P(0)$ to hold. $P(0)$ holds with no further conditions by Lemma~\ref{recursionbaseP}.

\item To establish $R(1)$, by Lemma~\ref{SimpliesR}, we need $S(1)$ and R(0) to hold. $S(1)$ is given by the statement above. $R(0)$ holds with no further conditions by Lemma~\ref{recursionbaseR}.

\item To establish $P(1)$, by Lemma~\ref{SQimpliesP}, we need $Q(0)$ and $S(1)$ to hold. $Q(0)$ holds, by Lemma~\ref{SPAimpliesQ} since $P(0)$ holds.

\item To establish $C(1)$, by Lemma~\ref{CQAimpliesC}, we need $C(0)$, $Q(0)$ and $A(0)$ to hold. $C(0)$ holds with no further conditions by Lemma~\ref{recursionbaseC}, similarly $A(0)$ by Lemma~\ref{CQAimpliesC}. $Q(0)$ holds by the  statement above.

\item To establish $A(1)$, by Lemma~\ref{SQAimpliesA}, we need $S(1)$ and $Q(0)$ and $A(0)$ to hold. The first two hold by the statements above, $A(0)$ holds by Lemma~\ref{SQAimpliesA}.

\item To establish $Q(1)$, by Lemma~\ref{SPAimpliesQ}, we need $S(1)$ and $P(1)$ and $A(0)$ to hold. The first two hold by the statements above, $A(0)$ holds by Lemma~\ref{SQAimpliesA}.
\end{itemize}

\item \emph{Inductive step}: For $n>1$, we assume that all $S(n{-}1)$, $R(n{-}1)$, $P(n{-}1)$, $C(n{-}1)$, $A(n{-}1)$, $Q(n{-}1)$ hold. 
\begin{itemize}
\item To establish $S(n)$, by Lemma~\ref{TPimpliesS}, we need $P(n{-}1)$, $C(n{-}1)$, $B(n{-}1)$  and $T(n)$ to hold. Since, $S(n{-}1)$ holds, by Lemma~\ref{SimpliesT} we obtain that $T(n)$ holds. Since $R(n{-}1)$ holds by Lemma~\ref{RimpliesL} we obtain that $B(n{-}1)$ holds. $P(n{-}1)$ and $C(n{-}1)$ hold by the inductive hypothesis.

\item To establish $R(n)$, by Lemma~\ref{SimpliesR}, we need $S(n)$, $T(n)$ and $R(n{-}1)$. We have that $S(n)$ and $T(n)$ hold by the statement above. $R(n{-}1)$ holds by the inductive hypothesis.

\item To establish $P(n)$, by Lemma~\ref{SQimpliesP}, we need $S(n)$ and $Q(n{-}1)$. $S(n)$ holds by the earlier part of this argument. $Q(n{-}1)$ holds by the inductive hypothesis.

\item To establish $Q(n{-}1)$, by Lemma~\ref{SPAimpliesQ}, we need $S(n{-}1)$ and $P(n{-}1)$ to hold, both hold by the inductive hypothesis.

\item To establish $C(n)$, by Lemma~\ref{CQAimpliesC}, we need $C(n{-}1)$, $A(n{-}1)$ and $Q(n{-}1)$ to hold. All hold by the inductive hypothesis.

\item To establish $A(n)$, by Lemma~\ref{SQAimpliesA}, we need $S(n)$, $Q(n{-}1)$ and $A(n{-}1)$ to hold. $S(n)$ holds by the statements above. $Q(n{-}1)$ and $A(n{-}1)$ hold by the inductive hypothesis.

\item To establish $Q(n)$, by Lemma~\ref{SPAimpliesQ}, we need $S(n)$, $P(n)$ and $A(n{-}1)$ to hold. $S(n)$ and $P(n)$ hold by the statements above. $A(n{-}1)$ holds by the inductive hypothesis.
\end{itemize}
Hence we obtain that $S(n)\wedge R(n)\wedge P(n)\wedge C(n)\wedge A(n)\wedge Q(n)$ holds.
\end{itemize}
By this inductive argument we established that the conjunction $S(n)\wedge R(n)\wedge P(n)\wedge C(n)\wedge A(n)\wedge Q(n)$ holds for all $n\geq 1$.
\end{proof}

We now bring all the above results together to obtain the desired result that the rewriting procedure produces well-defined $n$\-diagrams for any $n\geq 0$.

\begin{theorem}
\label{rewritewellfedined}
For any well-defined $n$-diagrams $D, S, T$ such that $S$ and $T$ are globular with respect to each other, and a well-defined embedding $e: S\hookrightarrow D$, the rewrite $\rewrite{D}{e}{S}{T}$ is a well-defined diagram.
\end{theorem}
\begin{proof} Two separate cases follow by different results proved previously:
\begin{itemize}
\item For $n=0$ follows by Lemma~\ref{recursionbaseR}.
\item For $n>0$, by Theorem~\ref{bigtheorem}, we have that $R(n)$ holds for any $n\geq 1$
\end{itemize}
By this we established that the rewrite $\rewrite{D}{e}{S}{T}$ of a well-defined $n$\-diagram $D$ is also well-defined.
\end{proof}

\begin{theorem}
Given well-defined $n$\-diagrams such that  $S,T$ are globular with respect to each other and for a well-defined embedding $e:S\hookrightarrow A$, the lifted embedding $\lift{e}{T}: T\hookrightarrow \rewrite{A}{e}{S}{T}$ is well-defined.
\end{theorem}
\begin{proof}
This result follows immediately by Lemma~\ref{RimpliesL}, since by Theorem~\ref{rewritewellfedined}, we have that  $R(n)$ is true for all $n\geq 0$.
\end{proof}

The final observation about the process of rewriting is that we may express it using higher level cells. The rewrite of an $n$\-diagram $D$ into $\rewrite{D}{e}{S}{T}$ gives rise to an $(n+1)$\-diagram $R$ of the following form:\JVcomm{I don't understand this} 
\begin{align*}
R = 
\begin{cases}
R.s = D \\
R[0].g = g & \text{ such that } g.s = S, g.t = T\\
R[0].e = e
\end{cases}
\end{align*}
This observation is crucial for our higher dimensional rewriting perspective and universal treatment of $(n+1)$\-cells as rewrites between composite $n$\-cells and gives us a wider expressivity to reason about higher\-level cells.

\subsection{Identity embeddings}
\label{sec:identityembeddings}

\begin{definition}
\label{identityembedding}
Given an  $n$\-diagram $D$, the \emph{identity embedding} $D.\id: D\hookrightarrow D$ is defined as follows:
\begin{itemize}
\item $D.\id.h = 0$
\item If $n>0$, then also: $D.\id.e = D.s.\id$
\end{itemize}
\end{definition}

\noindent
The intuition is that every diagram is trivially embedded in itself.

\begin{lemma}
\label{liftedidentityembedding}
Given well-defined $n$\-diagrams $S,T$ which are globular, the following holds:
\[
\lift{S.\id}{T} = T.\id\]
\end{lemma}
\begin{proof}
If $n=0$, by Definition~\ref{embeddingequiv}, to show that these embeddings are equivalent, we need to check that the domains and codomains of these embeddings are equivalent diagrams. The types of these embeddings are as follows:
\begin{align*}
&T.\id: T \hookrightarrow T & [\defref{identityembedding}] \\
&\lift{S.\id}{T} : T \hookrightarrow \rewrite{S}{S.\id}{S}{T} & [\defref{liftedembedding}] 
\end{align*}
By applying Definition~\ref{defrewrite} we obtain that $\rewrite{S}{S.\id}{S}{T} = T$, hence the domains and codomains of both embeddings match.

If $n>0$, we additionally need to check that:
\begin{itemize}

\item Component embeddings are equivalent:
\begin{align*}
&\lift{S.\id}{T}.e = S.\id.e & [\defref{liftedembedding}] 
\end{align*}

\item Heights are equal:
\begin{align*}
&\lift{S.\id}{T}.h = S.\id.h & [\defref{liftedembedding}] 
\end{align*}
\end{itemize}
By this argument both embeddings are equivalent, as required.
\end{proof}

A desired property of an identity embedding is that composing it with any other diagram embedding $e$, should have no effect on $e$.
\begin{lemma}
\label{identitycancellation}
Given a well-defined $n$\-diagram $D$, the identity embedding $D.\id: D\hookrightarrow D$ and a well-defined embeddings $e: A\hookrightarrow D$ and $f: D\hookrightarrow B$, the following holds:
\begin{align*}
&D.\id\circ e = e \\
&f\circ D.\id = f 
\end{align*}
\end{lemma}
\begin{proof}
Let us prove that $D.\id\circ e = e$ by induction on $n$. 

\begin{itemize}

\item \emph{Base case:} For $n=0$, by Definition~\ref{embeddingequiv} domains and codomains of both embeddings must be equal. The types are as follows:
\begin{align*}
&e : A \hookrightarrow D \\
&D.\id\circ e : A \hookrightarrow D \\
\end{align*}
Domains and codomains trivially are equivalent diagrams.

\item \emph{Inductive step:} For $n>0$, we assume that the result holds for $(n{-}1)$\-diagram embeddings \emph{(IH)}. Types are the same as in the base case, we additionally need to check that:
\begin{itemize}

\item Component embeddings are equivalent:
\begin{align*}
&(D.\id\circ e).e =  \\
&\quad= \lift{(D.\id).e}{D[e.h].d} \circ e.e & [\defref{defcomposition}] \\
&\quad= \lift{(D[0].d.\id)}{D[e.h].d} \circ e.e & [\defref{identityembedding}] \\
&\quad= D[e.h].d.\id \circ e.e & [\defref{liftedidentityembedding}] \\
&\quad= e.e & [\text{\emph{IH}}] 
\end{align*}

\item Heights are equal:
\begin{align*}
&(D.\id\circ e).h = \\
&\quad= (D.\id).h + e.h & [\defref{defcomposition}] \\
&\quad= e.h & [\defref{identityembedding}] 
\end{align*}
\end{itemize}
\end{itemize}

By this inductive argument both embeddings are equivalent for $n\geq 0$, as required. The proof of the other equivalence follows similarly.
\end{proof}

With the aid of these results we show that the identity embedding is well-defined.
\begin{lemma}
\label{identitywelldef}
Given a well-defined $n$\-diagram $D$, the identity embedding $D.\id: D\hookrightarrow D$ is well-defined.
\end{lemma}
\begin{proof}
$D.\id$ is an endomorphism of a well-defined diagram $D$, hence trivially both its domain and codomain are well-defined. We prove well-definedness of this embedding by induction on $n$.
\begin{itemize}
\item \emph{Base case:} For $n=0$, by Definition~\ref{defembwelldef} we just need the single cell in the domain being equal to the single cell in the codomain \emph{i.e.} $D[0].g = D[0].g$ which holds.

\item \emph{Inductive step:} For $n>0$, assume that all identity embeddings of $(n{-}1)$\-diagrams are well-defined.

By Definition~\ref{defembwelldef}  we need to check three conditions:

\begin{itemize}
\item The component embedding is well-defined. By Definition~\ref{identityembedding}, we have $D.\id.e = D.s.\id = D[0].d.\id$. Since $D[0].d$ is an $(n{-}1)$\-diagram, $D.\id.e$ is well-defined by the inductive hypothesis.

\item Corresponding generators are equal. For $0\leq i\leq |D|$ we have:
\begin{align*}
&(D[i + (D.\id).h].g = D[i].g & [\defref{identityembedding}] 
\end{align*}

\item Corresponding embeddings satisfy the following for $0\leq i\leq |D|$: 
\begin{align*}
&\lift{((D.\id).e)}{D[i].d}\circ D[i].e  \\
&\quad= \lift{(D.s.\id)}{D[i].d}\circ D[i].e & [\defref{identityembedding}] \\
&\quad= \lift{(D[0].d.\id)}{D[i].d}\circ D[i].e  \\
&\quad= D[i].d.\id \circ D[i].e  & [\defref{liftedidentityembedding}] \\
&\quad= D[i].e  & [\defref{identitycancellation}] 
\end{align*}
As required.
\end{itemize}
As all these conditions hold, $D.\id$ is well-defined.
\end{itemize}
By this inductive argument we established that $D.\id$ is well-defined for any $n$\-diagram $D$ for all $n\geq 0$.
\end{proof}

\subsection{Diagram composition}
\label{sec:diagramcomposition}

Here we prove correctness of statements concerning diagram composition. Intuitively, composition involves  gluing together diagrams along a common boundary. For instance for the following diagrams $S$ and~$D$:
\begin{align*}
S\quad=\quad
&
\begin{aligned}
\begin{tikzpicture}[thick, scale = 0.7]
\draw (-0.5,0) to (-0.5,0.5) to [out = up, in = left] (0,1)
to [out = right, in = up](0.5,0.5) to(0.5,0);
\draw (-0.5,3) to (-0.5,2.5) to [out = down, in = left] (0,2)
to [out = right, in = down](0.5,2.5) to(0.5,3);
\draw (0,1) to (0, 2);
\node [blob] at (0,1) {};
\node [blob] at (0,2) {};
\end{tikzpicture}
\end{aligned}
&
D\quad=\quad
&
\begin{aligned}
\begin{tikzpicture}[thick, scale = 0.7]
\draw (0.5,1) to(0.5,-2);
\draw (-0.5,1) to(-0.5,-2);
\node [blob] at (0.5,0) {};
\node [blob] at (-0.5,-1) {};
\end{tikzpicture}
\end{aligned}
\end{align*}
Their composite, which we denote by $S\circ D$ is:
\begin{align*}
S\circ D\quad=\quad
\begin{aligned}
\begin{tikzpicture}[thick, scale = 0.7]
\draw (-0.5,0) to (-0.5,0.5) to [out = up, in = left] (0,1)
to [out = right, in = up](0.5,0.5) to(0.5,0);
\draw (-0.5,5) to (-0.5,2.5) to [out = down, in = left] (0,2)
to [out = right, in = down](0.5,2.5) to(0.5,5);
\draw (0,1) to (0, 2);
\node [blob] at (0,1) {};
\node [blob] at (0,2) {};
\node [blob] at (0.5,4) {};
\node [blob] at (-0.5,3) {};
\end{tikzpicture}
\end{aligned}
\end{align*}
Note that we can also compose diagrams that are of different dimension, for instance given a 1\-diagram $M$:
\begin{align*}
M\quad=\quad
&
\begin{aligned}
\begin{tikzpicture}[thick, scale = 0.7]
\draw (-1,0) to (1, 0);
\node [blob] at (0,0) {};
\end{tikzpicture}
\end{aligned}
\end{align*}
The composite $M\circ D$ is:
\begin{align*}
\begin{aligned}
\begin{tikzpicture}[thick, scale = 0.7]
\draw (-1.5,1) to(-1.5,-2);
\draw (0.5,1) to(0.5,-2);
\draw (-0.5,1) to(-0.5,-2);
\node [blob] at (0.5,0) {};
\node [blob] at (-0.5,-1) {};
\end{tikzpicture}
\end{aligned}
\end{align*}
When we compose two diagrams, the higher-dimensional of the two---or both, if they have the same dimension---embed into the composite. We express this with \emph{inclusion} embeddings, the \emph{right inclusion} $\inc{S}{D}: D\hookrightarrow S\circ D$ and the \emph{left inclusion} $\increverse{S}{D}: D\hookrightarrow D\circ S$. We write $D.s^k$ and $D.t^k$ to denote taking the source or the target of the diagram $k$ times.

\begin{definition}
\label{definclusion}
Given an $n$\-diagram $D$ and an $m$\-diagram $S$:
\begin{itemize}
\item If $n\geq m$ and $S.t= D.s^{n-m+1}$, we define the \emph{inclusion} embedding $\inc{S}{D}: D\hookrightarrow S\circ D$ in the following way:

\begin{itemize}
\item If $n=m$, then:
\begin{align*}
\inc{S}{D}.h &= |S| \\
\inc{S}{D}.e &= D.\id.e
\end{align*}
\item If $n>m$, then:
\begin{align*}
\inc{S}{D}.h &= D.\id.h \\
\inc{S}{D}.e &= \inc{S}{D.s}
\end{align*}
\end{itemize}

\item If $n \leq m$ and $S.t^{m-n+1} = D.s$, we define the inclusion embedding $\increverse{S}{D}: S\hookrightarrow S\circ D$ in the following way:
\begin{align*}
\increverse{S}{D}.h &= S.\id.h \\
\increverse{S}{D}.e &= \inc{S.s}{D}
\end{align*}
\end{itemize}
\end{definition}
Note that in the latter case, the composition of $D$ with $S$ has no effect on the embedding data. However we do need to change the embedding codomain, so that the types match. 

Let us illustrate this with two examples, first consider the following diagrams $S$ and $D$ that are of the same dimension:
\begin{align*}
D\quad:=\quad
&
\begin{aligned}
\begin{tikzpicture}[thick, scale = 0.7]
\draw (0.5,1) to(0.5,-2);
\draw (-0.5,1) to(-0.5,-2);
\node [blob] at (0.5,0) {};
\node [blob] at (-0.5,-1) {};
\end{tikzpicture}
\end{aligned}
\qquad
\stackrel{\inc{S}{D}}{\hookrightarrow}
\qquad
\begin{aligned}
\begin{tikzpicture}[thick, scale = 0.7]
\draw (-0.5,0) to (-0.5,0.5) to [out = up, in = left] (0,1)
to [out = right, in = up](0.5,0.5) to(0.5,0);
\draw (-0.5,5) to (-0.5,2.5) to [out = down, in = left] (0,2)
to [out = right, in = down](0.5,2.5) to(0.5,5);
\draw (0,1) to (0, 2);
\node [blob] at (0,1) {};
\node [blob] at (0,2) {};
\node [blob] at (0.5,4) {};
\node [blob] at (-0.5,3) {};
\node [rectangle, minimum width = 30pt, minimum height = 40pt, draw, fill = none, color = blue, dashed] at (0,3.5) {};
\end{tikzpicture}
\end{aligned}
\quad:=\quad S\circ D
\end{align*}
In this picture, we indicate the inclusion embedding $\inc{S}{D}$ by placing a dashed rectangle over the instance of $D$ appearing in $S\circ D$. The diagram $S$ also embeds into the composite.

Now, consider $S$ and $D$ such that the dimension of $D$ is larger:
\begin{align*}
D\quad:=\quad
&
\begin{aligned}
\begin{tikzpicture}[thick, scale = 0.7]
\draw (0.5,1) to(0.5,-2);
\draw (-0.5,1) to(-0.5,-2);
\node [blob] at (0.5,0) {};
\node [blob] at (-0.5,-1) {};
\end{tikzpicture}
\end{aligned}
\qquad
\stackrel{\inc{S}{D}}{\hookrightarrow}
\qquad
\begin{aligned}
\begin{tikzpicture}[thick, scale = 0.7]
\draw (-1.5,1) to(-1.5,-2);
\draw (0.5,1) to(0.5,-2);
\draw (-0.5,1) to(-0.5,-2);
\node [blob] at (0.5,0) {};
\node [blob] at (-0.5,-1) {};
\node [rectangle, minimum width = 30pt, minimum height = 40pt, draw, fill = none, color = blue, dashed] at (0,-0.5) {};
\end{tikzpicture}
\end{aligned}
\quad:=\quad S\circ D
\end{align*}
Here again, we indicate $\inc{S}{D}$ by the dashed rectangle. In this case, $S$ does not embed into the composite, since it is of a lower dimension.

In this setup, an $n$\-diagram $D$ and an $m$\-diagram S can be composed in exactly one way, along their common $(\text{min}(m,n) {-} 1)$-dimensional boundary. In particular, two $n$\-diagrams can only be composed along their matching target and source, corresponding to vertical composition. This approach is still fully general and permits construction of generic perturbations of arbitrary composites.

We present a recursive definition of composition of an $n$\-diagram $D$ and an $m$\-diagram $S$. If $n>m$ we specify all the generators and embeddings of the composite and then we refer recursively to the $(n{-}1)$\-dimensional source of $D$ and to the $m$\-diagram $S$. That way with each recursive call $n$ decreases, hence $n-m$ decreases. Eventually, we decrease $n$ sufficiently that $n=m$, and the recursion terminates with the base clause. The case for $m>n$ is analogous. This ensures that the definition is well\-founded.
\begin{definition} 
\label{diagcompose}
Given an $n$\-diagram $D$ and an $m$-diagram $S$ such that $S.t= D.s^{n-m+1}$ if $m\leq n$ or $S.t^{m-n+1}=D.s$ otherwise, the \emph{composite} diagram $S\circ_{m, n} D$ is defined as follows:
\begin{itemize}
\item If $n=m$, the individual components of $S\circ_{m, n} D$ are as follows:\JVcomm{Needs aligning}
\begin{align}
&\begin{aligned}
\label{diagcomposesourcesnmequal}
&(S\circ_{m, n} D).s = S.s \\
\end{aligned} \\
&\begin{aligned}
\label{diagcomposesizesnmequal}
&|S\circ_{m, n} D| = |D|+|S|   \\
\end{aligned} \\
&\begin{aligned}
\label{diagcomposegeneratorsnmequal}
(S\circ_{m, n} D)[j].g =
\begin{cases}
S[j].g  & \text{for } 0\leq j < |S| \\
D[j - |S|].g  & \text{for } |S|\leq j < |D|+|S|  
\end{cases}
\end{aligned} \\
&\begin{aligned}
\label{diagcomposeembeddingsnmequal}
&(S\circ_{m, n} D)[j].e =
\begin{cases}
S[j].e  & \text{for } 0\leq j < |S| \\
D[j - |S|].e  & \text{for } |S|\leq j < |D|+|S|  
\end{cases}
\end{aligned} 
\end{align}
\item If $m<n$, the individual components of $S\circ_{m, n} D$ are as follows:
\begin{align}
&\begin{aligned}
\label{diagcomposesourcesnmgreater}
(S\circ_{m, n} D).s = S\circ_{m, n-1} D.s
\end{aligned} \\
&\begin{aligned}
\label{diagcomposesizesnmgreater}
|S\circ_{m, n} D)| = |D|
\end{aligned} \\
&\begin{aligned}
\label{diagcomposegeneratorsnmgreater}
(S\circ_{m, n} D)[j].g =
&D[j].g  & \quad\qquad\qquad\qquad\text{for } 0\leq j < |D| 
\end{aligned} \\
&\begin{aligned}
\label{diagcomposeembeddingsnmgreater}
&(S\circ_{m, n} D)[j].e =
\inc{S}{D[j].d}\circ D[j].d  & \text{for } 0\leq j < |D|  
\end{aligned} 
\end{align}

\item If $n<m$, the individual components of $S\circ_{m, n} D$ are as follows:
\begin{align}
&\begin{aligned}
\label{diagcomposesourcesnmsmaller}
(S\circ_{m, n} D).s = S.s\circ_{m-1, n} D\\
\end{aligned} \\
&\begin{aligned}
\label{diagcomposesizesnmsmaller}
|S\circ_{m, n} D| = |S|\\
\end{aligned} \\
&\begin{aligned}
\label{diagcomposegeneratorsnmsmaller}
(S\circ_{m, n} D)[j].g =
&S[j].g  & \text{for } 0\leq j < |S| \\
\end{aligned} \\
&\begin{aligned}
\label{diagcomposegeneratorsnmsmaller}
(S\circ_{m, n} D)[j].e =
&\increverse{S[j].d}{D}\circ S[i].e  & \text{for } 0\leq j < |S| 
\end{aligned} 
\end{align}
\end{itemize}
\end{definition}
As stated before, the type of composition is uniquely determined by the dimensions of the diagrams being composed. As a shorthand for a composite of an $n$\-diagram $D$ and an $m$\-diagram $S$, we label $S\circ_{m,n} D$ by $a=\min(n,m) -1$ instead. That way the terminology more closely matches the naming scheme for composition which is standard in category theory. This leads to a certain overloading of notation since we use the same label for composites $S\circ_{a} D$ and $D\circ_{a} S$. In some cases we omit the subscripts entirely and simply write $D\circ S$.

For composition of two $n$\-diagrams $S, D$, intuitively we just concatenate their lists of generators and embeddings. Due to that, the lists of slices also get concatenated. This is formalised by the following lemma.
\begin{lemma}
\label{composedslices}
Given well-defined $n$\-diagrams $S, D$, such that $S.t=D.s$, the following holds for slices of the composed diagram for $0\leq j\leq |S\circ D|$:
\begin{align*}
(S\circ D)[j].d=
\begin{cases}
S[j].d & 0\leq j \leq |S| \\
D[j - |S|].d & |S| \leq j \leq |S| + |D| 
\end{cases} 
\end{align*}
\end{lemma}
\begin{proof}
Note that, by~\ref{diagcomposesizesnmequal}, we obtain that: $|S\circ D|= |S|+|D|$, so let $0 \leq j \leq |S|+|D|$ and consider two ranges separately:

\begin{itemize}

\item For $0 \leq j \leq |S|$ we show the result by induction on $j$:

\begin{itemize}

\item \emph{Base case:}  For $j=0$, we have the following:
\begin{align*}
&(S\circ D)[0].d  \\
&\quad= (S\circ D).s &&[\defref{definitialterminalslice}]\\
&\quad= S.s &&[\myeqref{diagcomposesourcesnmequal}]\\
&\quad= S[0].d &&[\defref{definitialterminalslice}]
\end{align*}

\item \emph{Inductive step:} For $0 < j \leq |S|$, we assume that $(S\circ D)[j].d= S[j].d$. Let us consider $(S\circ D)[j+1].d$:
\begin{align*}
&(S\circ D)[j+1].d   \\
&\quad= \rewrite{(S\circ D)[j].d}{(S\circ D)[j].e}{(S\circ D)[j].g.s}{(S\circ D)[j].g.t} &&[\defref{defslice}]\\
&\quad= \rewrite{(S\circ D)[j].d}{(S\circ D)[j].e}{S[j].g.s}{S[j].g.t} &&[\myeqref{diagcomposegeneratorsnmequal}]\\
&\quad= \rewrite{(S\circ D)[j].d}{S[j].e}{S[j].g.s}{S[j].g.t} &&[\myeqref{diagcomposeembeddingsnmequal}]\\
&\quad= \rewrite{S[j].d}{S[j].e}{S[j].g.s}{S[j].g.t} &&[\text{\emph{IH}}]\\
&\quad= S[j+1].d &&[\defref{defslice}]
\end{align*}

\end{itemize}

\item For $|S| \leq j \leq |S| + |D|$, note that by the argument for $0 \leq j \leq |S|$ we obtain:
\begin{align*}
&(S\circ D)[|S|].d  \\
&\quad= S[|S|].d \\
&\quad= S.t &&[\defref{definitialterminalslice}]\\
&\quad= D.s &&[\text{assumption}]\\ 
&\quad= D[0].d &&[\defref{definitialterminalslice}]
\end{align*}
The result for $|S| \leq j \leq |S| + |D|$ is proved inductively in a similar way as the result for $0 \leq j \leq |S|$ .
\end{itemize}
By these two inductive arguments the result is shown for $0 \leq j \leq |S\circ D|$.
\end{proof}

\paragraph{Correctness.}
We now consider correctness of the composition construction. Again, due to the heavily recursive nature of the diagram and signature structures and their mutual references, we introduce several logical statements about properties of composite diagrams and inclusion embeddings. For composition of an $n$\-diagram and an $m$\-diagram, the statements depend on an integer $k\geq 0$ such that $k=|n-m|$. Hence the base case in our recursion covers the situation when the two diagrams composed are of the same dimension. Again, there is the main inductive proof that the conjunction of logical statements $K(k), L(k), M(k), N(k)$ holds for all $k\geq 0$. Since we make no other assumptions about $m,n$ this ensures that we show the results for any combination of dimensions. Again, the overall structure is that we first present the statements, then we show that they hold for $k=0$ to establish base cases. This is followed by proofs of a series of implications and finally, in the conclusion of this section, all the lemmas are brought together to establish the main result on the composite of two diagrams being well-defined.

Below, we define four logical statements on properties of diagram composition and inclusion embeddings:
\begin{definition}
For $k\geq 0$, $L(k)$ states that given an $n$\-diagram $D$ and an  $m$-diagram $S$, both well-defined, such that $|n-m|=k$, and such that  $S.t=s^{n-m+1}(D)$ if $m\leq n$ or $ t^{m-n+1}(S)=D.s$ otherwise, then $S\circ D$ is well-defined.
\end{definition}

\begin{definition}
For $k\geq 0$, $N(k)$ states that given an $n$\-diagram $D$ and  an $m$\-diagram $S$, both well-defined, such that $|n-m|=k$, then:
\begin{itemize}
\item if $n\geq m$ and $S.t=s^{n-m+1}(D)$, then $\inc{S}{D}: D\hookrightarrow S\circ D$   is well-defined; 
\item if $n< m$ and $S.t^{m-n+1}= D.s$, then $\increverse{S}{D}: D\hookrightarrow D\circ S$ is well-defined. 
\end{itemize}
\end{definition}

\noindent
A slice of the composite diagram can be written as a composite of the diagram of lower dimension and the appropriate slice of the diagram of higher dimension. 
\begin{definition}
For $k\geq 1$, $K(k)$ states that given an $n$\-diagram $D$ and an $m$\-diagram $S$, both well-defined, such that $|n-m|=k$ and such that the composite $S\circ D$ exists, the following equalities hold:
\begin{itemize}
\item if $n>m$, then $(S\circ D)[i].d= S\circ (D[i].d)$ for any $0\leq i < |D|$;
\item if $n<m$, then $(S\circ D)[i].d= (S[i].d) \circ D$ for any $0\leq i < |S|$.
\end{itemize}
\end{definition}

\noindent
The inclusion embedding can instead be expressed using the lifted embedding.
\begin{definition}
For $k\geq 1$, $M(k)$ states that given an $n$\-diagram $D$ and an  $m$\-diagram $S$, both well-defined, such that $|n-m|=k$, then for any $0\leq i < |D|$ the following equality holds:
\begin{itemize}
\item if $n>m$, then $\inc{S}{D[i].d}=\lift{\inc{S}{D}.e}{D[i].d}$ for any $0\leq i < |D|$;
\item if $n<m$, then$\increverse{S[i].d}{D}=\lift{\increverse{S}{D}.e}{S[i].d}$ for any $0\leq i < |S|$.
\end{itemize}
\end{definition}

\paragraph{Base cases.}
We now prove several lemmas that establish base cases for the main recursive proof.
\begin{lemma}
\label{recursionbaseN}
For $k=0$ the following holds: $L(0)\implies N(0)$
\end{lemma}
\begin{proof}
We assume that L(0) holds, \emph{i.e.} any two well-defined $n$\-diagrams $D,S$ such that $S.t=D.s$, the \emph{composite} diagram $S\circ D$ is well-defined.

As $k=0$, implies $n=m$, we only need to show that the inclusion embedding $\inc{S}{D}: D\hookrightarrow S\circ D$ is well-defined for $N(0)$ to hold.

The domain $D$ of the embedding $\inc{S}{D}$ is well-defined by assumption. The codomain $S\circ D$ is well-defined by  $L(0)$. Since $n,m\geq 1$, by Definition~\ref{defembwelldef} $\inc{S}{D}$ needs to satisfy three separate conditions:

\begin{itemize}

\item The component embedding $\inc{S}{D}.e$ is well-defined. As $n=m$, by Definition~\ref{definclusion}, $\inc{S}{D}.e = D.\id.e = D.s.\id$, which is well-defined by Lemma~\ref{identitywelldef}.

\item For every $0\leq j< |D|$ we have: 
\begin{align*}
&(S\circ D)[j+\inc{S}{D}.h].g \\
&\quad= (S\circ D)[j+|S|].g &&[\defref{definclusion}]\\
&\quad= D[(j+|S|)-|S|].g &&[\myeqref{diagcomposegeneratorsnmequal}]\\
&\quad= D[j].g 
\end{align*}
As required.

\item For every $0\leq j< |D|$ we have:  
\begin{align*}
& \lift{\inc{S}{D}.e}{D[j].d}\circ D[j].e  \\
&\quad= \lift{D.\id.e}{D[j].d}\circ D[j].e &&[\defref{definclusion}]\\
&\quad= \lift{D.s.\id}{D[j].d}\circ D[j].e &&[\defref{identityembedding}]\\
&\quad= \lift{D[0].d.\id}{D[j].d}\circ D[j].e &&[\defref{definitialterminalslice}]\\
&\quad= D[j].d.\id\circ D[j].e &&[\defref{liftedidentityembedding}]\\
&\quad= D[j].e &&[\defref{identitycancellation}]\\
&\quad= D[(j+|S|)-|S|].e \\
&\quad= (S\circ D)[j+|S|].e &&[\myeqref{diagcomposegeneratorsnmequal}]\\
&\quad=(S\circ D)[j+\inc{S}{D}.h].e &&[\defref{definclusion}]
\end{align*}
As required.

\end{itemize}
As all three conditions are satisfied, we conclude that for any two well-defined $n$\-diagrams $D, S$ such that  $S.t=D.s$, the inclusion embedding $\inc{S}{D}$ is well-defined.
\end{proof}

\begin{lemma}
\label{recursionbaseL}
$L(0)$.
\end{lemma}
\begin{proof}
For $k=0$, we have $n=m$, so to establish $L(0)$ we need to show that given two well-defined $n$\-diagrams $S,D$ such that $S.t=D.s$, the composite $S\circ D$ is well-defined.

Since $m,n\geq 1$, by Definition~\ref{defwellformed}, the diagram $S\circ D$ is well-defined if the following two conditions are satisfied:

\begin{itemize}
\item The source $(S\circ D).s = (S\circ D)[0].d$ is a well-defined diagram. 
\item For every $0 < j \leq |S\circ D|$ the slice $(S\circ D)[j].d$ exists and is well-defined.
\end{itemize}

As $n=m$, we apply Lemma~\ref{composedslices} to obtain that:
\begin{align*}
(S\circ D)[j].d=
\begin{cases}
S[j].d & 0\leq j \leq |S| \\
D[j - |S|].d & |S| \leq j \leq |S| + |D| 
\end{cases} 
\end{align*}
As both $S, D$ are well-defined diagrams, then all their slices are also well-defined. Since every slice of $S\circ D$ is equal to either a slice of $S$ or a slice of $D$, they are all well-defined, hence $S\circ D$ is also well-defined.
\end{proof}


\begin{lemma}
\label{recursionbaseK}
$K(1)$.
\end{lemma}
\begin{proof}

We need to show that for any $n$\-diagram $D$ and $m$\-diagram $S$ such that the composite $S\circ D$ exists and \begin{align*}
&\text{If } n>m& (S\circ D)[i].d&= S\circ (D[i].d)&&\text{for any } 0\leq i < |D| \\
& \text{If } n<m& (S\circ D)[i].d&= (S[i].d) \circ D&&\text{for any } 0\leq i < |S|
\end{align*}
We consider both these cases separately, first let $n>m$. We prove the result by induction on $0\leq j \leq |D|$.

\begin{itemize}

\item \emph{Base case:} For $j=0$, the result follows immediately form the definitions:
\begin{align*}
&(S\circ D)[0].d \\
&\quad=(S\circ D).s && [\defref{definitialterminalslice}]\\
&\quad= S\circ(D.s) && [\defref{diagcompose}] \\
&\quad= S\circ(D[0].d) && [\defref{definitialterminalslice}]
\end{align*}

\item \emph{Inductive step:} For $j>0$, assume that:
\begin{align*}
(S\circ D)[j].d = S\circ (D[j].d) &\qquad (\text{\emph{IH}})
\end{align*}
Let us consider $(S\circ D)[j+1].d$ and $S\circ (D[j+1].d)$, to show that two diagrams are equivalent, by Definition~\ref{diagequiv}, we need to check the following:

\begin{itemize}

\item Sources are equivalent diagrams:
\begin{align*}
&((S\circ D)[j+1].d).s \\
&\quad=  \rewrite{(S\circ D)[j].d}{(S\circ D)[j].e}{(S\circ D)[j].g.s}{(S\circ D)[j].g.t}).s && [\defref{defslice}]\\
&\quad=(S\circ D)[j].d.s && [\myeqref{defrewritesource}]\\
&\quad=(S\circ (D[j].d)).s && [\text{\emph{IH}}] \\
&\quad=S.s &&[\myeqref{diagcomposesourcesnmequal}] \\
&\quad=(S\circ (D[j+1].d)).s &&[\myeqref{diagcomposesourcesnmequal}] 
\end{align*}

\item Lengths of generator lists are equal:

Note that even though $m = n-1$, we still have $m<n$, hence the following is true for generators of the composed diagram for $0\leq j < |D|$:
\begin{align*}
(S\circ D)[j].g = D[j].g &\qquad [\defref{diagcomposegeneratorsnmgreater}]\\
\end{align*}
The comparison of lengths for the two diagrams is as follows:
\begin{align*}
&|(S\circ D)[j+1].d| \\
&\quad=  |\rewrite{(S\circ D)[j].d)}{(S\circ D)[j].e}{(S\circ D)[j].g.s}{(S\circ D)[j].g.t}| && [\defref{defslice}]\\
&\quad= |S\circ D)[j].d| - |(S\circ D)[j].g.s| + |(S\circ D)[j].g.t| && [\myeqref{defrewritesize}]\\
&\quad= |S\circ (D[j].d)| - |(S\circ D)[j].g.s| + |(S\circ D)[j].g.t| &&[IH] \\
&\quad= |S\circ (D[j].d)| - |D[j].g.s| + |D[j].g.t| &&[\myeqref{diagcomposegeneratorsnmgreater}]\\
&\quad= |S|+ |D[j].d| - |D[j].g.s| + |D[j].g.t| &&[\myeqref{diagcomposesizesnmequal}]\\
&\quad=|S|+|D[j+1].d| && [\myeqref{defrewritesize}]\\
&\quad=|S\circ (D[j+1].d)| &&[\myeqref{diagcomposesizesnmequal}]
\end{align*}

\item For generators and embeddings we need to show that for $0\leq k < |(S\circ D)[j+1].d| = |S\circ (D[j+1].d)|$, the $k$\-th generators in generator lists of both diagrams correspond and the same for $k$\-th embeddings. 


Since $m=n-1$, we can simplify the height of the $j$\-th embedding in $S\circ D$ in the following way. 
\begin{align*}
&((S\circ D)[j].e).h = \\ 
&\quad=(\inc{S}{D[j].d}\circ D[j].e).h  &&[\myeqref{diagcomposeembeddingsnmequal}]\\
&\quad= (\inc{S}{D[j].d}).h+ D[j].e.h &&[\defref{defcomposition}]\\
&\quad= |S| + D[j].e.h &&[\defref{definclusion}]
\end{align*}
Let us refer to this equality as $[*]$.

We show the necessary equivalences separately for four individual ranges:

\begin{itemize}

\item In the range: $0\leq k < |S|$, we have: 

Generators:
\begin{align*}
&((S\circ D)[j+1].d)[k].g \\
&\quad=(\rewrite{(S\circ D)[j].d)}{(S\circ D)[j].e}{(S\circ D)[j].g.s}{\\&\qquad\qquad (S\circ D)[j].g.t}).g && [\defref{defslice}]\\
&\quad=((S\circ D)[j].d)[k].g& & [\myeqref{defrewritegenerator}]\\
&\quad=(S\circ (D[j].d))[k].g && [IH]\\
&\quad= S[k].g &&[\myeqref{diagcomposegeneratorsnmequal}, k<|S|]\\
&\quad=(S\circ (D[j+1].d))[k].g &&[\myeqref{diagcomposegeneratorsnmequal}, k<|S|]
\end{align*}

Embeddings: 
\begin{align*}
&((S\circ D)[j+1].d)[k].e \\
&\quad=(\rewrite{(S\circ D)[j].d)}{(S\circ D)[j].e}{(S\circ D)[j].g.s}{\\&\qquad\qquad (S\circ D)[j].g.t})[k].e && [\defref{defslice}]\\
&\quad=((S\circ D)[j].d)[k].e && [\myeqref{defrewriteembedding}]\\
&\quad=(S\circ (D[j].d))[k].e && [IH]\\
&\quad= S[k].e &&[\myeqref{diagcomposeembeddingsnmequal}, k<|S|]\\
&\quad=(S\circ (D[j+1].d))[k].e &&[\myeqref{diagcomposegeneratorsnmequal}, k <|S|]
\end{align*}

\item In the range: $|S| \leq k < (S\circ D)[j].e.h$
\newline
We have: $(S\circ D)[j].e =D[i-|S|].e$, hence the range is equivalent to: 
\begin{align*}
|S| \leq k < D[i-|S|].e.h
\end{align*}
This means that when we apply Definition~\ref{defrewrite}, we use generators and embeddings from the original diagram, not from the target of the rewrite.

Generators:
\begin{align*}
&((S\circ D)[j+1].d)[k].g \\
&\quad=(\rewrite{(S\circ D)[j].d)}{(S\circ D)[j].e}{(S\circ D)[j].g.s}{\\&\qquad\qquad (S\circ D)[j].g.t}).g && [\defref{defslice}]\\
&\quad=((S\circ D)[j].d)[k].g && [\myeqref{defrewritegenerator}]\\
&\quad=(S\circ (D[j].d))[k].g && [IH]\\
&\quad= (D[j].d)[k - |S|].g &&[\myeqref{diagcomposegeneratorsnmequal}, |S|\leq k]\\
&\quad= (\rewrite{D[j].d}{D[j].d}{D[j].g.s}{D[j].g.t})[k-|S|].g && [\myeqref{defrewritegenerator}]\\
&\quad= (D[j+1])[(k- |S|].g && [\defref{defslice}]\\
&\quad=(S\circ (D[j+1].d))[k].g &&[\myeqref{diagcomposegeneratorsnmequal}, |S|\leq k]\\
\end{align*}
Embeddings:
\begin{align*}
&((S\circ D)[j+1].d)[k].e \\
&\quad=(\rewrite{(S\circ D)[j].d)}{(S\circ D)[j].e}{(S\circ D)[j].g.s}{\\&\qquad\qquad (S\circ D)[j].g.t}).e && [\defref{defslice}]\\
&\quad=((S\circ D)[j].d)[k].e && [\myeqref{defrewriteembedding}]\\
&\quad=(S\circ (D[j].d))[k].e && [IH]\\
&\quad= (D[j].d)[k - |S|].e &&[\myeqref{diagcomposeembeddingsnmequal}, |S|\leq k]\\
&\quad= (\rewrite{D[j].d}{D[j].d}{D[j].g.s}{D[j].g.t})[k-|S|].e && [\myeqref{defrewriteembedding}]\\
&\quad= (D[j+1])[(k- |S|].e && [\defref{defslice}]\\
&\quad=(S\circ (D[j+1].d))[k].e &&[\myeqref{diagcomposeembeddingsnmequal}, |S|\leq k]
\end{align*}

\item In the range: $(S\circ D)[j].e.h\leq k < (S\circ D)[j].e.h + |(S\circ D)[j].g.t|$
\newline
We have: $(S\circ D)[j].g =D[j].g$, hence the range is equivalent to: 
\begin{align*}
|S| + D[j].e.h\leq k < |S| + D[j].e.h + |D[j].g.t|
\end{align*}

Generators:
\begin{align*}
&((S\circ D)[j+1].d)[k].g \\
&\quad=(\rewrite{(S\circ D)[j].d)}{(S\circ D)[j].e}{(S\circ D)[j].g.s}{\\&\qquad\qquad (S\circ D)[j].g.t}).g && [\defref{defslice}]\\
&\quad= (S\circ D)[j].g.t[k - (S\circ D)[j].e.h] && [\myeqref{defrewritegenerator}]\\
&\quad= D[j].g.t[k - (S\circ D)[j].e.h].g &&[\myeqref{diagcomposegeneratorsnmgreater}]\\
&\quad= D[j].g.t[(k- (D[j].e.h + |S|)].g &&[*] \\
&\quad= D[j].g.t[(k- |S|) - D[j].e.h].g \\
&\quad= (D[j+1].d)[(k- |S|].g && [\myeqref{defrewritegenerator}]\\
&\quad=(S\circ (D[j+1].d))[k].g &&[\myeqref{diagcomposegeneratorsnmequal}, |S|\leq k]
\end{align*}

Embeddings:
\begin{align*}
&((S\circ D)[j+1].d)[k].e \\
&\quad=(\rewrite{(S\circ D)[j].d)}{(S\circ D)[j].e}{(S\circ D)[j].g.s}{(S\circ D)[j].g.t}).e && [\defref{defslice}]\\
&\quad= \lift{(((S\circ D)[j].e).e)}{(S\circ D)[j].g.t[k - (S\circ D)[j].e.h].d} \\&\qquad\qquad \circ (S\circ D)[j].g.t[k - (S\circ D)[j].e.h].e && [\myeqref{defrewritegenerator}]\\
&\quad= \lift{((\inc{S}{D[j].d}\circ D[j].e).e)}{\\&\qquad\qquad (S\circ D)[j].g.t[k - (S\circ D)[j].e.h].d} \\&\qquad\qquad\qquad\qquad \circ (S\circ D)[j].g.t[k - (S\circ D)[j].e.h].e && [\myeqref{diagcomposeembeddingsnmgreater}]\\
&\quad= \lift{((\inc{S}{D[j].d}\circ D[j].e).e)}{D[j].g.t[k - (S\circ D)[j].e.h].d} \\&\qquad\qquad \circ D[j].g.t[k - (S\circ D)[j].e.h].e &&[\myeqref{diagcomposegeneratorsnmgreater}]\\
&\quad= \lift{((\inc{S}{D[j].d}\circ D[j].e).e)}{D[j].g.t[k - (|S| + D[j].e.h)} \\&\qquad\qquad \circ D[j].g.t[k - (|S| + D[j].e.h)].e &&[*]\\
&\quad= \lift{(\lift{(\inc{S}{D[j].d}.e)}{(D[j].d)[D[j].e.h].d}\\&\qquad\qquad \circ D[j].e.e}{D[j].g.t[k - (|S| + D[j].e.h)} \\&\qquad\qquad\qquad\qquad \circ D[j].g.t[k - (|S| + D[j].e.h)].e &&[\defref{defcomposition}]\\
&\quad= \lift{\lift{D[j].d.\id.e}{(D[j].d)[D[j].e.h].d}\\&\qquad\qquad \circ D[j].e.e}{D[j].g.t[k - (|S| + D[j].e.h)} \\&\qquad\qquad\qquad\qquad \circ D[j].g.t[k - (|S| + D[j].e.h)].e &&[\defref{definclusion}]\\
&\quad= \lift{(\lift{D[j].d[0].d.\id}{(D[j].d)[D[j].e.h].d}\\&\qquad\qquad \circ D[j].e.e}{D[j].g.t[k - (|S| + D[j].e.h)} \\&\qquad\qquad\qquad\qquad \circ D[j].g.t[k - (|S| + D[j].e.h)].e &&[\defref{identityembedding}]\\
&\quad= \lift{({D[j].d[D[j].e.h].d.\id})\\&\qquad\qquad \circ D[j].e.e}{D[j].g.t[k - (|S| + D[j].e.h)} \\&\qquad\qquad\qquad\qquad \circ D[j].g.t[k - (|S| + D[j].e.h)].e &&[\defref{liftedidentityembedding}]\\
&\quad= \lift{D[j].e.e}{D[j].g.t[k - (|S| + D[j].e.h)} \\&\qquad\qquad \circ D[j].g.t[k - (|S| + D[j].e.h)].e &&[\defref{identitycancellation}]\\ 
&\quad= \lift{D[j].e.e}{D[j].g.t[(k- |S|) - D[j].e.h].d}\\&\qquad\qquad \circ D[j].g.t[(k- |S|) - D[j].e.h].e \\
&\quad= \rewrite{D[j].d}{D[j].e}{D[j].g.s}{D[j].g.t}[k- |S|].e \\
&\quad= D[j+1].d[k- |S|].e && [\myeqref{defrewritegenerator}]\\
&\quad=(S\circ D[j+1].d)[k].e &&[\myeqref{diagcomposegeneratorsnmequal}]
\end{align*}

\item In the range: 
\begin{align*}
&(S\circ D)[j].e.h + |(S\circ D)[j].g.t|\leq k \\
&\qquad < |(S\circ D)[j].d| -  |(S\circ D)[j].g.s| +  |(S\circ D)[j].g.t|
\end{align*}
We have: $(S\circ D)[j].e =D[j-|S|].e$, hence the range is equivalent to: 
\begin{align*}
&D[j-|S|].e.h + |D[j-|S|].g.t|\leq k  \\
&\qquad <| S| + |D[j-|S|]| -  |D[j-|S|].g.s| +  |D[j-|S|].g.t|
\end{align*}

Generators:
\begin{align*}
&((S\circ D)[j+1].d)[k].g \\
&\quad=(\rewrite{(S\circ D)[j].d)}{(S\circ D)[j].e}{(S\circ D)[j].g.s}{\\&\qquad\qquad (S\circ D)[j].g.t}).g && [\defref{defslice}] \\
&\quad= ((S\circ D)[j].d)[k - |(S\circ D)[j].g.t| \\&\qquad\qquad + |(S\circ D)[j].g.s|].g &&[\myeqref{defrewritegenerator}]\\
&\quad= (S\circ (D[j].d))[k - |(S\circ D)[j].g.t| \\&\qquad\qquad + |(S\circ D)[j].g.s|].g && [IH] \\
&\quad= (S\circ (D[j].d))[(k - |D[k].g.t| + |D[k].g.s|)].g &&[\myeqref{diagcomposegeneratorsnmgreater}, |S|\leq k]\\
&\quad= (D[j].d)[(k - |D[k].g.t| + |D[k].g.s|) - |S|].g &&[\myeqref{diagcomposegeneratorsnmequal}]\\
&\quad= (D[j].d)[(k - |S|) - |D[k].g.t| + |D[k].g.s|].g \\
&\quad= (\rewrite{D[j].d}{D[j].e}{D[j].g.s}{D[j].g.t})[k- |S|].e &&[\myeqref{defrewritegenerator}]\\
&\quad= (D[j+1].d)[(k- |S|].g &&[\defref{defslice}] \\
&\quad=(S\circ (D[j+1].d))[k].g  &&[\myeqref{diagcomposegeneratorsnmgreater}, |S|\leq k]
\end{align*}

Embeddings:
\begin{align*}
&((S\circ D)[j+1].d)[k].e \\
&\quad=(\rewrite{(S\circ D)[j].d)}{(S\circ D)[j].e}{(S\circ D)[j].g.s}{\\&\qquad\qquad (S\circ D)[j].g.t}).e && [\defref{defslice}] \\
&\quad= ((S\circ D)[j].d)[k - |(S\circ D)[j].g.t| \\&\qquad\qquad + |(S\circ D)[j].g.s|].e &&\myeqref{defrewriteembedding}]\\
&\quad= (S\circ (D[j].d))[k - |(S\circ D)[j].g.t| \\&\qquad\qquad + |(S\circ D)[j].g.s|].e && [IH] \\
&\quad= (S\circ (D[j].d))[(k - |D[k].g.t| + |D[k].g.s|)].e &&[\myeqref{diagcomposeembeddingsnmequal}, |S|\leq k]\\
&\quad= (D[j].d)[(k - |D[k].g.t| + |D[k].g.s|) - |S|].e &&[\myeqref{diagcomposeembeddingsnmequal}]\\
&\quad= (D[j].d)[(k - |S|) - |D[k].g.t| + |D[k].g.s|].e \\
&\quad= \rewrite{(D[j].d)}{D[j].e}{D[j].g.s}{D[j].g.t})[k- |S|].e &&[\myeqref{defrewriteembedding}]\\
&\quad= (D[j+1].d)[(k- |S|].g &&[\defref{defslice}] \\
&\quad=(S\circ (D[j+1].d))[k].e &&[\myeqref{diagcomposeembeddingsnmequal}, |S|\leq k]
\end{align*}
\end{itemize}
\end{itemize}
\end{itemize}
All embeddings and generators in both diagrams correspond. Hence, as all these conditions are satisfied, the two diagrams are equivalent by Definition~\ref{diagequiv}. The argument for $n<m$ is analogous. By this we established that $K(1)$ holds.
\end{proof}

\paragraph{Inductive steps.}
With all the base cases established, we prove a series of implications between the logical statements defined earlier in this section. Again, for each implication we only take the minimal subset of expressions that implies the given statement for $n$.
\begin{lemma}
\label{NLimpliesL}
For $k\geq 1$ the following holds: $N(k-1)\wedge L(k-1) \implies L(k)$.
\end{lemma}
\begin{proof}
We assume that $N(k-1)$ holds, \emph{i.e.} that for any well-defined $x$\-diagram $B$ and a well-defined $y$-diagram $A$, such that $|x-y|=k-1$ and $A.t= s^{x-y+1}(B)$, their composite $A\circ B$ is well-defined. We also assume that $L(k-1)$ holds.

Now consider two well-defined diagrams: an $n$\-diagram $D$ and an $m$\-diagram $S$, such that $S.t=s^{n-m+1}(D)$ and $|n-m|=k$. We need to show that $S\circ D$ is well-defined.  We consider two cases separately, first let $m\geq n$.  By Definition~\ref{defwellformed}, $S\circ D$ is well-defined if for $0\leq j\leq |S\circ D|$ all the slices $(S\circ D)[j].d$ are well-defined.  We prove this result by induction on $0\leq j\leq |S\circ D|$:
\begin{itemize}

\item \emph{Base case:} For $j=0$, we need to show that the source $(S\circ D)[0].d=(S\circ D).s$ is a well-defined diagram. As $m\leq n$, by the clause~\myeqref{diagcomposesourcesnmgreater} in Definition~\ref{diagcompose} we obtain the following: 
\begin{align*}
&(S\circ D).s = S\circ (D.s)
\end{align*}
The dimension of $D.s$ is $n-1$, hence we get $(n{-}1)-m= k-1$ and since $L(k-1)$ holds we obtain that $(S\circ D).s $ is well-defined as the composite of two well-defined diagrams whose difference in dimensions is $k-1$.

\item \emph{Inductive step:}  For $0 < j \leq |S\circ D|$, we assume that the slice $(S\circ D)[j].d$ exists and is well-defined. Let us consider the subsequent slice $(S\circ D)[j+1].d$, then we have the following:
\begin{align*}
&(S\circ D)[j+1].d =  \\
&\quad= \rewrite{((S\circ D)[j].d)}{(S\circ D)[j].e}{(S\circ D)[j].g.s}{(S\circ D)[j].g.t} &&[\defref{defslice}]\\
&\quad= \rewrite{((S\circ D)[j].d)}{(S\circ D)[j].e}{D[j].g.s}{D[j].g.t} &&[\myeqref{diagcomposegeneratorsnmgreater}]\\
&\quad= \rewrite{((S\circ D)[j].d)}{\inc{S}{D[j].d}\circ D[j].e}{D[j].g.s}{D[j].g.t} &&[\myeqref{diagcomposeembeddingsnmgreater}]
\end{align*}

The following hold:
\begin{itemize}
\item $(S\circ D)[j].d$ is well-defined by the inductive hypothesis.
\item $D[j].g.s$ and $D[j].g.t$ are globular with respect to each other by Definition~\ref{defsignature}.
\item $D[j].e$ is well-defined, since $D$ is well-defined.
\item $\inc{S}{D[j].d}$ is well-defined, by application of $N(k-1)$, since the dimension of $D[j].d$ is $n-1$.
\item $\inc{S}{D[j].d}\circ D[j].e$ is well-defined as the composite of two well-defined embeddings by $C(n)$~\ref{defC} which holds by Theorem~\ref{bigtheorem}.
\end{itemize}

Hence, we apply Theorem~\ref{rewritewellfedined} to conclude that $(S\circ D)[j+1].d$ is a well-defined diagram as the rewrite of a well-defined diagram $(S\circ D)[j].d$.
\end{itemize}
By this inductive argument all slices $(S\circ D)[j].d$ are well-defined for $0\leq j\leq |S\circ D|$, hence $S\circ D$ is well-defined.

The proof for $m>n$ is analogous. This establishes that $L(k)$ holds.
\end{proof}


\begin{lemma}
\label{LMNimpliesN}
For $k\geq 1$ the following holds: $L(k)\wedge M(k)\wedge N(k-1)\implies N(k)$.
\end{lemma}
\begin{proof}
We make three assumptions:
\begin{itemize}

\item $L(k)$ holds: that is, for any well-defined $n$\-diagram $D$ and a well-defined $m$-diagram $S$ such that $|n-m|=k$ and $S.t=s^{n-m+1}(D)$ if $m\leq n$ or $ t^{m-n+1}(S)=D.s$ otherwise, then the \emph{composite} diagram $S\circ D$ is well-defined.

\item $M(k)$ holds: that is, for any well-defined $n$\-diagram $D$ and a well-defined diagram $m$\-diagram $S$ such that $|n-m|=k$ the following equalities hold:
\begin{align*}
&\text{If } n>m& \inc{S}{D[i].d}&=\lift{(\inc{S}{D}.e)}{D[i].d}&&\text{for any } 0\leq i < |D| \\
& \text{If } n<m& \increverse{S[i].d}{D}&=\lift{(\increverse{S}{D}.e)}{S[i].d}&&\text{for any } 0\leq i < |S|
\end{align*}

\item $N(k-1)$ holds.
\end{itemize}
We need to show that given an $n$\-diagram $D$ and an $m$\-diagram $S$, such that $1\leq m, n$:
\begin{itemize}
\item If $n>m$ and $S.t=s^{n-m+1}(D)$ and the composite $S\circ D$ exists, the inclusion embedding $\inc{S}{D}: D\hookrightarrow S\circ D$   is well-defined.
\item If $n<m$ and  $t^{m-n+1}(S)\equiv D.s$ and the composite $S\circ D$ exists, the inclusion embedding $\increverse{S}{D}: D\hookrightarrow D\circ S$ is well-defined.
\end{itemize}
The case for $n=m$ cannot happen as $|n-m|=k\geq 1$. We consider the two above cases separately, first let $n>m$ and consider $\inc{S}{D}$:

The domain $D$ of the embedding $\inc{S}{D}$ is well-defined by assumption. The codomain $S\circ D$ is well-defined by  $L(k)$. By Definition~\ref{defembwelldef} the inclusion embedding $\inc{S}{D}: D\hookrightarrow S\circ D$ is well-defined if the following three conditions are satisfied:

\begin{itemize}

\item The component embedding $\inc{S}{D}.e$ is well-defined. As $n>m$, by Definition~\ref{definclusion}, we have $\inc{S}{D}.e=\inc{S}{D.s}$

The dimension of $D.s$ is $n-1$, hence we get $(n{-}1)-m= k-1$ and by $N(k-1)$ we obtain that $\inc{S}{D.s}$ is well-defined as the inclusion for two well-defined diagrams whose difference in dimensions is $k-1$.

\item For every $0\leq j< |D|$ we have: 
\begin{align*}
&(S\circ D)[j+\inc{S}{D}.h].g \\
&\quad= (S\circ D)[j+D.\id.h].g &&[\defref{definclusion}]\\
&\quad= (S\circ D)[j].g &&[\defref{identityembedding}]\\
&\quad= D[j].g &&[\myeqref{diagcomposegeneratorsnmgreater}]
\end{align*}
As required.

\item For every $0\leq j< |D|$ we have:  
\begin{align*}
& \lift{(\inc{S}{D}.e)}{D[j].d}\circ D[j].e = \\
&\quad=\inc{S}{D[j].d}\circ D[j].e &&[M(k)] \\
&\quad=(S\circ D)[j].e &&[\myeqref{diagcomposeembeddingsnmgreater}]\\
&\quad=(S\circ D)[j+\inc{S}{D}.h].e &&[\defref{definclusion}]
\end{align*}
As required.

\end{itemize}
As all conditions are satisfied, we can conclude that, hence $\inc{S}{D}$ is well-defined.

The argument for $n<m$ that $\increverse{S}{D}$ is well-defined is analogous. This establishes that $N(k)$ holds and the implication is true. 
\end{proof}


\begin{lemma}
\label{KimpliesM}
For $n\geq 1$ the following holds $K(k)\implies M(k)$
\end{lemma}
\begin{proof}

We assume that $K(n)$ holds \emph{i.e.} that for any $n$\-diagram $D$ and any $m$\-diagram $S$ such that $|n-m|=k$ and such that the composite $S\circ D$ exists, the following equalities hold:
\begin{align*}
&\text{If } n>m& (S\circ D)[j].d&= S\circ (D[j].d)&&\text{for any } 0\leq j < |D| \\
& \text{If } n<m& (S\circ D)[j].d&= (S[j].d) \circ D&&\text{for any } 0\leq j < |S|
\end{align*}
We need to show that for any well-defined $n$\-diagram $D$ and a well-defined diagram $m$\-diagram $S$ such that $|n-m|=k$, the following equalities hold:
\begin{align*}
&\text{If } n>m& \inc{S}{D[j].d}&=\lift{(\inc{S}{D}.e)}{D[j].d}&&\text{for any } 0\leq j < |D| \\
&\text{If } n<m& \increverse{S[j].d}{D}&=\lift{(\increverse{S}{D}.e)}{S[j].d}&&\text{for any } 0\leq j < |S|
\end{align*}
We consider both cases above separately, first let $n>m$. As $n\geq 1$, to show that these two embeddings are equivalent, by Definition~\ref{embeddingequiv}, we need to check three conditions:

\begin{itemize}

\item Domains and codomains are equivalent diagrams:

\begin{itemize}

\item By Definitions~\ref{liftedembedding},~\ref{definclusion} the type of $\inc{S}{D[i].d}$ is as follows:
\begin{align*}
\inc{S}{D[i].d}: D[i].d\hookrightarrow S\circ (D[i].d)
\end{align*}

\item By Definitions~\ref{liftedembedding},~\ref{definclusion} the type of $\lift{(\inc{S}{D}.e)}{D[i].d}$ is as follows:
\begin{align*}
&\lift{(\inc{S}{D}.e)}{D[i].d}: \\&\qquad D[i].d\hookrightarrow \rewrite{((S\circ D)[\inc{S}{D}.h])}{\inc{S}{D}.e}{D.s}{D[i].d}
\end{align*}
We see that the domains of both are immediately equivalent. For codomains, we need to simplify first. In this derivation we make use of statement $S(n)$ defined in Definition~\ref{defS}, which holds by Theorem~\ref{bigtheorem}.
\begin{align*}
&\rewrite{((S\circ D)[\inc{S}{D}.h])}{\inc{S}{D}.e}{D.s}{D[i].d}= \\
&\quad=\rewrite{((S\circ D)[\inc{S}{D}.h])}{\inc{S}{D}.e}{D.s}{\\&\qquad\qquad D[i - \inc{S}{D}.h].d} &&[\defref{definclusion}]\\
&\quad= \rewrite{(S\circ D)}{(\inc{S}{D})}{D}{D}[i].d &&[S(n)]\\
&\quad= (S\circ D)[i].d &&[\text{Lemma }~\ref{vacuousrewrite}]
\end{align*}
Need $S\circ D$ well-defined for this.

Since $K(n)$ holds, we now obtain that the codomains are equivalent diagrams.

\end{itemize}

\item Component embeddings are equivalent. In this derivation we make use of statement $T(n)$ defined in Definition~\ref{defT}, which holds by Theorem~\ref{bigtheorem}.
\begin{align*}
&(\lift{(\inc{S}{D}.e)}{D[i].d}).e= \\
&\quad= (\inc{S}{D}.e).e &&[\defref{liftedembedding}] \\
&\quad= (\inc{S}{D.s}).e &&[\defref{definclusion}] \\
&\quad= \inc{S}{D.s.s}&&[\defref{definclusion}] \\
&\quad= \inc{S}{(D[i].d).s}&&[T(n)] \\
&\quad=(\inc{S}{D[i].d}).e&&[\defref{definclusion}] 
\end{align*}

\item Heights are equal:
\begin{align*}
&(\lift{(\inc{S}{D}.e)}{D[i].d}).h = \\
&\quad= (\inc{S}{D}.e).h &&[\defref{liftedembedding}] \\
&\quad= (\inc{S}{D.s}).h &&[\defref{definclusion}] \\
&\quad= D.s.\id.h &&[\defref{definclusion}] \\
&\quad=0\\
&\quad= {\id}_{D[i].d}.h &&[\defref{definclusion}] \\
&\quad=(\inc{S}{D[i].d}).h&&[\defref{definclusion}]
\end{align*}
\end{itemize}
As all these conditions are fulfilled, the two embeddings are equivalent. The argument for $n<m$ is analogous. This establishes that $M(k)$ holds and the implication is true.
\end{proof}

\begin{lemma}
\label{MimpliesK}
For $n\geq 1$ the following holds: $M(k-1) \implies K(k)$
\end{lemma}
\begin{proof}
Assume that $M(k-1)$ holds, \emph{i.e.} for any well-defined $n$\-diagram $D$ and a well-defined diagram $m$\-diagram $S$ such that $|n-m|=k$ the following equalities hold:
\begin{align*}
&\text{If } n>m& \inc{S}{D[j].d}&=\lift{(\inc{S}{D}.e)}{D[j].d}&&\text{for any } 0\leq j < |D| \\
&\text{If } n<m& \increverse{S[j].d}{D}&=\lift{(\increverse{S}{D}.e)}{S[j].d}&&\text{for any } 0\leq j < |S|
\end{align*}
We need to show that for any $n$\-diagram $D$ and $m$\-diagram $S$ such that the composite $S\circ D$ exists and $|n-m|=k$ and such that the composite $S\circ D$ exists, the following equalities hold:
\begin{align*}
&\text{If } n>m& (S\circ D)[j].d&= S\circ (D[j].d)&&\text{for any } 0\leq j < |D| \\
& \text{If } n<m& (S\circ D)[j].d&= (S[j].d) \circ D&&\text{for any } 0\leq j < |S|
\end{align*}
We consider both cases above separately, first let $n>m$. We prove this result by induction on $0\leq j \leq |D|$.

\paragraph{Base case.} For $j=0$, the result follows immediately from the definitions:
\begin{align*}
&(S\circ D)[0].d \\
&\quad=(S\circ D).s && [\defref{definitialterminalslice}]\\
&\quad= S\circ(D.s) && [\defref{diagcompose}] \\
&\quad= S\circ(D[0].d) &&[\defref{definitialterminalslice}]
\end{align*}

\paragraph{Inductive step.} For $j>0$, let us assume the following inductive hypothesis:
\begin{align*}
(S\circ D)[j].d = S\circ (D[j].d) &\qquad (\text{I.H.})
\end{align*}
To prove that $(S\circ D)[j+1].d=S\circ (D[j+1].d)$, by Definition~\ref{diagequiv}, we need to check the following:

\begin{itemize}

\item Sources are equivalent diagrams:
\begin{align*}
&((S\circ D)[j+1].d).s \\
&\quad=  (\rewrite{(S\circ D)[j].d)}{(S\circ D)[j].e}{(S\circ D)[j].g.s}{(S\circ D)[j].g.t}).s && [\defref{defslice}]\\
&\quad=((S\circ D)[j].d).s && [\myeqref{defrewritesource}]\\
&\quad=(S\circ (D[j].d)).s && [\text{\emph{IH}}] \\
&\quad=S\circ (D[j].d).s &&[\myeqref{diagcomposesourcesnmgreater}] \\
&\quad=S\circ (\rewrite{(D[j].d)}{D[j].e}{D[j].g.s}{D[j].g.t}).s && [\myeqref{defrewritesource}]\\
&\quad=S\circ (D[j+1].d).s && [\defref{defslice}]\\
&\quad=(S\circ (D[j+1].d)).s  &&[\myeqref{diagcomposesourcesnmgreater}] 
\end{align*}

\item Lengths of generator lists are equal:
\begin{align*}
&|(S\circ D)[j+1].d| \\
&\quad=  |\rewrite{(S\circ D)[j].d)}{(S\circ D)[j].e}{(S\circ D)[j].g.s}{(S\circ D)[j].g.t}| && [\defref{defslice}]\\
&\quad= |S\circ D)[j].d| - |(S\circ D)[j].g.s| + |(S\circ D)[j].g.t| && [\myeqref{defrewritesize}]\\
&\quad= |S\circ (D[j].d)| - |(S\circ D)[j].g.s| + |(S\circ D)[j].g.t|  &&[IH] \\
&\quad= |S\circ (D[j].d)| - |D[j].g.s| + |D[j].g.t| &&[\myeqref{diagcomposegeneratorsnmgreater}]\\
&\quad= |D[j].d| - |D[j].g.s| + |D[j].g.t|  &&[\myeqref{diagcomposesizesnmgreater}]\\
&\quad=|D[j+1].d| && [\myeqref{defrewritesize}]\\
&\quad=|S\circ (D[j+1].d)| &&[\myeqref{diagcomposesizesnmgreater}]
\end{align*}
\end{itemize}

\item For generators and embeddings we need to show that for $0\leq k < |(S\circ D)[j+1].d| = |S\circ (D[j+1].d)|$, the $k$th generators in generator lists of both diagrams correspond and the same for $k$th embeddings. Firstly, we distinguish between cases for $m=n-1$ and $m<n-1$.

Since $m=n-1$, we can simplify the height of the $j$\-th embedding in $S\circ D$ in the following way. 
\begin{align*}
&(S\circ D)[j].e.h \\
&\quad= (\inc{S}{D[j].d}\circ D[j].e).h &&[\myeqref{diagcomposeembeddingsnmequal}]\\
&\quad= (\inc{S}{D[j].d}).h+ D[j].e.h &&[\defref{defcomposition}] \\
&\quad= D.\id.h + D[j].e.h &&[\defref{definclusion}] \\
&\quad= D[j].e.h &&[\defref{identityembedding}] 
\end{align*}
Let us refer to this equality as $[*]$.

We consider these generators and embeddings in three separate ranges. The first range we consider is4$0\leq k < (S\circ D)[j].e.h$, which by $[*]$ is equivalent to $0\leq k < D[j].e.h$. We show that the generators are the same:
\begin{align*}
&((S\circ D)[j+1].d)[k].g \\
&\quad=(\rewrite{(S\circ D)[j].d)}{(S\circ D)[j].e}{(S\circ D)[j].g.s}{(S\circ D)[j].g.t})[k].g && [\defref{defslice}] \\
&\quad=((S\circ D)[j].d)[k].g &&[\myeqref{defrewritegenerator}]\\
&\quad=(S\circ (D[j].d))[k].g && [IH]\\
&\quad= (D[j].d)[k].g &&[\myeqref{diagcomposegeneratorsnmgreater}]\\
&\quad= \rewrite{(D[j].d)}{D[j].d}{D[j].g.s}{D[j].g.t})[k].g &&[\myeqref{defrewritegenerator}]\\
&\quad= (D[j+1])[k].g && [\defref{defslice}] \\
&\quad=(S\circ (D[j+1].d))[k].g &&[\myeqref{diagcomposegeneratorsnmgreater}]
\end{align*}

Embeddings:
\begin{align*}
&((S\circ D)[j+1].d)[k].e \\
&\quad=(\rewrite{(S\circ D)[j].d)}{(S\circ D)[j].e}{(S\circ D)[j].g.s}{(S\circ D)[j].g.t})[k].e && [\defref{defslice}] \\
&\quad=((S\circ D)[j].d)[k].e &&[\myeqref{defrewriteembedding}]\\
&\quad=(S\circ (D[j].d))[k].e && [IH]\\
&\quad= \inc{S}{(D[j].d)[k].d}\circ(D[j].d)[k].e &&[\myeqref{diagcomposeembeddingsnmgreater}]\\
&\quad= \inc{S}{(\rewrite{D[j].d}{D[j].e}{D[j].g.s}{D[j].g.t})[k].d}\circ(D[j].d)[k].e &&[S(n{-}1)]\\
&\quad= \inc{S}{(\rewrite{D[j].d}{D[j].e}{D[j].g.s}{D[j].g.t})[k].d}\\&\qquad\qquad\circ((\rewrite{D[j].d}{D[j].e}{D[j].g.s}{D[j].g.t})[k].d)[k].e &&[\myeqref{defrewriteembedding}]\\
&\quad= \inc{S}{(D[j+1].d)[k].d}\circ(D[j+1].d)[k].e && [\defref{defslice}] \\
&\quad=(S\circ (D[j+1].d))[k].e &&[\myeqref{diagcomposeembeddingsnmgreater}]
\end{align*}

\item Range: 
\begin{align*}
(S\circ D)[j].e.h\leq k < (S\circ D)[j].e.h + |(S\circ D)[j].g.t|
\end{align*} 
By $[*]$ this is equivalent to 
\begin{align*}
D[j].e.h \leq k < D[j].e.h + |D[j].g.t|
\end{align*}
Generators:
\begin{align*}
&((S\circ D)[j+1].d)[k].g \\
&\quad=(\rewrite{(S\circ D)[j].d)}{(S\circ D)[j].e}{(S\circ D)[j].g.s}{(S\circ D)[j].g.t})[k].g && [\defref{defslice}] \\
&\quad= (S\circ D)[j].g.t[k - (S\circ D)[j].e.h].g &&[\myeqref{defrewritegenerator}]\\
&\quad= D[j].g.t[k - (S\circ D)[j].e.h].g &&[\myeqref{diagcomposegeneratorsnmgreater}]\\
&\quad= D[j].g.t[(k - D[j].e.h].g &&[*]\\
&\quad= \rewrite{(D[j].d)}{D[j].d}{D[j].g.s}{D[j].g.t})[k].g &&[\myeqref{defrewritegenerator}]\\
&\quad= (D[j+1])[k].g && [\defref{defslice}] \\
&\quad=(S\circ (D[j+1].d))[k].g &&[\myeqref{diagcomposegeneratorsnmgreater}]
\end{align*}
Embeddings: In this derivation we make use of statement $Q(n)$ described in Definition~\ref{defQ}, whose correctness is proved in Theorem~\ref{bigtheorem}. 
\begin{align*}
\lift{(\lift{f}{A}\circ e)}{T} = \lift{f}{\rewrite{A}{e}{S}{T}}\circ\lift{e}{T} &\qquad (Q(n))
\end{align*}

Here we instantiate $Q(n)$ for use in this particular application:
\begin{align*}
&f=\inc{S}{D[j].d}.e &&e=D[j].e.e\\
&A=(D[j].d)[D[j].e.h].d &&T= D[j].g.t[k - D[j].e.h].d 
\end{align*}

Then, the following equality is given by $[Q(n{-}1)]$:
\begin{align*} 
&\lift{((\lift{\inc{S}{D[j].d}.e)}{(D[j].d)[D[j].e.h].d}\circ D[j].e.e)}{\\&\qquad\qquad D[j].g.t[k - D[j].e.h].d} \\
&= \lift{(\inc{S}{D[j].d}.e)}{(\rewrite{(D[j].d)[D[j].e.h]}{D[j].e.e}{X}{\\&\qquad\qquad D[j].g.t[k-D[j].e.h].d})}\circ \lift{(D[j].e.e)}{D[j].g.t[k - D[j].e.h].d} 
\end{align*}
We also make use of the statement $S(n)$ defined in Definition~\ref{defS} and which holds by Theorem~\ref{bigtheorem}. 
\begin{align*}
&((S\circ D)[j+1].d)[k].e \\
&\quad=(\rewrite{(S\circ D)[j].d)}{(S\circ D)[j].e}{(S\circ D)[j].g.s}{(S\circ D)[j].g.t})[k].e && [\defref{defslice}] \\
&\quad= \lift{(((S\circ D)[j].e).e)}{(S\circ D)[j].g.t[k - (S\circ D)[j].e.h].d}\circ \\&\qquad\qquad (S\circ D)[j].g.t[k - (S\circ D)[j].e.h].e &&[\myeqref{defrewriteembedding}]\\
&\quad= \lift{(((S\circ D)[j].e).e)}{D[j].g.t[k - (S\circ D)[j].e.h].d}\circ\\&\qquad\qquad  D[j].g.t[k - (S\circ D)[j].e.h].e  &&[\myeqref{diagcomposegeneratorsnmgreater}]\\
&\quad= \lift{(((S\circ D)[j].e).e)}{D[j].g.t[k - D[j].e.h].d}\circ \\&\qquad\qquad  D[j].g.t[k - D[j].e.h].e &&[*]\\
&\quad= \lift{((\inc{S}{D[j].d}\circ D[j].e).e)}{D[j].g.t[k - D[j].e.h].d}\circ\\&\qquad\qquad  D[j].g.t[k - D[j].e.h].e &&[\myeqref{diagcomposeembeddingsnmgreater}]\\
&\quad= \lift{((\lift{\inc{S}{D[j].d}.e)}{(D[j].d)[D[j].e.h].d}\circ D[j].e.e)}{\\&\qquad\qquad D[j].g.t[k - D[j].e.h].d}\circ D[j].g.t[k - D[j].e.h].e && [\myeqref{defcomposition}]]
\\
&\quad= \lift{(\inc{S}{D[j].d}.e)}{(\rewrite{(D[j].d)[D[j].e.h]}{D[j].e.e}{X}{\\&\qquad\qquad D[j].g.t[k-D[j].e.h].d})}\circ \lift{(D[j].e.e)}{D[j].g.t[k - D[j].e.h].d}\\&\qquad\qquad\qquad\qquad\circ D[j].g.t[k - D[j].e.h].e &&[Q(n{-}1)]\\
&\quad= \lift{(\inc{S}{D[j].d}.e)}{(\rewrite{(D[j].d)}{D[j].e}{D[j].g.s}{\\&\qquad\qquad D[j].g.t})[k].d}\circ \lift{(D[j].e.e)}{D[j].g.t[k - D[j].e.h].d}\\&\qquad\qquad\qquad\qquad\circ D[j].g.t[k - D[j].e.h].e &&[S(n{-}1)]\\
&\quad= \lift{(\inc{S}{(D[j].d).s})}{(\rewrite{(D[j].d)}{D[j].e}{D[j].g.s}{\\&\qquad\qquad D[j].g.t})[k].d}\circ \lift{(D[j].e.e)}{D[j].g.t[k - D[j].e.h].d}\\&\qquad\qquad\qquad\qquad\circ D[j].g.t[k - D[j].e.h].e &&[\defref{definclusion}]\\
&\quad= \lift{(\inc{S}{\rewrite{((D[j].d)}{D[j].e}{D[j].g.s}{D[j].g.t}).s})}{(\rewrite{(D[j].d)}{D[j].e}{D[j].g.s}{\\&\qquad\qquad D[j].g.t})[k].d}\circ \lift{(D[j].e.e)}{D[j].g.t[k - D[j].e.h].d}\\&\qquad\qquad\qquad\qquad\circ D[j].g.t[k - D[j].e.h].e &&[\myeqref{defrewritesource}]\\
&\quad= \lift{(\inc{S}{\rewrite{(D[j].d)}{D[j].e}{D[j].g.s}{D[j].g.t}}.e)}{(\rewrite{(D[j].d)}{D[j].e}{D[j].g.s}{\\&\qquad\qquad D[j].g.t})[k].d}\circ \lift{(D[j].e.e)}{D[j].g.t[k - D[j].e.h].d}\\&\qquad\qquad\qquad\qquad\circ D[j].g.t[k - D[j].e.h].e &&[\defref{definclusion}]\\
&\quad= \inc{S}{(\rewrite{(D[j].d)}{D[j].e}{D[j].g.s}{D[j].g.t})[k].d}\\&\qquad\qquad\circ \lift{(D[j].e.e)}{D[j].g.t[k - D[j].e.h].d}\circ D[j].g.t[k - D[j].e.h].e &&[M(k-1)]\\
&\quad= \inc{S}{(\rewrite{(D[j].d)}{D[j].e}{D[j].g.s}{D[j].g.t})[k].d}\\&\qquad\qquad\circ((\rewrite{(D[j].d)}{D[j].e}{D[j].g.s}{D[j].g.t})[k].d)[k].e &&[\myeqref{defrewriteembedding}]\\
&\quad= \inc{S}{(D[j+1].d)[k].d}\circ(D[j+1].d)[k].e && [\defref{defslice}] \\
&\quad=(S\circ (D[j+1].d))[k].e &&[\myeqref{diagcomposeembeddingsnmgreater}]
\end{align*}

\item Range: 
\begin{align*}
&(S\circ D)[j].e.h + |(S\circ D)[j].g.t|\leq k < 
\\&\qquad |(S\circ D)[j].d| -  |(S\circ D)[j].g.s| +  |(S\circ D)[j].g.t|
\end{align*}
By $[*]$ this is equivalent to: 
\begin{align*}
D[j].e.h + |D[j].g.t|\leq k < |D[j].d| -  |D[j].g.s| +  |D[j].g.t| 
\end{align*}

Generators:
\begin{align*}
&((S\circ D)[j+1].d)[k].g \\
&\quad=(\rewrite{(S\circ D)[j].d)}{(S\circ D)[j].e}{(S\circ D)[j].g.s}{\\&\qquad\qquad (S\circ D)[j].g.t})[k].g && [\defref{defslice}] \\
&\quad= ((S\circ D)[j].d)[k - |(S\circ D)[j].g.t| + |(S\circ D)[j].g.s|].g &&[\myeqref{defrewritegenerator}]\\
&\quad= (S\circ (D[j].d))[k - |(S\circ D)[j].g.t| + |(S\circ D)[j].g.s|].g && [IH] \\
&\quad= (S\circ (D[j].d))[(k - |D[k].g.t| + |D[k].g.s|)].g &&[\myeqref{diagcomposegeneratorsnmgreater}]\\
&\quad= (D[j].d)[(k - |D[k].g.t| + |D[k].g.s|)].g &&[\myeqref{diagcomposegeneratorsnmgreater}]\\
&\quad= (D[j].d)[k - |D[k].g.t| + |D[k].g.s|].g \\
&\quad= (\rewrite{D[j].d}{D[j].d}{D[j].g.s}{D[j].g.t})[k].g &&[\defref{defrewritegenerator}]\\
&\quad= (D[j+1].d)[(k- |S|].g && [\defref{defslice}] \\
&\quad=(S\circ (D[j+1].d))[k].g &&[\myeqref{diagcomposegeneratorsnmgreater}]
\end{align*}

Embeddings:
\begin{align*}
&((S\circ D)[j+1].d)[k].e \\
&\quad=(\rewrite{(S\circ D)[j].d)}{(S\circ D)[j].e}{(S\circ D)[j].g.s}{\\&\qquad\qquad (S\circ D)[j].g.t})[k].e && [\defref{defslice}] \\
&\quad= ((S\circ D)[j].d)[k - |(S\circ D)[j].g.t| + |(S\circ D)[j].g.s|].e &&[\myeqref{defrewriteembedding}]\\
&\quad= (S\circ (D[j].d))[k - |(S\circ D)[j].g.t| + |(S\circ D)[j].g.s|].e && [IH] \\
&\quad= (S\circ (D[j].d))[(k - |D[k].g.t| + |D[k].g.s|)].e &&[\myeqref{diagcomposegeneratorsnmgreater}]\\
&\quad= \inc{S}{(D[j].d)[k + |D[j].g.s| - |D[j].g.t|]}\\&\qquad\qquad\circ(D[j].d)[k - |D[k].g.t| + |D[k].g.s|)].e &&[\myeqref{diagcomposeembeddingsnmgreater}]\\
&\quad= \inc{S}{(\rewrite{(D[j].d)}{D[j].e}{D[j].g.s}{D[j].g.t})[k].d}\\&\qquad\qquad\circ(D[j].d)[k - |D[k].g.t| + |D[k].g.s|].e &&[S(n{-}1)]\\
&\quad= \inc{S}{(\rewrite{(D[j].d)}{D[j].e}{D[j].g.s}{D[j].g.t})[k].d}\\&\qquad\qquad\circ((\rewrite{(D[j].d)}{D[j].e}{D[j].g.s}{D[j].g.t})[k].d)[k].e &&[\myeqref{defrewriteembedding}]\\
&\quad= \inc{S}{(D[j+1].d)[k].d}\circ(D[j+1].d)[k].e && [\defref{defslice}] \\
&\quad=(S\circ (D[j+1].d))[k].e &&[\myeqref{diagcomposeembeddingsnmgreater}]
\end{align*}
All embeddings and generators in both diagrams correspond. Hence, as all these conditions are satisfied, the two diagrams are equivalent by Definition~\ref{diagequiv}. The argument for $n<m$ is analogous. Altogether we established that $K(k)$ holds, hence the implication is true.
\end{proof}

Finally, we bring all these lemmas together to prove the following.
\begin{theorem}
\label{maintheoremcomposition}
For $k\geq 1$ the following logical statement holds: $K(k)\wedge L(k)\wedge M(k)\wedge N(k)$
\end{theorem}
\begin{proof}
We prove this by induction on $k$

\begin{itemize}
\item \emph{Base case}: For $k=1$:

\begin{itemize}
\item $K(1)$, holds with no further conditions by Lemma~\ref{recursionbaseK}.

\item To establish $M(1)$, by Lemma~\ref{KimpliesM}, we need $K(1)$. This holds by the argument above.

\item To establish $L(1)$, by Lemma~\ref{NLimpliesL}, we need $N(0)$ and L(0) to hold.  $L(0)$ holds by the argument above. $N(0)$ holds by Lemma~\ref{recursionbaseN}, since $L(0)$ holds.

\item To establish $N(1)$, by Lemma~\ref{LMNimpliesN}, we need all $L(1)$, $M(1)$, $N(0)$ to hold. All hold by the argument above.

\end{itemize}

\item \emph{Inductive step:} For $k>1$, we assume that all $K(k-1)$, $L(k-1)$, $M(k-1)$, $N(k-1)$ hold.

\begin{itemize}
\item To establish $K(k)$, by Lemma~\ref{MimpliesK}, we need $M(k-1)$. This holds by the inductive hypothesis.

\item To establish $M(k)$, by Lemma~\ref{KimpliesM}, we need $K(k)$. This holds by the argument above.

\item To establish $L(k)$, by Lemma~\ref{NLimpliesL}, we need $N(k-1)$ and L(k-1) to hold.  Both hold by the inductive hypothesis.

\item To establish $N(k)$, by Lemma~\ref{LMNimpliesN}, we need all $L(k)$, $M(k)$, $N(k-1)$ to hold. The initial two statements hold by the argument above. $N(k-1)$ holds by the inductive hypothesis.
\end{itemize}
As all statements $K(k)$, $L(k)$, $M(k)$, $N(k)$ hold, this establishes that their conjunction is true.
\end{itemize}

By this inductive argument the logical statement: $K(k)\wedge L(k)\wedge M(k)\wedge N(k)$ holds for $k\geq 1$.

\end{proof}

The two main results on the composite of two diagrams being well-defined and the inclusion embedding being well-defined follow immediately:
\begin{theorem}
for any well-defined $n$\-diagram $D$ and a well-defined diagram $m$\-diagram $S$ such that $|n-m|=k$:
\begin{itemize}
\item If $n\geq m$ and $S.t=s^{n-m+1}(D)$ and the composite $S\circ D$ exists, the inclusion embedding $\inc{S}{D}: D\hookrightarrow S\circ D$   is well-defined. 
\item If $n< m$ and  $t^{m-n+1}(S)\equiv D.s$ and the composite $S\circ D$ exists, the inclusion embedding $\inc{S}{D}: D\hookrightarrow D\circ S$ is well-defined. 
\end{itemize}
\end{theorem}
\begin{proof}
By Lemma~\ref{recursionbaseN} this holds for $n=m$ and by Theorem~\ref{maintheoremcomposition} for $n\neq m$.
\end{proof}

\begin{theorem}
For any well-defined $n$\-diagram $D$ and a well-defined $m$-diagram $S$ such that $|n-m|=k$ and $S.t=s^{n-m+1}(D)$ if $m\leq n$ or $t^{m-n+1}(S)=D.s$ otherwise, then the \emph{composite} diagram $S\circ D$ is well-defined
\end{theorem}
\begin{proof}
By Lemma~\ref{recursionbaseL} this holds for $n=m$ and by Theorem~\ref{maintheoremcomposition} for $n\neq m$.
\end{proof}

\subsection{Identity diagrams}
\label{sec:identitydiagrams}

\begin{definition}
\label{defidiagid}
Given an $n$\-diagram $D$, the \emph{identity diagram} $\diagid{D}$ on $D$ is the $(n{+}1)$\-diagram defined as follows:
\begin{align*}
\diagid{D}.s &= D \\
|\diagid{D}| &= 0 \\
\end{align*}
\end{definition}

\begin{lemma}
Given a well-defined $n$\-diagram $D$ the identity $\diagid{D}$ on the diagram $D$ is well-defined.
\end{lemma}
\begin{proof}
Since $\diagid{D}.s= \diagid{D}[0].d$ it is the only and final slice, there is no signature element and no embedding associated with it.  $\diagid{D}.s=D$ is well-defined as $D$ is well-defined, hence by Definition~\ref{defwellformed} $\diagid{D}$ is also well-defined. 
\end{proof}
We also refer to this operation as \emph{boosting} a diagram $D$. As expected, composing an $n$\-diagram $D$ with an identity on an $m$\-diagram $S$ such that $n>m$ leaves $D$ unaltered. This is because of the requirements on matches between sources and targets of the diagrams being composed.
\begin{lemma}
\label{compositionidempotent}
Given a well-defined $n$\-diagram $D$ and a well-defined $m$\-diagram $S$ such that $m<n$, the following holds:
\begin{itemize}
\item If $S = \diagid{S}.t =s^{n-m+1}(D)$: 
\begin{align*}
\diagid{S}\circ D = D
\end{align*}

\item If $t^{n-m+1}(D)= \diagid{S}.s  = S$: 
\begin{align*}
D\circ\diagid{S}= D
\end{align*}

\end{itemize}

\end{lemma}
\begin{proof}
First let us assume that $S=s^{n-m+1}(D)$. We show that $\diagid{S}\circ D$ and $D$ are equivalent diagrams by induction on $k=(n{-}1)-m$. 

\begin{itemize}

\item \emph{Base case}: For $k=0$, we have $n-1=m$. By Definition~\ref{diagequiv}, we need to show four separate conditions:

\begin{itemize}
\item Sources are equivalent diagrams:
\begin{align*}
&(\diagid{S}\circ D).s \\
&\quad= \diagid{S}.s &[\myeqref{diagcomposesourcesnmequal}]\\
&\quad= S &[\defref{defidiagid}]\\
&\quad= \diagid{S}.t &[\defref{defidiagid}]\\
&\quad= D.s &[\text{assumption}]
\end{align*}

\item Sizes of generator lists are equal:
\begin{align*}
&|\diagid{S}\circ D| \\
&\quad= |\diagid{S}| + |D| &[\myeqref{diagcomposesizesnmequal}]\\
&\quad= |D| &[\defref{defidiagid}]
\end{align*}

\item Corresponding generators are equal, for $0\leq j\leq |D|$:
\begin{align*}
&(\diagid{S}\circ D)[j].g \\
&\quad= D[j].g &[\myeqref{diagcomposegeneratorsnmequal}, |\diagid{S}|=0]
\end{align*}

\item Corresponding embeddings are equivalent, for $0\leq j\leq |D|$:
\begin{align*}
&(\diagid{S}\circ D)[j].e \\
&\quad= D[j].e &[\myeqref{diagcomposeembeddingsnmequal}, |\diagid{S}|=0]
\end{align*}
\end{itemize}

\item \emph{Inductive step}: For $k>0$, we assume that the result holds, \emph{i.e.} for all $x$\-diagrams $M$ and $y$\-diagrams $M$ such that $k=(x-1)-y$, we have $\diagid{N}\circ M = M$ (\emph{IH}).

Now consider an $n$\-diagram $D$ and an $m$\-diagram $S$, such that $k+1=(n{-}1)-m$, then for $\diagid{S}\circ D = D$ to hold, by Definition~\ref{diagequiv}, we need to show four separate conditions:

\begin{itemize}
\item Sources are equivalent diagrams:
\begin{align*}
&(\diagid{S}\circ D).s \\
&\quad= \diagid{S} \circ D.s&&[\myeqref{diagcomposesourcesnmgreater}]\\
&\quad=  D.s &&[\text{\emph{IH}}]
\end{align*}

\item Sizes of generator lists are equal:
\begin{align*}
&|\diagid{S}\circ D| \\
&\quad= |D|&[\myeqref{diagcomposesizesnmgreater}]
\end{align*}

\item Corresponding generators are equal, for $0\leq j\leq |D|$:
\begin{align*}
&(\diagid{S}\circ D)[j].g \\
&\quad= D[j].g &[\myeqref{diagcomposegeneratorsnmgreater}]
\end{align*}

\item Corresponding embeddings are equivalent, for $0\leq j\leq |D|$:
\begin{align*}
&(\diagid{S}\circ D)[j].e \\
&\quad= \inc{\diagid{S}}{D[j].d}\circ D[j].e &&[\myeqref{diagcomposeembeddingsnmgreater}]\\
&\quad= D[j].e &&[\text{\emph{IH}}]
\end{align*}
\end{itemize}
\end{itemize}
By this inductive argument, we established that for any $n$\-diagram $D$ and any $m$\-diagram $S$ such that $m<n$, such that $S=s^{n-m+1}(D)$ the following holds: $\diagid{S}\circ D = D$.

The argument for $S=t^{n-m+1}((D)$ is analogous, so the entire result holds.
\end{proof}

However, if the diagram $S$ that the identity operation acts on is of the same dimension $n$ as $D$, we get slightly different behaviour.
\begin{lemma}
\label{boostingcomposition}
Given a well-defined $n$\-diagram $D$ and a well-defined diagram $S$ such that $m,n>0$, the following holds:

\begin{itemize}
\item If $S=s^{n-m+1}(D)$: 
\begin{align*}
S\circ \diagid{D} = \diagid{S\circ D}
\end{align*}

\item If $t^{m-n+1}(S)=D$: 
\begin{align*}
\diagid{S}\circ D = \diagid{S\circ D}
\end{align*}

\end{itemize}

\end{lemma}
\begin{proof}

We prove the result for both cases separately. If $S=s^{n-m+1}(D)$, we have $n\geq m$.

Since $n,m >0$ to show that these two diagrams are equivalent, by Definition~\ref{diagequiv}, we need to check four separate conditions:

\begin{itemize}
\item Sources are equivalent diagrams:\JVcomm{Don't follow this}
\begin{align*}
&(S\circ \diagid{D} = \diagid{S\circ D}).s \\
&\quad= S \circ \diagid{D}.s &&[\myeqref{diagcomposesourcesnmsmaller}]\\
&\quad=  S\circ D &&[\defref{defidiagid}]\\
&\quad= \diagid{S\circ D}.s &&[\defref{defidiagid}]
\end{align*}

\item Sizes of generator lists are equal:
\begin{align*}
&|\diagid{S}\circ D| \\
&\quad= |\diagid{S}|&&[\myeqref{diagcomposesizesnmsmaller}]\\
&\quad= 0 &&[\defref{defidiagid}]\\
&\quad= |\diagid{S\circ D}| &&[\defref{defidiagid}]
\end{align*}

\item Since $|\diagid{S}\circ D| = 0$, we do not need to show anything further for generators and embeddings and the remaining two conditions are vacuously true.

\end{itemize}
This establishes that the two diagrams are equivalent, as required. The argument for $D.t^{n-m+1}=S$ is analogous.

\end{proof}

\subsection{Associativity and distributivity}
\label{secassocdist}

With the goal of modelling quasistrict $n$\-categories in mind and taking into account that the only non\-trivial morphisms that we want to keep are the interchange law and coherences derived from it, we do not want to include associator morphisms. For this reason, we need to show that certain associativity and distributivity results are built-in properties of diagram composition.

Let us consider different possible ways in which three diagrams can be composed. For any three well-defined diagrams: an $n$\-diagram $D$, $m$\-diagram $S$ and an $l$\-diagram $M$ the form of the composite depends on the order of binary compositions (bracketing) and the ordering of natural numbers $m,n,l$. Certain combinations allow for associativity or distributivity rules, others do not yield any interesting behaviour. Before we proceed to listing these formally, we give several examples to illustrate how associativity or distributivity of diagram composition arises.

Firstly, consider the following 1\-diagrams $S$, $D$, and a 2\-diagram $M$, note that $S$ and $D$ are of the same dimension:
\begin{align*}
S\quad=\quad
&
\begin{aligned}
\begin{tikzpicture}[thick, scale = 0.7]
\draw (-3.5,0) to (-3.5, 3);
\end{tikzpicture}
\end{aligned}
&
D\quad=\quad
&
\begin{aligned}
\begin{tikzpicture}[thick, scale = 0.7]
\draw (-1.5,0) to (-1.5, 3);
\draw (-2.5,0) to (-2.5, 3);
\end{tikzpicture}
\end{aligned}
&
M\quad=\quad
&
\begin{aligned}
\begin{tikzpicture}[thick, scale = 0.7]
\draw (-0.5,0) to (-0.5,0.5) to [out = up, in = left] (0,1)
to [out = right, in = up](0.5,0.5) to(0.5,0);
\draw (-0.5,3) to (-0.5,2.5) to [out = down, in = left] (0,2)
to [out = right, in = down](0.5,2.5) to(0.5,3);
\draw (0,1) to (0, 2);
\node [blob] at (0,1) {};
\node [blob] at (0,2) {};
\end{tikzpicture}
\end{aligned}
\end{align*}
We could then form the following composites:
\begin{align*}
S\circ D\quad=\quad
\begin{aligned}
\begin{tikzpicture}[thick, scale = 0.7]
\draw (-3.5,0) to (-3.5, 3);
\draw (-1.5,0) to (-1.5, 3);
\draw (-2.5,0) to (-2.5, 3);
\end{tikzpicture}
\end{aligned}
\qquad
D\circ M\quad=\quad
\begin{aligned}
\begin{tikzpicture}[thick, scale = 0.7]
\draw (-0.5,0) to (-0.5,0.5) to [out = up, in = left] (0,1)
to [out = right, in = up](0.5,0.5) to(0.5,0);
\draw (-0.5,3) to (-0.5,2.5) to [out = down, in = left] (0,2)
to [out = right, in = down](0.5,2.5) to(0.5,3);
\draw (0,1) to (0, 2);
\node [blob] at (0,1) {};
\node [blob] at (0,2) {};
\draw (-1.5,0) to (-1.5, 3);
\draw (-2.5,0) to (-2.5, 3);
\end{tikzpicture}
\end{aligned}
\end{align*}

Now note that the order (bracketing) in which we decided to perform the binary compositions does not have any effect on the final result:
\begin{align*}
S\circ(D\circ M)=(S\circ D)\circ M =\quad
\begin{aligned}
\begin{tikzpicture}[thick, scale = 0.7]
\draw (-0.5,0) to (-0.5,0.5) to [out = up, in = left] (0,1)
to [out = right, in = up](0.5,0.5) to(0.5,0);
\draw (-0.5,3) to (-0.5,2.5) to [out = down, in = left] (0,2)
to [out = right, in = down](0.5,2.5) to(0.5,3);
\draw (0,1) to (0, 2);
\node [blob] at (0,1) {};
\node [blob] at (0,2) {};
\draw (-1.5,0) to (-1.5, 3);
\draw (-2.5,0) to (-2.5, 3);
\draw (-3.5,0) to (-3.5, 3);
\end{tikzpicture}
\end{aligned}
\end{align*}

Secondly, let $D$ and $M$ have the same dimensions which are different than the dimension of $S$:
\begin{align*}
S\quad=\quad
&
\begin{aligned}
\begin{tikzpicture}[thick, scale = 0.7]
\draw (-3.5,0) to (-3.5, 3);
\end{tikzpicture}
\end{aligned}
&
D\quad=\quad
&
\begin{aligned}
\begin{tikzpicture}[thick, scale = 0.7]
\draw (-0.5,0) to (-0.5,0.5) to [out = up, in = left] (0,1)
to [out = right, in = up](0.5,0.5) to(0.5,0);
\draw (-0.5,3) to (-0.5,2.5) to [out = down, in = left] (0,2)
to [out = right, in = down](0.5,2.5) to(0.5,3);
\draw (0,1) to (0, 2);
\node [blob] at (0,1) {};
\node [blob] at (0,2) {};
\end{tikzpicture}
\end{aligned}
&
M\quad=\quad
&
\begin{aligned}
\begin{tikzpicture}[thick, scale = 0.7]
\draw (0.5,1) to(0.5,-2);
\draw (-0.5,1) to(-0.5,-2);
\node [blob] at (0.5,0) {};
\node [blob] at (-0.5,-1) {};
\end{tikzpicture}
\end{aligned}
\end{align*}

This results in different behaviour, as composition with $S$ distributes over composition of diagrams $D$ and $M$, \emph{i.e.} we could first separately compose $S$ with $D$ and with $M$ and then compose the resulting diagrams vertically, or alternatively we could first vertically compose $D$ with $M$ and then compose the result with $S$:
\begin{align*}
S\circ D\quad=\quad
\begin{aligned}
\begin{tikzpicture}[thick, scale = 0.7]
\draw (-0.5,0) to (-0.5,0.5) to [out = up, in = left] (0,1)
to [out = right, in = up](0.5,0.5) to(0.5,0);
\draw (-0.5,3) to (-0.5,2.5) to [out = down, in = left] (0,2)
to [out = right, in = down](0.5,2.5) to(0.5,3);
\draw (0,1) to (0, 2);
\node [blob] at (0,1) {};
\node [blob] at (0,2) {};
\draw (-1.5,0) to (-1.5, 3);
\end{tikzpicture}
\end{aligned}
\qquad
S\circ M\quad=\quad
\begin{aligned}
\begin{tikzpicture}[thick, scale = 0.7]
\draw (0.5,1) to(0.5,-2);
\draw (-0.5,1) to(-0.5,-2);
\draw (-1.5,1) to(-1.5,-2);
\node [blob] at (0.5,0) {};
\node [blob] at (-0.5,-1) {};
\end{tikzpicture}
\end{aligned}
\end{align*}

\begin{align*}
S\circ(D\circ M)=(S\circ D)\circ (S\circ M)\quad=\quad
\begin{aligned}
\begin{tikzpicture}[thick, scale = 0.7]
\draw (-0.5,0) to (-0.5,0.5) to [out = up, in = left] (0,1)
to [out = right, in = up](0.5,0.5) to(0.5,0);
\draw (-0.5,5) to (-0.5,2.5) to [out = down, in = left] (0,2)
to [out = right, in = down](0.5,2.5) to(0.5,5);
\draw (0,1) to (0, 2);
\draw (-1.5,0) to (-1.5, 5);
\node [blob] at (0,1) {};
\node [blob] at (0,2) {};
\node [blob] at (0.5,4) {};
\node [blob] at (-0.5,3) {};
\end{tikzpicture}
\end{aligned}
\end{align*}

At this point let us remind ourselves that each composite can either be denoted explicitly by the dimensions of the diagrams involved, such as $S\circ_{m,n}D$ or using an overloaded notation $S\circ_{\mathrm{min}(m,n) - 1}D$. Below we present a theorem that summarises all the interesting associativity and distributivity laws for composition of three diagrams. The most interesting clauses are proved in the later part of this section, through the familiar technique of making several logical statements depending on a natural number $k$ (in this instance, the difference between diagram dimensions) and then proving a conjunction of all these statements by induction on $k$.
\begin{theorem}
\label{assocdisttheorem}
Given three well-defined diagrams: an $n$\-diagram $D$, $m$\-diagram $S$ and an $l$\-diagram $M$, let $a=\mathrm{min}(n,l)-1$, $b=\mathrm{min}(m, \mathrm{max}(n,l)) - 1$, $c=\mathrm{min}(m,n) -1$ and $d=\mathrm{min}(\mathrm{max}(m,n), l) - 1$, then, provided that these composites exist, the following hold:
\begin{align}
\label{eqassoc}
S\circ_{a}(D\circ_{a} M) &= (S\circ_{a} D)\circ_{a} M
&\text{if } a=b
\\
\label{eqdistleft}
S\circ_{b}(D\circ_{a} M) &= (S\circ_{b} D)\circ_{a}(S\circ_{b} M)
&\text{if } b<a
\\
\label{eqdistright}
(S\circ_{c} D)\circ_{d} M &= (S\circ_{d} M)\circ_{c}(D\circ_{d} M)
&\text{if } d<c
\end{align}
\end{theorem}
\begin{proof}
We consider these three equalities separately:

\begin{itemize}

\item For equality~\myeqref{eqassoc}, we have $a=b$, this implies $\mathrm{min}(n,l)=\mathrm{min}(m, \mathrm{max}(n,l))$, which in turn forces one of the following three options:
\begin{itemize}

\item $n=m\leq l$

\item $n=l \leq m$ 

\item $m=l\leq n$

\end{itemize}
We prove the equality for the first of these cases in Lemma~\ref{EFhold}. The setup is prepared by Definitions~\ref{defE} and~\ref{defF}. The remaining two cases follow by an analogous argument.

\item For equality~\myeqref{eqdistleft}, we have $a>b$, this implies $\mathrm{min}(n,l)>\mathrm{min}(m, \mathrm{max}(n,l))$, which in turn forces $n,l>m$. We prove this in Lemma~\ref{GHhold}, after preparing the setup by Definitions~\ref{defG},~\ref{defH}.

\item For equality~\myeqref{eqdistright}, similarly we have $c>d$, this implies $\mathrm{min}(n,m)>\mathrm{min}(l, \mathrm{max}(n,m))$, which in turn forces $n,m>l$. The argument for this case is analogous to the proof for equality~\myeqref{eqdistleft}.
\end{itemize}
This completes the proof.
\end{proof}

Instances for $b>a$ or $d>c$ are not included, as they do not give rise to any associativity or distributivity laws. As an example, consider the composite $S\circ_{b}(D\circ_{a} M) $ such that $m>n>l$, then $a=l-1$, $b=n-1$ and we have $b>a$. Let us consider when such a composite exists. We need the following:
\begin{itemize}
\item $t^{n-l+1}(D)=M.s$
\item $t^{m-n+1}(S)=(D\circ_{a} M).s$
\end{itemize}
Since $n>l$, by equation~\ref{diagcomposesourcesnmgreater}, we obtain that: $(D\circ M).s = D.s \circ M$. By this, we see that the composite of the source of $D$ with $M$ must match the appropriate target boundary of $S$, so we cannot compose them in any other order, as the relevant sources and targets would not match.

We now prove the results stated in Theorem~\ref{assocdisttheorem}. First, we make two logical statements which, when established for all $k\geq 0$, prove that equality~\myeqref{eqassoc} holds. Apart from the main result on associativity of composition, we additionally need a statement on composition of inclusion embeddings.
\begin{definition}
\label{defE}
For $k\geq 0$, $E(k)$ states that for two well-defined $n$\-diagrams $D, S$, and a well-defined $l$\-diagram $M$ such that $l>n>0$ and $l-n=k$, the following holds:
\[
S\circ(D\circ M) = (S\circ D) \circ M
\]
\end{definition}

\begin{definition}
\label{defF}
For $k, n\geq 0$, $F(k)$ states that for two well-defined $n$\-diagrams $D, S$, and a well-defined $l$\-diagram $M$ such that $l>n$ and $l-n=k$, the following holds:
\[
\inc{S}{D\circ M}\circ \inc{D}{M} = \inc{S\circ D}{M}
\]
\end{definition}

First, in a separate lemma, we establish the base for the recursive proof for the Lemma~\ref{EFhold}, which comes later in the section.
\begin{lemma}
\label{recursionbaseE}
The statement $E(0)$ holds with no further assumptions
\end{lemma}
\begin{proof}
We need to show that given three well-defined diagrams: an $n$\-diagram $D$, and $m$\-diagram $S$ and an $l$\-diagram $M$, the following holds:
\begin{align*}
S\circ(D\circ M) = (S\circ D) \circ M
\end{align*}
As $k=0$, this gives us $n=l$. Since $n>1$
we need to check the following four conditions:
\begin{itemize}
\item Sources are equivalent diagrams, the derivation follows by~\myeqref{diagcomposesourcesnmequal}:
\begin{align*}
&(S\circ(D\circ M)).s \\
&\quad= S.s  \\
&\quad= (S\circ D).s  \\
&\quad=((S\circ D) \circ M).s 
\end{align*}

\item Sizes of generator lists are equal, the derivation follows by~\myeqref{diagcomposesizesnmequal}:
\begin{align*}
&|S\circ(D\circ M) |\\
&\quad= |S| + |D\circ M| \\
&\quad=|S| + |D| + |M| \\
&\quad=|(S\circ D)| + |M|  \\
&\quad=|(S\circ D) \circ M| 
\end{align*}

\item Generators are equal for $0\leq j < |S| + |D| + |M|$, we show this for three separate ranges. First let $0\leq j< |S|$, the derivation follows by~\myeqref{diagcomposegeneratorsnmequal}:
\begin{align*}
&(S\circ(D\circ M))[j].g\\
&\quad= S[j].g \\
&\quad=(S\circ D)[j].g \\
&\quad=((S\circ D) \circ M)[j].g
\end{align*}
The argument is analogous for the remaining two ranges $|S|\leq j < |S| + |D|$ and $|S|+|D|\leq j < |S| + |D| + |M|$.

\item Embeddings are equivalent for $0\leq j < |S| + |D| + |M|$, we show this for three separate ranges. First let us show this for $|S|\leq j < |S| + |D|$ , the derivation follows by~\myeqref{diagcomposeembeddingsnmequal}:
\begin{align*}
&(S\circ(D\circ M))[j].e\\
&\quad= (D\circ M)[j-|S|].e \\
&\quad=D[j-|S|].e\\
&\quad=((S\circ D)[j].e\\
&\quad=((S\circ D)\circ M)[j].e
\end{align*}
The argument is analogous for the remaining two ranges $0\leq j< |S|$ and $|S|+|D|\leq j < |S| + |D| + |M|$.
\end{itemize}

By this argument, the diagrams $S\circ(D\circ M)$ and $(S\circ D) \circ M$ are equivalent and the statement $E(0)$ holds, as required.
\end{proof}

The following two implications establish $E(k)$ and $F(k)$:
\begin{lemma}
\label{FEimpliesE}
For $k\geq 1$ the following statement holds: $F(k-1) \wedge E(k-1)\implies E(k)$
\end{lemma}
\begin{proof}
Let us assume that both $F(k-1)$ and $E(k-1)$ hold. We need to show that given three well-defined diagrams: an $n$\-diagram $D$, and $m$\-diagram $S$ and an $l$\-diagram $M$, the following holds:
\begin{align*}
S\circ(D\circ M) = (S\circ D) \circ M
\end{align*}
We need to show this result for all the possible orderings of $n$, $m$, $l$. First, assume $m<n<l$:

Since $n>1$, by Definition~\ref{diagequiv}, we need to check the following four conditions:

\begin{itemize}
\item Sources are equivalent diagrams, the derivation follows by~\myeqref{diagcomposesourcesnmgreater}:
\begin{align*}
&(S\circ(D\circ M)).s\\
&\quad= S\circ (D\circ M).s \\
&\quad= S\circ (D\circ M.s)  \\
&\quad= (S\circ D)\circ M.s) &[E(k-1)]\\
&\quad=((S\circ D) \circ M).s
\end{align*}

\item Sizes of generator lists are equal, the derivation follows by~\myeqref{diagcomposesizesnmgreater}:
\begin{align*}
&|S\circ(D\circ M)|\\
&\quad= |D\circ M| \\
&\quad= |M| \\
&\quad=|(S\circ D) \circ M| 
\end{align*}

\item Generators are equal for $0\leq j\leq |S\circ(D\circ M)|$, the derivation follows by~\myeqref{diagcomposegeneratorsnmgreater}:
\begin{align*}
&(S\circ(D\circ M))[j].g\\
&\quad= (D\circ M)[j].g \\
&\quad= M[j].g\\
&\quad=((S\circ D) \circ M)[j].g 
\end{align*}

\item Embeddings are equivalent for $0\leq j\leq |S\circ(D\circ M)|$, the derivation follows by~\myeqref{diagcomposeembeddingsnmgreater}:
\begin{align*}
&(S\circ(D\circ M))[j].e\\
&\quad= \inc{S}{(D\circ M)[j].d}\circ (D\circ M)[j].e \\
&\quad= \inc{S}{(D\circ M)[j].d}\circ (\inc{D}{M[j].d}\circ M[j].e)  \\
&\quad= (\inc{S}{(D\circ M)[j].d}\circ \inc{D}{M[j].d})\circ M[j].e &&[A(n), \defref{bigtheorem}] \\
&\quad= \inc{S\circ D}{M[j].d}\circ M[j].e && [F(n{-}1)] \\
&\quad=((S\circ D) \circ M)[j].e 
\end{align*}
We already showed that the statement $A(n)$ holds for all $n\geq 1$ in the proof of the Theorem~\ref{bigtheorem}, so there is no need to include it separately in the conjunction of logical statements being proved here.
\end{itemize}
As all these conditions are satisfied, we can conclude that both diagrams are equivalent, as required. Arguments for other orderings of $n$, $m$ and $l$ are analogous. By this, we established that $E(k)$ holds, hence the implication is true.
\end{proof}

\begin{lemma}
\label{FEimpliesF}
For $k\geq 1$ the following statement holds: $F(k-1) \wedge E(k)\implies F(k)$, additionally F(0) holds with no further assumptions.
\end{lemma}
\begin{proof}

Let us assume that both $F(k-1)$ and $E(k)$ hold. We need to show that given two well-defined $n$\-diagrams $D, S$, and a well-defined $l$\-diagram $M$ such that $l>n>0$ and $l-n=k$, the following holds:
\begin{align*}
\inc{S}{D\circ M}\circ \inc{D}{M} = \inc{S\circ D}{M}
\end{align*}

Since $l>n>0$ by Lemma~\ref{embeddingequiv},  we need to check the following conditions:

\begin{itemize}

\item Types are the same. By Definition~\ref{definclusion}, these types are as follows:
\begin{align*}
\inc{D}{M} &:  M\hookrightarrow D\circ M\\
\inc{S}{(D\circ M)} &: D\circ M\hookrightarrow S\circ(D\circ M)\\
\inc{S}{D\circ M}\circ \inc{D}{M}  &: M\hookrightarrow S\circ(D\circ M) \\
\inc{S\circ D}{M} &: M\hookrightarrow (S\circ D)\circ M
\end{align*}
The domains are immediately equivalent, codomains are equivalent by $E(n)$.

\item Heights are equal. For this condition, we consider two scenarios:
\begin{itemize}
\item For $k=0$:
\begin{align*}
&(\inc{S}{D\circ M}\circ \inc{D}{M}).h \\
&\quad= \inc{S}{D\circ M}.h + \inc{D}{M}.h &&[\defref{defcomposition}]\\
&\quad= |S|+ |D| &&[\defref{definclusion}]\\
&\quad= |S\circ D| &&[\myeqref{diagcomposesizesnmequal}]\\
&\quad=(\inc{S\circ D}{M}).h &&[\defref{definclusion}]
\end{align*}

\item For $k>0$:
\begin{align*}
&(\inc{S}{D\circ M}\circ \inc{D}{M}).h \\
&\quad= \inc{S}{D\circ M}.h + \inc{D}{M}.h &&[\defref{defcomposition}]\\
&\quad= (D \circ M).\id.h + M.\id.h &&[\myeqref{definclusion}]\\
&\quad=0 \\
&\quad= M.\id.h &&[\defref{definclusion}]\\
&\quad=(\inc{S\circ D}{M}).h &&[\defref{definclusion}]
\end{align*}

\end{itemize}

\item Component embeddings are equivalent. Here, again we consider two scenarios:

\begin{itemize}
\item For $k=0$:
\begin{align*}
&(\inc{S}{D\circ M}\circ \inc{D}{M} ).e \\
&\quad=\lift{(\inc{S}{D\circ M}.e)}{(D\circ M)[\inc{D}{M}.h].d} \circ  \inc{D}{M}.e &&[\defref{defcomposition}]\\
&\quad=\lift{(\inc{S}{D\circ M}.e)}{(D\circ M)[\inc{D}{M}.h].d} \circ  M.\id.e&&[\defref{definclusion}]\\
&\quad=\lift{((D \circ M).\id.e)}{(D\circ M)[\inc{D}{M}.h].d} \circ  M.\id.e&&[\defref{definclusion}]\\
&\quad={(D\circ M)[\inc{D}{M}.h].d.\id} \circ  M.\id.e&&[\defref{liftedidentityembedding}]\\
&\quad={(D\circ M)[|D|].d.\id} \circ  M.\id.e&&[\defref{definclusion}]\\
&\quad=M.s.\id \circ  M.\id.e&&[\defref{defcomposition}]\\
&\quad=M.s.\id \circ  M.s.\id&&[\defref{identityembedding}]\\
&\quad=M.s.\id&&[\defref{identitycancellation}]\\
&\quad=M.\id.e &&[\defref{identityembedding}]\\
&\quad=(\inc{S\circ D}{M}).e &&[\defref{definclusion}]
\end{align*}

\item For $k>0$:
\begin{align*}
&(\inc{S}{D\circ M}\circ \inc{D}{M} ).e \\
&\quad=\lift{(\inc{S}{D\circ M}.e)}{(D\circ M)[\inc{D}{M}.h].d} \circ  \inc{D}{M}.e &&[\defref{defcomposition}]\\
&\quad=\lift{(\inc{S}{D\circ M}.e)}{(D\circ M)[0].d} \circ  \inc{D}{M}.e &&[\defref{definclusion}]\\
&\quad=\inc{S}{(D\circ M)[0].d} \circ \inc{D}{M}.e &&[M(n)] \\
&\quad=\inc{S}{(D\circ M).s} \circ \inc{D}{M}.e &&[\defref{defcomposition}]\\
&\quad=\inc{S}{(D\circ M).s} \circ \inc{D}{M.s} &&[\defref{defcomposition}]\\
&\quad=\inc{S}{D\circ M.s} \circ \inc{D}{M.s} &&[\myeqref{diagcomposesourcesnmgreater}]\\
&\quad= \inc{S\circ D}{M.s}  &&[F(n{-}1)]\\
&\quad=(\inc{S\circ D}{M}).e &&[\defref{definclusion}]
\end{align*}
As with the proof of Lemma~\ref{FEimpliesE},  already proved that the statement $M(n)$ holds for all $n\geq 1$ by Theorem~\ref{bigtheorem}, so there is no need to include it separately in the conjunction of logical statements being proved here.
\end{itemize}
\end{itemize}
As all these conditions are satisfied both embeddings are equivalent. By this, we established that $F(k)$ holds, hence the implication is true. Additionally, since we used the assumptions $E(k)$ and $F(k-1)$ only for $n>0$, F(0) holds with no further assumptions.
\end{proof}

These three lemmas allow us to prove that the conjunction of statements $E$ and $F$ holds for all $k\geq 0$, therefore proving that equality~\myeqref{eqassoc} holds.
\begin{lemma}
\label{EFhold}
For $k\geq 0$, $E(k)\wedge F(k)$ holds.
\end{lemma}
\begin{proof}
We prove this by induction on $k$:

\begin{itemize}
\item \emph{Base case}: For $k=0$

\begin{itemize}

\item $F(0)$ holds by Lemma~\ref{FEimpliesF}.

\item $E(0)$ holds by LemmA~\ref{recursionbaseE}.

\end{itemize}

\item \emph{Inductive step}: For  $k>0$ we assume that both $F(k-1)$ and $E(k-1)$ hold.
\begin{itemize}
\item To establish $E(k)$ by Lemma~\ref{FEimpliesE} we need both $E(k-1)$ and $F(k-1)$ to hold. Both hold by the inductive hypothesis.

\item To establish $F(k)$ by Lemma~\ref{FEimpliesF} we need both $E(k)$ and $F(k-1)$ to hold. Both hold by the statement above.
\end{itemize}
\end{itemize}
By this inductive argument the statement $E(k)\wedge F(k)$ holds for $k\geq 0$.
\end{proof}

Below, in a similar way to Definitions~\ref{defE},~\ref{defF}, we make two logical statements such that when their conjunction is shown for all $k\geq 0$, equality~\myeqref{eqdistleft} is proved. Here $k$ is the difference between dimensions of the two diagrams in the bracket.
\begin{definition}
\label{defG}
For $k\geq 0$, $G(k)$ states that for an $n$\-diagram $D$, an $m$\-diagram $S$ and an $l$\-diagram $M$, all well-defined, such that $l,n>m>0$ and $|l-n|=k$, the following holds:
\begin{align*}
& S\circ_{b}(D\circ_{a} M)  = (S\circ_{b} D)\circ_{a}(S\circ_{b} M) &\text{if } b<a
\end{align*}
Here, we have $a=\mathrm{min}(n,l)-1$, $b=\mathrm{min}(m, \mathrm{max}(n,l)) - 1$.
\end{definition}

\begin{definition}
\label{defH}
For $k\geq 0$, $H(k)$ states that for an $n$\-diagram $D$, an $m$\-diagram $S$ and an $l$\-diagram $M$, all well-defined, such that $l,n>m>0$ and $|l-n|=k$, then, provided that these composites exist, the following holds:
\begin{align*}
\inc{S}{D\circ M}\circ \inc{D}{M} = \inc{S\circ D}{S\circ M}\circ \inc{S}{M}
\end{align*}
Here, we have $a=\mathrm{min}(n,l)$, $b=\mathrm{min}(m, \mathrm{max}(n,l)) - 1$.
\end{definition}

Again, in a similar fashion to Lemma~\ref{recursionbaseE}, we first establish the statement $G$ for $k=0$.
\begin{lemma}
\label{recursionbaseG}
The statement $G(0)$ holds without any further assumptions.
\end{lemma}
\begin{proof}

For $G(0)$ to hold, we need to show that for any three well-defined diagrams: an $n$\-diagrams $D$, an $m$\-diagram $S$ and an $l$\-diagram $M$, such that $l,n>m>0$ and $|l-n|=k$, the following holds:
\begin{align*}
&S\circ_{b}(D\circ_{a} M)  = (S\circ_{b} D)\circ_{a}(S\circ_{b} M)
\end{align*}
Since $k=|l-n|=0$, we have  $l=n>m$. Bearing that in mind we drop the composition indices. As $m>0$, by Definition~\ref{diagequiv}, to show equivalence of these two diagrams, we need to check the following four conditions:

\begin{itemize}

\item Sources are equivalent diagrams:
\begin{align*}
&(S\circ (D\circ M)).s \\
&\quad= S\circ (D\circ M).s &&[\myeqref{diagcomposesourcesnmgreater}]\\
&\quad= S\circ D.s &&[\myeqref{diagcomposesourcesnmequal}]\\
&\quad= (S\circ D).s &&[\myeqref{diagcomposesourcesnmgreater}]\\
&\quad= ((S\circ D)\circ (S\circ M)).s &&[\myeqref{diagcomposesourcesnmequal}]
\end{align*}

\item Sizes of generator and embedding lists for both diagrams are equal:
\begin{align*}
&|S\circ (D\circ M)| \\
&\quad= |D\circ M| &&[\myeqref{diagcomposesizesnmgreater}]\\
&\quad= |D| + |M|&&[\myeqref{diagcomposesourcesnmequal}]\\
&\quad= |(S\circ D)| + |(S\circ M)| &&[\myeqref{diagcomposesizesnmgreater}]\\
&\quad= |(S\circ D)\circ (S\circ M)| &&[\myeqref{diagcomposesourcesnmequal}]
\end{align*}

\item Corresponding generators are equal for $0\leq j\leq |M| + |D|$, we show this for two separate ranges, first assume $|D|\leq j < |D| + |M|$:
\begin{align*}
&(S\circ (D\circ M))[j].g \\
&\quad= (D\circ M)[j].g &&[\myeqref{diagcomposegeneratorsnmgreater}]\\
&\quad= M[j-|D|].g &&[\myeqref{diagcomposegeneratorsnmequal}]\\
&\quad= (S\circ M)[j-|D|].g &&[\myeqref{diagcomposegeneratorsnmgreater}]\\
&\quad= (S\circ M)[j-|S\circ D|].g &&[\myeqref{diagcomposesizesnmgreater}]\\
&\quad= ((S\circ D)\circ (S\circ M))[j].g &&[\myeqref{diagcomposegeneratorsnmequal}]
\end{align*}
The argument for $0\leq j < |D|$ is analogous.

\item Corresponding embeddings are equivalent for $0\leq j\leq |M| + |D|$, we show this for two separate ranges, first assume $|D|\leq j < |D| + |M|$:
\begin{align*}
&(S\circ (D\circ M))[j].e \\
&\quad= \inc{S}{(D\circ M)[j].d}\circ (D\circ M)[j].e &&[\myeqref{diagcomposeembeddingsnmgreater}]\\
&\quad= \inc{S}{(D\circ M)[j].d}\circ M[j-|D|].e &&[\myeqref{diagcomposeembeddingsnmequal}]\\
&\quad= \inc{S}{M[j-|D|.d]}\circ M[j-|D|].e &&[\defref{composedslices}]\\
&\quad= (S\circ M)[j - |D|].e &&[\myeqref{diagcomposesizesnmequal}]\\
&\quad= (S\circ M)[j - |S\circ D|].e &&[\myeqref{diagcomposesizesnmgreater}]\\
&\quad= ((S\circ D)\circ (S\circ M))[j].e &&[\myeqref{diagcomposeembeddingsnmequal}]
\end{align*}
The argument for $0\leq j < |D|$ is analogous.
\end{itemize}
Since all these conditions are satisfied, we established that $S\circ (D\circ M)=(S\circ D)\circ (S\circ M)$ for an $m$\-diagram $S$ and an $l$\-diagram $M$, such that $l=n>m>0$. Hence, the statement $G(0)$ holds, as required.
\end{proof}

This is followed by implications establishing $G(k)$ and $H(k)$.
\begin{lemma}
\label{GHimpliesG}
For $k\geq1$ the following statement holds: $G(k-1)\wedge H(k-1)\implies G(k)$
\end{lemma}
\begin{proof}
Assume that both $G(k-1)$ and $H(k-1)$ hold. For $G(k)$ to hold, we need to show that for any three well-defined diagrams: an $n$\-diagrams $D$, an $m$\-diagram $S$ and an $l$\-diagram $M$, such that $l,n>m>0$ and $|l-n|=k$, the following holds:
\begin{align*}
&S\circ_{b}(D\circ_{a} M)  = (S\circ_{b} D)\circ_{a}(S\circ_{b} M)
\end{align*}
Since $k=|l-n|>0$, we have  $l>n>m$ or $n>l>m$. First, assume $l>n>m$, then bearing that in mind we drop the composition indices. As $m>0$, by Definition~\ref{diagequiv}, to show equivalence of these two diagrams, we need to check the following four conditions:

\begin{itemize}
\item Sources are equivalent diagrams:
\begin{align*}
&(S\circ (D\circ M)).s \\
&\quad= S\circ (D\circ M).s &&[\myeqref{diagcomposesourcesnmgreater}]\\
&\quad= S\circ (D\circ M.s) &&[\myeqref{diagcomposesourcesnmgreater}]\\
&\quad= (S\circ D)\circ (S\circ M.s) &&[G(k-1)]\\
&\quad= (S\circ D)\circ (S\circ M).s &&[\myeqref{diagcomposesourcesnmgreater}]\\
&\quad= ((S\circ D)\circ (S\circ M)).s &&[\myeqref{diagcomposesourcesnmgreater}]
\end{align*}

\item Sizes of generator and embedding lists for both diagrams are equal:
\begin{align*}
&|S\circ (D\circ M)| \\
&\quad= |D\circ M| &&[\myeqref{diagcomposesizesnmgreater}]\\
&\quad= |M|&&[\myeqref{diagcomposesizesnmgreater}]\\
&\quad= |(S\circ M)| &&[\myeqref{diagcomposesizesnmgreater}]\\
&\quad= |(S\circ D)\circ (S\circ M)|&&[\myeqref{diagcomposesizesnmgreater}]
\end{align*}

\item Corresponding generators are equal for $0\leq j\leq |M|$:
\begin{align*}
&(S\circ (D\circ M))[j].g \\
&\quad= M[j].g &&[\myeqref{diagcomposegeneratorsnmgreater}]\\
&\quad= (S\circ M)[j].g &&[\myeqref{diagcomposegeneratorsnmgreater}]\\
&\quad= ((S\circ D)\circ (S\circ M))[j].g &&[\myeqref{diagcomposegeneratorsnmgreater}]
\end{align*}

\item Corresponding embeddings are equal for $0\leq j\leq |M|$:
\begin{align*}
&(S\circ (D\circ M))[j].e \\
&\quad= \inc{S}{(D\circ M)[j].d}\circ (D\circ M)[j].e &&[\myeqref{diagcomposeembeddingsnmgreater}]\\
&\quad= \inc{S}{(D\circ M)[j].d}\circ (\inc{D}{M[j].d}\circ M[j].e) &&[\myeqref{diagcomposeembeddingsnmgreater}]\\
&\quad= (\inc{S}{(D\circ M)[j].d}\circ \inc{D}{M[j].d})\circ M[j].e &&[A(n)]\\
&\quad= (\inc{S\circ D}{(S\circ M)[j].d}\circ \inc{S}{M[j].d})\circ M[j].e &&[H(k-1)]\\
&\quad= \inc{S\circ D}{(S\circ M)[j].d}\circ (\inc{S}{M[j].d}\circ M[j].e) &&[A(n)]\\
&\quad= \inc{S\circ D}{(S\circ M)[j].d}\circ (S\circ M)[j].e &&[\myeqref{diagcomposeembeddingsnmgreater}]\\
&\quad= ((S\circ D)\circ (S\circ M))[j].e &&[\myeqref{diagcomposeembeddingsnmgreater}]
\end{align*}
\end{itemize}
Here $A(n)$ is the statement on associativity of embedding composition proved in Theorem~\ref{bigtheorem}. Since all these conditions are satisfied, we established that $S\circ (D\circ M)=(S\circ D)\circ (S\circ M)$ for an $m$\-diagram $S$ and an $l$\-diagram $M$, such that $l>n>m>0$.

The argument for $n>l>m$ is analogous. Hence, the statement $G(k)$ holds and the implication is true.
\end{proof}

\begin{lemma}
\label{GHimpliesH}
For $k\geq 1$ the following statement holds: $H(k-1) \wedge G(k)\implies H(k)$, additionally H(0) holds with no further assumptions.
\end{lemma}
\begin{proof}

Let us assume that both $H(k-1)$ and $G(k)$ hold. 

For $H(k)$ to hold, we need to show that for any three well-defined diagrams: an $n$\-diagrams $D$, an $m$\-diagram $S$ and an $l$\-diagram $M$, such that $l,n>m>0$ and $|l-n|=k$, the following holds:
\begin{align*}
\inc{S}{D\circ M}\circ \inc{D}{M} = \inc{S\circ D}{S\circ M}\circ \inc{S}{M}
\end{align*}

\begin{itemize}
\item For $k=0$, we have $l=n>m$
\item For $k=|l-n|>0$, we have  $l>n>m$ or $n>l>m$
\end{itemize}
First, assume $l\geq n>m$, to consider the cases $l=n>m$ and $l>n>m$ simultaneously. Since $l>n>0$, to establish that these two embeddings are equivalent, by Lemma~\ref{embeddingequiv}, we need to check the following conditions.

\paragraph{Types are the same.} By Definition~\ref{definclusion} these types are as follows:
\begin{align*}
&\inc{D}{M}: M\hookrightarrow D\circ M \\
&\inc{S}{M}: M\hookrightarrow S\circ M \\
&\inc{S}{D\circ M}: D\circ M\hookrightarrow S\circ (D\circ M) \\
&\inc{S\circ D}{S\circ M}: S\circ M\hookrightarrow (S\circ D)\circ (S\circ M) 
\end{align*}
We could see that the domains are immediately equivalent, the codomains are equivalent by G(k).

\paragraph{Heights are equal.} For this condition, we consider two scenarios:
\begin{itemize}
\item For $l=n>m$:
\begin{align*}
&(\inc{S}{D\circ M}\circ \inc{D}{M}).h = \\
&\quad= (\inc{S}{D\circ M}).h + (\inc{D}{M}).h &&[\defref{defcomposition}] \\
&\quad= (D \circ M).\id.h + |D| &&[\defref{definclusion}] \\
&\quad= |D| &&[\defref{identityembedding}] \\
&\quad= |S\circ D| &&[\myeqref{diagcomposesizesnmgreater}] \\
&\quad= (\inc{S\circ D}{S\circ M}).h + (\inc{S}{M}).h &&[\defref{definclusion}]
\\
&\quad= (\inc{S\circ D}{S\circ M}\circ \inc{S}{M}).h &&[\defref{defcomposition}]
\end{align*}
\item For $l>n>m$:
\begin{align*}
&(\inc{S}{D\circ M}\circ \inc{D}{M}).h = \\
&\quad= (\inc{S}{D\circ M}).h + (\inc{D}{M}).h   &&[\defref{defcomposition}]
\\
&\quad= (D \circ M).\id.h + M.\id.h  &&[\defref{definclusion}] \\
&\quad= M.\id.h  &&[\defref{identityembedding}] \\
&\quad= (S\circ M).\id.h + M.\id.h  &&[\defref{identityembedding}] \\
&\quad= (\inc{S\circ D}{S\circ M}).h + (\inc{S}{M}).h  &&[\defref{definclusion}]
\\
&\quad= (\inc{S\circ D}{S\circ M}\circ \inc{S}{M}).h  &&[\defref{defcomposition}]
\end{align*}
\end{itemize}

\paragraph{Component embeddings are equivalent.} Again, we consider two scenarios:
\begin{itemize}
\item For $l=n>m$:
\begin{align*}
&(\inc{S}{D\circ M}\circ \inc{D}{M}).e = \\
&\quad= \lift{(\inc{S}{D\circ M}).e}{(D\circ M)[\inc{D}{M}.h].d}\circ \inc{D}{M}.e  &&[\defref{defcomposition}] \\
&\quad= \lift{(\inc{S}{D\circ M}).e}{(D\circ M)[|D|].d}\circ \inc{D}{M}.e  &&[\defref{definclusion}] \\
&\quad=\inc{S}{(D\circ M)[|D|].d}\circ \inc{D}{M}.e &&[M(k)] \\
&\quad=\inc{S}{M[|D| - |D|].d}\circ \inc{D}{M}.e &&[\defref{composedslices}] \\
&\quad=\inc{S}{M[0].d}\circ \inc{D}{M}.e &&[\defref{definitialterminalslice}] \\
&\quad=\inc{S}{M.s}\circ \inc{D}{M}.e &&[\defref{definitialterminalslice}] \\
&\quad=\inc{S}{M.s}\circ M.\id.e  &&[\defref{definclusion}] \\
&\quad=\inc{S}{M.s}\circ M.s.\id  &&[\defref{identityembedding}] \\
&\quad=\inc{S}{M.s} &&[\defref{identitycancellation}] \\
&\quad=S\circ M.s.\id \circ \inc{S}{M.s}& &[\defref{identitycancellation}] \\
&\quad=S\circ M.s.\id \circ \inc{S}{M}.e &&[\defref{definclusion}] \\
&\quad=(S\circ M).s.\id \circ \inc{S}{M}.e &&[\myeqref{diagcomposesourcesnmgreater}] \\
&\quad=(S\circ M)[0].d.\id \circ \inc{S}{M}.e &&[\defref{definitialterminalslice}] \\
&\quad= \lift{((S\circ M).\id.e)}{(S\circ M)[0].d}\circ \inc{S}{M}.e &&[M(k)] \\
&\quad= \lift{((\inc{S\circ D}{S\circ M}).e)}{(S\circ M)[0].d}\circ \inc{S}{M}.e &&[\defref{identityembedding}] \\
&\quad= \lift{((\inc{S\circ D}{S\circ M}).e)}{(S\circ M)[\inc{S}{M}.h].d}\circ \inc{S}{M}.e &&[\defref{definclusion}] \\
&\quad= (\inc{S\circ D}{S\circ M}\circ \inc{S}{M}).e &&[\defref{defcomposition}]
\end{align*}

\item For $l>n>m$:
\begin{align*}
&(\inc{S}{D\circ M}\circ \inc{D}{M}).e = \\
&\quad= \lift{(\inc{S}{D\circ M}).e}{(D\circ M)[\inc{D}{M}.h].d}\circ \inc{D}{M}.e  &&[\defref{defcomposition}] \\
&\quad= \lift{(\inc{S}{D\circ M}).e}{(D\circ M)[0].d}\circ \inc{D}{M}.e  &&[\defref{definclusion}] \\
&\quad=\inc{S}{(D\circ M)[0].d}\circ \inc{D}{M}.e &&[M(k)] \\
&\quad=\inc{S}{(D\circ M).s}\circ \inc{D}{M}.e &&[\defref{definitialterminalslice}] \\
&\quad=\inc{S}{(D\circ M.s)}\circ \inc{D}{M.s} &&[\myeqref{diagcomposesourcesnmgreater}] \\
&\quad=\inc{S\circ D}{S\circ M.s}\circ \inc{S}{M.s} &&[IH]\\
&\quad=\inc{S\circ D}{(S\circ M).s}\circ \inc{S}{M}.e &&[\myeqref{diagcomposesourcesnmgreater}] \\
&\quad=\inc{S\circ D}{(S\circ M)[0].d}\circ \inc{S}{M}.e &&[\defref{definitialterminalslice}] \\
&\quad= \lift{(\inc{S\circ D}{S\circ M}).e}{(S\circ M)[0].d}\circ \inc{S}{M}.e &&[M(k)] \\
&\quad= \lift{(\inc{S\circ D}{S\circ M}).e}{(S\circ M)[\inc{S}{M}.h].d}\circ \inc{S}{M}.e &&[\defref{definclusion}] \\
&\quad= (\inc{S\circ D}{S\circ M}\circ \inc{S}{M}).e &&[\defref{defcomposition}]
\end{align*}
\end{itemize}
As all these conditions are satisfied both embeddings are equivalent. The argument for $n>l>m$ is analogous. By this, we established that $H(k)$ holds, hence the implication is true. Additionally, since we used the assumptions $G(k)$ and $H(k-1)$ only for $k>0$, H(0) holds with no further assumptions.
\end{proof}

Finally, we can put all the pieces together to establish that the conjunction of $G$ and $H$ holds for all $k\geq 0$:
\begin{lemma}
\label{GHhold}
For $k\geq 0$: the following holds: $G(k)\wedge H(k)$
\end{lemma}
\begin{proof}
We prove this by induction on $k$.

\paragraph{Base case.} For $k=0$:
\begin{itemize}
\item $G(0)$ holds by Lemma~\ref{recursionbaseG}.
\item $H(0)$ holds by Lemma~\ref{GHimpliesH}.
\end{itemize}

\paragraph{Inductive step.} For  $k>0$ we assume that both $F(k-1)$ and $E(k-1)$ hold.
\begin{itemize}
\item To establish $G(k)$ by Lemma~\ref{GHimpliesG} we need both $G(k-1)$ and $H(k-1)$ to hold. Both hold by the inductive hypothesis.
\item To establish $H(k)$ by Lemma~\ref{GHimpliesH} we need both $G(k)$ and $H(k-1)$ to hold. Both hold by the statement above.
\end{itemize}
By this inductive argument the statement $G(k)\wedge H(k)$ holds for $k\geq 0$.
\end{proof}

\section{Coherent adjunctions}
\label{chapterbutterfly}

The purpose of this section is to prove the following theorem.

\begin{theorem}
\label{thm:butterfly}
In a quasistrict 4\-category, an adjunction of 1\-morphisms gives rise to a coherent adjunction satisfying the butterfly equations.
\end{theorem}

\noindent
The butterfly equations are equations holding between specified 4\-cells. By general results of Riehl and Verity~\cite{Riehl_2016}, this is expected to hold an any correct algebraic model of 4\-categories. We believe that our proof is the first to be given; indeed, we believe that it is the first substantial proof of any sort to be presented in an algebraic 4\-category. The full proof itself is formalized in the \emph{Globular} proof assistant at \href{http://globular.science/1605.002}{globular.science/1605.002}, and given in substantial detail later in this section.

There is a corresponding result at the level of 3\-categories, which reads as follows.
\begin{theorem}
\label{thm:swallowtail}
In a tricategory, an adjunction of 1\-morphisms gives rise to a coherent adjunction satisfying the swallowtail equations.
\end{theorem}

\noindent
The swallowtail equations are equations holding between specified 3\-cells. This was first established by Verity~\cite{VeritySwallowtail}, and later discussed in depth by Gurski~\cite{GurskiThesis, GurskiTricategories} and Pstragowski~\cite{piotr-thesis}. Given an adjunction in a quasistrict 4\-category, one can apply the proof of Theorem~\ref{thm:swallowtail} (suitably expressed in a Gray category) to obtain invertible 4\-morphisms called \emph{swallowtailators}, which witness the satisfaction up to isomorphism of the swallowtail equations. The butterfly equations are certain equations that these swallowtailators may, or may not, satisfy.

Analogously to the proof of Theorem~\ref{thm:swallowtail}, the proof of Theorem~\ref{thm:butterfly} proceeds by redefining one of the swallowtailators in a certain way, and then demonstrating explicitly that a certain equation holds. There are two butterfly equations; the proofs of each are in fact rather different, but for reasons of length we only present one of them here.

This section is structured as follows. In Section~\ref{sec:butterflyhighlights}, we examine some interesting steps of the proof, and see how they arise as instances of our homotopy generators, by viewing their local sources and targets as movies of 2\-projected 3\-diagrams. In Section~\ref{sec:butterflyproject3} we illustrate the entire proof as a single 2\-projected 5\-diagram of height 140. In Section~\ref{sec:butterflyproject2} we illustrate the entire proof in more detail, as a movie of 2\-projected 4\-diagrams; each of these steps is itself a movie of movies of 2\-diagrams.\JVcomm{It would be nice to get a measure of the total number of 2d frames in the entire fully-expanded proof.}

\subsection{Proof highlights}
\label{sec:butterflyhighlights}

Here we visualize some of the steps of the proof sequence given in Section~\ref{sec:butterflyproject2}. In each case, it is interesting to see how it arises as a special case of various homotopy generators.

\paragraph{Step 17.}
This step is an instance of a type $\III'$ homotopy generator. To visualise is, we isolate the subdiagram of the previous movie frame indicated by the yellow highlight box, and represent it as a movie of 2\-projected 3\-diagrams:
\def\imgscale{0.4}
\begin{align*}
\begin{aligned}
\includegraphics[scale=\imgscale]{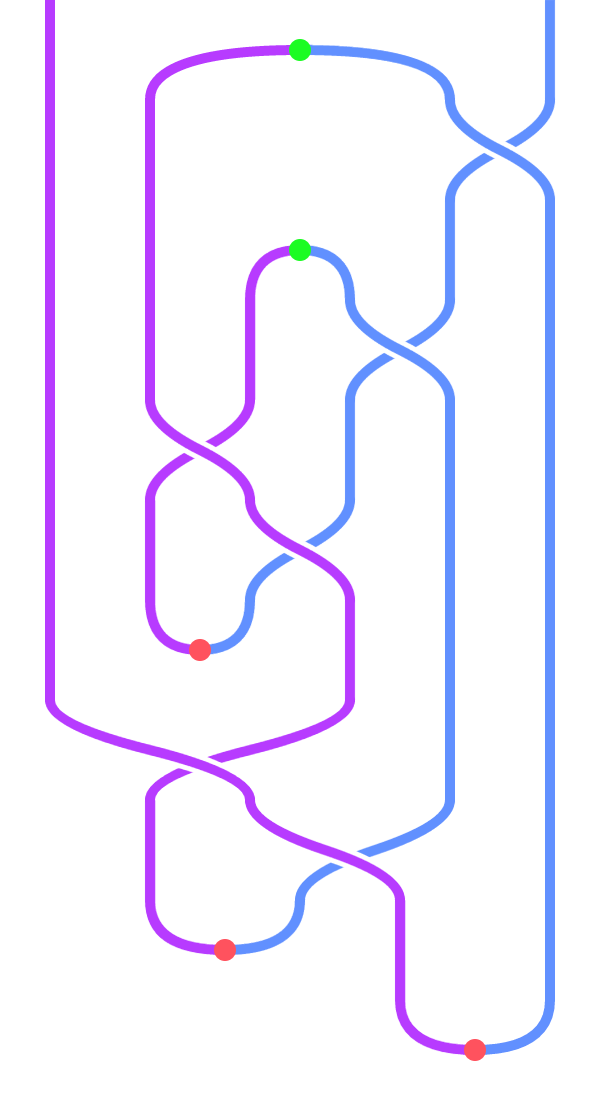}
\end{aligned}
\quad\stackrel{\II_3}{\rightarrow}\quad
\begin{aligned}
\includegraphics[scale=\imgscale]{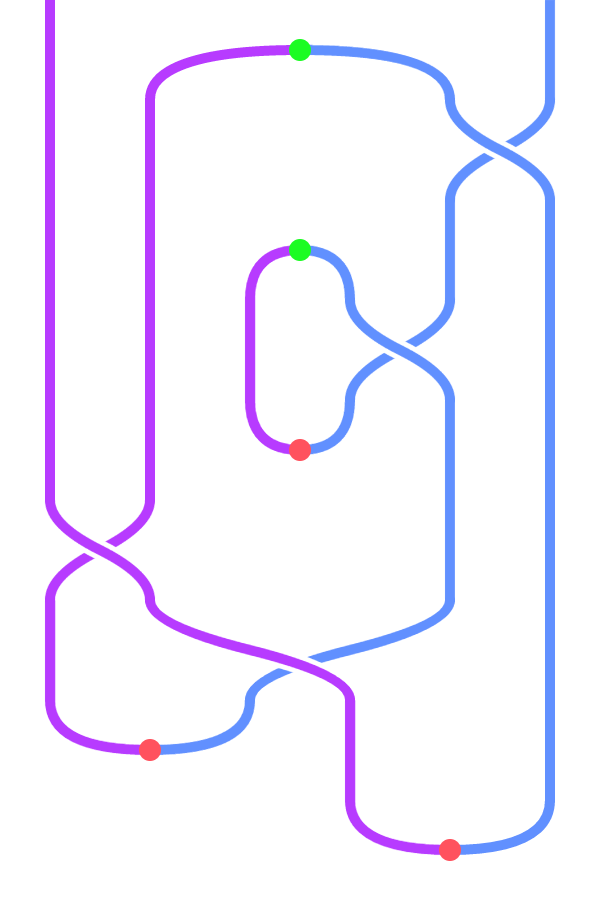}
\end{aligned}
\quad\stackrel{\I_3}{\rightarrow}\quad
\begin{aligned}
\includegraphics[scale=\imgscale]{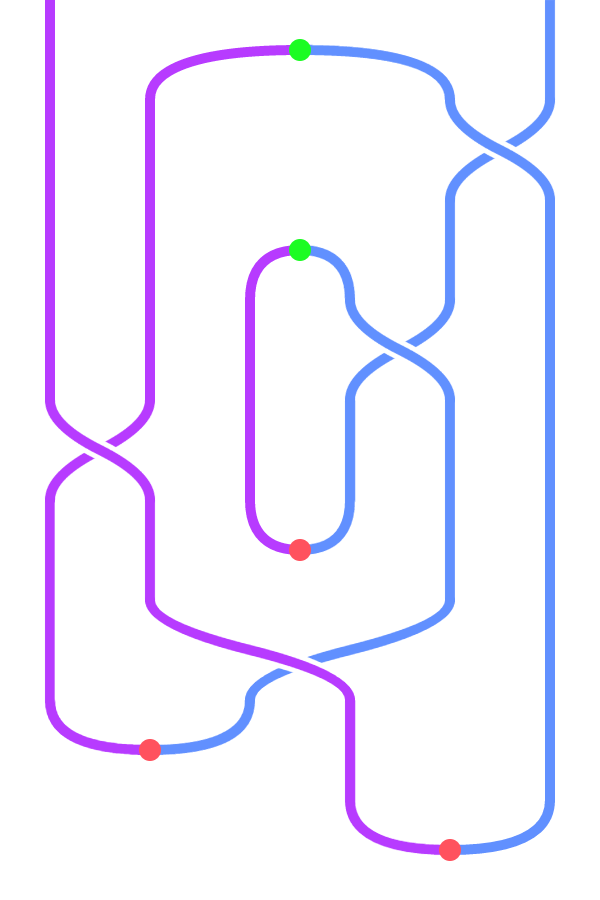}
\end{aligned}
\quad\stackrel{\II_3}{\rightarrow}\quad
\begin{aligned}
\includegraphics[scale=\imgscale]{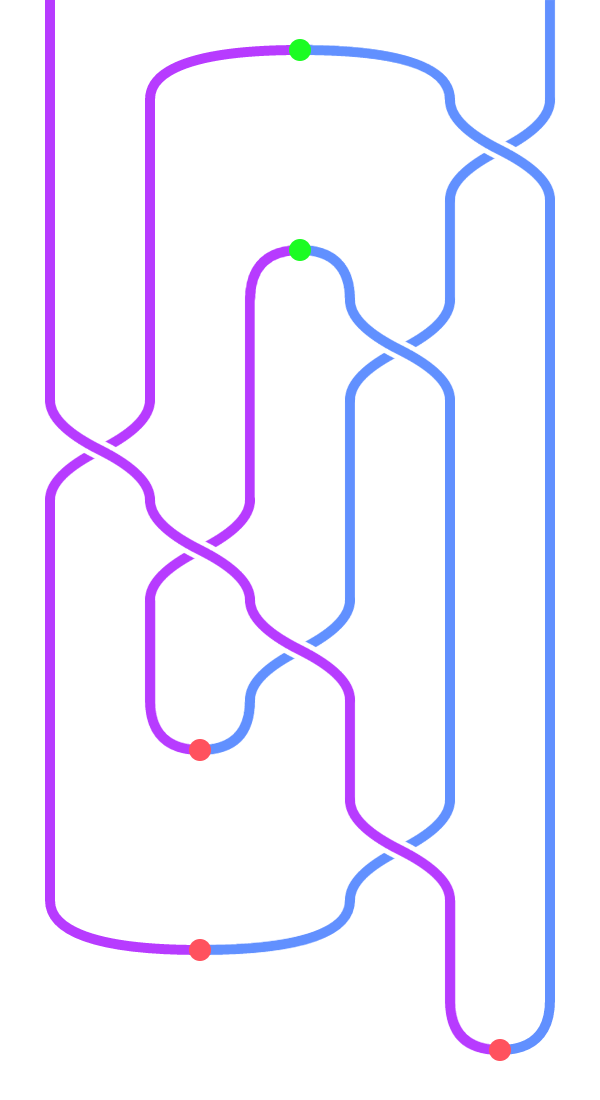}
\end{aligned} 
\end{align*}

\noindent
Step 17 replaces this with the following movie:
\begin{align*}
\begin{aligned}
\includegraphics[scale=\imgscale]{graphics/ButterflySub1.png}
\end{aligned}
&\quad\stackrel{\II_3}{\rightarrow}\quad&
\begin{aligned}
\includegraphics[scale=\imgscale]{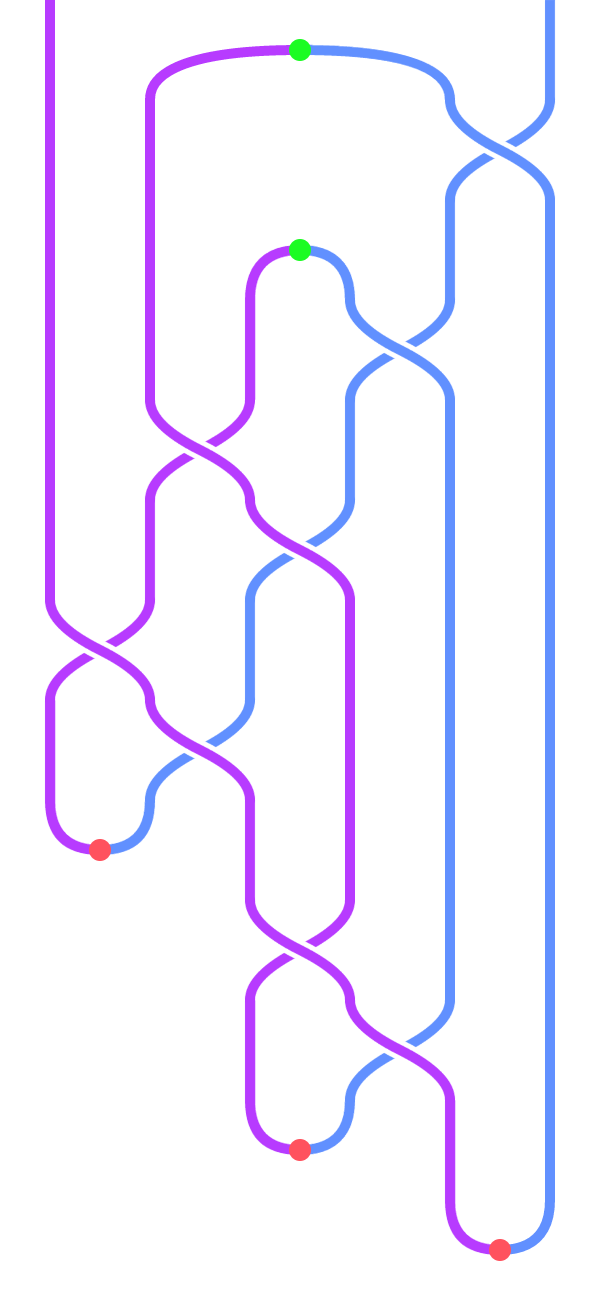}
\end{aligned}
&\quad\stackrel{\I_3}{\rightarrow}\quad&
\begin{aligned}
\includegraphics[scale=\imgscale]{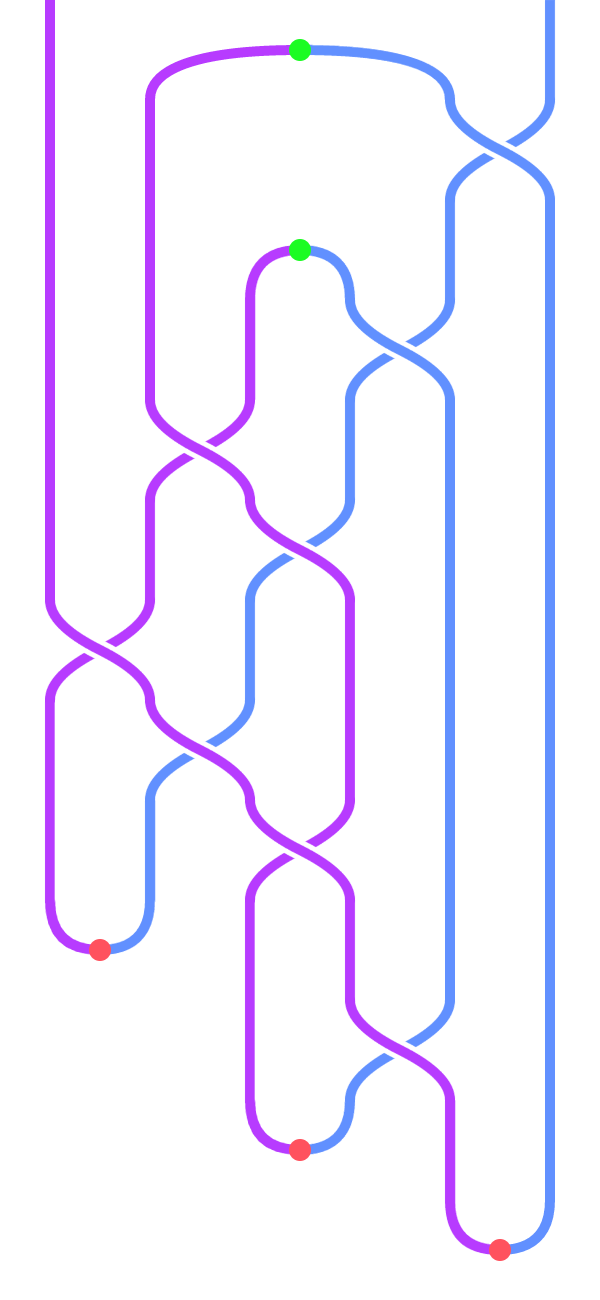}
\end{aligned}
&\quad\stackrel{\I_3}{\rightarrow}\quad&
\begin{aligned}
\includegraphics[scale=\imgscale]{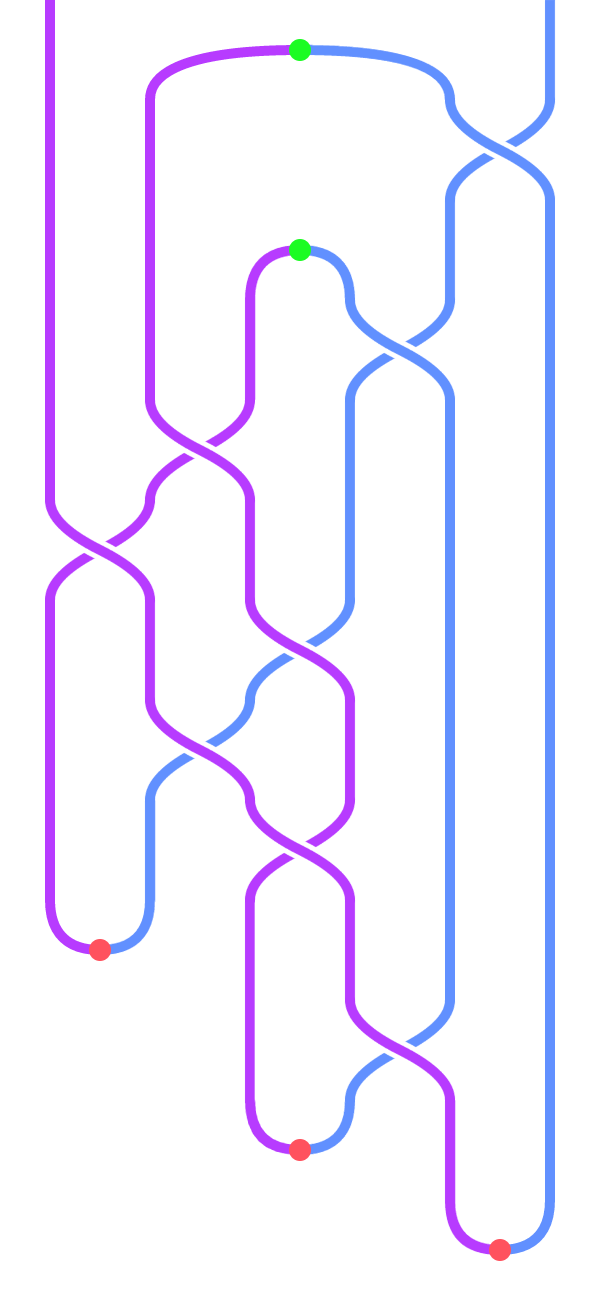}
\end{aligned}
&\quad\stackrel{\II_3}{\rightarrow}\quad& \\
\begin{aligned}
\includegraphics[scale=\imgscale]{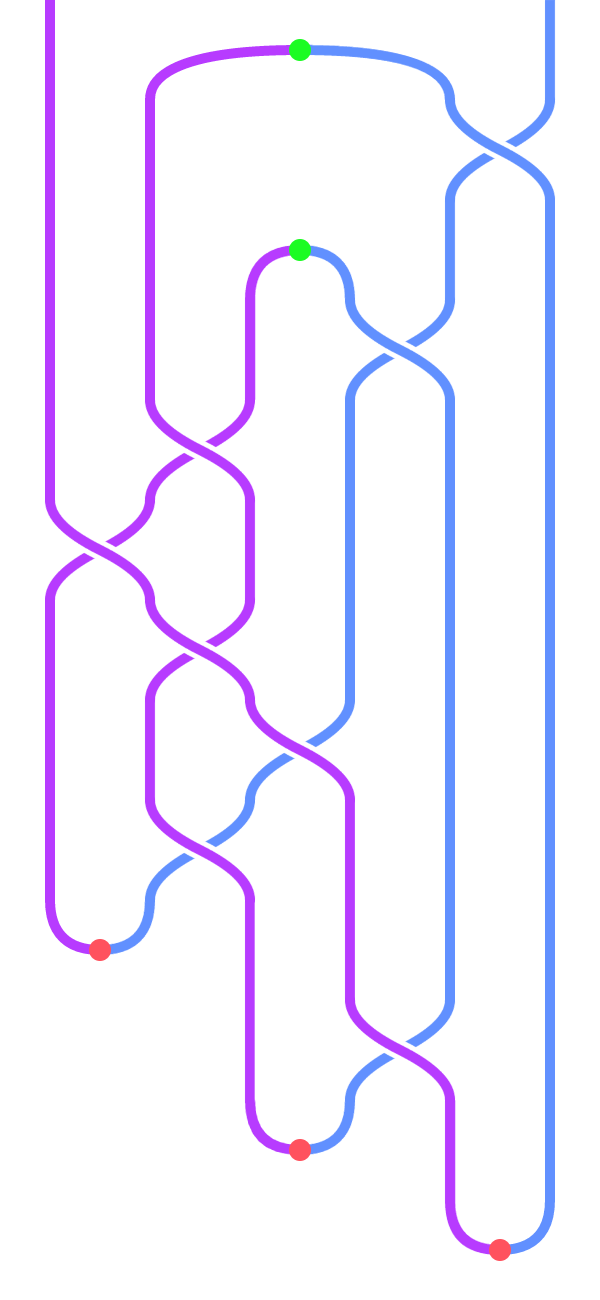}
\end{aligned}
&\quad\stackrel{\II_3}{\rightarrow}\quad&
\begin{aligned}
\includegraphics[scale=\imgscale]{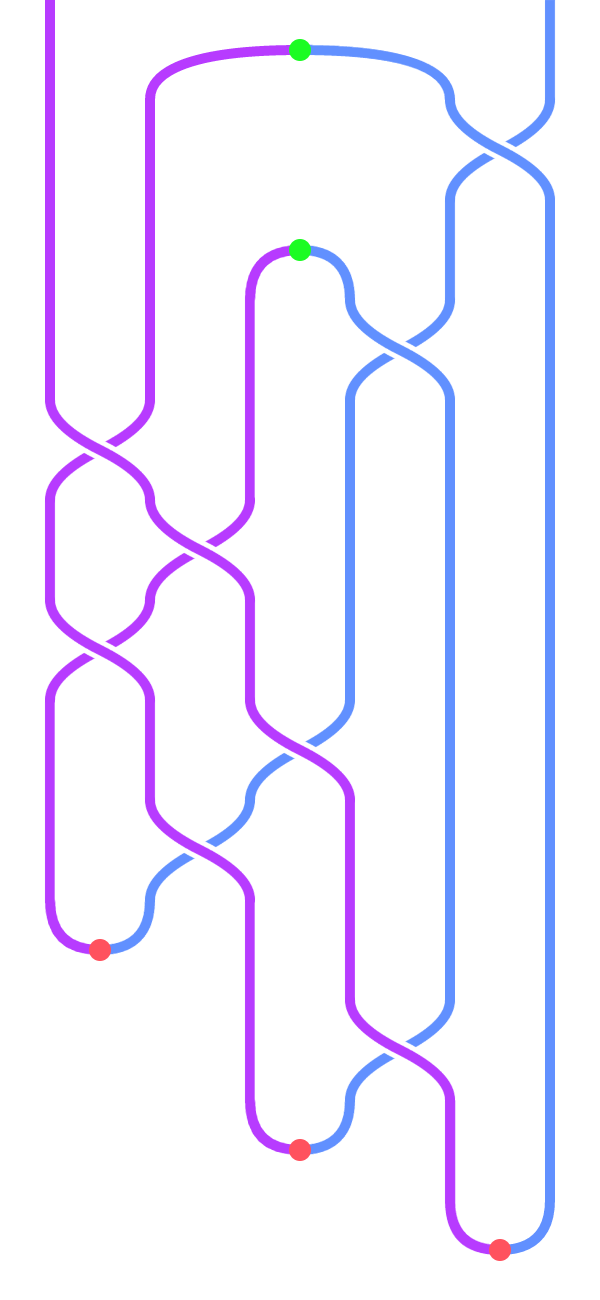}
\end{aligned}
&\quad\stackrel{\I_3}{\rightarrow}\quad&
\begin{aligned}
\includegraphics[scale=\imgscale]{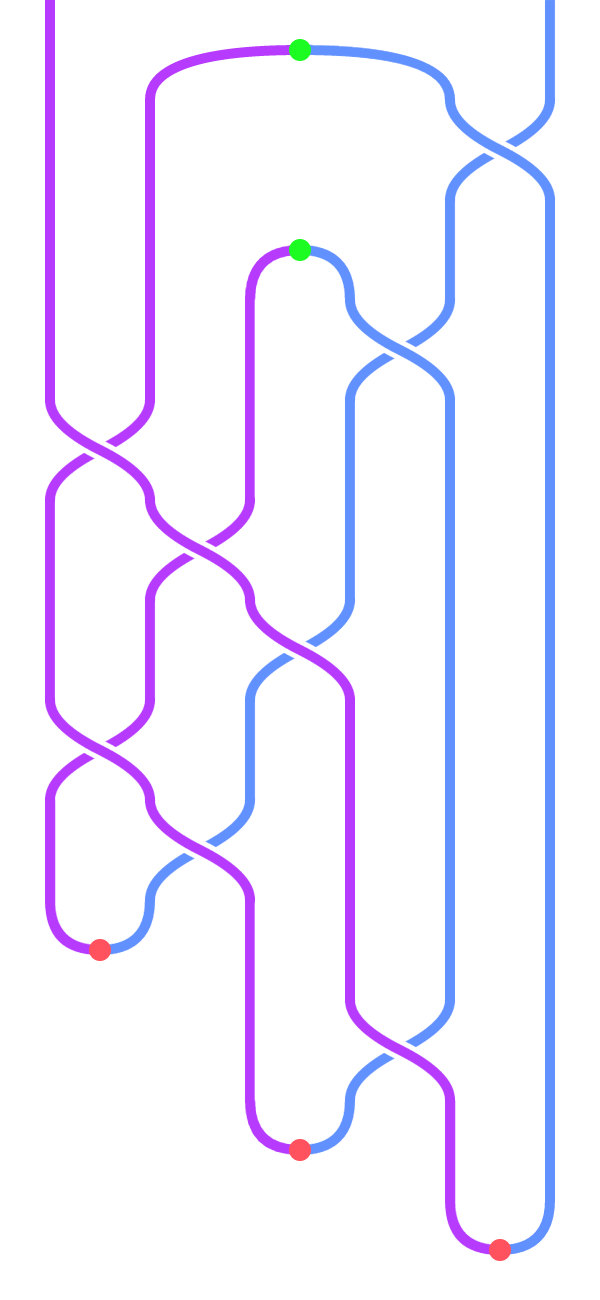}
\end{aligned}
&\quad\stackrel{\II_3}{\rightarrow}\quad&
\begin{aligned}
\includegraphics[scale=\imgscale]{graphics/ButterflySub4.png}
\end{aligned} 
\end{align*}

\noindent

\paragraph{Step 83.}
This is another instance of a type $\III'$ homotopy generator. The highlighted box of the initial slice has the following representation as a movie of 2\-projected 3\-diagrams:
\begin{align*}
&\begin{aligned}
\includegraphics[scale=\imgscale]{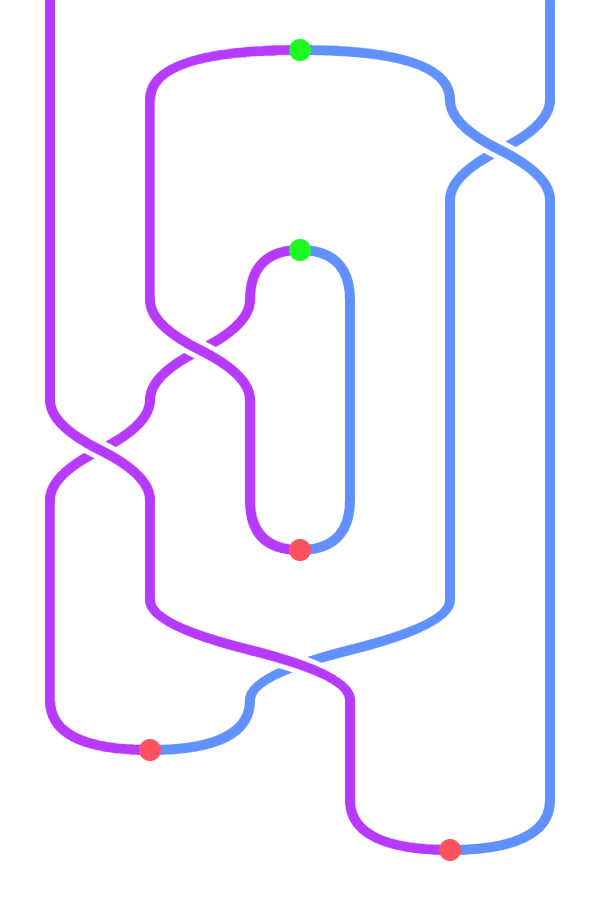}
\end{aligned}
\quad\stackrel{\II_3}{\rightarrow}\quad
\begin{aligned}
\includegraphics[scale=\imgscale]{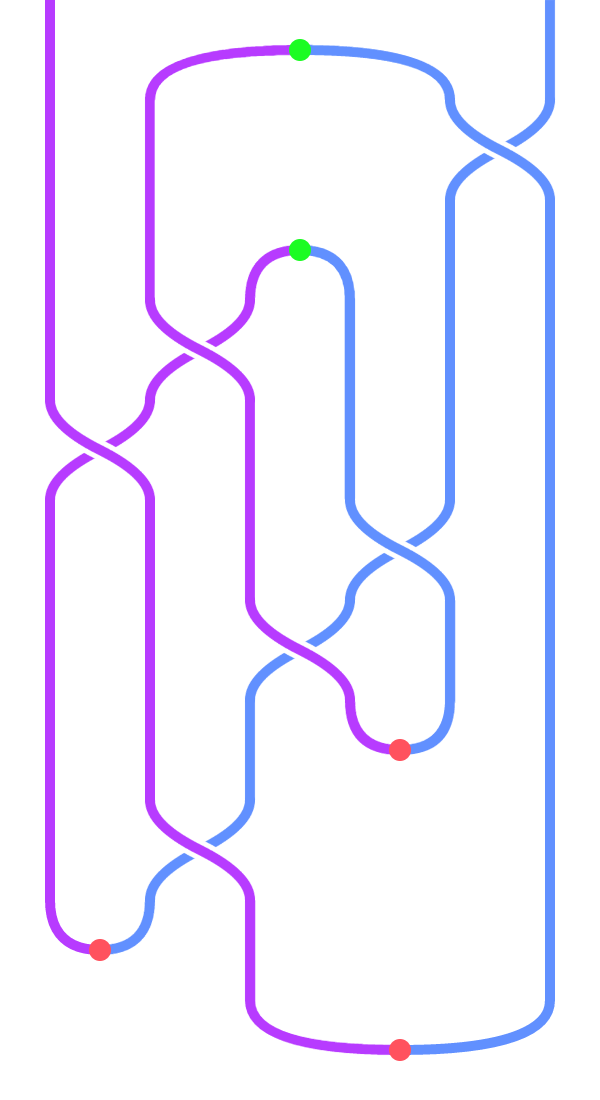}
\end{aligned}
\quad\stackrel{\I_3}{\rightarrow}\quad
\begin{aligned}
\includegraphics[scale=\imgscale]{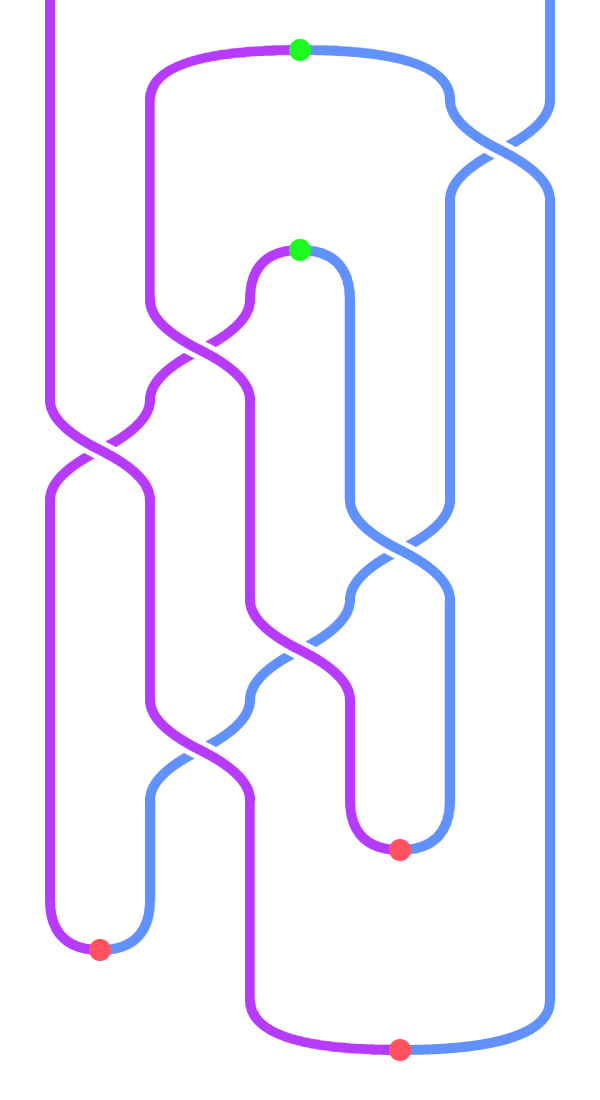}
\end{aligned}
\quad\stackrel{\II_3}{\rightarrow}\quad
\begin{aligned}
\includegraphics[scale=\imgscale]{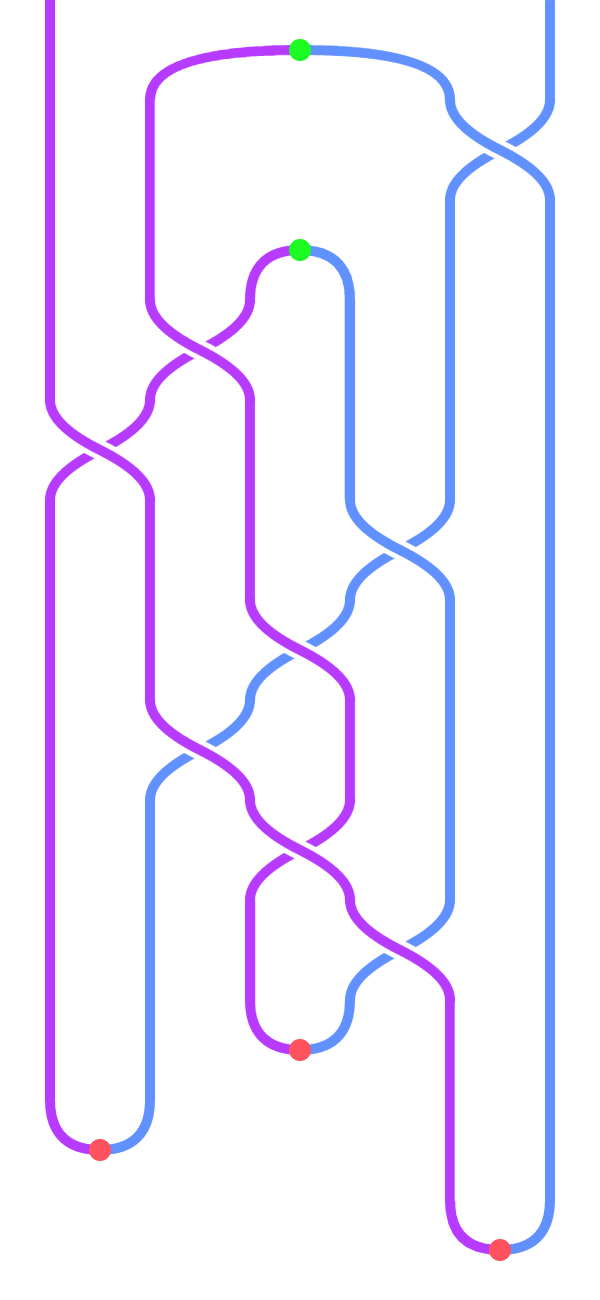}
\end{aligned}
\quad\stackrel{\II_3}{\rightarrow}\quad \\ 
&\begin{aligned}
\includegraphics[scale=\imgscale]{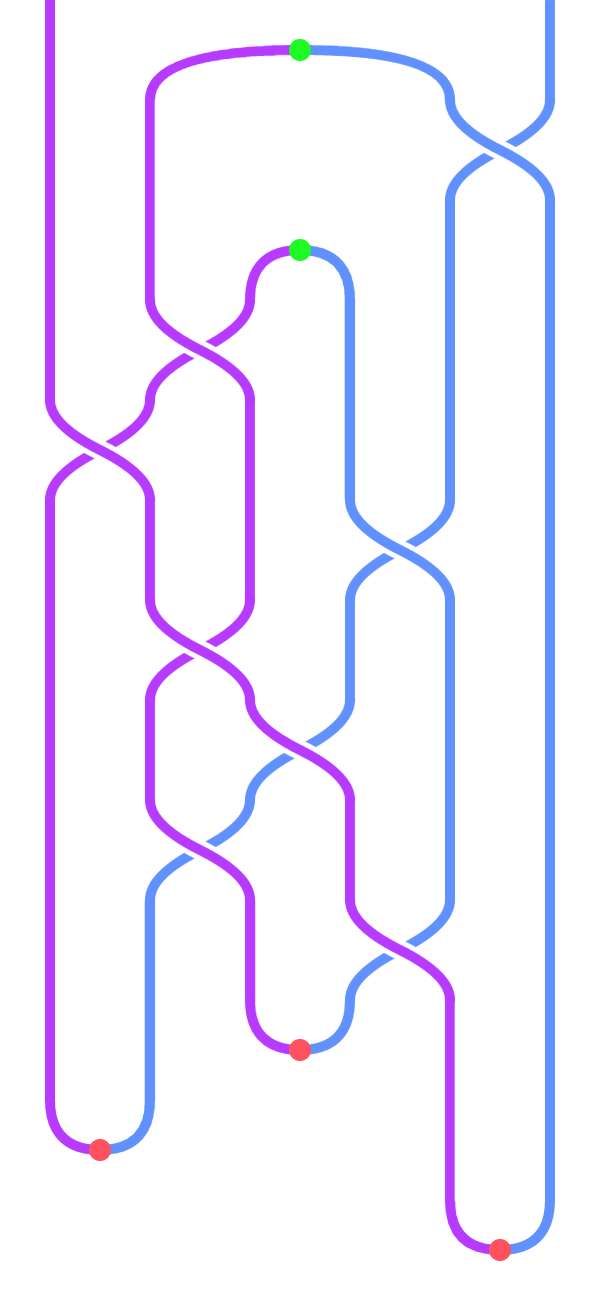}
\end{aligned}
\quad\stackrel{\I_3}{\rightarrow}\quad
\begin{aligned}
\includegraphics[scale=\imgscale]{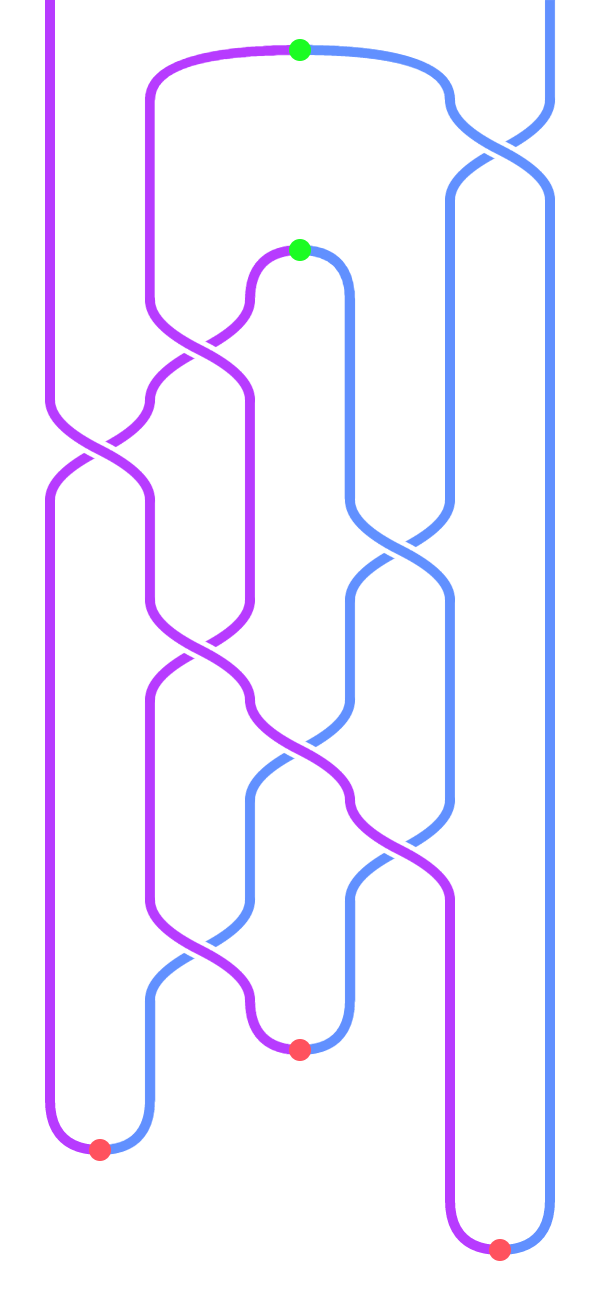}
\end{aligned}
\quad\stackrel{\I_3}{\rightarrow}\quad
\begin{aligned}
\includegraphics[scale=\imgscale]{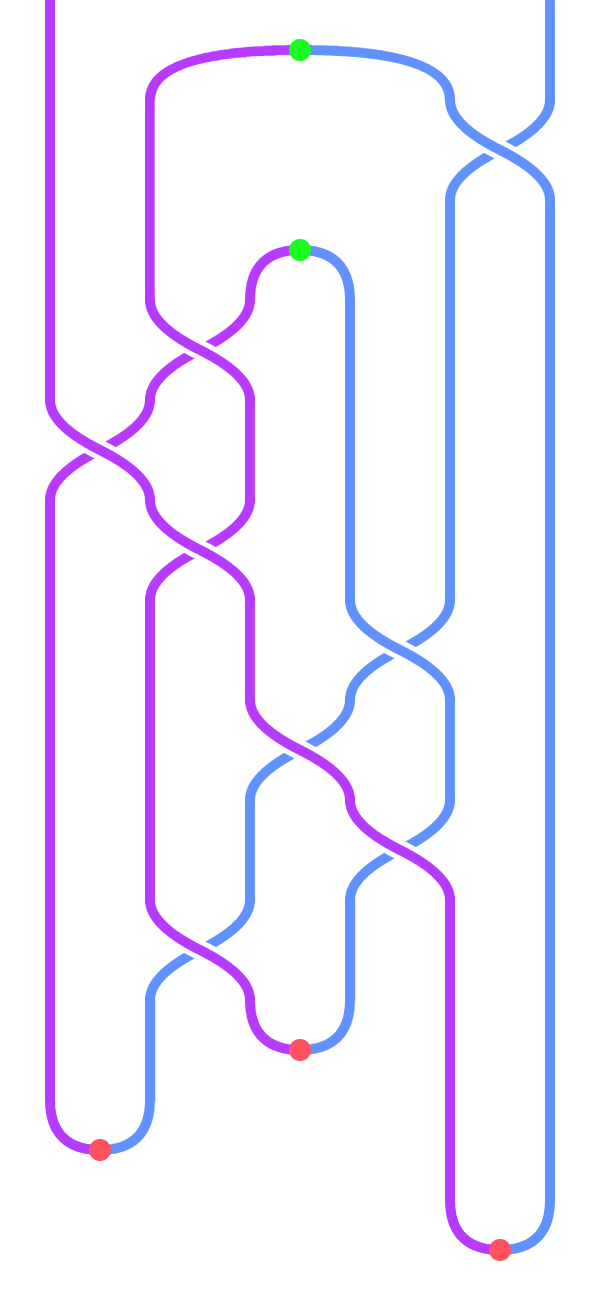}
\end{aligned}
\quad\stackrel{\II_3}{\rightarrow}\quad
\begin{aligned}
\includegraphics[scale=\imgscale]{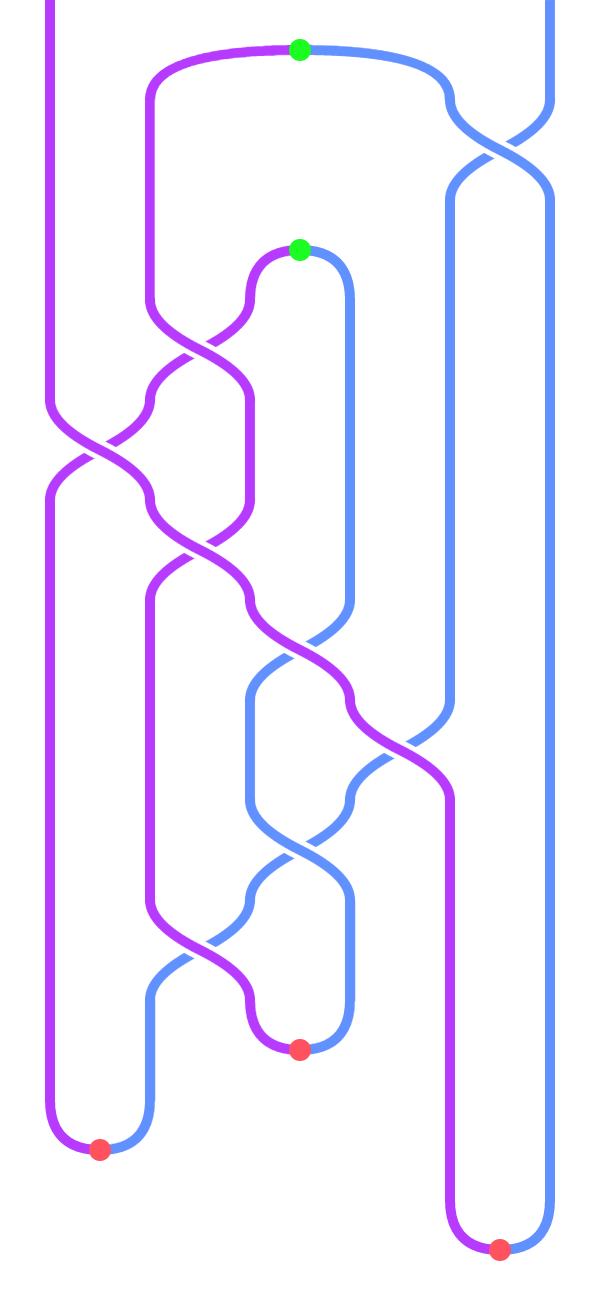}
\end{aligned} 
\end{align*}
This is rewritten as follows:
\begin{align*}
\begin{aligned}
\includegraphics[scale=\imgscale]{graphics/ButterflySubB1.png}
\end{aligned}
\quad\stackrel{\II_3}{\rightarrow}\quad
\begin{aligned}
\includegraphics[scale=\imgscale]{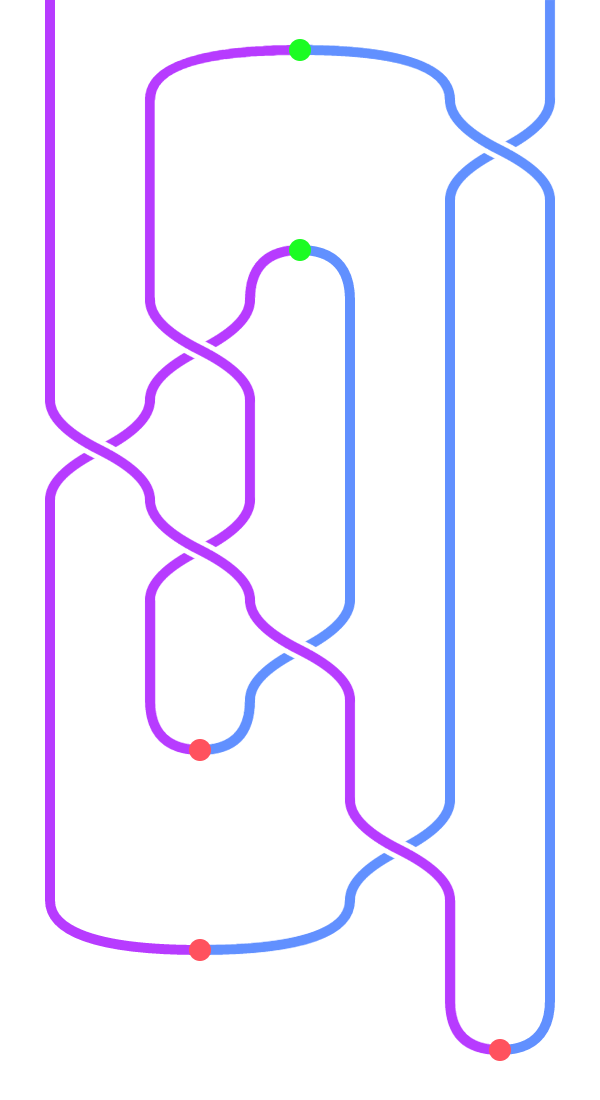}
\end{aligned}
\quad\stackrel{\I_3}{\rightarrow}\quad
\begin{aligned}
\includegraphics[scale=\imgscale]{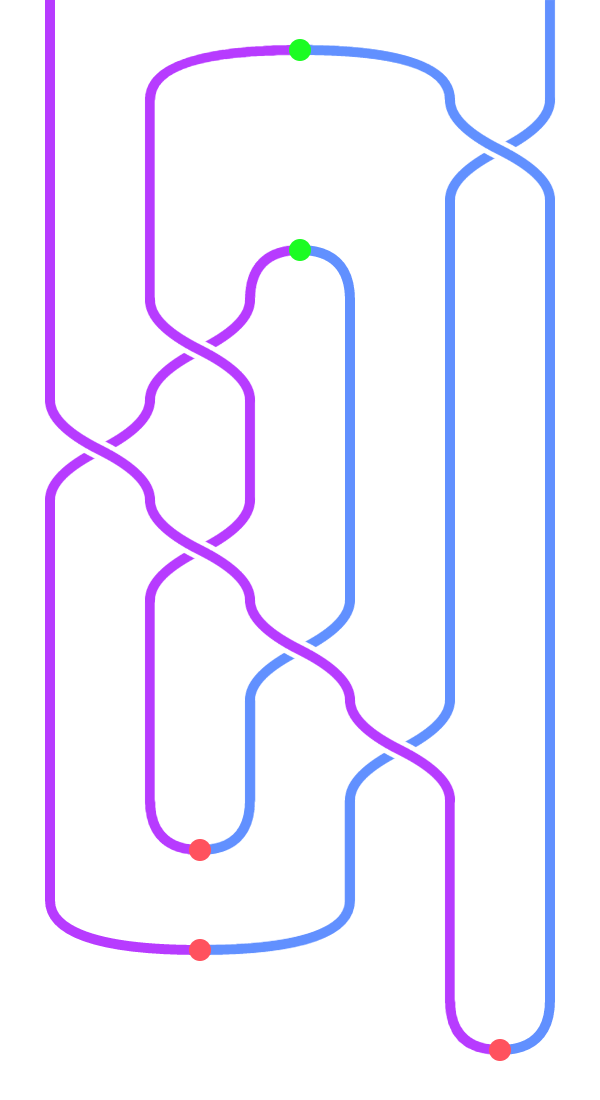}
\end{aligned}
\quad\stackrel{\II_3}{\rightarrow}\quad
\begin{aligned}
\includegraphics[scale=\imgscale]{graphics/ButterflySubB8.png}
\end{aligned} 
\end{align*}
\noindent
Note how the 4\-cell of type $\II_3$, where the red node is pulled\-through the blue wire, gets executed right away in the first sequence and then, in the second sequence, gets executed at the end after being pulled\-through the purple wire.

\paragraph{Step 97.}
Once again we draw the highlighted region of the previous slice as a movie of 2\-projected 3\-diagrams:
\def\imgscale{0.5}
\begin{align*}
\begin{aligned}
\includegraphics[scale=\imgscale]{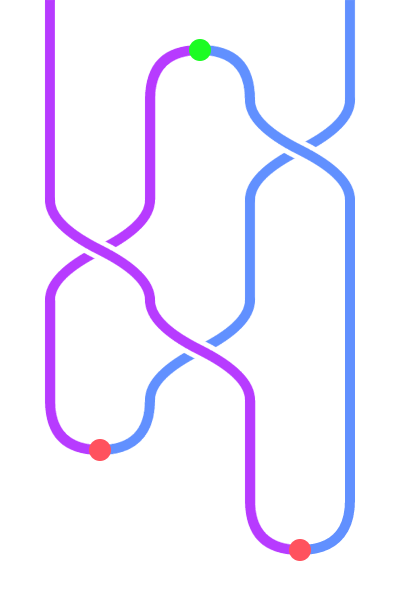}
\end{aligned}
&\quad\stackrel{{\zeta}^{-1}}{\rightarrow}\quad&
\begin{aligned}
\includegraphics[scale=\imgscale]{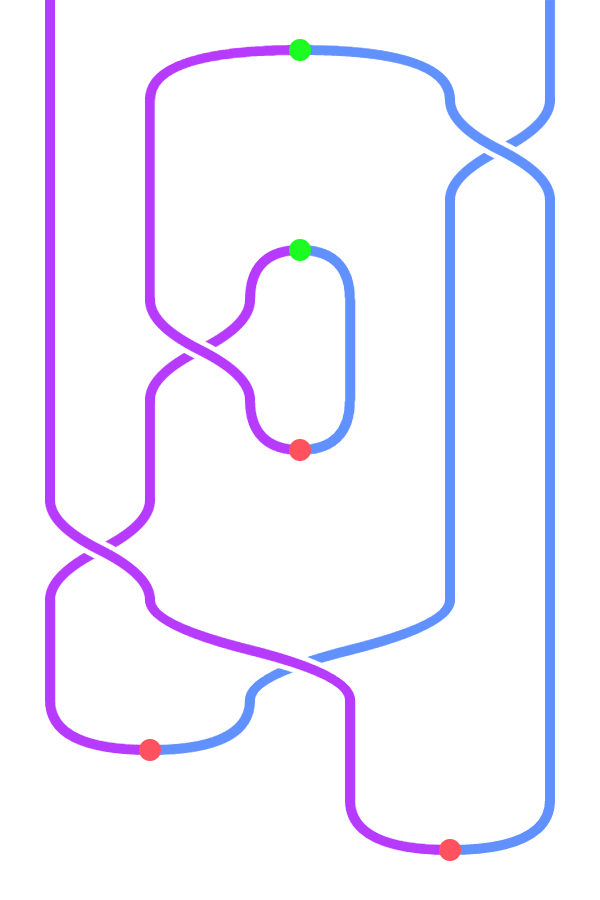}
\end{aligned}
&\quad\stackrel{\I_3}{\rightarrow}\quad&
\begin{aligned}
\includegraphics[scale=\imgscale]{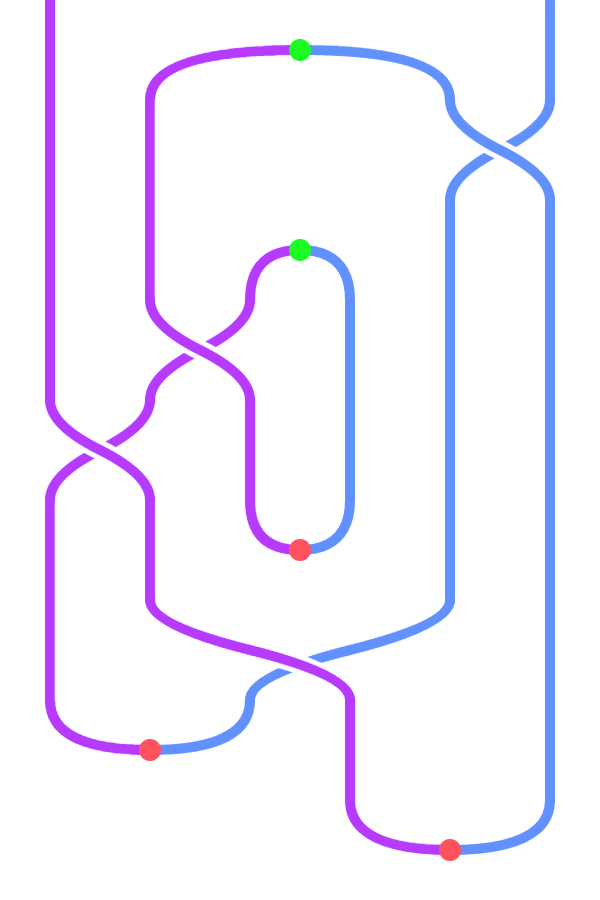}
\end{aligned}
&\quad\stackrel{\II_3}{\rightarrow}\quad&
\begin{aligned}
\includegraphics[scale=\imgscale]{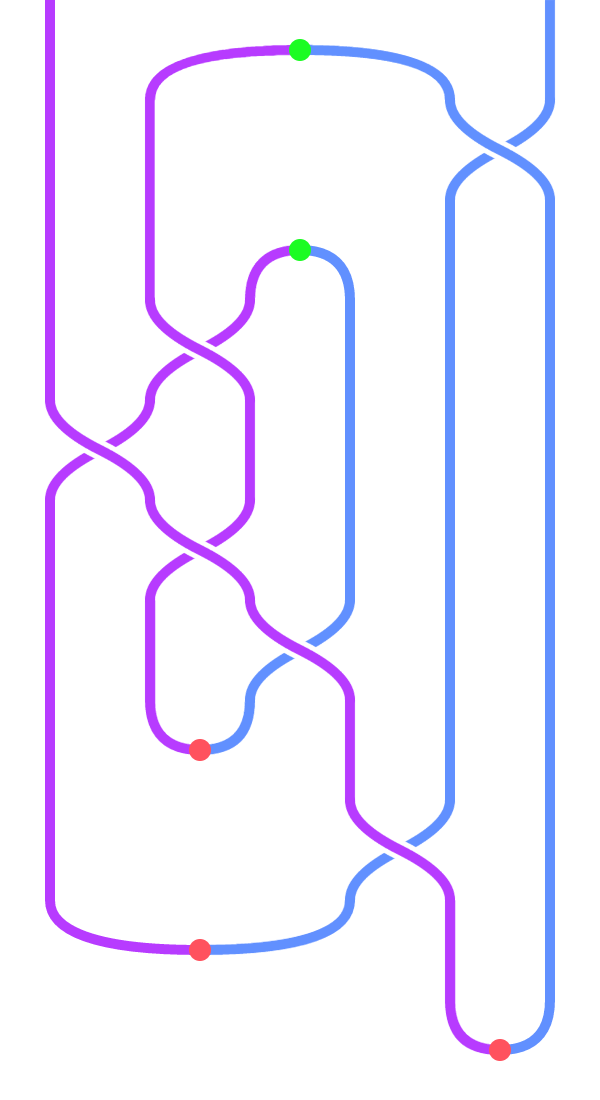}
\end{aligned}
&\quad\stackrel{\II_3}{\rightarrow}\quad& \\
\begin{aligned}
\includegraphics[scale=\imgscale]{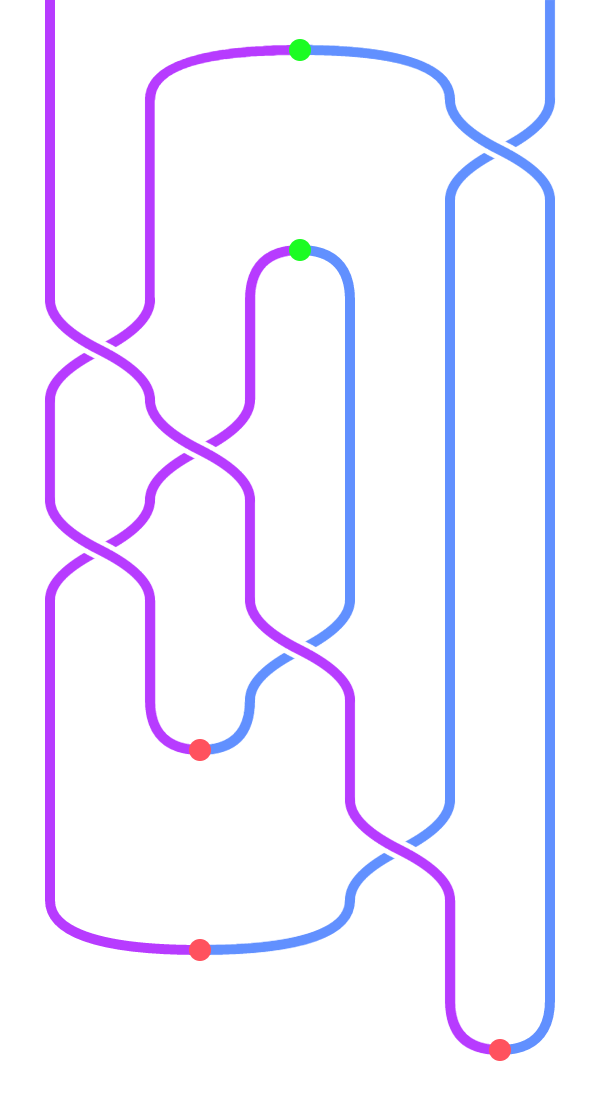}
\end{aligned}
&\quad\stackrel{\I_3}{\rightarrow}\quad&
\begin{aligned}
\includegraphics[scale=\imgscale]{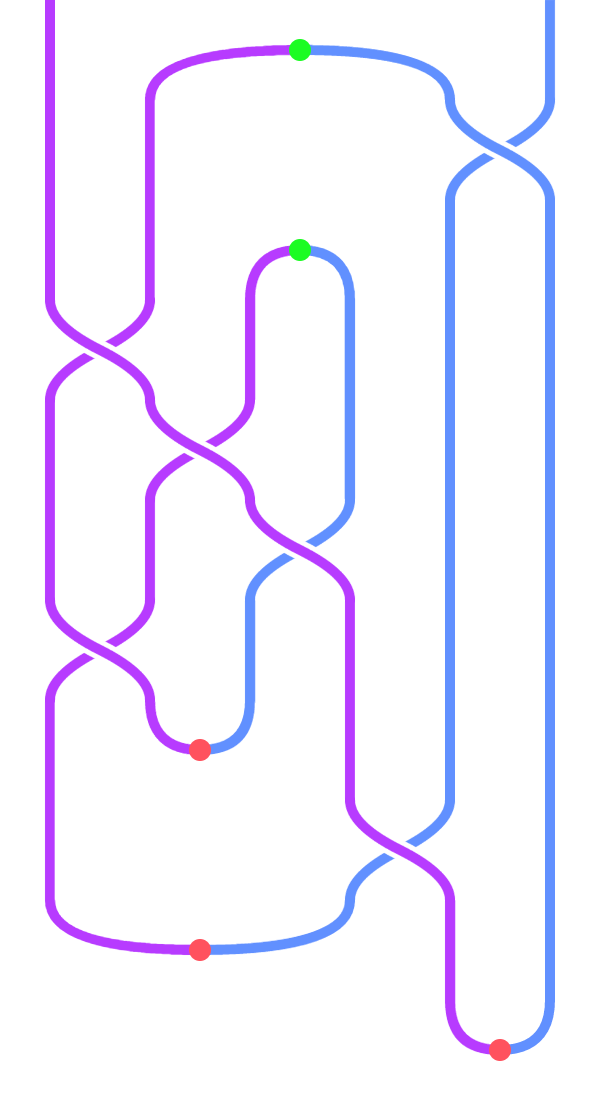}
\end{aligned}
&\quad\stackrel{\I_3}{\rightarrow}\quad&
\begin{aligned}
\includegraphics[scale=\imgscale]{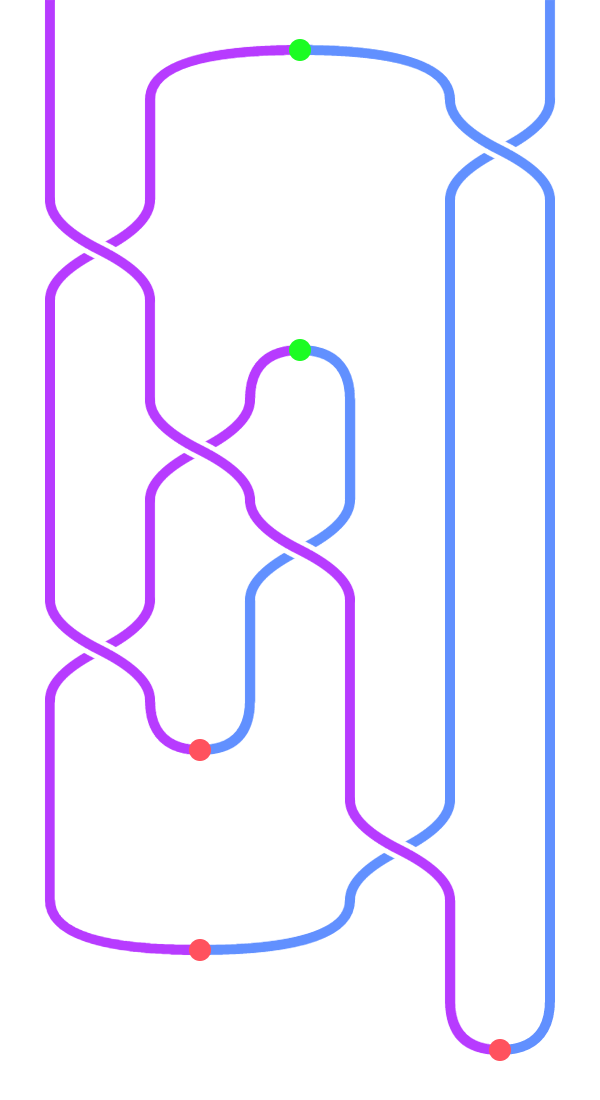}
\end{aligned}
&\quad\stackrel{\II_3}{\rightarrow}\quad&
\begin{aligned}
\includegraphics[scale=\imgscale]{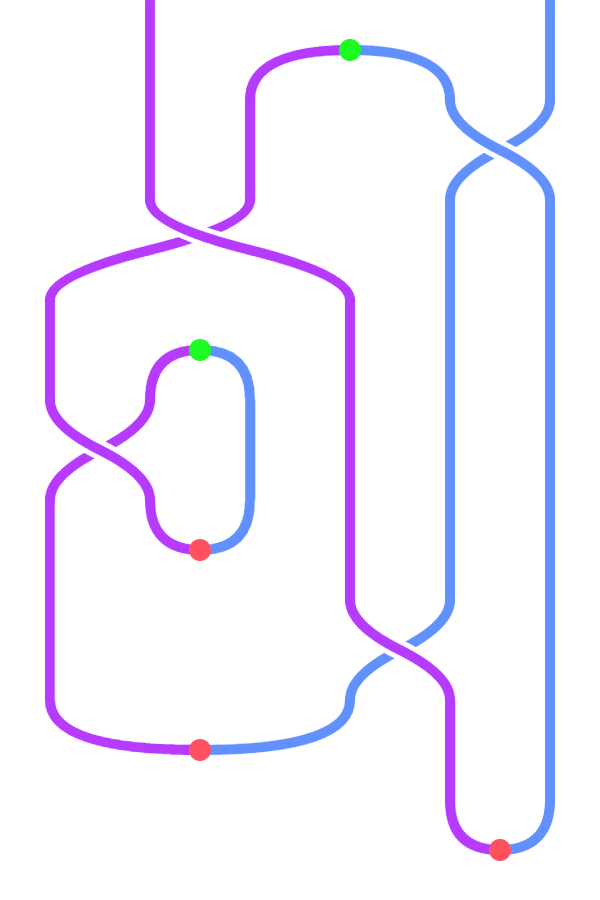}
\end{aligned} 
\end{align*}

\noindent
This is rewritten into the following:
\begin{align*}
\begin{aligned}
\includegraphics[scale=\imgscale]{graphics/ButterflySubC1.png}
\end{aligned}
\quad\stackrel{{\zeta}^{-1}}{\rightarrow}\quad
\begin{aligned}
\includegraphics[scale=\imgscale]{graphics/ButterflySubC8.png}
\end{aligned}
\end{align*}
\noindent
This is another instance of an interchanger of type $\III'$. Here, the 4\-cell ${\zeta}^{-1}$ gets executed right away in the first sequence and then, in the second sequence gets executed after being pulled-through the purple wire.

\subsection{Projected view}
\label{sec:butterflyproject3}

\newcounter{butterflycounter}
\def\butterflyscale{0.9}
\newcommand\butterflyimage[1]{\ensuremath{\begin{aligned}\includegraphics[scale=\butterflyscale]{graphics/butterfly_30px_#1}\hspace{0pt}\end{aligned}}\xspace}
\newcommand\butterflystep[2]{\stepcounter{butterflycounter}\mbox{\ensuremath{\substack{#1 \\ \to \\ \arabic{butterflycounter}}\hspace{-0.5pt}{\butterflyimage{#2}}}}\xspace}

\begin{figure}
$$\includegraphics[scale=0.75]{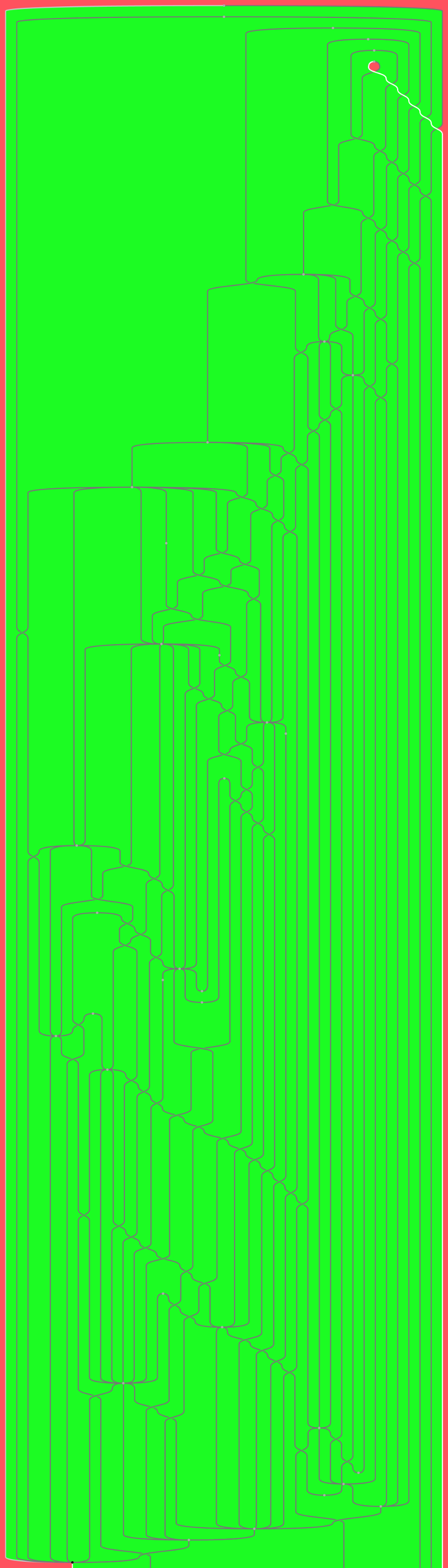}$$
\caption{\label{fig:wholeproof}The 2-projection of the proof 5-diagram.}
\end{figure}

Figure~\ref{fig:wholeproof} gives a 2\-projection of the proof 5\-diagram.

\subsection{Movie view}
\label{sec:butterflyproject2}

The following is a movie of 2\-projection of 4\-diagram slices of the proof, in the manner of the `Movies' subsection within Section~\ref{secgraphicalformalism}. Each frame is only a partial description of that proof step, and to be fully specified would require a movie of movies of 2\-dimensional images.

\raggedright
\butterflyimage{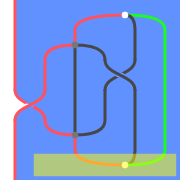}
\butterflystep{\zeta}{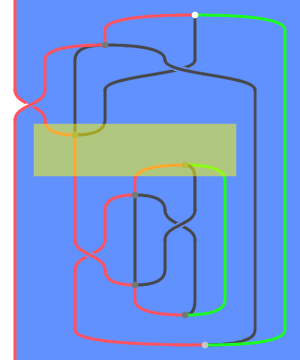}
\butterflystep{\I_4}{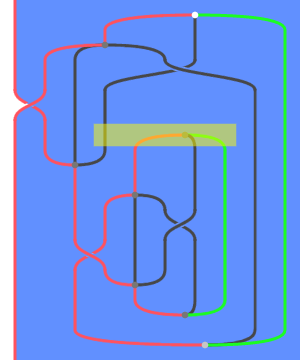}
\butterflystep{\II_4}{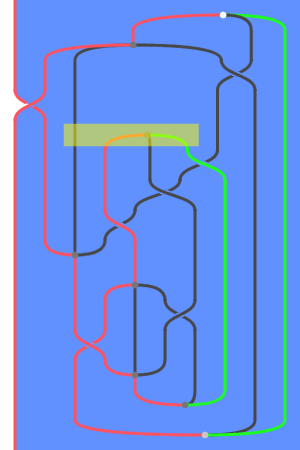}
\butterflystep{\II_4}{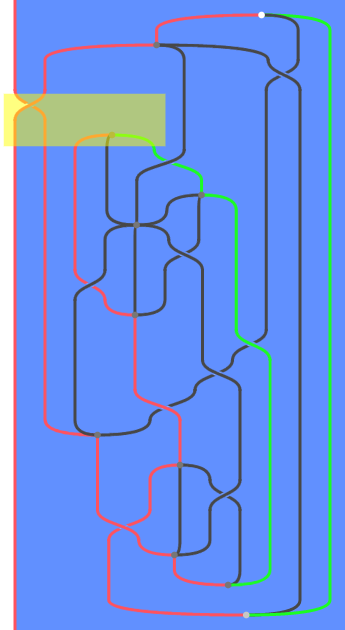}
\butterflystep{\I_4}{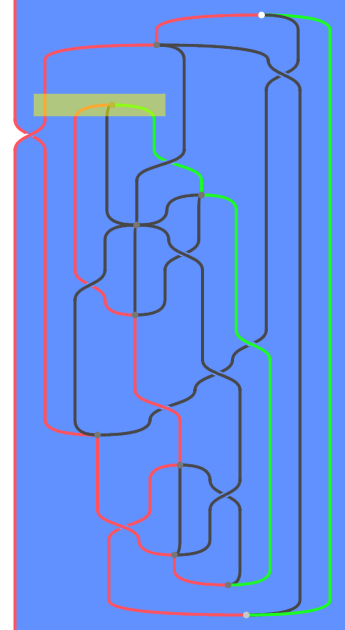}
\butterflystep{\II_4}{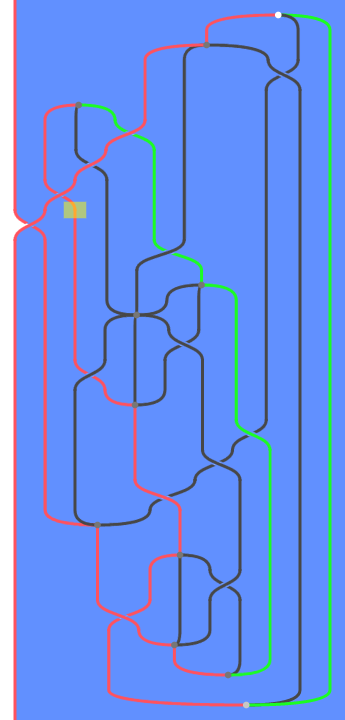}
\butterflystep{\II_3'}{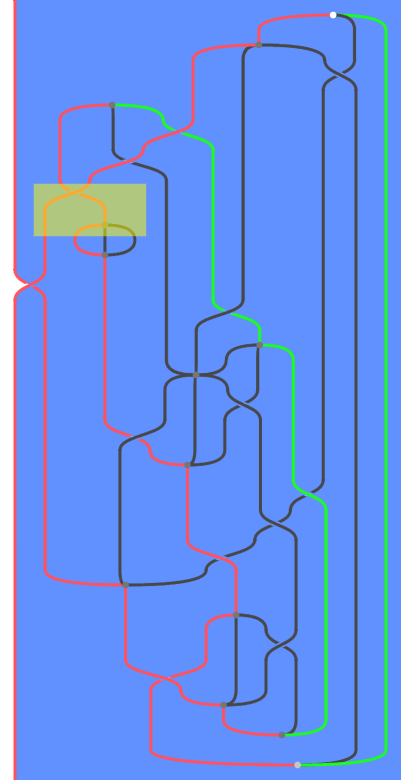}
\butterflystep{\II_4}{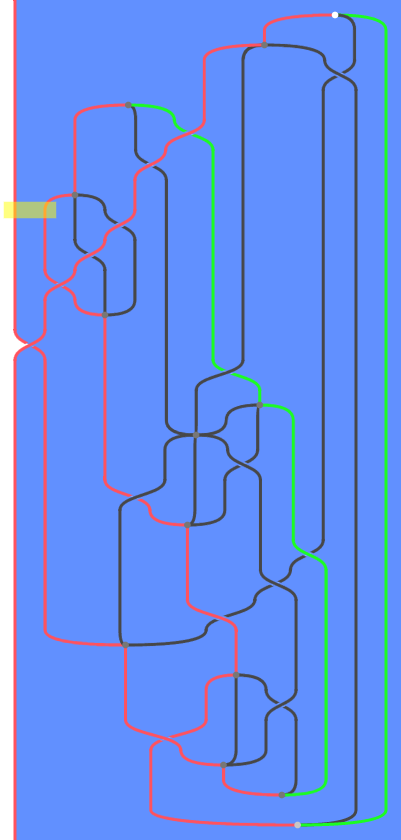}
\butterflystep{\I_3'}{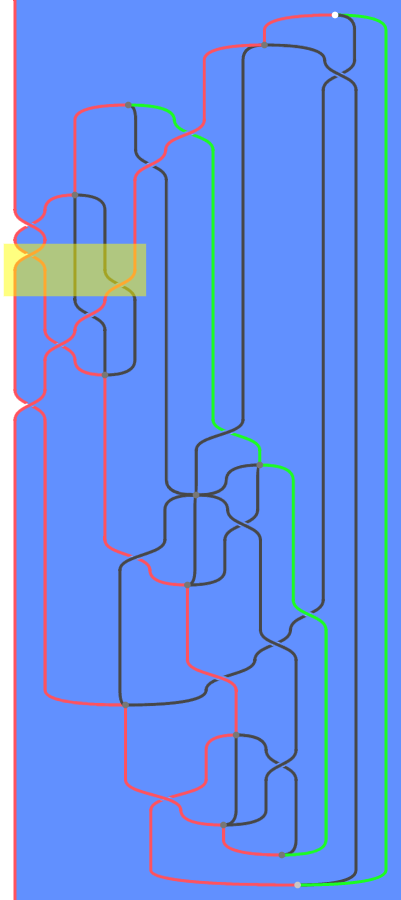}
\butterflystep{\I_4}{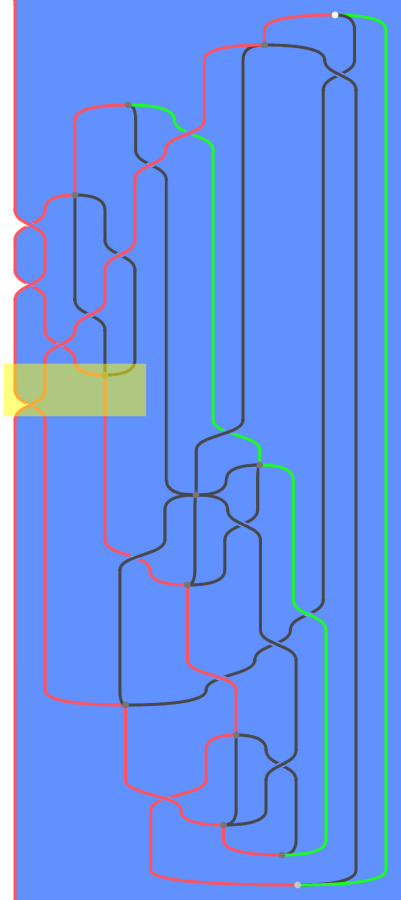}
\butterflystep{\I_4}{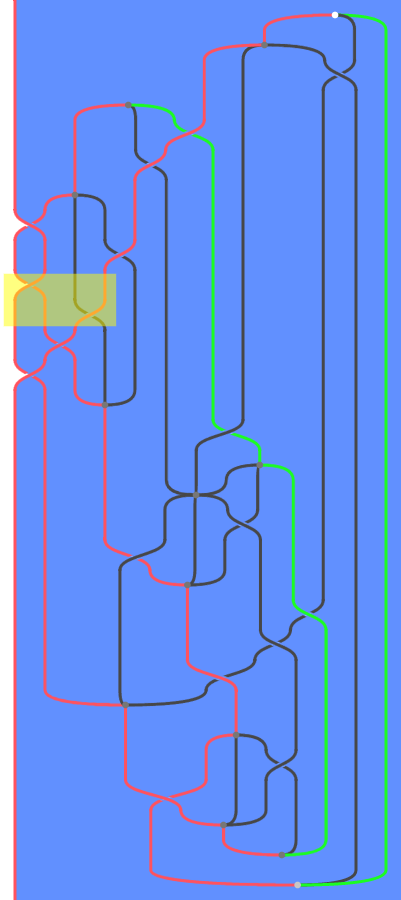}
\butterflystep{\I_4}{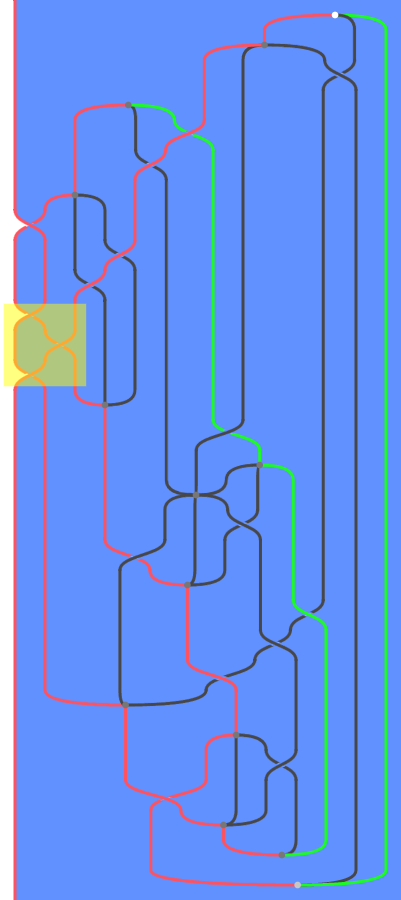}
\butterflystep{\II_4}{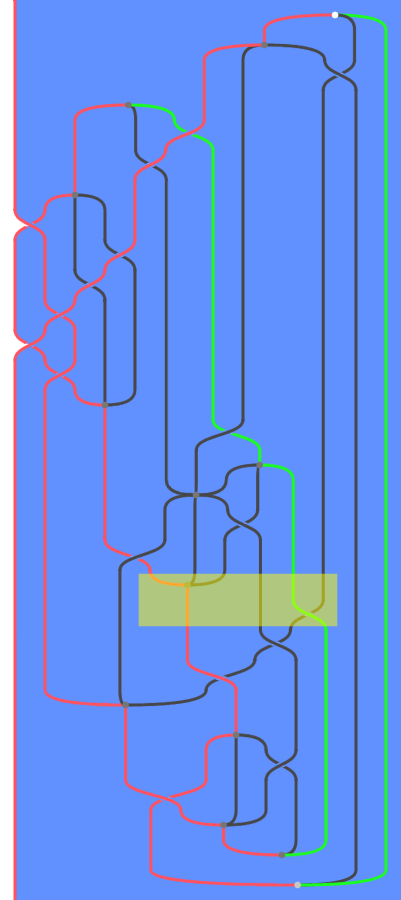}
\butterflystep{\I_4}{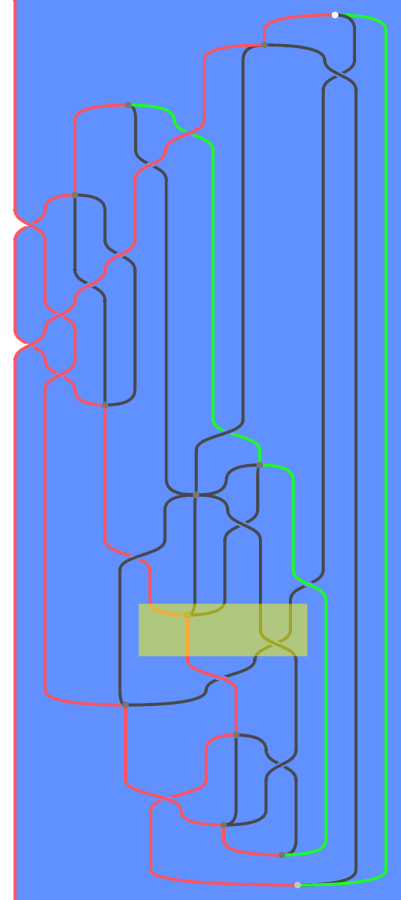}
\butterflystep{\I_4}{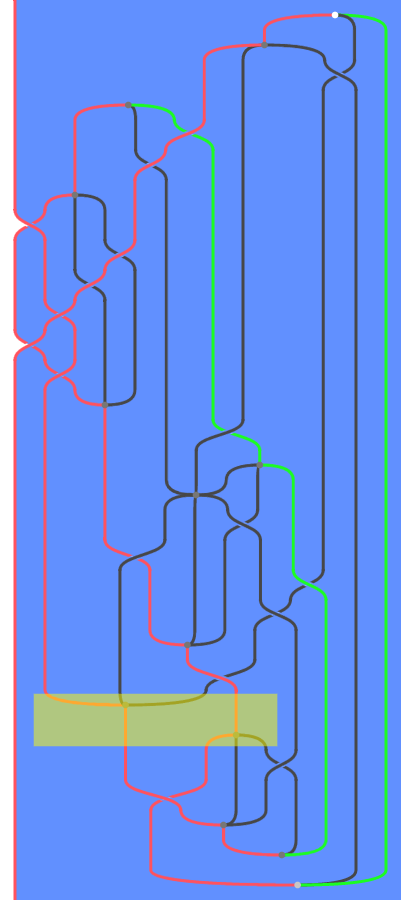}
\butterflystep{\I_4}{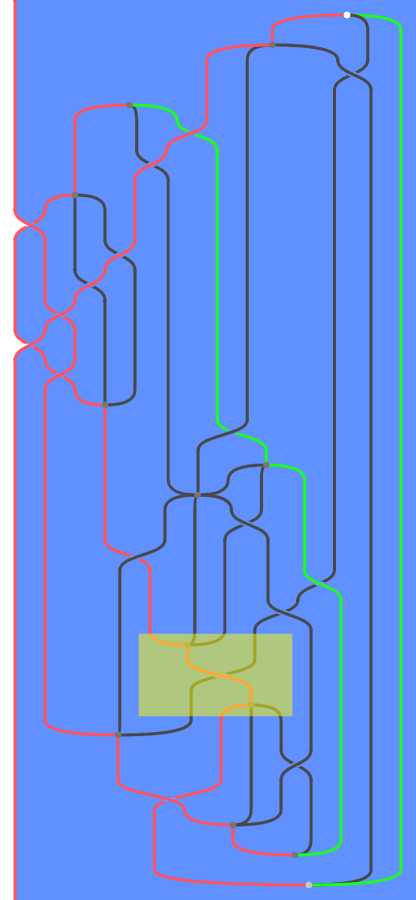}
\butterflystep{\III_4}{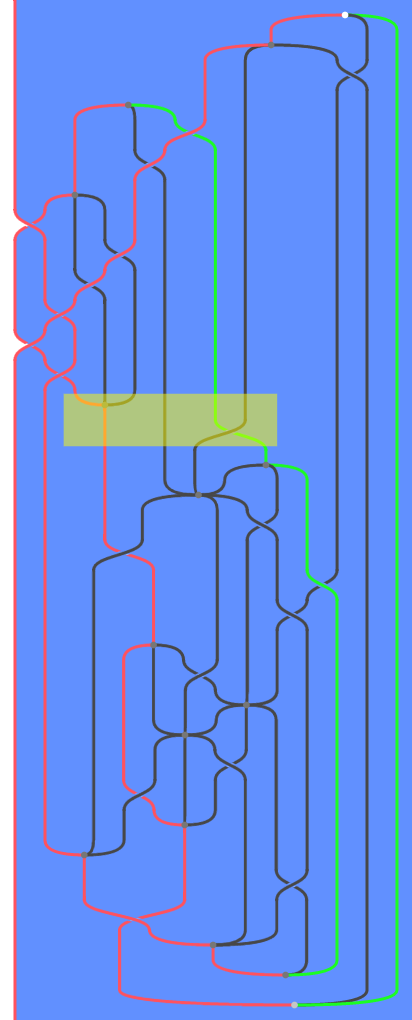}
\butterflystep{\I_4}{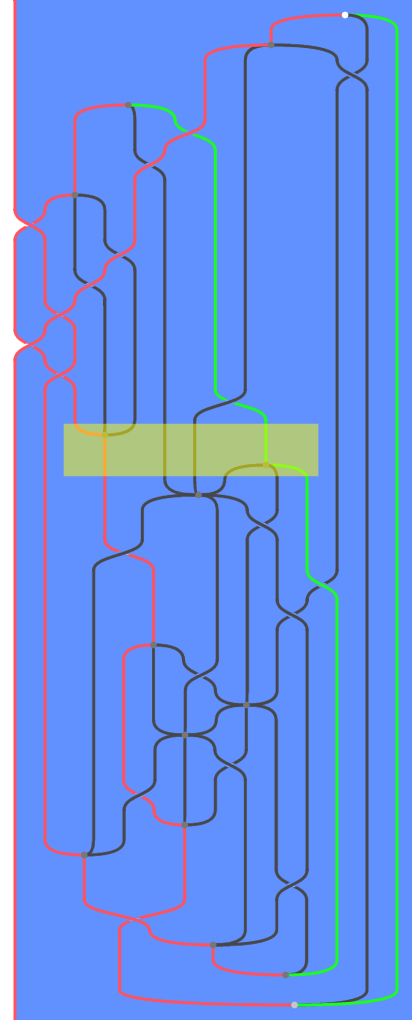}
\butterflystep{\I_4}{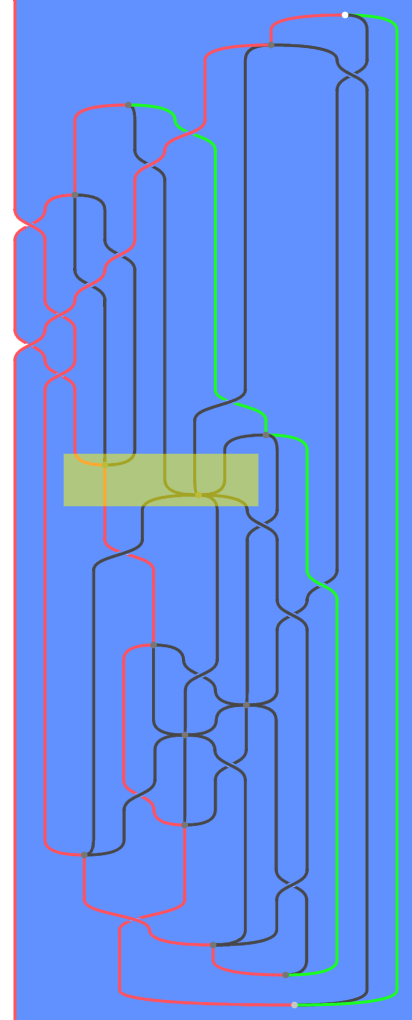}
\butterflystep{\I_4}{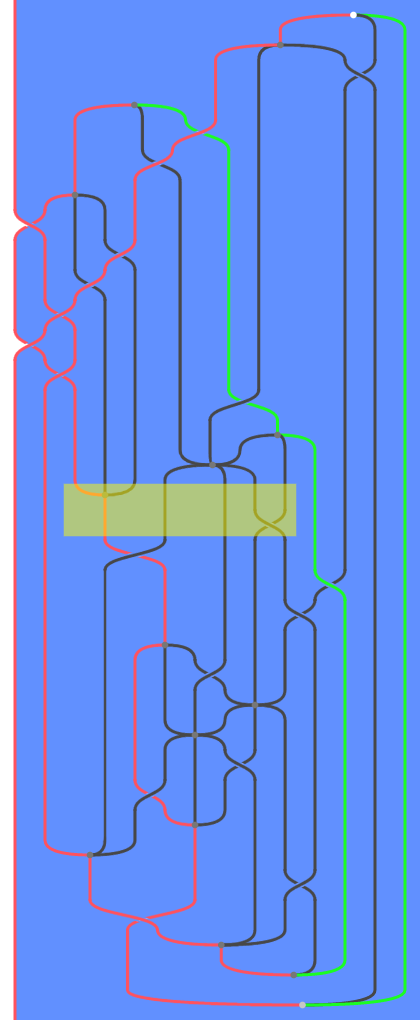}
\butterflystep{\I_4}{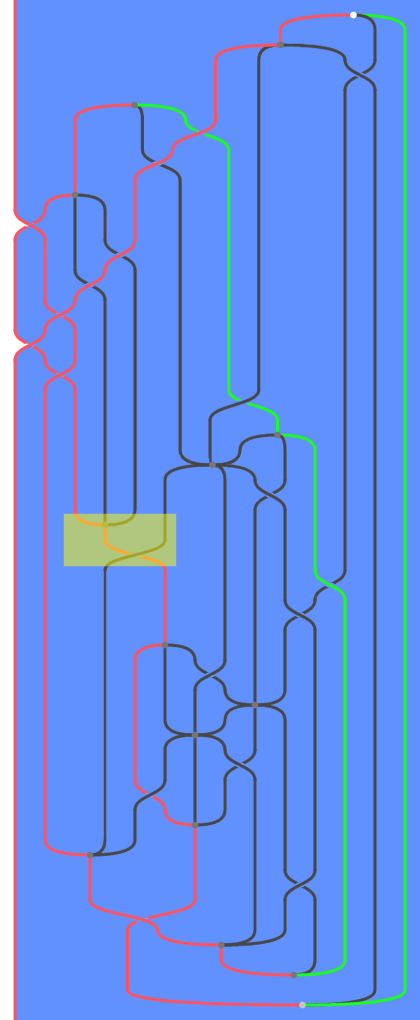}
\butterflystep{\I_4}{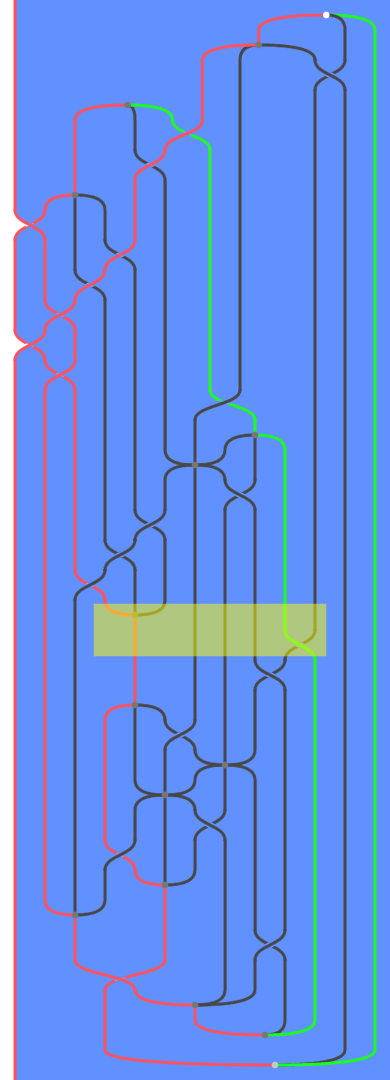}
\butterflystep{\I_4}{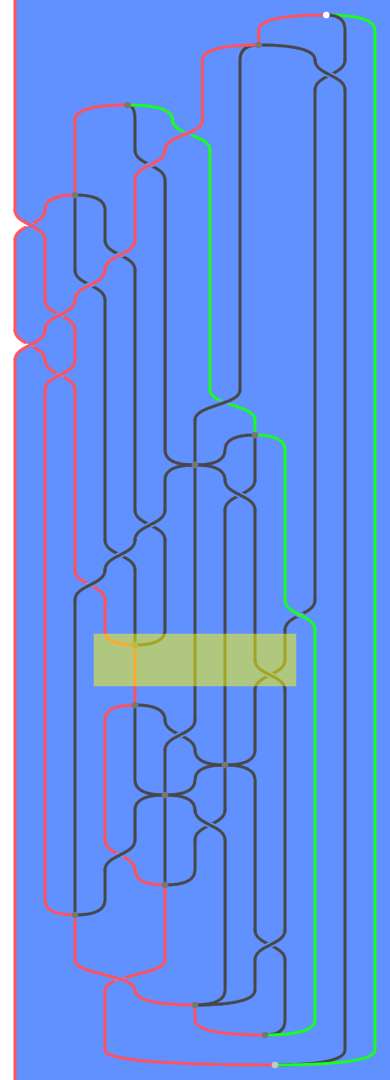}
\butterflystep{\I_4}{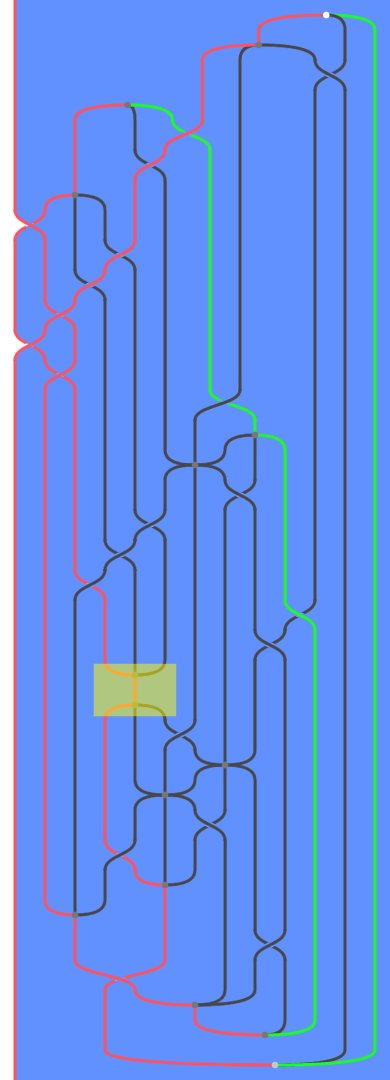}
\butterflystep{\I_4}{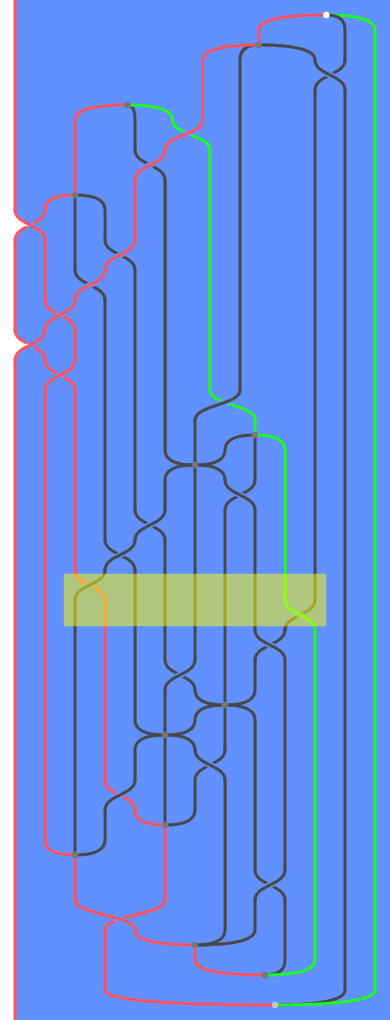}
\butterflystep{\I_4}{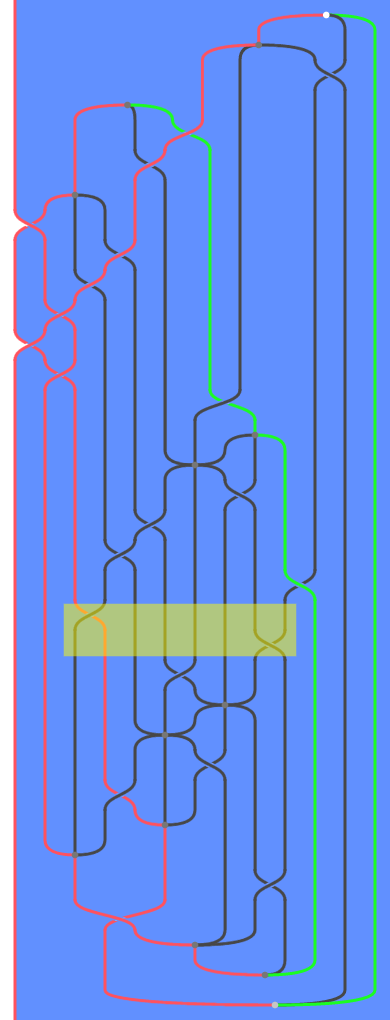}
\butterflystep{\I_4}{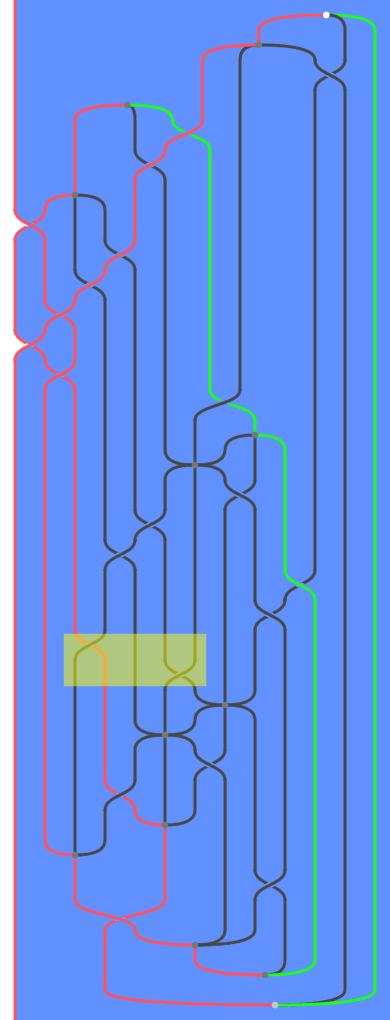}
\butterflystep{\I_4}{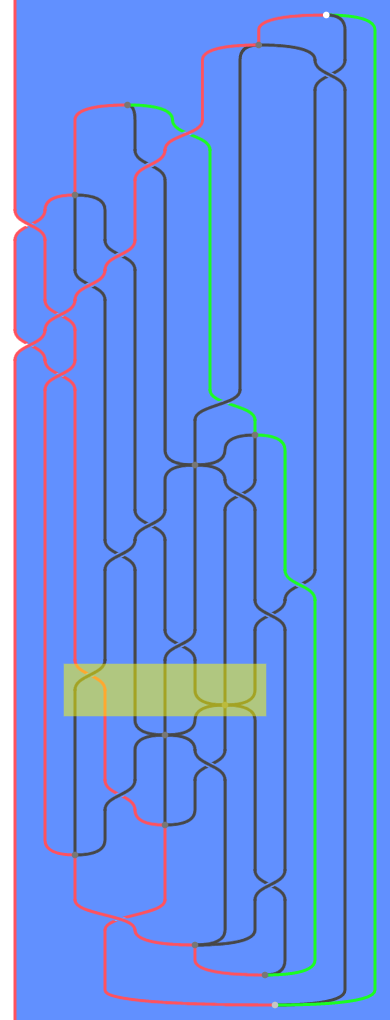}
\butterflystep{\I_4}{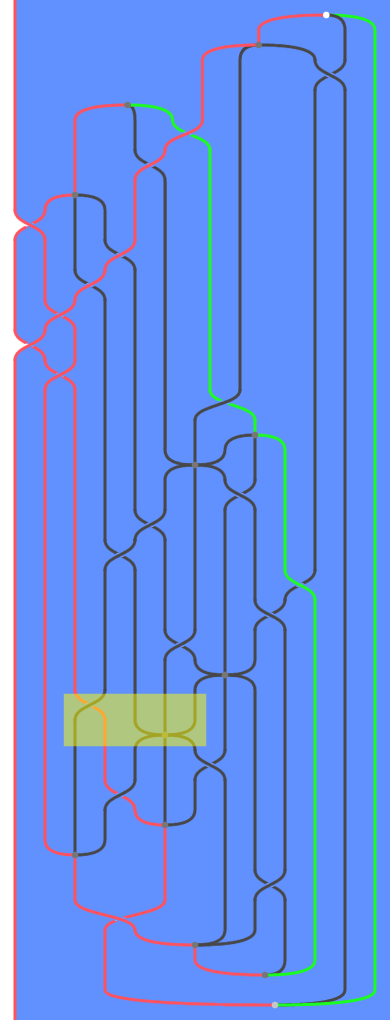}
\butterflystep{\I_4}{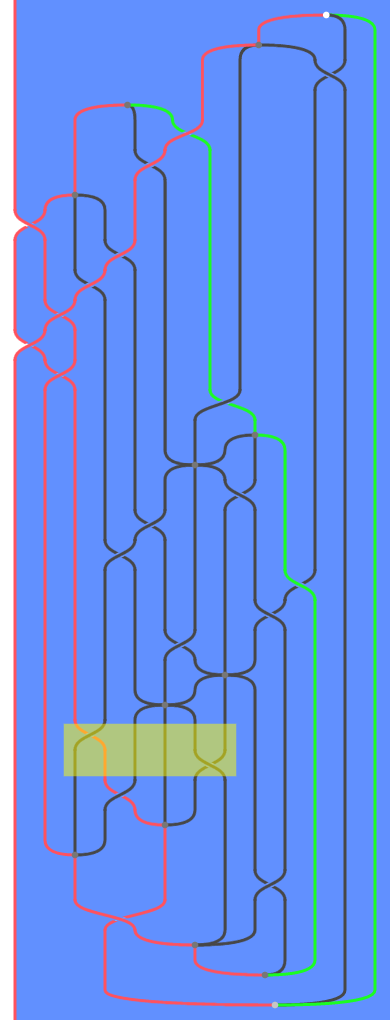}
\butterflystep{\I_4}{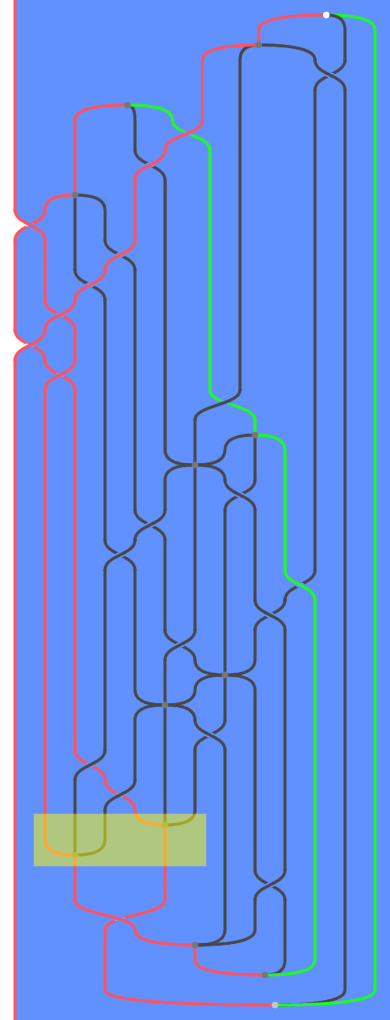}
\butterflystep{\I_4}{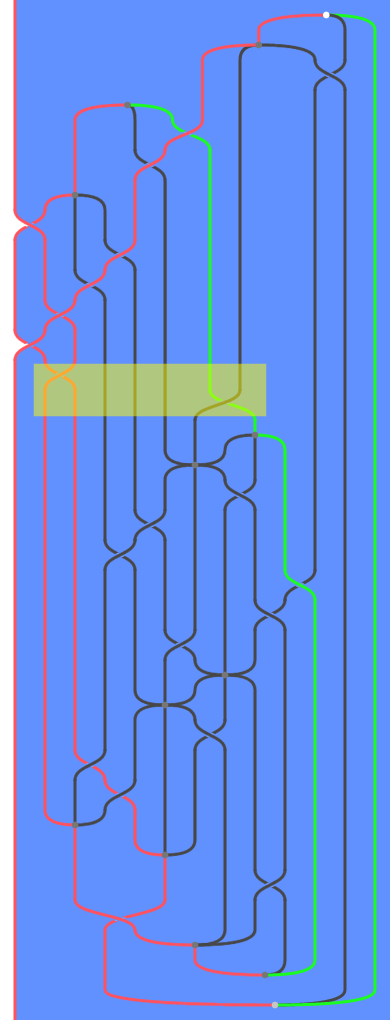}
\butterflystep{\I_4}{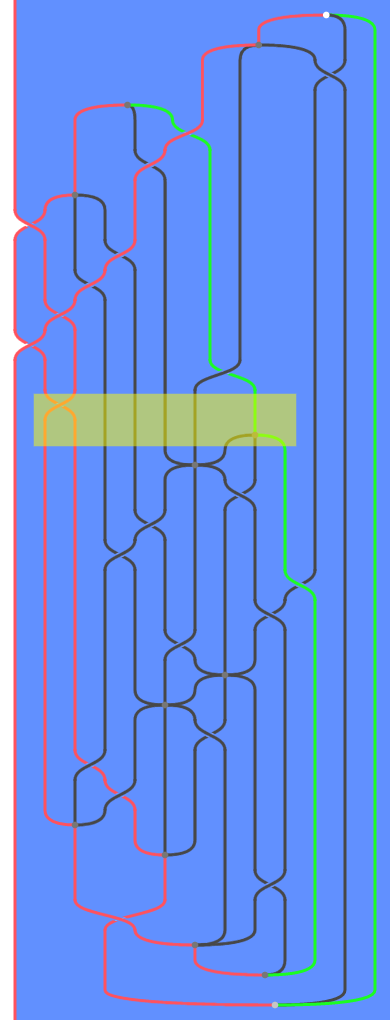}
\butterflystep{\I_4}{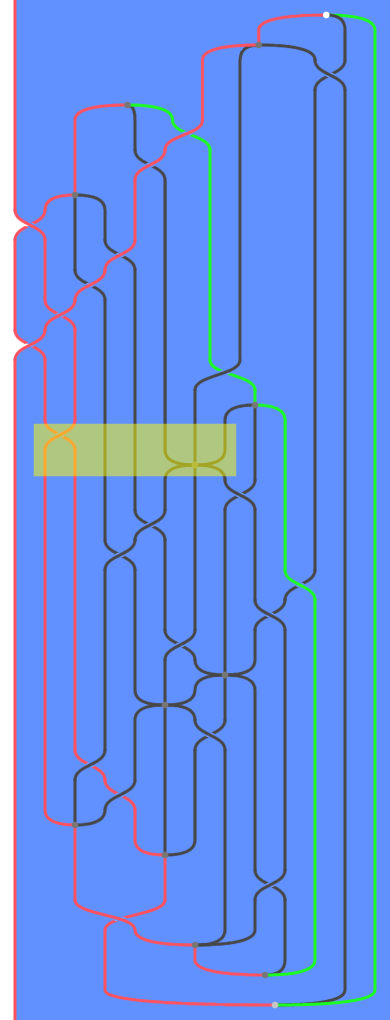}
\butterflystep{\I_4}{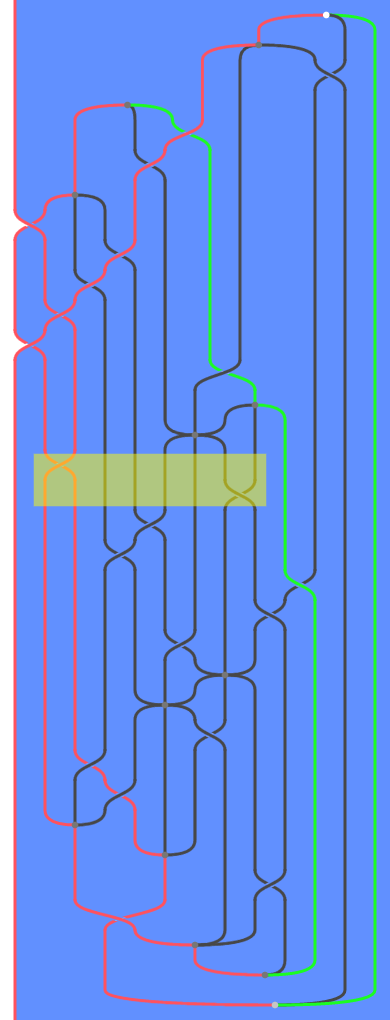}
\butterflystep{\I_4}{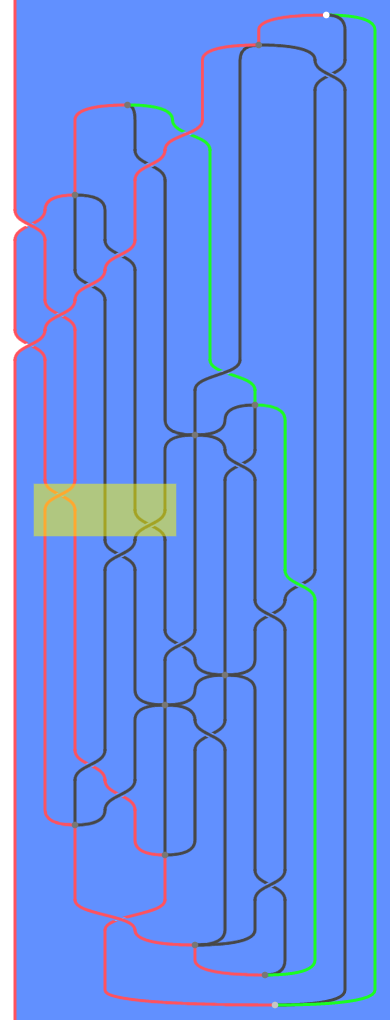}
\butterflystep{\I_4}{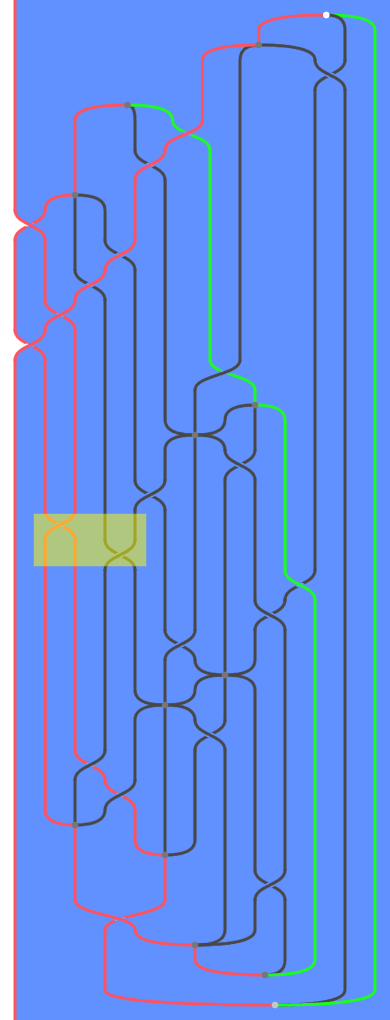}
\butterflystep{\I_4}{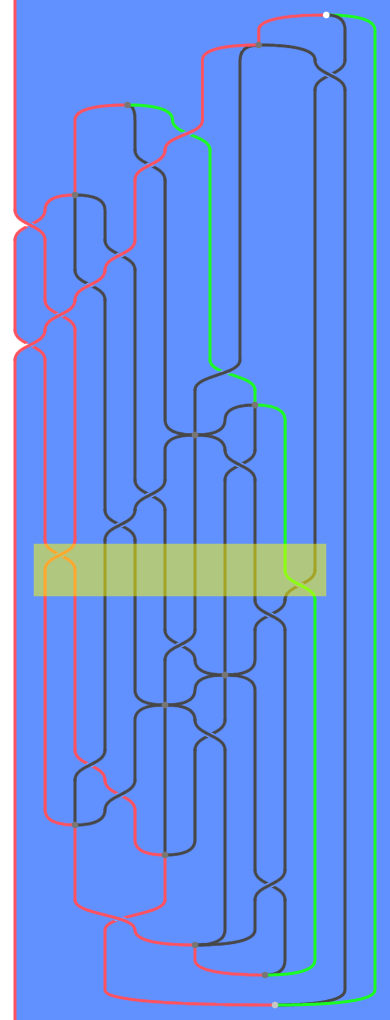}
\butterflystep{\I_4}{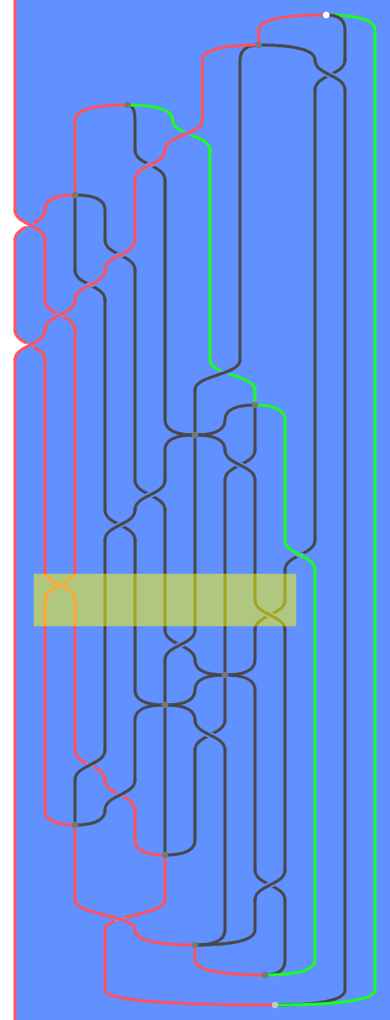}
\butterflystep{\I_4}{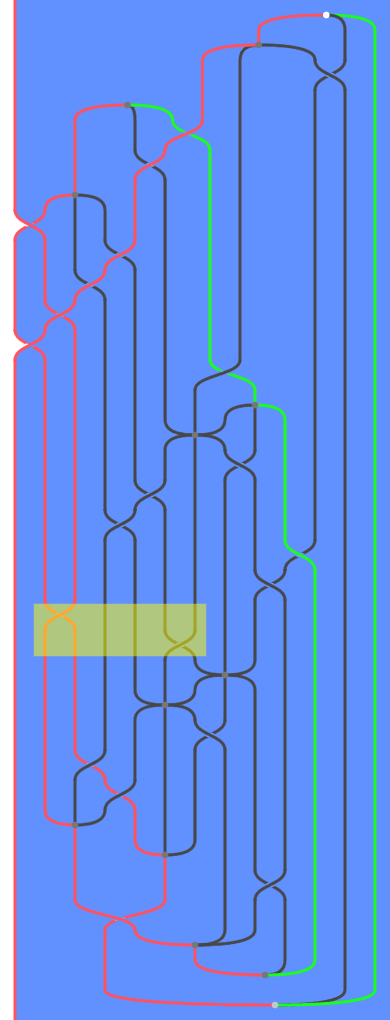}
\butterflystep{\I_4}{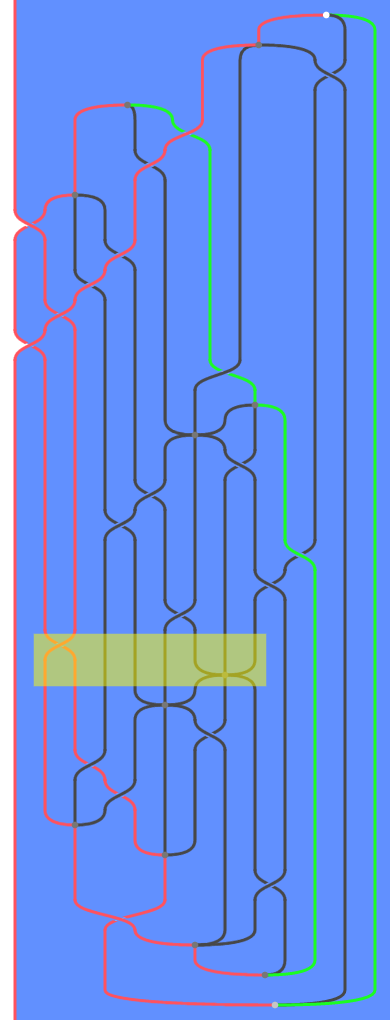}
\butterflystep{\I_4}{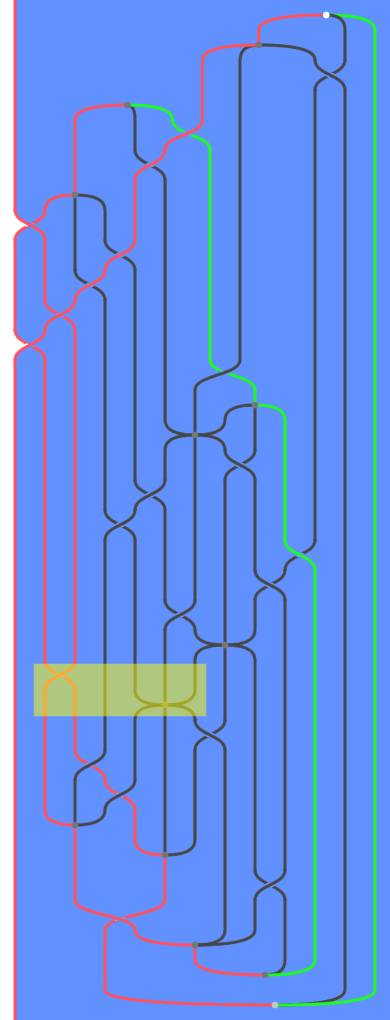}
\butterflystep{\I_4}{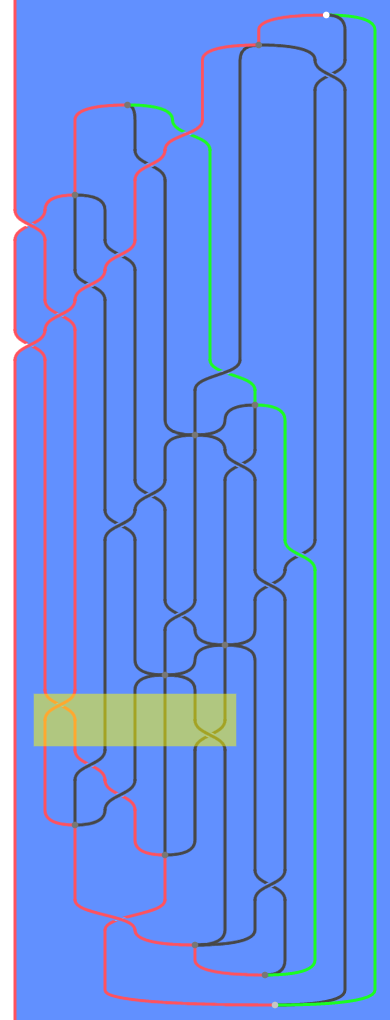}
\butterflystep{\I_4}{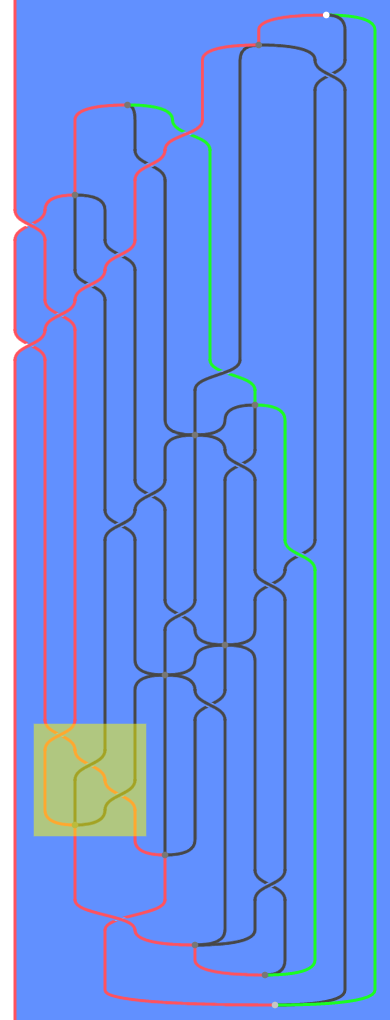}
\butterflystep{\I_4}{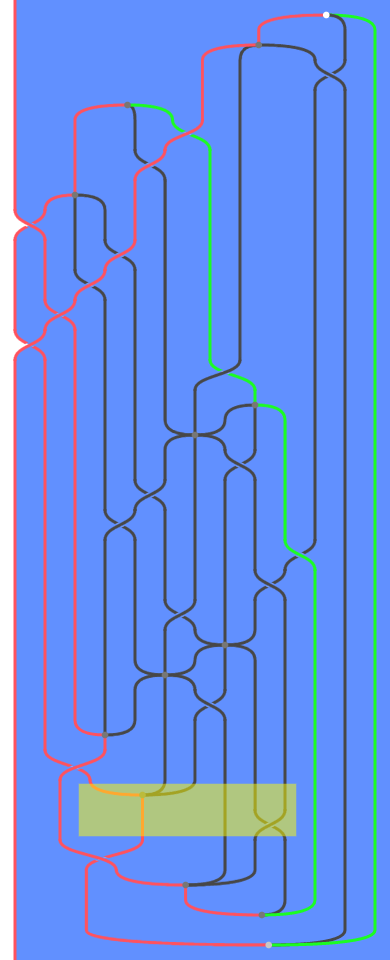}
\butterflystep{\I_4}{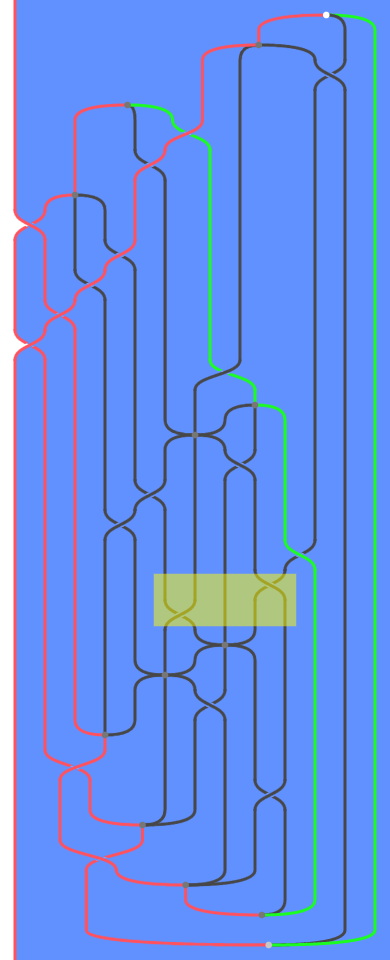}
\butterflystep{\I_4}{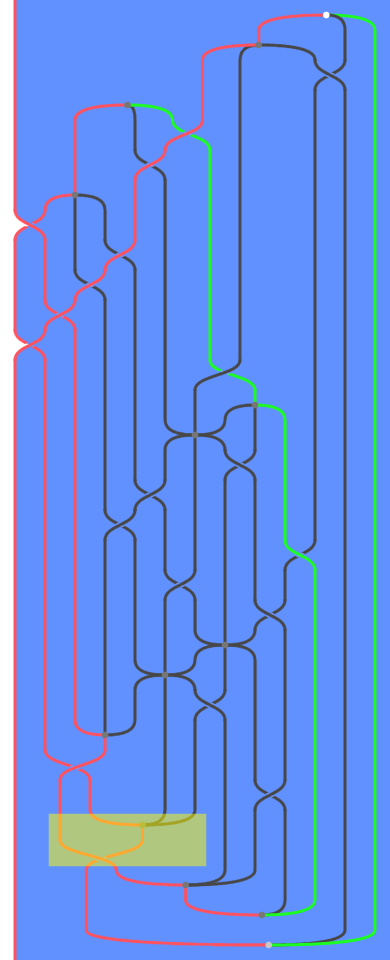}
\butterflystep{\I_4}{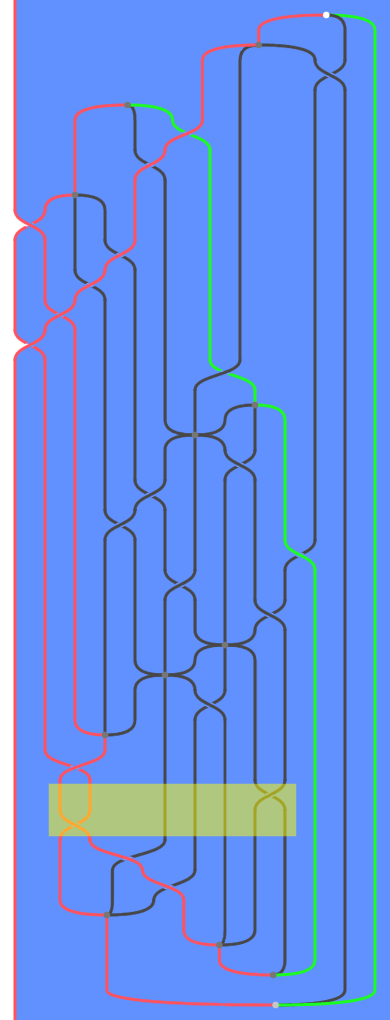}
\butterflystep{\I_4}{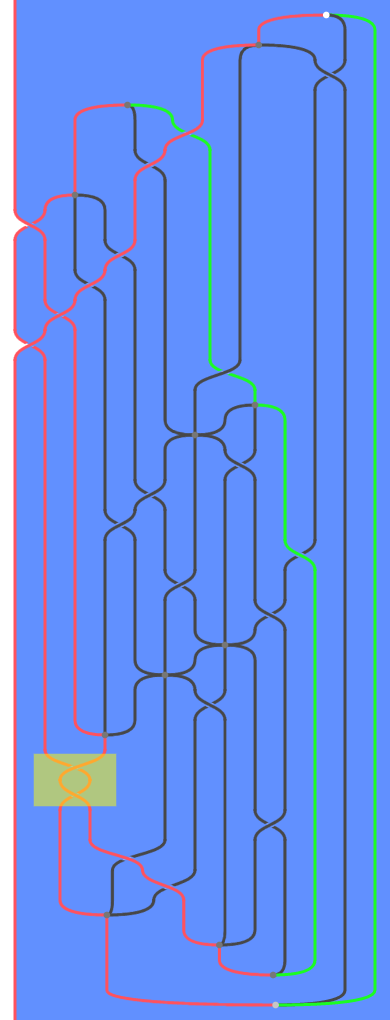}
\butterflystep{\I_3'}{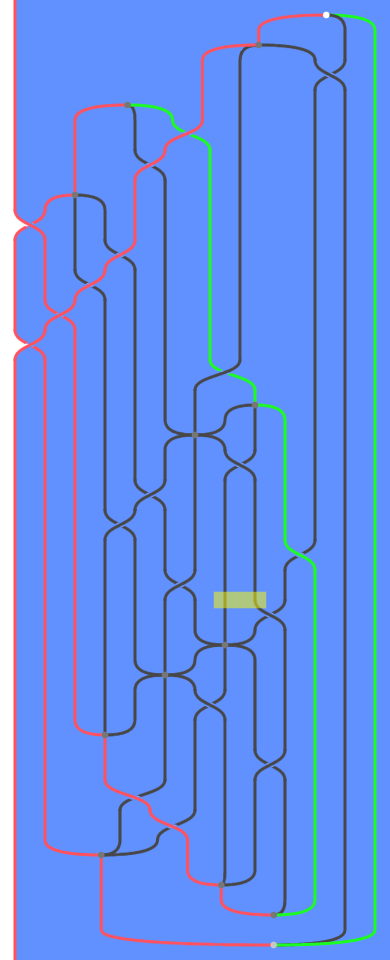}
\butterflystep{\I_3'}{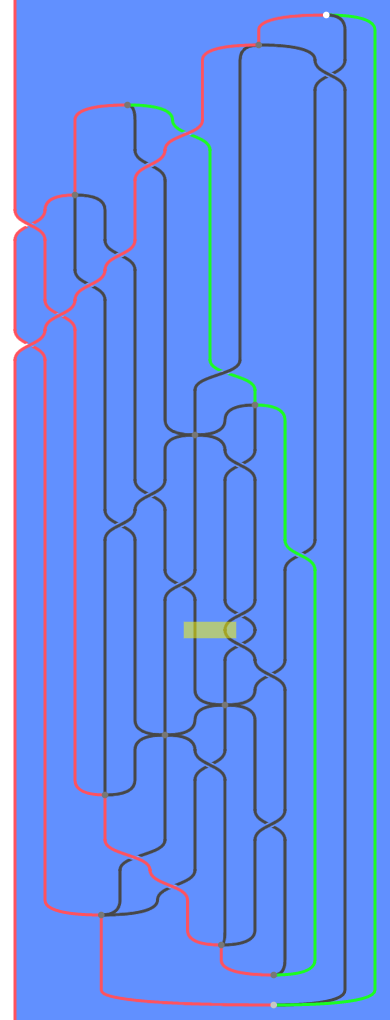}
\butterflystep{\I_3'}{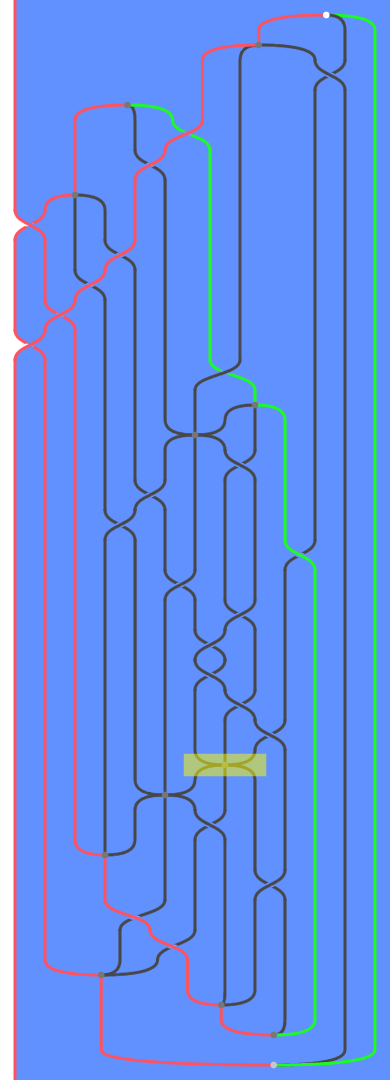}
\butterflystep{\VI_4}{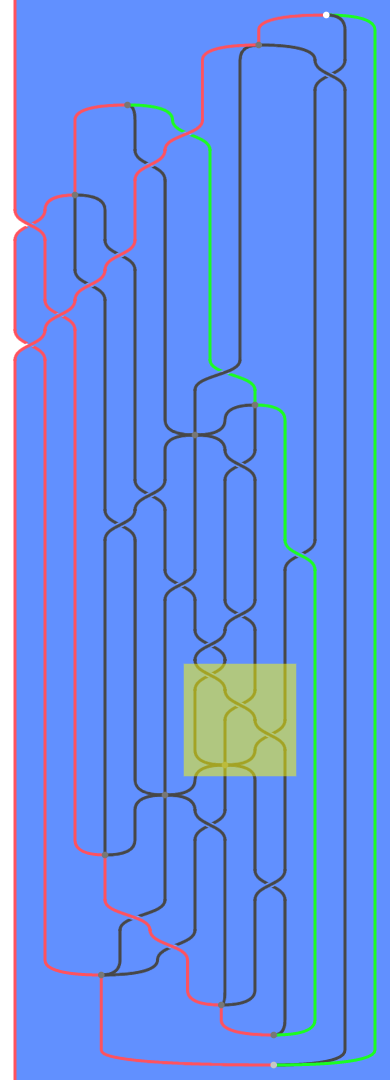}
\butterflystep{\II_4}{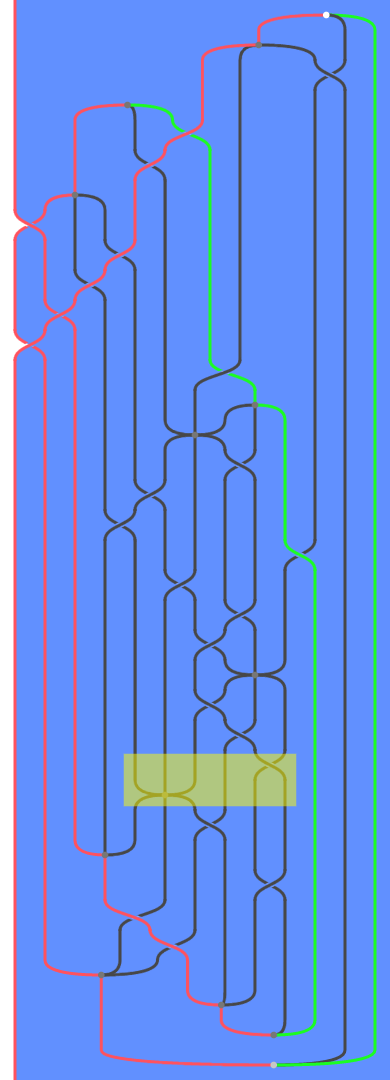}
\butterflystep{\I_4}{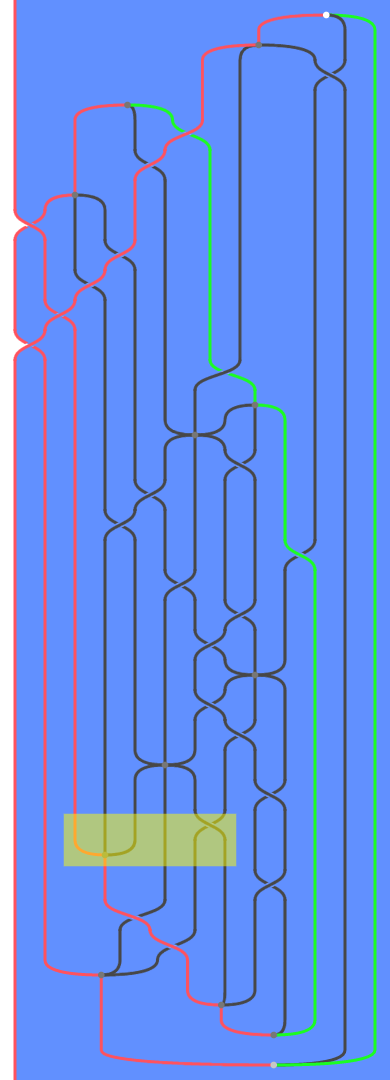}
\butterflystep{\I_4}{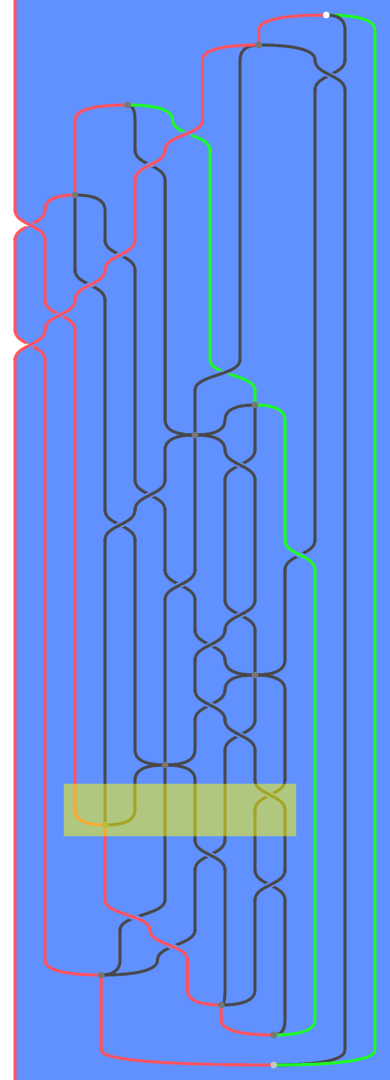}
\butterflystep{\I_4}{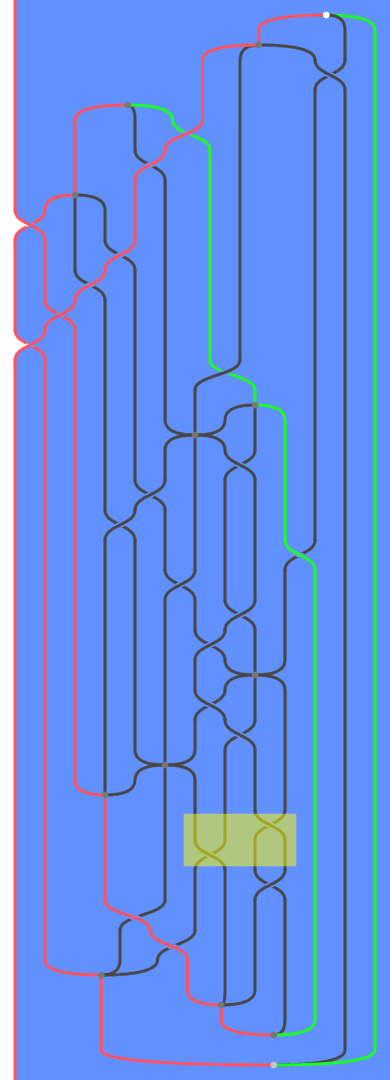}
\butterflystep{\I_4}{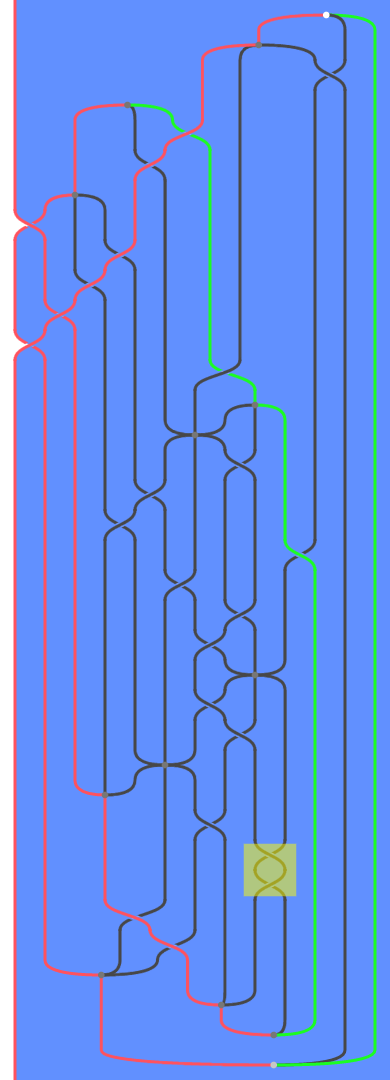}
\butterflystep{\I_3'}{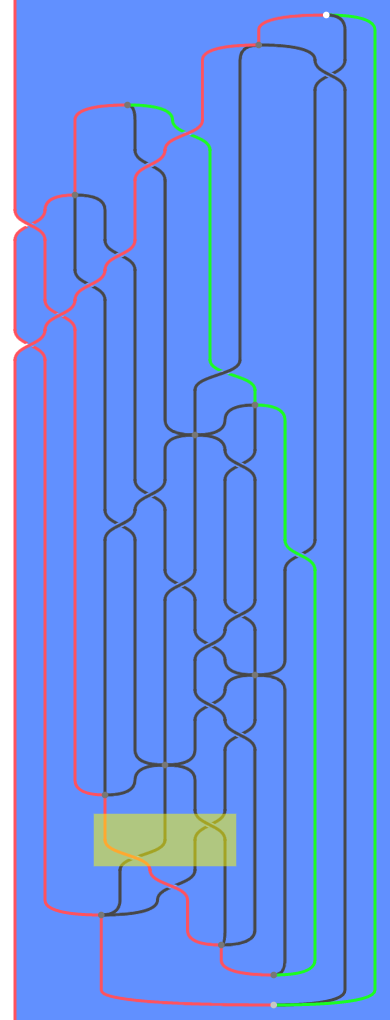}
\butterflystep{\I_4}{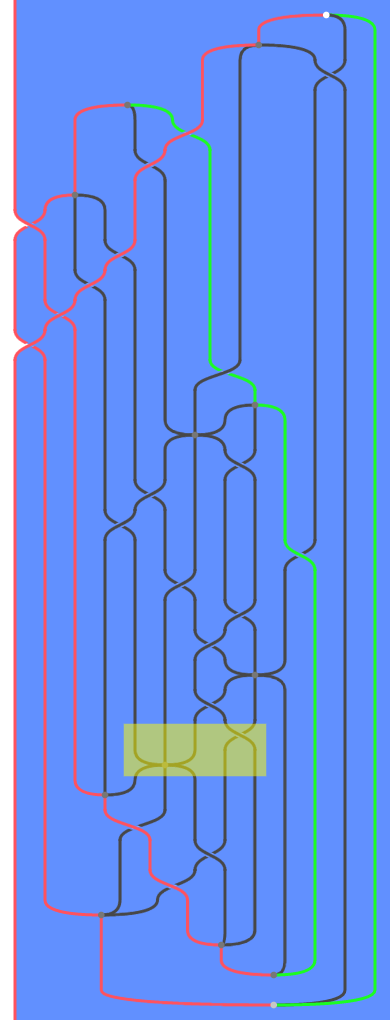}
\butterflystep{\I_4}{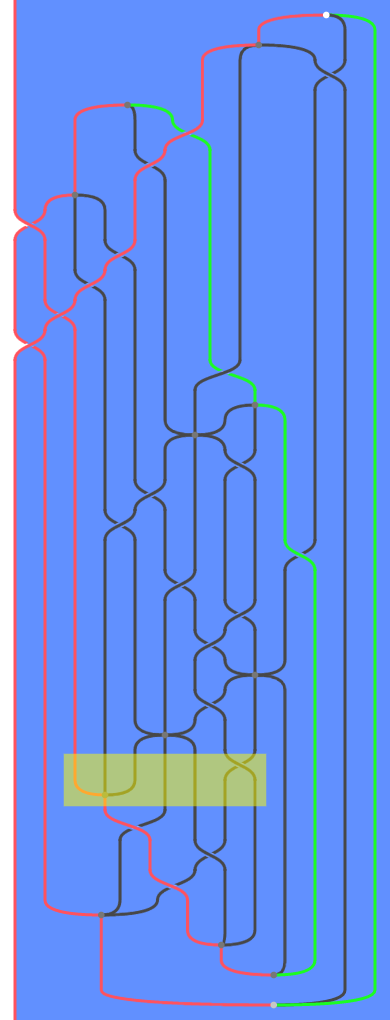}
\butterflystep{\I_4}{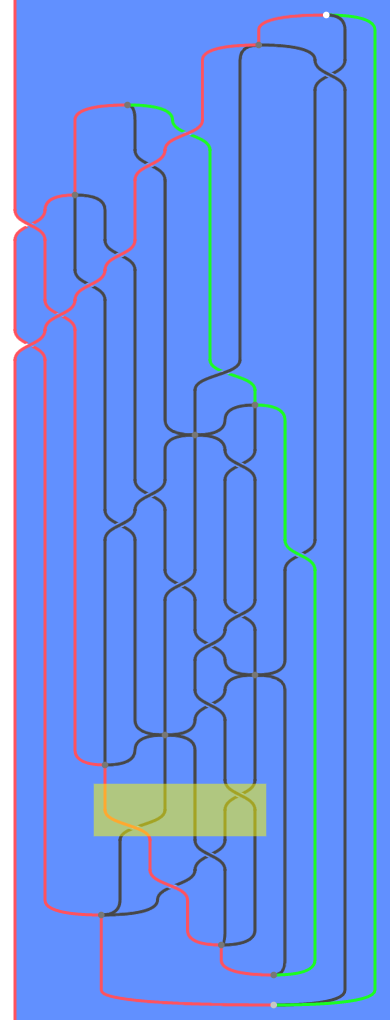}
\butterflystep{\I_4}{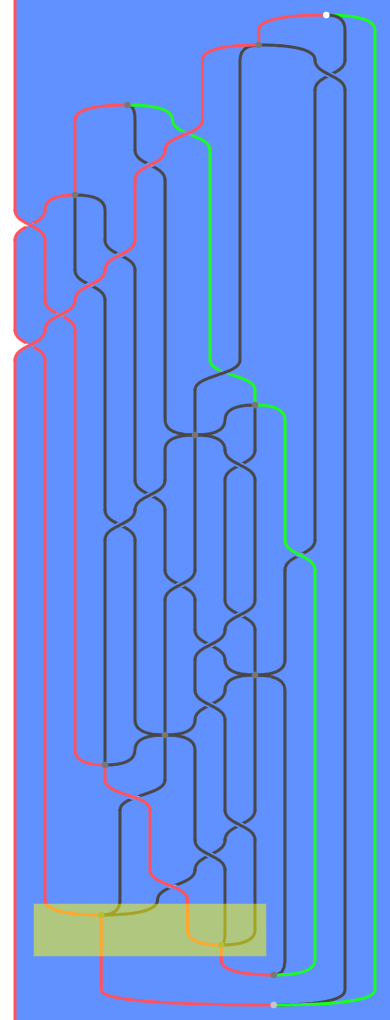}
\butterflystep{\I_4}{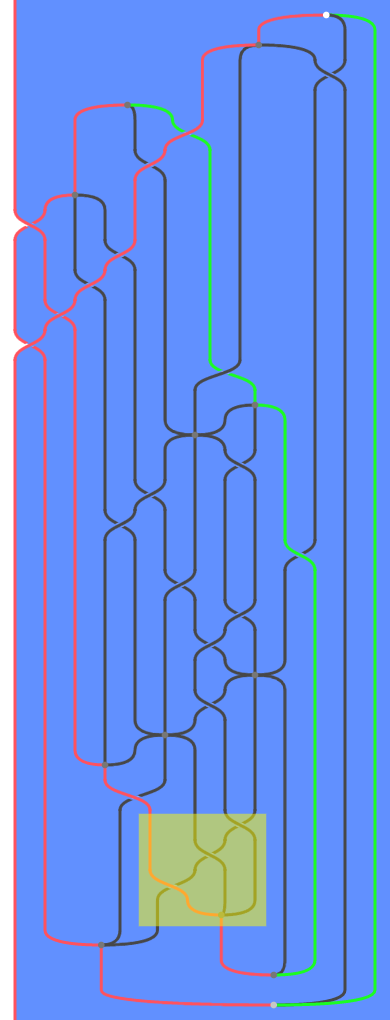}
\butterflystep{\II_4}{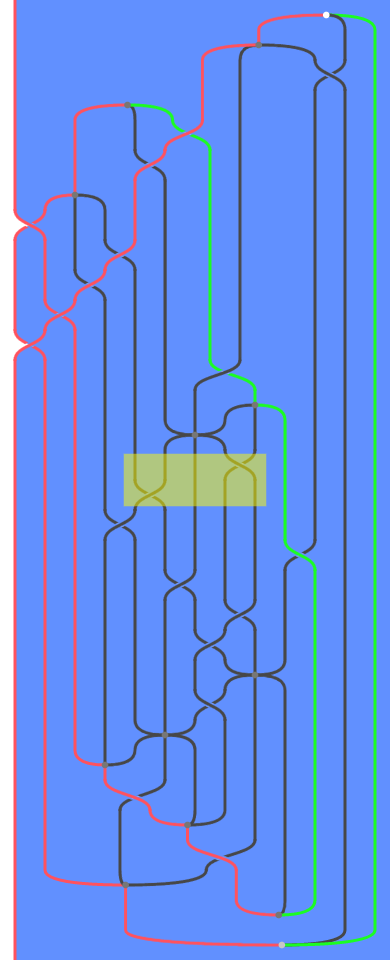}
\butterflystep{\I_4}{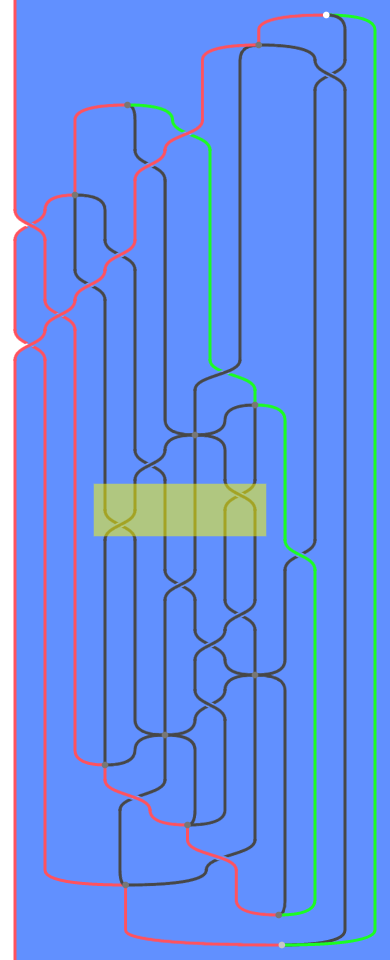}
\butterflystep{\I_4}{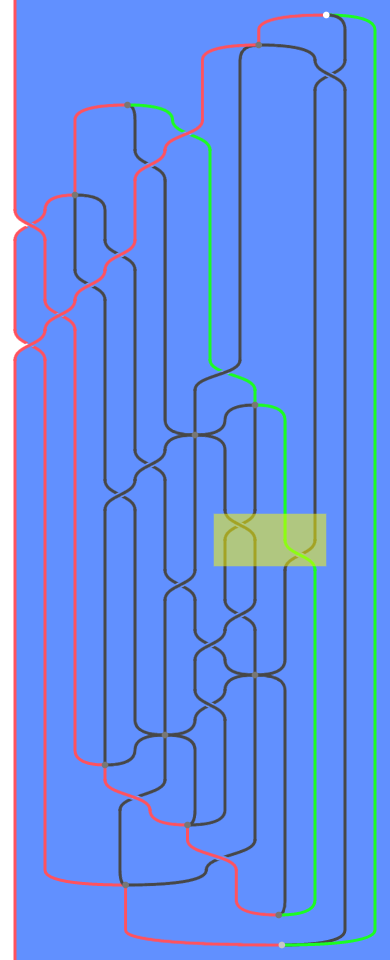}
\butterflystep{\I_4}{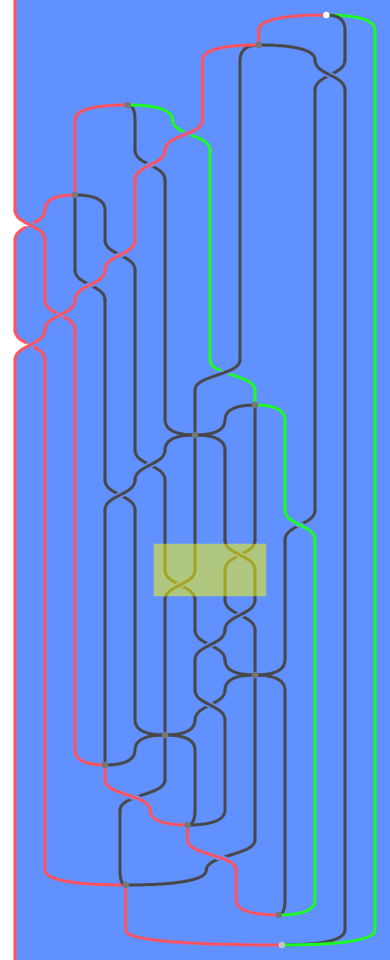}
\butterflystep{\I_4}{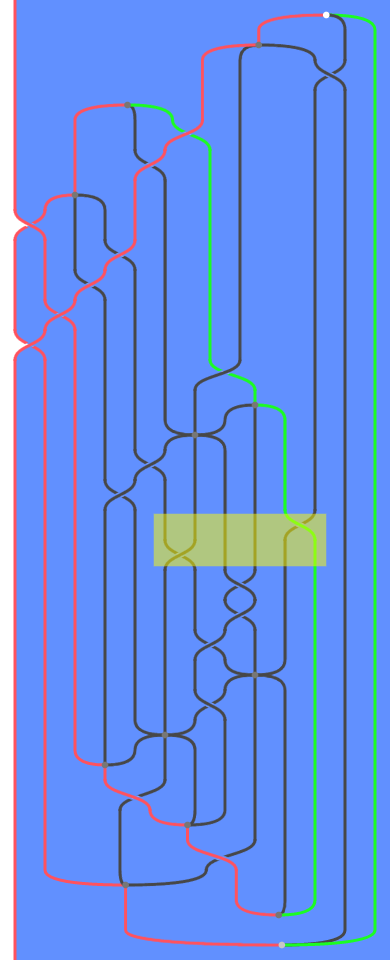}
\butterflystep{\I_4}{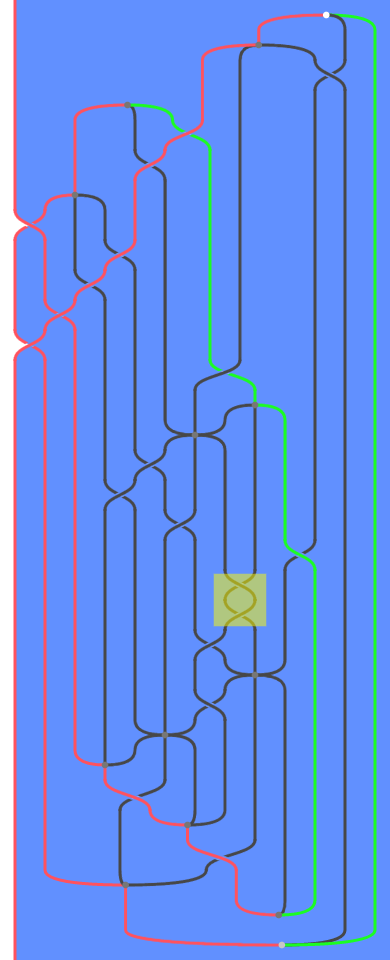}
\butterflystep{\I_3'}{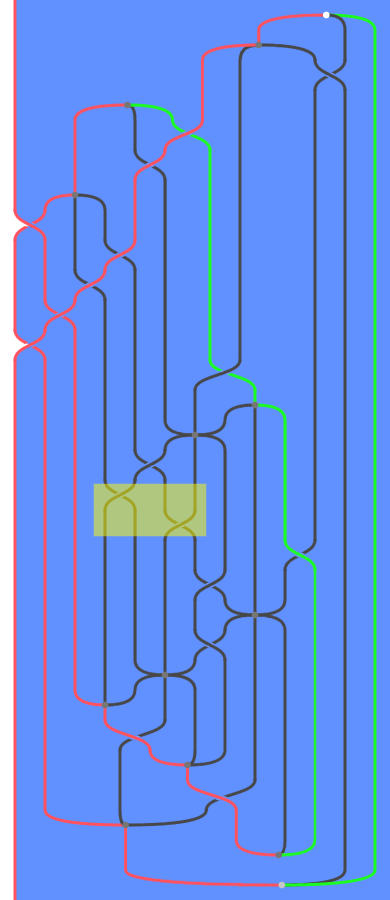}
\butterflystep{\I_4}{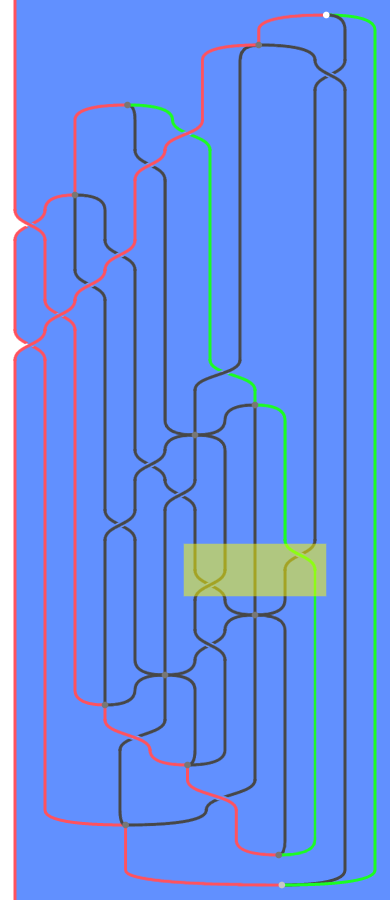}
\butterflystep{\I_4}{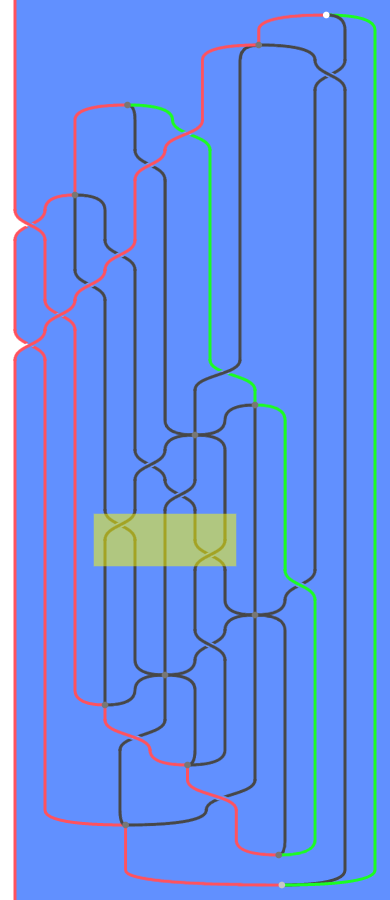}
\butterflystep{\I_4}{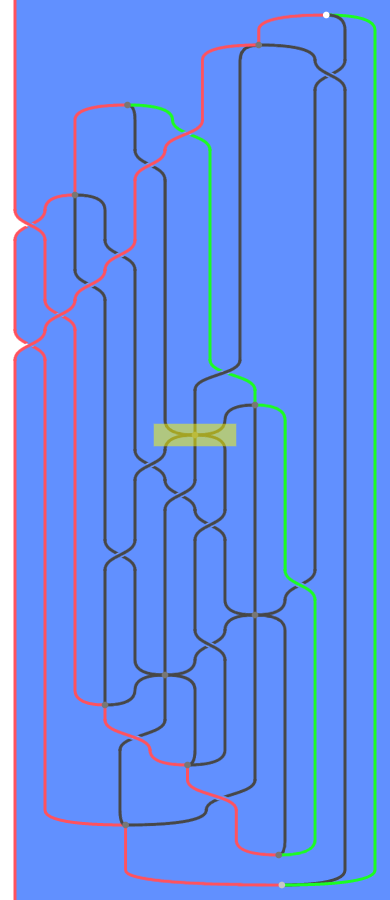}
\butterflystep{\VI_4}{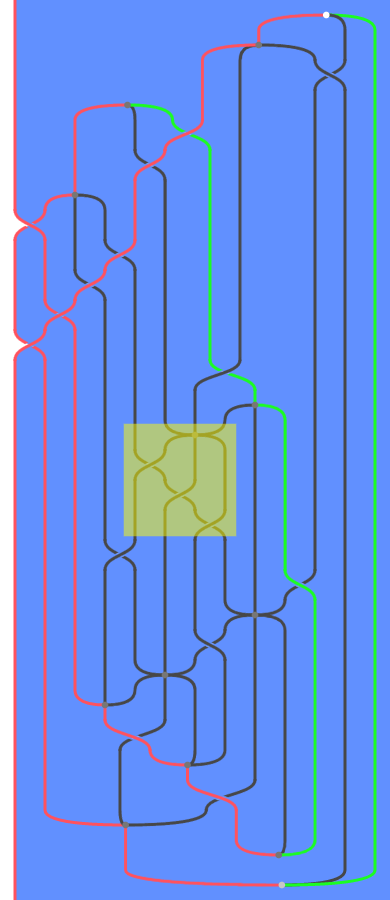}
\butterflystep{\II_4}{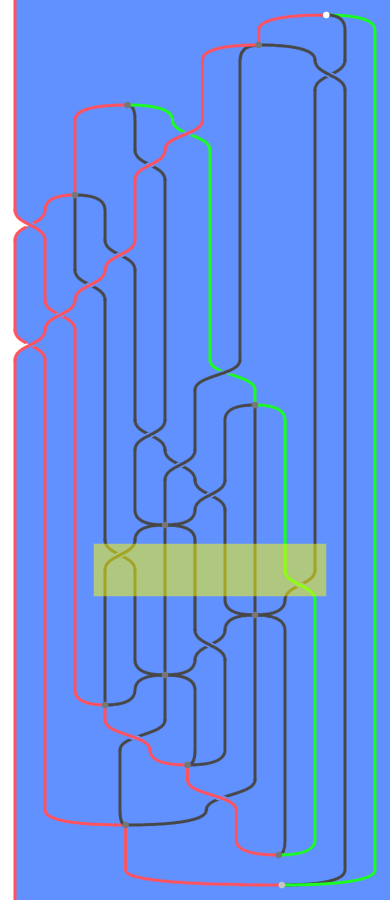}
\butterflystep{\I_4}{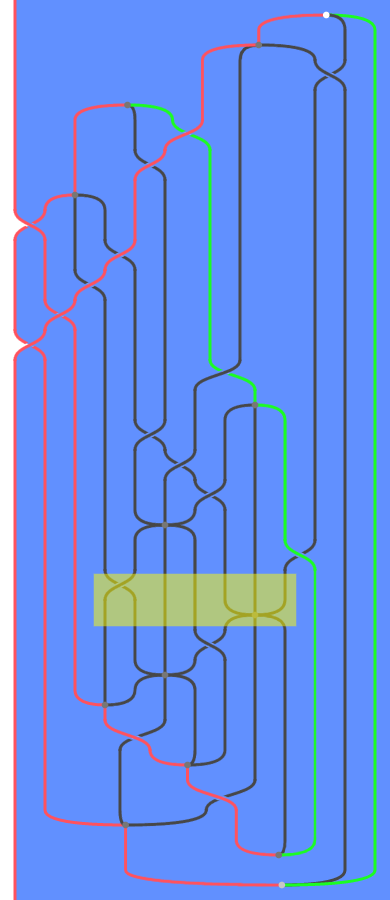}
\butterflystep{\I_4}{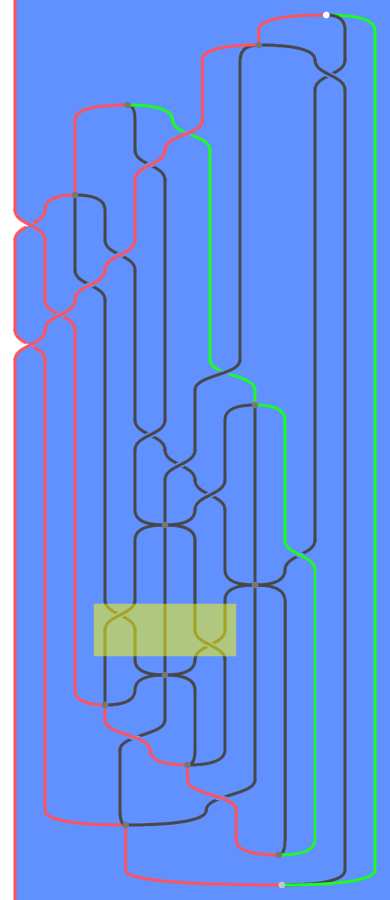}
\butterflystep{\I_4}{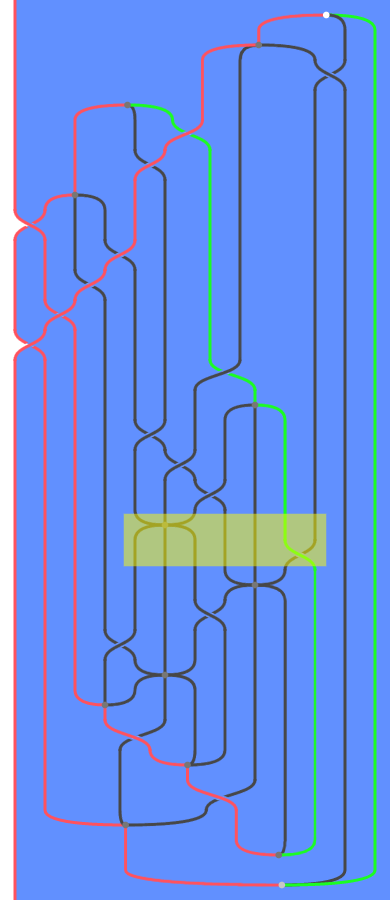}
\butterflystep{\I_4}{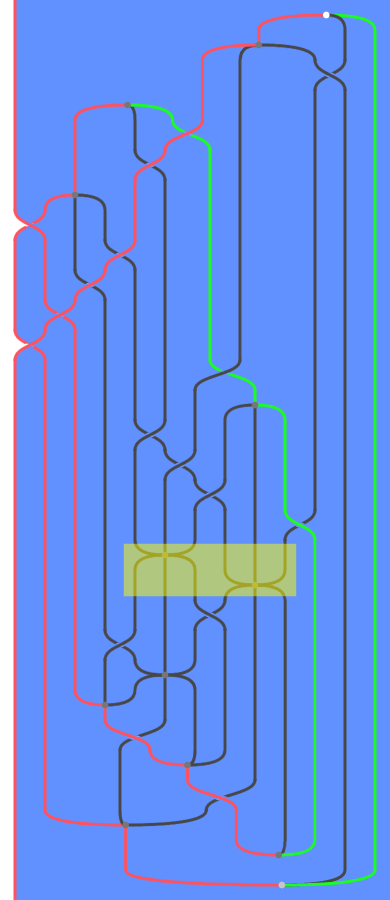}
\butterflystep{\I_4}{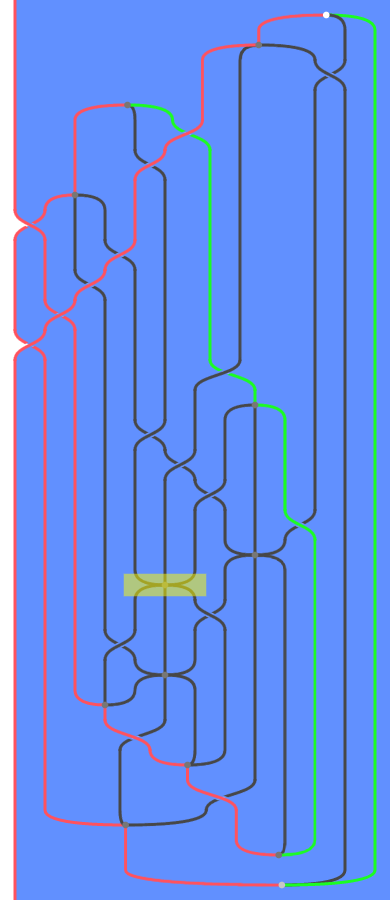}
\butterflystep{\VI_4}{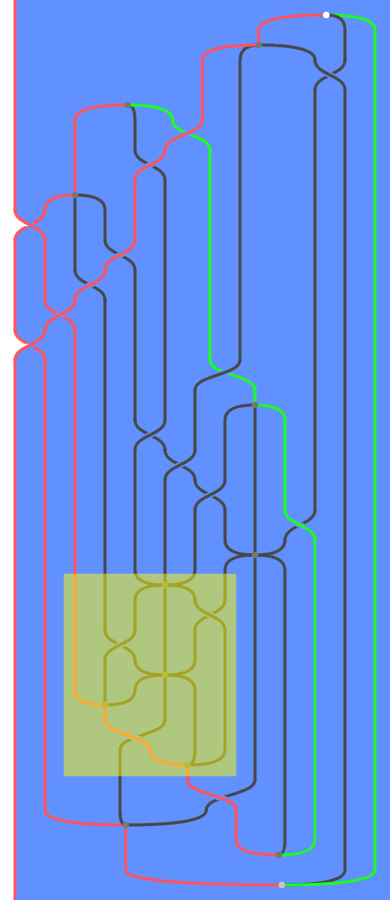}
\butterflystep{\III_4}{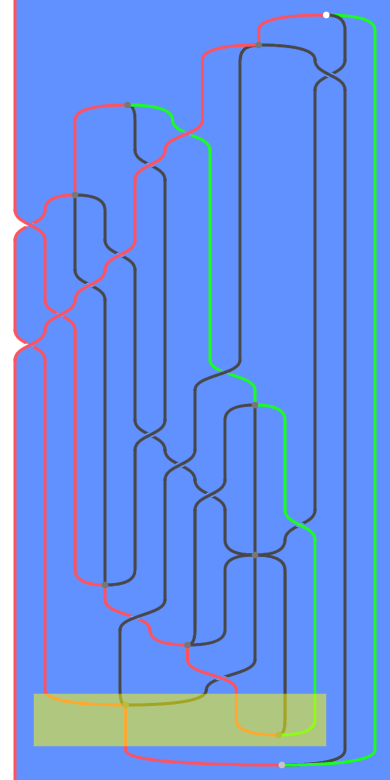}
\butterflystep{\I_4}{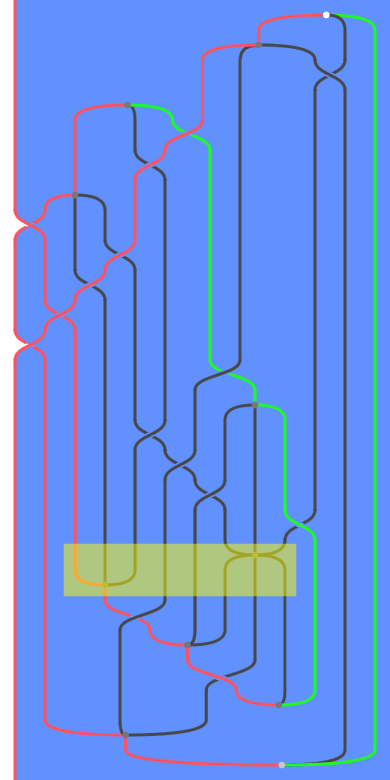}
\butterflystep{\I_4}{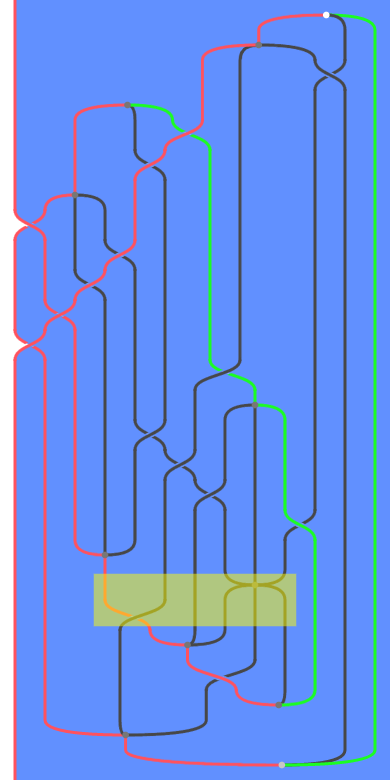}
\butterflystep{\I_4}{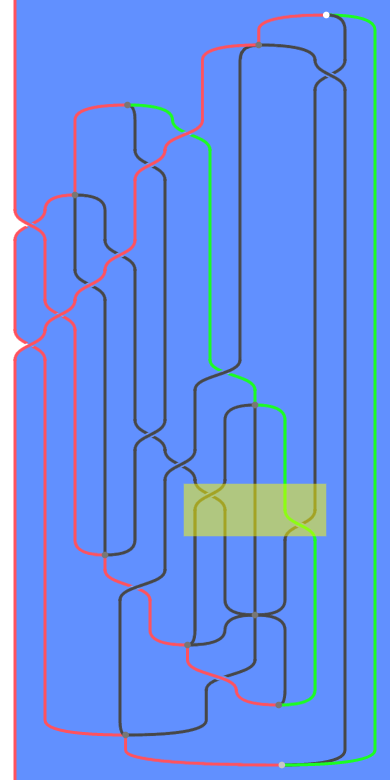}
\butterflystep{\I_4}{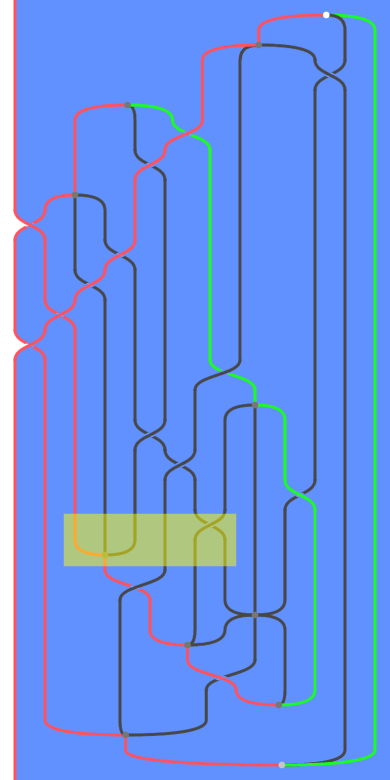}
\butterflystep{\I_4}{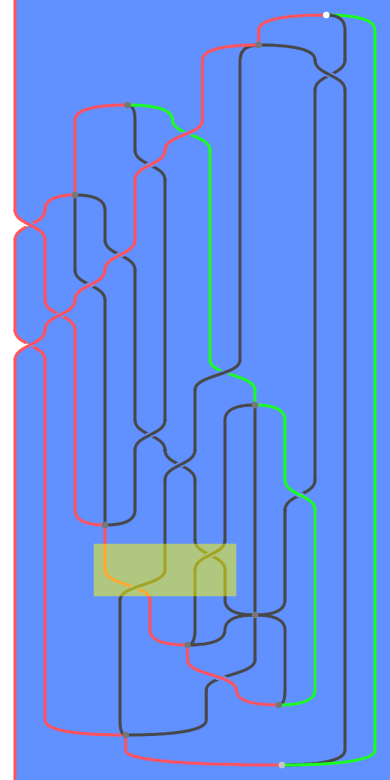}
\butterflystep{\I_4}{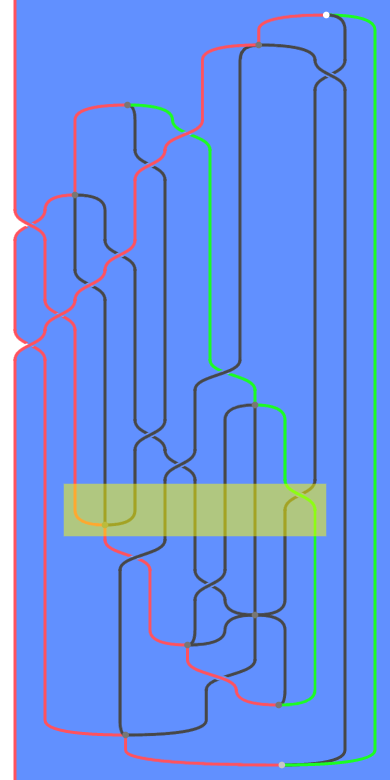}
\butterflystep{\I_4}{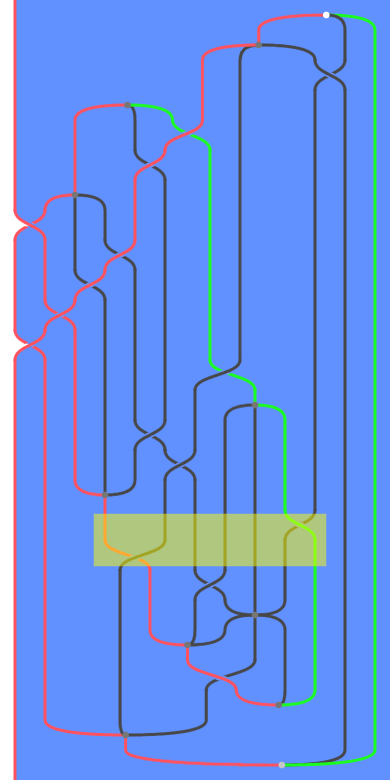}
\butterflystep{\I_4}{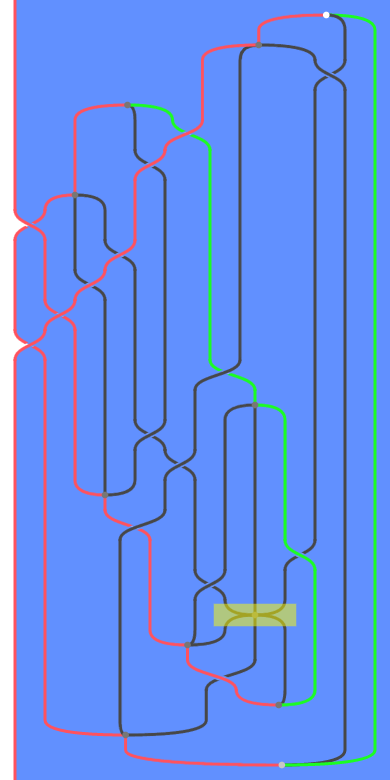}
\butterflystep{\VI_4}{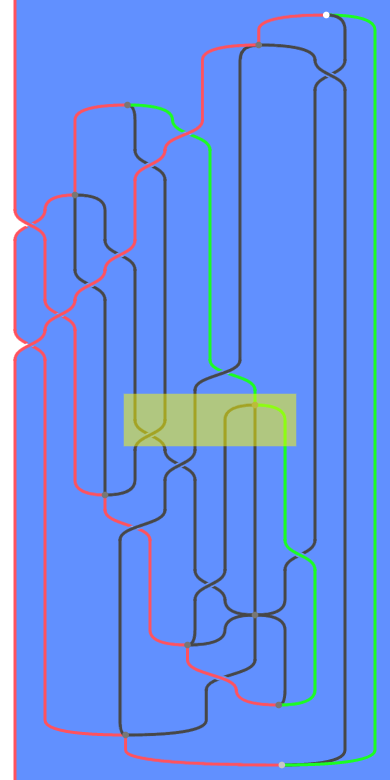}
\butterflystep{\I_4}{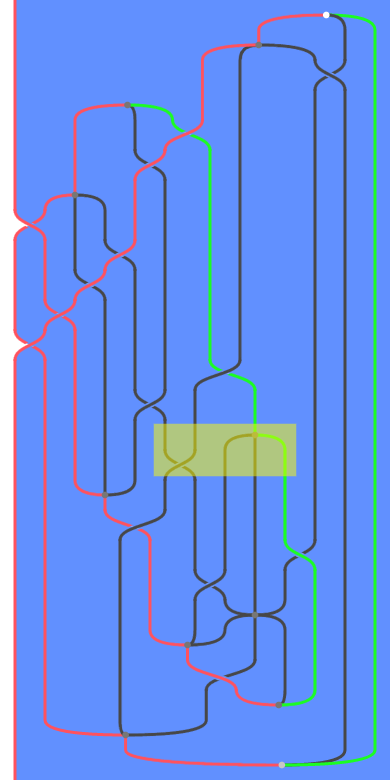}
\butterflystep{\I_4}{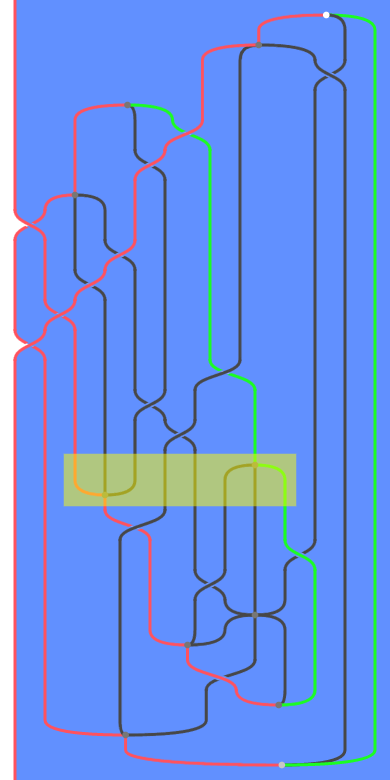}
\butterflystep{\I_4}{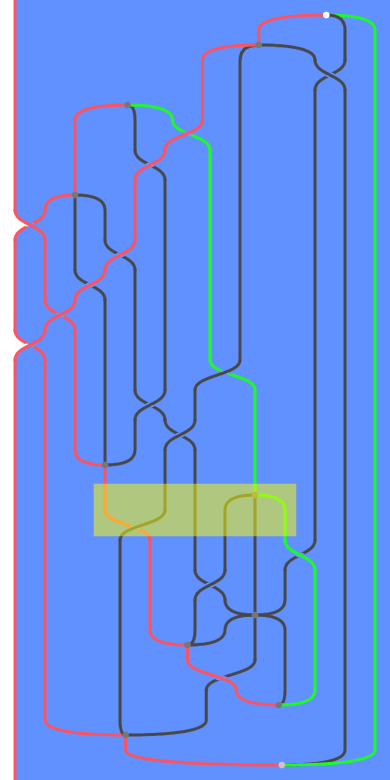}
\butterflystep{\I_4}{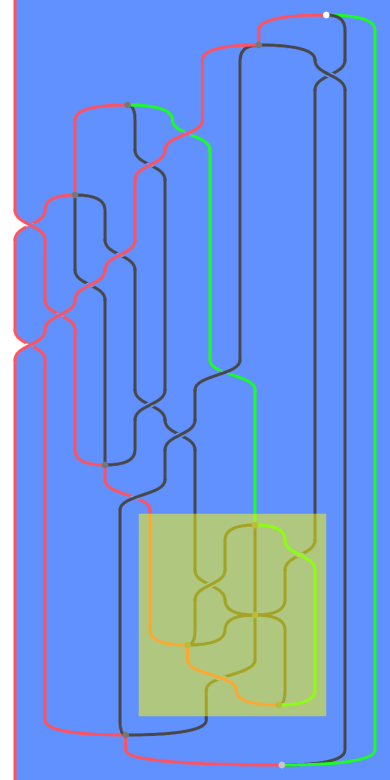}
\butterflystep{\I_2}{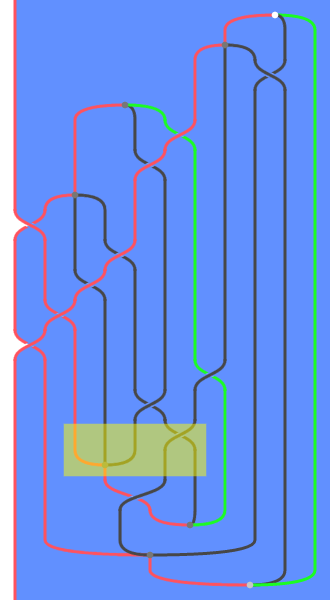}
\butterflystep{\I_4}{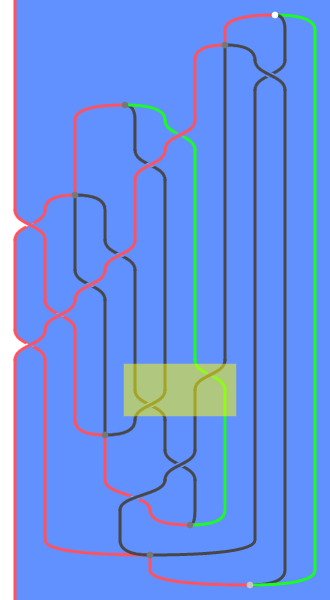}
\butterflystep{\I_4}{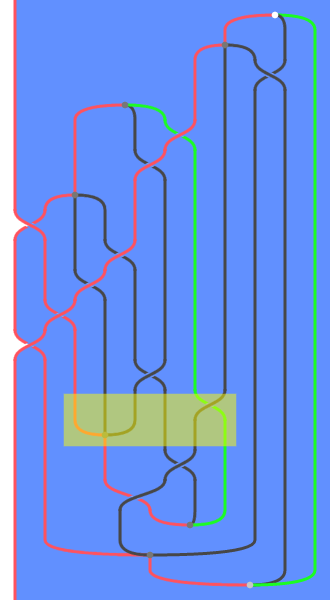}
\butterflystep{\I_4}{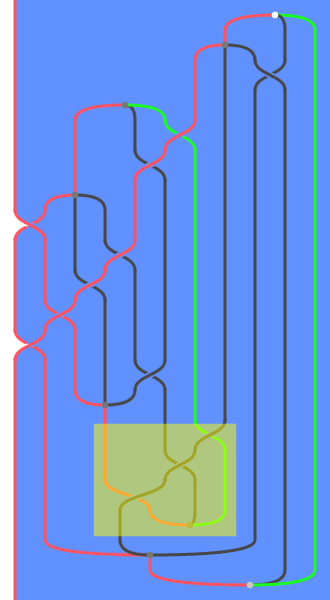}
\butterflystep{\II_4}{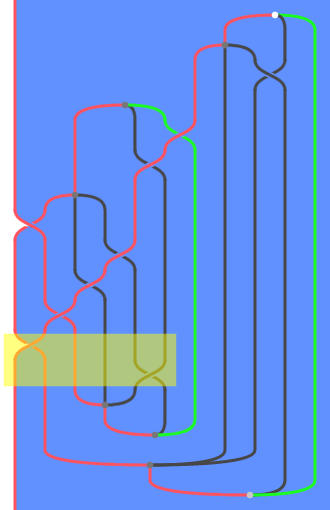}
\butterflystep{\I_4}{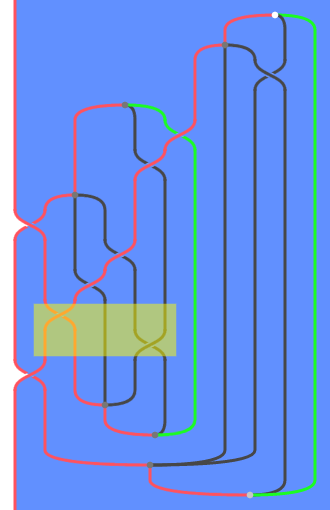}
\butterflystep{\I_4}{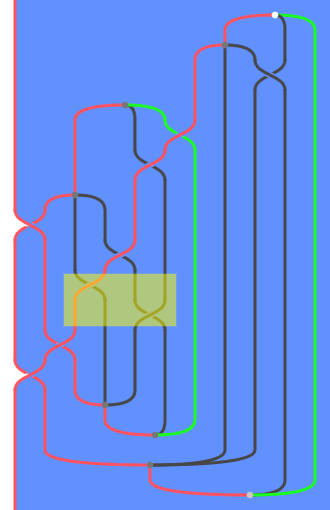}
\butterflystep{\I_4}{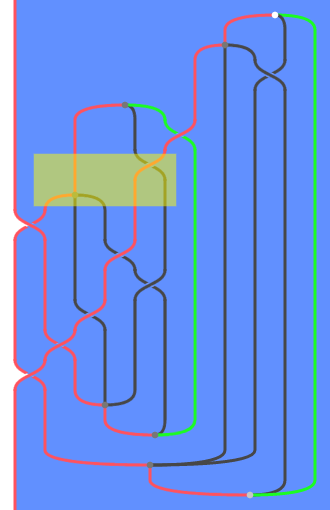}
\butterflystep{\I_4}{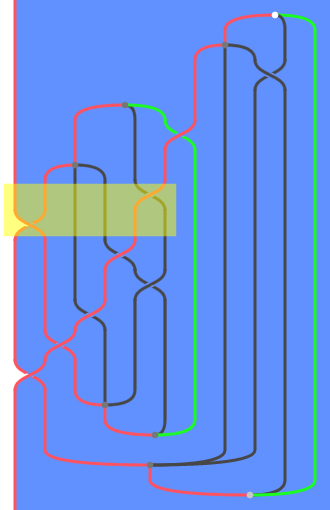}
\butterflystep{\I_4}{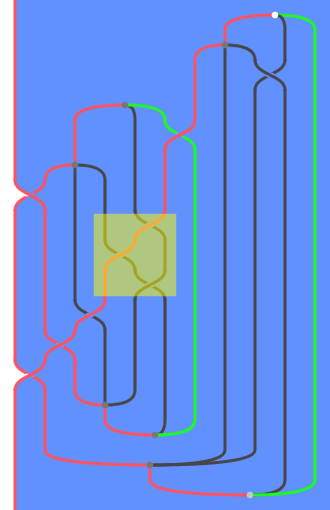}
\butterflystep{\II_4}{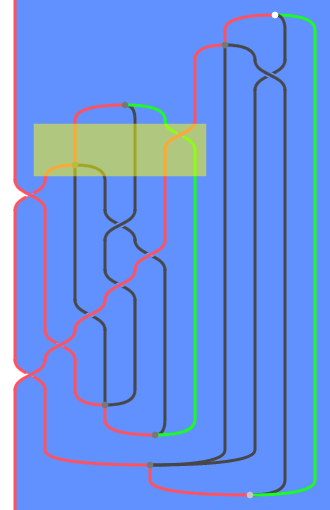}
\butterflystep{\I_4}{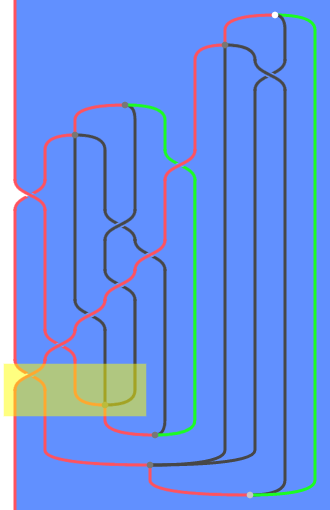}
\butterflystep{\I_4}{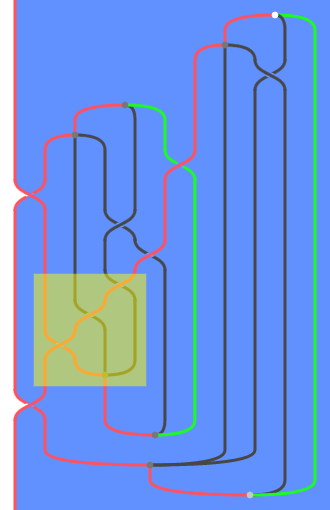}
\butterflystep{\II_4}{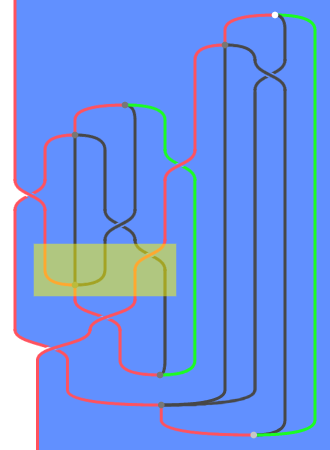}
\butterflystep{\I_4}{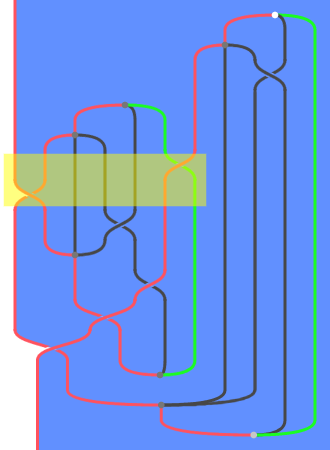}
\butterflystep{\I_4}{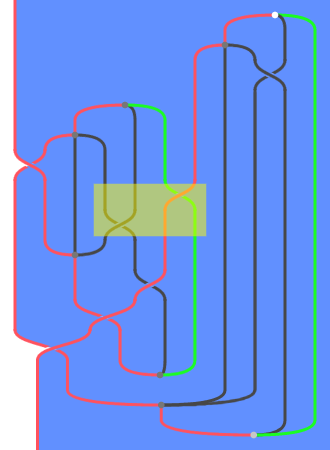}
\butterflystep{\I_4}{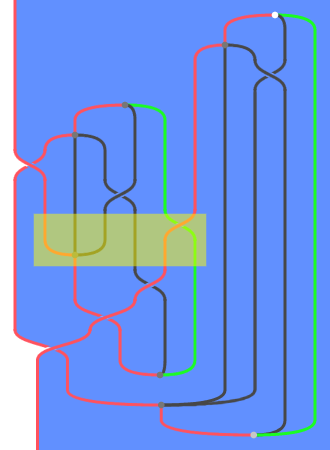}
\butterflystep{\I_4}{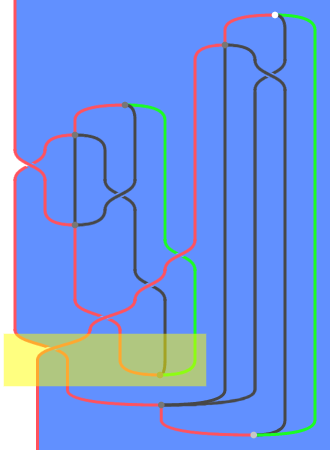}
\butterflystep{\I_4}{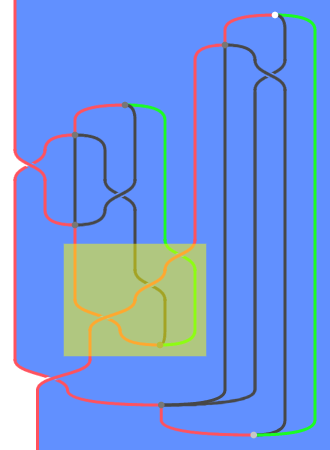}
\butterflystep{\II_4}{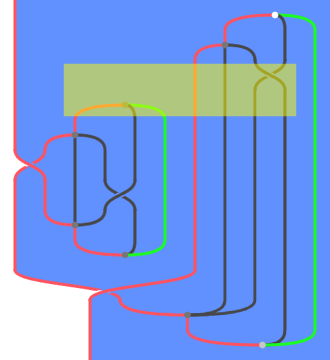}
\butterflystep{\I_4}{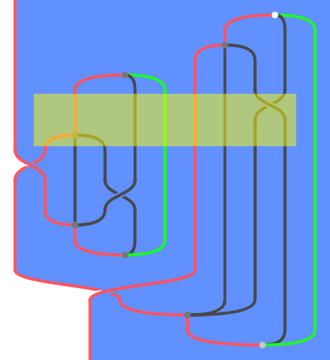}
\butterflystep{\I_4}{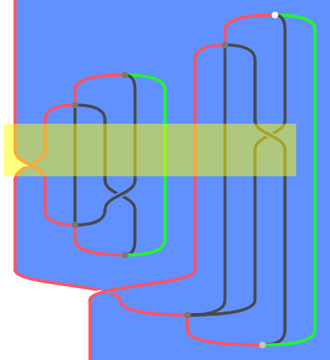}
\butterflystep{\I_4}{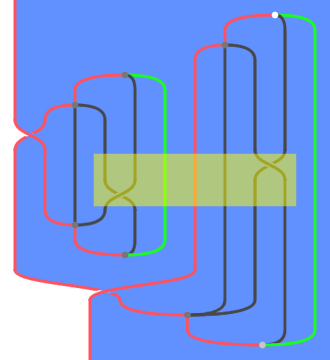}
\butterflystep{\I_4}{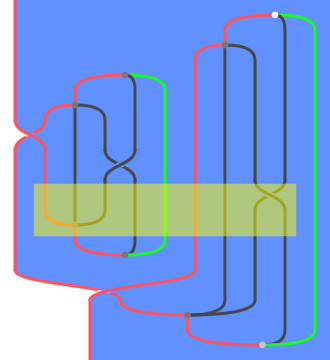}
\butterflystep{\I_4}{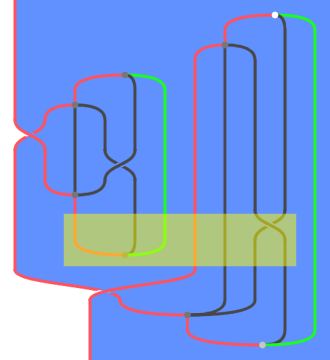}
\butterflystep{\I_4}{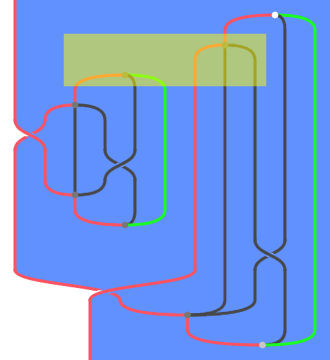}
\butterflystep{\I_4}{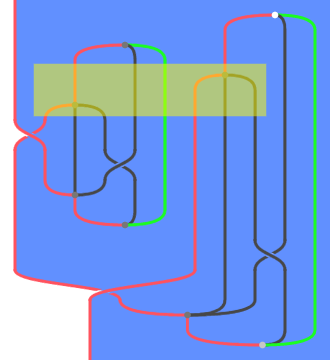}
\butterflystep{\I_4}{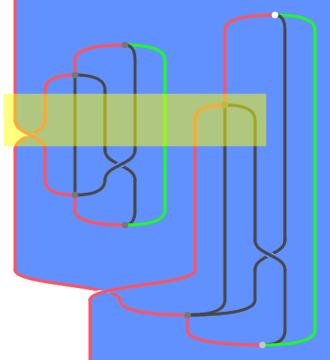}
\butterflystep{\I_4}{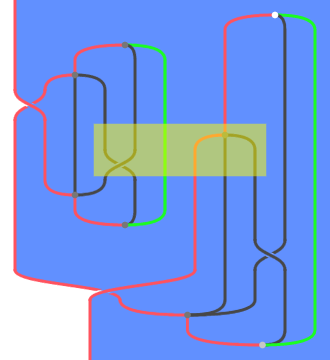}
\butterflystep{\I_4}{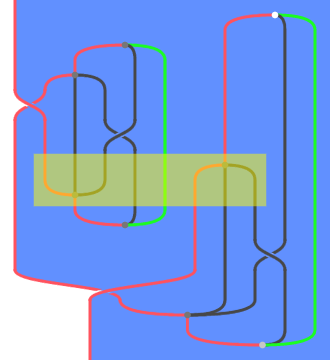}
\butterflystep{\I_4}{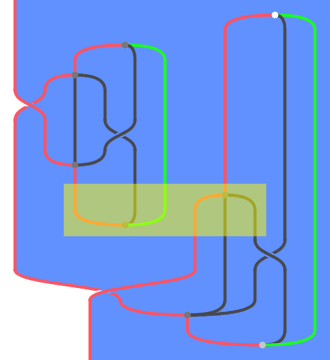}
\butterflystep{\I_4}{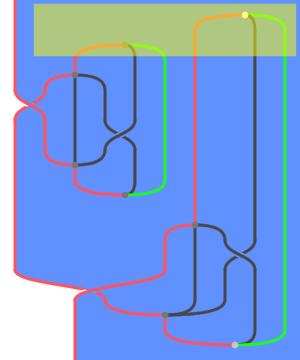}
\butterflystep{\I_4}{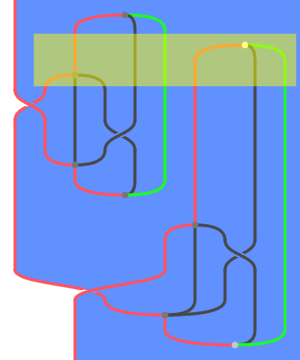}
\butterflystep{\I_4}{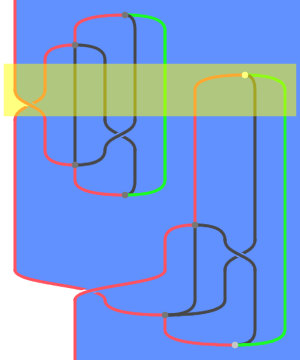}
\butterflystep{\I_4}{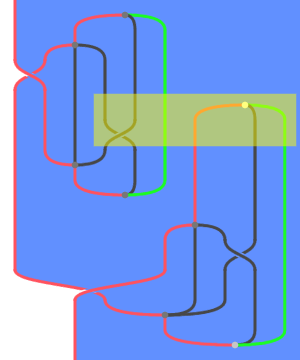}
\butterflystep{\I_4}{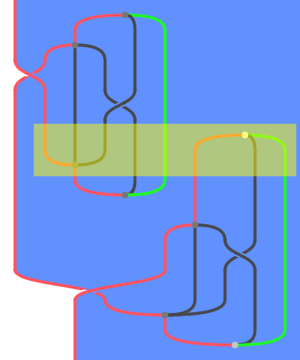}
\butterflystep{\I_4}{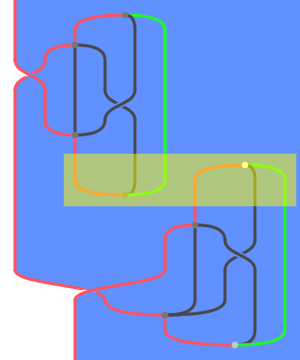}
\butterflystep{\I_4}{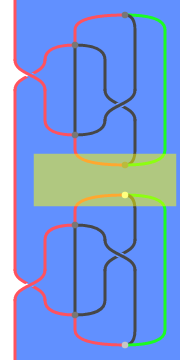}
\butterflystep{}{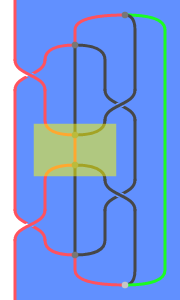}
\butterflystep{}{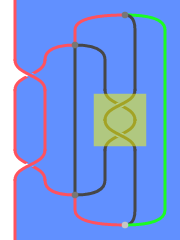}
\butterflystep{}{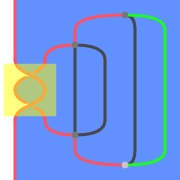}
\butterflystep{}{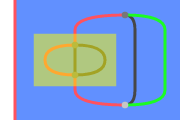}
\butterflystep{}{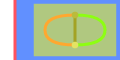}
\butterflystep{}{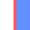}

\noindent
The 6 final steps are applications of the coinductive definition of invertibility.

\bibliographystyle{plainurl} 
\bibliography{references}

\end{document}